\documentclass[10pt]{amsart} 
\usepackage[utf8]{inputenc}
\usepackage [english]{babel}
\usepackage{a4wide}
\usepackage{amssymb}
\usepackage{amsmath}
\usepackage{amsthm}
\usepackage{amsbsy}
\usepackage{mathrsfs}
\usepackage{amsfonts}
\usepackage{mathtools}
\usepackage{verbatim}
\usepackage{tikz}
\usepackage{tikz-cd}
\usetikzlibrary{cd}
\usepackage{amscd}
\usepackage{hyperref}
\usepackage{xcolor,pifont}
\hypersetup{
	colorlinks,
	linkcolor={red!50!black},
	citecolor={blue!50!black},
	urlcolor={blue!80!black}
}
\usepackage[normalem]{ulem}
\usepackage{cleveref}   
\usepackage[all]{xy}
\usepackage{xcolor}
\theoremstyle{plain}
\usepackage[margin=1 in, a4paper]{geometry}
\usepackage{relsize}
\usepackage{lipsum}
\usepackage[autostyle]{csquotes}
\usepackage{relsize}
\usepackage{enumitem}  
\usepackage{titlesec}
\usepackage{setspace}
\usepackage{stmaryrd}
\allowdisplaybreaks
\titleformat{\section}
{\centering\large\sc}{\thesection.}{1em}{}
\titleformat{\subsection}
{\centering\normalfont\bfseries}{\thesubsection.}{1em}{}
\titleformat{\subsubsection}
{\normalfont\normalsize\bfseries}{\thesubsubsection.}{1em}{}
\titleformat{\subsubsection}
{\normalfont\normalsize\bfseries}{\thesubsubsection.}{1em}{}

\newcommand{\gl}{\operatorname{GL}}

\newcommand{\ind}{\operatorname{ind}}

\newcommand{\N}{\mathbb{N}}
\newcommand{\Z}{\mathbb{Z}}
\newcommand{\Q}{\mathbb{Q}}
\newcommand{\f}{\mathbb{F}_{p}}

\newcommand{\cmark}{\textcolor{red}{\ding{51}}} 
\newcommand{\xmark}{$\times$}    
\newcommand{\stkout}[1]{\ifmmode\text{\sout{\ensuremath{#1}}}\else\sout{#1}\fi}
\renewcommand{\mod}[1]{~\mathrm{mod}~ #1}

\newtheorem{theorem}{Theorem}[section]
\newtheorem*{theorem*}{Theorem}
\newtheorem{corollary}[theorem]{Corollary}
\newtheorem{lemma}[theorem]{Lemma}
\newtheorem*{lemma*}{Lemma}

\theoremstyle{definition}
\newtheorem{remark}[theorem]{Remark}
\newtheorem{definition}[theorem]{Definition}

\numberwithin{equation}{section}		
\numberwithin{figure}{section}			
\numberwithin{table}{section}				
\makeatletter

\renewcommand{\l@section}[2]{%
  \@dottedtocline{1}{0em}{2.3em}{#1}{#2}
}

\renewcommand{\l@subsection}[2]{%
  \@dottedtocline{2}{1.5em}{2.5em}{#1}{#2}
}

\def\@dotsep{10000}

\makeatother
\begin{document}
\title[REDUCTIONS OF CRYSTALLINE REPRESENTATIONS  OF FRACTIONAL SLOPE $<p-1$]{\textbf{REDUCTIONS OF CRYSTALLINE REPRESENTATIONS \\ OF FRACTIONAL SLOPE $< {p-1}$}}
\author{Shalini Bhattacharya, Eknath Ghate and Ravitheja Vangala}
\address{S. Bhattacharya, School of Mathematics \& Statistics, University of Hyderabad,  Hyderabad 500046, India}\email{shalinib@uohyd.ac.in} 
\address{E. Ghate, School of Mathematics, TIFR,  Mumbai 400005, India}\email{eghate@math.tifr.res.in}
\address{R. Vangala, Department of Mathematics, IISc,  Bangalore 560012, India}\email{ravithejav@iisc.ac.in}
\date{\today}

\begingroup
\def\uppercasenonmath#1{} 

\begin{abstract}
    Let $p$ be an odd prime and let $V_{k,a_p}$ be the two-dimensional crystalline representation of the Galois group of ${\mathbb Q}_p$ of weight $k \geq 2$ and parameter $a_p \in \bar{\mathbb Q}_p$. We study the semi-simplification $\bar{V}_{k,a_p}$ of the mod $p$ reduction of $V_{k,a_p}$ when the slope (valuation of $a_p$) is a positive fraction $< p-1$ using the mod $p$ local Langlands correspondence. We describe the {\it exact shape} of $\bar{V}_{k,a_p}$ for all such slopes and all (sufficiently large, depending on the slope) weights $k$, as long as certain Jordan-H\"older factors of dimension $p-1$ do not intervene in the computation (when $k$ is odd), though we also provide some criteria which further determine the shape of $\bar{V}_{k,a_p}$ in some of these exceptional cases.
    To keep this paper a reasonable length, we assume that for certain bad congruence classes of $k$ mod $p$, the slope is less than 
    the representative -  taken in the range $[1,p-1]$ - of the congruence class of 
    $k-2$ mod $(p-1)$,  which is generically the case if the slope is small. 
    Finally, a folklore conjecture predicts that the reduction $\bar{V}_{k,a_p}$ is {\it irreducible} for fractional slopes if $k$ is even. We deduce this conjecture for all fractional slopes $< p-2$ and all (sufficiently large, even)
    weights $k$ under the aforementioned slope assumption.
\end{abstract}

\maketitle
\endgroup
\tableofcontents

\section{Introduction}

Let $p$ be an odd prime. In this article, we compute the mod $p$ reductions of  certain two-dimensional crystalline representations of the Galois group $G_{\Q_p}$ of $\Q_p$.

Let $E$ be a finite extension field of $\Q_p$. There is an equivalence of categories between $n$-dimensional crystalline representations of $G_{\Q_p}$ defined over $E$ and the category of $n$-dimensional admissible filtered $\varphi$-modules over $E$ induced by Fontaine's functor $D_\mathrm{cris}$.  Let $a_p \in E$ with $v(a_p) > 0$ where $v$ is the valuation of $\bar\Q_p$ normalized so that $v(p) = 1$ and let $k \geq 2$. Let $V_{k,a_p}$ be the irreducible two-dimensional crystalline representation of $G_{\Q_p}$ over $E$ with Hodge-Tate weights $(0,k - 1)$ such that $D_{\mathrm{cris}}(V^\ast_{k,a_p} ) = D_{k,a_p}$, where $(\>)^\ast$ denotes contragredient and $D_{k,a_p}$ is the filtered $\varphi$-module defined in \cite[\S 2.3]{Ber11}. The local $p$-adic Galois representation at $p$ of a normalized cuspidal newform of weight $k \geq 2$, level coprime to $p$, trivial nebentypus at $p$ and with $p$-th Fourier coefficient $a_p$, is known to be $V_{k,a_p}$.

Let $\bar{V}_{k,a_p}$ be the semi-simplification of the reduction of a lattice in $V_{k,a_p}$ modulo the maximal ideal of the ring of integers of $E$. It is well-known that $\bar{V}_{k,a_p}$ is independent of the choice of lattice. It is an open problem to determine whether $\bar{V}_{k,a_p}$ is irreducible or reducible, let alone its exact shape. 

In recent decades, $\bar{V}_{k,a_p}$ has been studied extensively. The exact shape of $\bar{V}_{k,a_p}$ for $k \leq 2p+1$ was determined by the work of Fontaine-Edixhoven and Breuil  in \cite{Edi92} and 
\cite{Breuil}. The shape  of $\bar{V}_{k,a_p}$  when the slope $v(a_p)$ is greater than $\lfloor \frac{k-2}{p-1} \rfloor$ was determined by  Berger, Li and Zhu,  in \cite{BLZ}. A variant of this bound,  
namely $\lfloor \frac{k-1}{p} \rfloor$, was provided by Bergdall and Levin \cite{BL22}. 

 The reduction $\bar{V}_{k,a_p}$ has also been completely determined when $0< v(a_p) < 2$ by the work of Buzzard, Gee, Ganguli, Ghate, Bhattacharya, Rozensztajn and Rai in \cite{BG09}, \cite{BG13}, \cite{GG15}, \cite{BG15}, \cite{BGR18}, \cite{GR} (some of these results assume $p \geq 5$; more recently \cite{GM} have treated the case of $v(a_p) = 2$ and Nagel-Pande \cite{NP24} have partially treated the case $v(a_p) \in (2,3)$). 
 In this article, we determine the reduction $\bar{V}_{k,a_p}$ for $v(a_p) < p-1$
 for all {\it fractional slopes} less than $p-1$, i.e., for 
 $$v(a_p) \in (i,i+1)$$
 for $0 \leq i < p-1$.

 To state our results, we recall some notation. On the Galois side, let $\omega = \omega_1$ and $\omega_2$ denote the fundamental characters of level $1$ and $2$ respectively and let $\ind(\omega_2^c)$ be the unique irreducible representation of $G_{\Q_p}$  whose determinant is $\omega^c$ with whose restriction to the inertia group $I_p$ equals $\omega_2^c \oplus \omega_2^{pc}$ (for $p + 1 \nmid c$). Let $\mu_{\lambda}$ be the unramified character of $G_{\Q_p}$ mapping a (geometric) Frobenius at $p$ to $\lambda \in \bar{\mathbb{F }}_p^\times$.

On the automorphic side, we let $G=\mathrm{GL}_2(\Q_p)$, and write $K=\mathrm{GL}_2(\Z_p)$ and $Z=\Q_p^\times$ for its maximal compact subgroup and center respectively. For a $KZ$-representation $V$, let  $\mathrm{ind}_{KZ}^{G}(V)$ denote  the compact induction of $V$ from $KZ$ to $G$.  For a ${\mathbb Z}_p$-algebra $R$, let $\mathrm{Sym}^{r}R^2$ be the $r$-th symmetric power representation of $\mathrm{GL}_{2}(R)$. Let $V_r:=\mathrm{Sym}^r\bar{\mathbb{F}}^2_p$ be the base change to $\bar{\mathbb F}_p$ of the $r$-th symmetric power representation of $\Gamma:=\mathrm{GL}_{2}(\f)$. Let $D$ denote the determinant character $\Gamma$.
Both of these representations of $\Gamma$ are thought of as  representations of $K$ by inflation, and as representations of $KZ$ by making $p \in Z$ act trivially. 
Let
 	$$
 	\Pi_{k,a_p} := \frac{\ind_{KZ}^{G}(\mathrm{Sym}^{k-2} \bar{\Q}_p^2)}{T-a_p}
 	$$
    be the usual locally algebraic representation of $G$,
 where $T$ is the standard Hecke operator acting on $\mathrm{Sym}^{k-2}\bar{\Q}_p^2$
 and $a_p \in \bar{\Q}_p$. Let $\Theta_{k,a_p}$ denote the standard lattice in $\Pi_{k,a_p}$ which is the image of $\mathrm{ind}_{KZ}^{G}(\mathrm{Sym}^{k-2}\bar{\Z}_p^2)$ in $\Pi_{k,a_p}$, and $\bar{\Theta}_{k,a_p}$ be the semi-simplification of its reduction modulo $p$. 
 
 Recall that $\bar{\Theta}_{k,a_p}$ corresponds to $\bar{V}_{k,a_p}$ under the mod $p$ local Langlands correspondence.

   With this we can now state the first main result of this article:
 \begin{theorem}\label{main theorem good}
 	Let $p \geq 3$ and   $ v(a_p) \in(i, i+1)$ for some $0 \leq i < p-1$.  Suppose $k-2 =: r \geq i(p+1)+p$ and $r \equiv a ~\mathrm{mod}~(p-1)$ with $1 \leq a \leq p-1$.  Set
 	\begin{align*}
 	    b = \begin{cases}
 	     	a \quad &\mathrm{if}~i<a, \\
 	     	p-1+a & \mathrm{if}~i \geq a.
 	    \end{cases}
 	\end{align*}
    \begin{enumerate}
 	\item[$(i)$]If $b>2i+1$ and $ r\not \equiv b-i+1, \ldots, b~\mathrm{mod}~p$, then
 	    $$\bar{V}_{k,a_p} \simeq \ind(\omega_2^{b-i+1+ip}).$$
 	\item[$(ii)$] If $b=2i+1$ and $ r\not \equiv b-i+1, \ldots, b~\mathrm{mod}~p$, then $$\mathrm{ind}_{KZ}^{G}(V_{p-2} \otimes D^{i+1}) \twoheadrightarrow \bar{\Theta}_{k,a_p}.$$
 	However if $i<a$, and when $v(a_p)=i+\frac{1}{2}$, we further assume that
    \begin{eqnarray}
    \label{minimal val}
        v\left(a_p^2 - \binom{r-i}{i+1} \binom{r-i-1}{i}p^{2i+1}\right) = 2i+1,
    \end{eqnarray}
 	     then  $$\bar{V}_{k,a_p} \simeq \ind(\omega_2^{i+2+ip}).$$

 	\item[$(iii)$] If $b \leq 2i$ and $ r\not \equiv b-i-1, \ldots, b~\mathrm{mod}~p$, then $$\bar{V}_{k,a_p} \simeq \ind(\omega_2^{b-i-1+(i+2)p}).$$ 
    \end{enumerate}
 \end{theorem}

The representation $V_r$ is equipped with the theta filtration $V_r \supset V_r^{(1)}\cdots \supset V_r^{(m)} \supset \cdots$. If $r = k - 2$ and $v(a_p) \in (i,i+1)$, then there is a surjection $\ind_{KZ}^G (V_r/V_{r}^{(i+1)}) \twoheadrightarrow \bar{\Theta}_{k,a_p}$.  Using  past work of Ghate-Ravitheja \cite{GR19} and the first main technical result of this paper, namely Theorem~\ref{Elimination i < a and not in interval} below, we eliminate  all but one Jordan-H\"older (JH) factor of $V_r/V_{r}^{(i+1)}$ on the left. 
 See Corollary~\ref{Shape of Theta i < a not in interval} for the case $i < a$ and Corollary \ref{Shape of Theta i > a not in interval} for the
case $i \geq a$. Theorem~\ref{main theorem good} then follows using the mod $p$ 
local Langlands correspondence (see Section~\ref{sec: mod p LLC}) when $b \neq 2i+1$.

In the case
$b = 2i+1$, the remaining
JH factor is a twist of $V_{p-2}$ so that  
$\bar{V}_{k,a_p}$ may be reducible.
However, when $i < a$, we show in Theorem~\ref{thm: a=2i+1 good cases and non-half integer slopes},~\ref{a=2i+1 good cases and half integer slopes} that only the irreducible possibility occurs, under the additional condition \eqref{minimal val} which says that if $v(a_p)$ is at the mid-point of the interval $(i,i+1)$, then the  expression in \eqref{minimal val} has minimal valuation.  (We do not work out what happens when $b = 2i+1$ and $i \geq a$.)

 Theorem~\ref{main theorem good} treats the so called {\it good} congruence classes of $r$ mod $p$ in that it avoids certain trickier {\it bad} congruence classes of $r$ mod $p$. A large portion of this paper is devoted to treating these bad congruence classes. To keep this paper within reasonable limits, we assume that $i<a$ so that $b = a$. This is not a serious assumption if  $i$ is small (which was the case for many of the papers mentioned earlier in this introduction which this paper attempts to generalize). 
 
 When $a\geq 2i+1$,  the bad congruence classes are $r \equiv a-i+1,\ldots,a~\mathrm{mod}~ p$, whereas for $i<a\leq 2i$, they are  $r \equiv a-i-1,\ldots,a~\mathrm{mod}~ p$. Equivalently, these congruence classes are $$r\equiv a-i+n~\mathrm{mod}~ p$$ for some $1 \leq n \leq i$ (if $a\geq 2i+1$) or $-1\leq n \leq i$ (if $i<a\leq 2i$). For each such value of $n$, we define
\begin{align}
    s = a-i+n+(i-n)p 
\end{align}
and set
\[
 t=v(r-s).
\]
These numbers play a crucial role in determining the shape of $\bar{V}_{k,a_p}$. Indeed, by past work on the local constancy of reductions \cite{Ber12} and the work of \cite{BLZ}, \cite{BL22}  one has that in many cases $\bar{V}_{k,a_p} \simeq \bar{V}_{s+2,a_p} \simeq \bar{V}_{s+2, 0}$ for $t$ and $v(a_p)$ is sufficiently large. The last reduction has been completely determined by Breuil \cite{Breuil}. This shows the importance of the parameters $s$ and $t$. In this paper, we go much further and try to determine which sub-quotient in the theta filtration surjects onto $\bar{\Theta}_{k,a_p}$ for each value of $t \geq 1$. This usually determines the reduction $\bar{V}_{k,a_p}$ using the mod $p$ local Langlands correspondence.

With the notation as above, we prove the following theorems for the bad congruence classes of $r$ mod $p$. First, we consider the case $a > 2i$.
   \begin{theorem}\label{a>2i bad}
 	Let $p \geq 3$ and   $v(a_p) \in (i,i+1)$ for some $0 \leq i < p-1$.  Suppose $k-2 =: r \geq i(p+1)+p$ and $r \equiv a ~\mathrm{mod}~(p-1)$ with $1 \leq a \leq p-1$. Let $a \geq 2i+1$ and $r \equiv a-i+1,\ldots,a~\mathrm{mod}~ p$. If $a = 2i+1$ and $r\equiv a~\mathrm{mod}~p $, we further assume that  $v(a_p) \neq i+\frac{1}{2}$. If $ r \equiv a-i+n ~\mathrm{mod}~p$ for some $1 \leq n \leq i$ and $t = v(r-s) \geq 1$, then
		\begin{align*}
		    \bar{V}_{k,a_p} \simeq
		     \begin{cases}
		          \ind(\omega_2^{a-i+t+1+(i-t)p}) ~ & \text{ if } 1 \leq t \leq n,\\
		          \ind(\omega_2^{a-i+n+1+(i-n)p}) ~ & \text{ if } t \geq n.
		     \end{cases}
		\end{align*}
 \end{theorem}
This result is a consequence of Theorem~\ref{Shape theta a>2i bad congruence class} (for $n<i$) and Theorem~\ref{Shape theta a>2i bad n=i} (for $n=i$). 
We remark that the excluded case $a = 2i+1$, $r \equiv a \mod p$ and $v(a_p) = i + \frac{1}{2}$ is  trickier to handle and comes under the purview of the zig-zag conjecture (now proved in \cite{Gha22}, though the statement only treats $t \gg 0$). 

Next, we consider the case $a = 2i$.

\begin{theorem}\label{a=2i bad}
 	Let $p \geq 3$ and  $v(a_p) \in (i,i+1)$ for some $0 \leq i < p-1$.  Suppose $k-2 =: r \geq i(p+1)+p$ and $r \equiv a ~\mathrm{mod}~(p-1)$ with $1 \leq a \leq p-1$. Let $a =2i$ and $r \equiv a-i-1,\ldots,a~\mathrm{mod}~ p$.
	\begin{enumerate}
	 	\item[$(i)$] If $ r \equiv a-i+n ~\mathrm{mod}~p$ for some $-1 \leq n \leq i-1$ and $t=v(r-s)\geq 1$, then
		\begin{align*}
		    \bar{V}_{k,a_p} \simeq
		     \begin{cases}
		          \ind(\omega_2^{a-i+t+1+(i-t)p}) ~ & \text{ if } 1 \leq t \leq n,\\
		          \ind(\omega_2^{a-i+n+1+(i-n)p}) ~ & \text{ if } t \geq n.
		     \end{cases}
		\end{align*}
        \item[$(ii)$] If $ r \equiv a-i+n ~\mathrm{mod}~p$ with  $n=i$ and $t=v(r-s)\geq 1$, then
		\begin{align*}
		    \bar{V}_{k,a_p} \simeq
		     \begin{cases}
		          \ind(\omega_2^{a-i+t+2+(i-t-1)p}) ~ & \text{ if } 1 \leq t < n,\\
		          \ind(\omega_2^{a+1}) ~ & \text{ if } t \geq n.
		     \end{cases}
		\end{align*}
    \end{enumerate}
 \end{theorem}
 This result is a consequence of Theorem~\ref{Shape theta a=2i bad n=-1,0} (for $n=-1,0$), Theorem~\ref{Shape theta a>2i bad congruence class} (for $1\leq n<i$) and Theorem~\ref{Shape theta a<2i bad n=i} (for $n=i$).

 Finally, consider the case $a < 2i$.
 \begin{theorem}\label{a<2i bad n<i}
 	Let $p \geq 3$ and   $v(a_p) \in (i,i+1)$ for some $0 \leq i < p-1$.  Suppose $k-2 =: r \geq i(p+1)+p$ and $r \equiv a ~\mathrm{mod}~(p-1)$ with $1 \leq a \leq p-1$. Let $i< a  < 2i$. Assume $r\equiv a-i-1,\ldots,a-1 ~\mathrm{mod}~p$, that is, $r \equiv a-i+n ~\mathrm{mod}~p$ for some $-1 \leq n \leq i-1$ and let $t=v(r-s) \geq 1$. 
	\begin{enumerate}
	 	\item[$(i)$] For $a<2i-2n-1$,  we have 
        \begin{align*}
		    \bar{V}_{k,a_p} \simeq
		     \begin{cases}
		           \ind(\omega_2^{a-i+t-1+(i-t+2)p})~ & \text{ if } 1 \leq t \leq n+1,\\
		          \ind(\omega_2^{a-i+n+1+(i-n)p}) ~ & \text{ if } t \geq n+2.
		     \end{cases}
		\end{align*}
        \item[$(ii)$] For $a=2i-2n-1$, we have
        \begin{enumerate}
            \item[$(a)$] {\makebox[6.5 cm]{$\bar{V}_{k,a_p} \simeq \ind(\omega_2^{a-i+t-1+(i-t+2)p})$ \hfill} if $1 \leq t \leq n$.}
            \item[$(b)$] {\makebox[6.5 cm]{$\mathrm{ind}_{KZ}^{G}(V_{r}^{(i-n)}/V_r^{(i-n+1)}) \twoheadrightarrow \bar{\Theta}_{k,a_p}$ \hfill} if $t\geq n+1$.} 
        \end{enumerate}
        \item[$(iii)$] For $a>2i-2n-1$, we have
        \begin{enumerate}
	        \item[$(a)$] {\makebox[6.5 cm]{$ \bar{V}_{k,a_p} \simeq \ind(\omega_2^{a-i+t-1+(i-t+2)p})$ \hfill} for $1\leq t \leq i- \frac{a}{2}$} 
	       \item[$(b)$] {\makebox[6.5 cm]{$ \mathrm{ind}_{KZ}^{G}(V_r^{(a-i+t-1)}/V_r^{(a-i+t)}) \twoheadrightarrow \bar{\Theta}_{k,a_p} $ \hfill} for $ t = i- \frac{a-1}{2}$}
	       \item[$(c)$] {\makebox[6.5 cm]{$ \bar{V}_{k,a_p} \simeq \ind(\omega_2^{a-i+t+1+(i-t)p})$ \hfill}  for $ i- \frac{a-1}{2} < t \leq n$}
	       \item[$(d)$] {\makebox[6.5 cm]{$ \bar{V}_{k,a_p} \simeq \ind(\omega_2^{a-i+n+1+(i-n)p})$  \hfill} for $t\geq n+1$ and $a\neq 2i-2n+1$}
	     \item[$(e)$] {\makebox[6.5 cm]{$\mathrm{ind}_{KZ}^{G}(V_{p-2}\otimes D^{i-n+1}) \twoheadrightarrow \bar{\Theta}_{k,a_p}$ \hfill} for $t\geq n+1$ and $a=2i-2n+1$,} 
\end{enumerate}
        where we assume $i=p-2 \Longrightarrow a \neq p-1$ if $2i-n-a < t \leq n$.
    \end{enumerate}
 \end{theorem}
 Part $(i)$ follows from  Theorem~\ref{a<2i bad n=-1} and Theorem~\ref{diagonal conj a-i<i}. Part $(ii)$ follows from  Theorem~\ref{diagonal conj a-i<i}. Part $(iii)$ follows from Theorems~\ref{thm: Shape theta a<2i-n} and \ref{Shape theta 2i-n<a}. 

In part $(ii)$ $(b)$ and $(iii)$ $(b)$
 of Theorem~\ref{a<2i bad n<i}, there is the possibility that a twist of the JH factor $V_{p-2}$ contributes to the reduction, and it becomes difficult to specify the reduction without more work. In the former case, we note that, in fact, a local constancy argument guarantees that this factor appears (and contributes reducibly) for $t \gg 0$. 
 In the latter case, taking $n = 1$ and $a = 2i-1$ so that $t =1$, it turns out that again a twist of $V_{p-2}$ appears and  contributes both irreducibly and reducibly,  e.g., when $i = 2$ (see Remark~\ref{r = 27}). 
 
 In part $(iii)$ $(e)$ of Theorem~\ref{a<2i bad n<i}, it is proved that $V_{p-2}$ definitely contributes 
 if $a=2i-2n+1$ and $t \geq n+1$. Thus, $\bar{V}_{k,a_p}$ could again possibly be reducible. In particular, if $n=1$, the reducibility question arises again when $a = 2i-1$ and $t \geq 2$. 
 
 We offer the following rather delicate theorem when $a = 2i-1$, $n=1$ and $t \geq 1$, which shows that the reduction is generically irreducible.

       \begin{theorem}
             Let $p \geq 3$ and   $v(a_p) \in (i,i+1)$ for some $0 \leq i < p-1$.  Suppose $k-2=:r \geq i(p+1)+p $,  $r \equiv 2i-1 \mod{(p-1)}$ and $r\equiv i \mod p$ with $ 2 \leq i \leq \frac{p-1}{2} $. Let $$-d = \frac{1}{p} \binom{r-i+1}{i} + \frac{(-1)^{i+1}}{i}.$$ Furthermore, when $ v(a_{p})  = i + \frac{1}{2} $,  assume that $$v(a_p^2 - i d^2 p^{2i+1}) = 2i+1.$$ Then $\bar{V}_{k,a_p} \simeq \ind(\omega_2^{i-1+(i+1)p})$ is irreducible.
        \end{theorem}
        
        This result follows from Theorems~\ref{a=2i-1, left and right half, irred}, \ref{a=2i-1, right half, p not divides d} and \ref{a=2i-1, mid-point, irred}. Note that the binomial coefficient is divisible by $p$ by Lucas' theorem. Note also the similarity between the minimal valuation assumption being made at the mid-point here with condition \eqref{minimal val}. 

        In Theorem~\ref{a<2i bad n<i} the case $n=i$, that is $r\equiv a \mod p$ was left untreated. We prove the following.
  \begin{theorem}
\label{intro: n = i}
 	Let $p \geq 3$ and   $v(a_p) \in (i,i+1)$ for some $0 \leq i < p-1$.  Suppose $k-2 =: r \geq i(p+1)+p$ and $r \equiv a ~\mathrm{mod}~(p-1)$ with $1 \leq a \leq p-1$. Let $i< a  < 2i$. Assume $r\equiv a-i+n~\mathrm{mod}~p$ with $n=i$ and let $t=v(r-a)$.
	\begin{enumerate}
	 	\item[$(i)$] If  $t < a-i$, then $ \mathrm{ind}_{KZ}^{G}(V_r^{(a-i-t-1)}/V_r^{(a-i-t)}) \twoheadrightarrow \bar{\Theta}_{k,a_p} $. In addition, if $2i+2t+2-a \neq p-2$ then $ \bar{V}_{k,a_p} \simeq \ind(\omega_2^{a-i-t-1+(i+t+2)p})$.
	    \item[$(ii)$] If $t\geq a-i$, then $ \bar{V}_{k,a_p} \simeq \ind(\omega_2^{a+1})$.
    \end{enumerate}
 \end{theorem}
 This result follows from  Theorem~\ref{Shape theta a<2i bad n=i} and the fact that in part $(i)$,  the sub-quotient $V_r^{(a-i-t-1)}/V_r^{(a-i-t)}$ does not have a JH factor that is twist of $V_{p-2}$ if $2i+2t+2-a \neq p-2$.

It is a folklore conjecture, attributed to Breuil, Buzzard, and Emerton, that if $k$ is even and $v(a_p)$ is not an integer, then
$\bar{V}_{k,a_p}$ is irreducible.
The results above allow us to make considerable headway on this conjecture for slopes up to $p-2$.

\begin{corollary}
Let $p \geq 3$ and $v(a_p) \in (i,i+1)$ for $0 \leq i < p-2$ and let  $k-2 =:r \geq 
4i+4$
with $r \equiv a ~\mathrm{mod}~(p-1)$ with $1 \leq a \leq p-1$. 
If $k$ is even, then
$\bar{V}_{k,a_p}$ is irreducible if either
\begin{itemize}
    \item[$(i)$] $i < a$, or \item[$(ii)$] $i \geq a$, and $r$ is in a good congruence class mod $p$.
\end{itemize}
\end{corollary}

\begin{proof}
Assume first that $r \geq i(p+1)+p$. An inspection of the (proofs of the) theorems above shows that 
there is a surjective map  
$\mathrm{ind}_{KZ}^{G}(V_{r}^{(m)}/V_r^{(m+1)}) \twoheadrightarrow \bar{\Theta}_{k,a_p}$
for some $0 \leq m \leq i$. Moreover, if $r$ is even, then $a-2m$ is even, so it cannot be congruent modulo $(p-1)$ to $1$ or $p-2$. This shows that $V_{r}^{(m)}/V_r^{(m+1)}$ does not contain (a twist of) $V_{p-2}$ as a JH factor. This forces $\bar{V}_{k,a_p}$ to be irreducible by the mod $p$ local Langlands correspondence.

Now assume that $4i+4\leq r < i(p+1)+p$. Since $i\leq p-3$, we get $r+1 \leq p^2-p-3$. In the notation of \cite[Theorem B]{Ber12}, we have $\alpha(r+1) = \lfloor \frac{r+1}{p-1} \rfloor \leq i+2$. Hence $r+2 > 3 v(a_p)+\alpha(r+1)+1$ and it follows from \cite[Theorem B]{Ber12} that $\bar{V}_{k,a_p} \simeq \bar{V}_{k',a_p}$ whenever $k'-k \equiv 0 \mod p^m(p-1)$ and $m \gg 0$. Taking $k'$ large enough we are reduced to the previous paragraph. 
\end{proof}

\begin{remark}
We note:
    \begin{enumerate}
        \item \cite{Breuil1} has computed  $\bar{V}_{k,a_p}$ completely when $r<2p$. If $i\leq \frac{p-2}{2}$, then $4i+4 \leq 2p$, so the corollary holds for all $r \geq 0$, for such $i$.
        This extends the work of \cite[Theorem 1.1]{Ars21} who proved 
        the conjecture for slopes  $< \frac{p-1}{2}$ and $a > 2i +1$.
        See also the forthcoming work \cite{LTXZ} which proves the conjecture under a `very generic' hypothesis.
       \item The trick in the corollary of using a local constancy argument to compute the reduction for $r$ smaller than $i(p+1)+p$ can also be used to extend the range of $r$ in Theorems~\ref{main theorem good} - \ref{intro: n = i} above to some smaller $r$. We leave the details to the reader.
   \end{enumerate}
\end{remark}

A final word about the proofs. It turns out that
if  $r = k - 2$ and $v(a_p) \in (i,i+1)$, then there is a further surjection $\ind_{KZ}^G Q(i) \twoheadrightarrow \bar{\Theta}_{k,a_p}$
where $$Q(i) = \dfrac{V_r}{X_{r-i} + V_{r}^{(i+1)}}, $$
for $X_{r-i}$ the 
submodule of $V_r$ generated by the
$i$-th monomial in 
a polynomial model of $V_r$. An exhaustive study of $X_{r-i}$ and the quotient
$Q(i)$ was made in \cite{GR19} for $0 \leq i \leq p-1$. Indeed, 
in that paper the cases for which 
$Q(i)$ is irreducible were described completely, allowing the authors to write down the structure of $\bar{V}_{k,a_p}$
in these cases (with the usual exception when the dimension of the JH factor is $p-1$), see \cite[Corollary 1.12]{GR19}.
In this paper, we go much further. In 
Chapter~\ref{sec: JH factors}, we derive 
a complete list of all JH factors of $Q(i)$ using the results from \cite{GR19}. In  Chapter~\ref{sec: good}, we consider the good congruence classes of $r $ mod $p$. We use explicit computations with the Hecke operator $T$ to eliminate all but one JH factor
of $Q(i)$ using the important Theorem~\ref{Elimination i < a and not in interval}. 

In Chapter~\ref{sec: bad}, we turn to the bad congruence classes of $r$ mod $p$, though as mentioned above we assume that $i < a$. Assume momentarily that $n \lneq i$. The cases where $a \leq 2i-2n -1$ are treated in Section~\ref{sec: bad small}. We first eliminate the shallow JH factors in $Q(i)$ using Theorem~\ref{Elimination i < a and not in interval}. Then we eliminate all but
one of the remaining deep JH factors so that  a certain \textit{diagonal} pattern of JH factors survives. This is described in the pictures Figures~\ref{fig:diagonal conj. a=2i-2n-1}, \ref{fig:diagonal conj. a<2i-2n-1}. In Section~\ref{sec: bad large}, we treat the cases $a \geq 2i$ following a similar strategy. This time, the deep JH factors of $Q(i)$ that survive are explained by a \textit{superdiagonal} pattern, see Figure~\ref{fig:super diagonal conj. a>2i}. In the intermediate region $2i-2n -1 < a < 2i$, the surviving JH factors are determined in Section~\ref{sec: bad medium} and are explained by a striking \textit{hybrid} version of the previous two patterns (with an additional line segment parallel to the \textit{anti-diagonal}!). See Figure~\ref{fig: hybrid a odd} for the
case of $a$ odd and 
Figure~\ref{fig: hybrid a even} for the case of $a$ even. The patterns in
these figures took several years to determine and the final results caught us by surprise. The proofs of these patterns involve the construction of certain master functions which contain certain constants $\beta_l$ whose existence depends on the evaluation of certain large determinants built out of double binomial sums. We decided to give complete details about the computation of these determinants, partly explaining the size of this paper. Finally, in the boundary case $n = i$, the JH factors of $Q(i)$ that survive are determined in Section~\ref{sec: n=i}.
 
  All the results in Theorems~\ref{main theorem good} through \ref{intro: n = i} now follow from the mod $p$ local Langlands correspondence.

\vspace{.3cm}

{\bf \noindent Acknowledgments.} The first author is supported by ARG-Matrics grant from ANRF. The second author 
is supported by project 1303/9/2025-R\&D-II-DAE/TIFR-17312. The third author is supported by NBHM Fellowship 0204/27/(23)/2023-R{\&}D-II/11914.

\section{Preliminaries}
In this chapter, we recall and prove some combinatorial results, some facts about the theta filtration and some facts about Hecke operators.

\subsection{Congruences for binomial coefficients and sums}
In this section, we study congruences between binomial coefficients and  sums of binomial coefficients. These congruences will be used in the later chapters in the Hecke operator computations.

Let $n \geq 0$. Recall that $\binom{n}{0} =1$ and $\binom{n}{m} =0$ for $m<0$.

The following result describes binomial coefficients mod $p$.
\begin{lemma}[Lucas' theorem]\label{Lucas}
	For any prime $p$, let $m$ and $n$ be two non-negative integers with base $p$ expansions given by $m = m_k p^k + m_{k-1} p^{k-1} + \cdots + m_0$ and $n = n_k p^k + n_{k-1}p^{k-1} + \cdots + n_0$ respectively. Then $\binom{m}{n} \equiv \binom{m_{k}}{n_{k}} \cdots \binom{m_{0}}{n_{0}} \mod p$. 
\end{lemma}

However, we will need deeper congruences for binomial coefficients. To this end, we first study congruences between values of polynomials at two integers that are $p$-adically  close. 

\begin{lemma}\label{polynomial under congruences}
    Let $p$ be an odd prime. Let $r,s$ be positive integers such that $r \equiv s ~\mathrm{mod}~ p^t$ for some $t \geq 1$. Let $\alpha_1,\ldots,\alpha_N \in \Z_p$ and $f(X) = (X-\alpha_1) \cdots (X-\alpha_N) $ be a polynomial. 
    \begin{enumerate}
    	\item[$(i)$] If $s \not\equiv \alpha_1,\ldots,\alpha_N ~\mathrm{mod}~p$, then $f(r)- f(s) \equiv (r-s) f(s) \left( \sum\limits_{i=1}^N \frac{1}{s-\alpha_i}\right)~\mathrm{mod}~p^{t+1}$ .
    	  \item[$(ii)$] If $s \equiv \alpha_n ~\mathrm{mod}~ p$ for some $ n $, then $f(r)- f(s) \equiv (r-s)  \prod\limits_{\substack{1\leq j \leq N \\ j \neq n}} (s-\alpha_j) \mod p^{t+1}$.
    \end{enumerate} 
 \end{lemma}
 \begin{proof}
     Let $g(X) = f(X)-f(s)$. As $g(s)=0$,  we get $g(X) = (X-s)h(X)$ for some $h(X) \in \Z_p[X]$. Further, note that $h(s) = g'(s)$, where $g'$ denotes the derivative of $g$. Thus
     \begin{align}\label{eq: cong. poly}
         f(r)-f(s) = g(r) = (r-s) h(r) \equiv (r-s) h(s) \equiv (r-s) g'(s) \text{ mod } p^{t+1}.
     \end{align}
     To prove the lemma, it is enough to determine $g'(s)$ modulo $p$. Clearly, we have 
     \begin{align*}
        g'(s) = \sum_{j=1}^N (s-\alpha_1)\cdots \widehat{(s-\alpha_j)}\cdots (s-\alpha_N),
     \end{align*}
     where $\widehat{\quad}$ denotes that the term is omitted from the product. Now $(i)$ follows immediately. If $s \equiv \alpha_n \mod p$ for some $n$, then $g'(s) \equiv (s-\alpha_1)\cdots \widehat{(s-\alpha_n)}\cdots (s-\alpha_N) \mod p$ and $(ii)$ follows.    
 \end{proof}
 
We now apply the previous result to study the congruences between the binomial coefficients of two integers congruent modulo $p^t$. 
    
\begin{lemma}\label{binomial coefficient under congruences}
    Let $p$ be an odd prime. Let $r,s$ be positive integers such that $r \equiv s ~\mathrm{mod}~ p^t$ for some $t \geq 1$. Then for  $0 < m < \min\{r,s\} $, we have
    \begin{enumerate}
    	  \item[$(i)$] $\binom{r}{m} \equiv \binom{s}{m} \mod p^{t-v(m!)}\mathbb{Z}_p$.
    	  \item[$(ii)$] If $p\nmid(s-n)$ for all $0 \leq n < m$, then $\binom{r}{m} - \binom{s}{m} \equiv (r-s) \binom{s}{m} (H_{s}-H_{s-m})\mod p^{t+1-v(m!)}\mathbb{Z}_p$, where $H_{n}$ denotes the $n^{th}$-harmonic number.
    	  \item[$(iii)$] If $p\mid(s-n)$ for some $0 \leq n < m$, then $\binom{r}{m} - \binom{s}{m} \equiv \frac{(r-s)}{s-n} \binom{s}{m} \mod p^{t+1-v(m!)}\mathbb{Z}_p$.
    \end{enumerate} 
\end{lemma}
\begin{proof}
    For every integer $m \geq 0$, let $f_m(X) := X(X-1)\cdots (X-m+1)$. Clearly $\binom{r}{m} = \frac{1}{m!} f_m(r)$ and $\binom{s}{m} = \frac{1}{m!} f_m(s)$ for every $m \geq 0$.
    \begin{enumerate}
        \item[$(i)$] Follows from the congruence $f_m(r) \equiv f_m(s) \mod p^t \mathbb{Z}_p$. Multiplying both sides by $1/m!$ we obtain the result. 
    	\item[$(ii)$]  Applying Lemma~\ref{polynomial under congruences} $(i)$  with $f$ there equal to $f_m$, we have
    		\begin{align*}
    			f_{m}(r) - f_m(s) \equiv  (r-s)f_{m} (s) \left( H_{s} - H_{s-m}\right) \mod p^{t+1}\Z_p.
    		\end{align*}
    		 Multiplying both sides by $1/m!$ we obtain the result. 
    	\item [$(iii)$] Applying Lemma~\ref{polynomial under congruences} $(ii)$  with $f$ there equal to $f_m$, we have
    		\begin{align*}
    			   f_{m}(r) - f_m(s) \equiv (r-s)\frac{f_m(s)}{(s-n)}  \mod p^{t+1} \mathbb{Z}_p.
    		\end{align*}
    	    Multiplying both sides of the congruence by $1/m!$ we obtain $(iii)$. \qedhere
    \end{enumerate}
\end{proof}
 We now recall an identity regarding the sums of powers of roots of unity, which we need to use frequently:
    \begin{align}\label{sums of roots of unity}
        \sum_{\xi \in \mu_{p-1}} \xi^{n} = 
        \begin{cases}
            p-1 &\text{ if } (p-1) \mid n,\\
            0 & \text{ otherwise. }
        \end{cases}
    \end{align}

In the next lemma, we show that if two positive integers are congruent modulo $p^t(p-1)$, then the corresponding binomial sums are congruent modulo $p^{t+1}$. 
\begin{lemma}\label{binomial sums under congruences}
    Let $p$ be an odd prime. Let $r,s$ be positive integers such that $r \equiv s ~\mathrm{mod}~ p^t(p-1)$ for some $t \geq 0$. For every positive integer $b$, we have 
    \begin{align*}
        \sum_{\substack{ 0 \leq j \leq r \\ j \equiv b ~\mathrm{mod}~(p-1)}} \binom{r}{j}  \equiv   \sum_{\substack{ 0 \leq j \leq s \\ j \equiv b ~\mathrm{mod}~(p-1)}} \binom{s}{j} \mod p^{t+1}.
    \end{align*}
\end{lemma}
\begin{proof}
    For every $\xi \in \mu_{p-1}$ and $\xi \neq -1$, we have $(1+\xi)^{(p-1)p^t} = (1+p z_{\xi})^{p^t} \equiv 1 \mod p^{t+1}$ for some $z_{\xi} \in \Z_p$. For $\xi = -1$, we have $ (1+\xi)^r = 0 = (1+\xi)^s $. Thus, we get $ (1+\xi)^r \equiv (1+\xi)^s \mod p^{t+1}$ for all $\xi \in \mu_{p-1}$. Hence, by \eqref{sums of roots of unity}, we have
    \begin{small}
        \begin{align*}
            \sum_{\substack{ 0 \leq j \leq r \\ j \equiv b ~\mathrm{mod}~(p-1)}} \binom{r}{j} = \frac{1}{(p-1)} \sum_{\xi \in \mu_{p-1}} \xi^{-b} (1+\xi)^{r}  &\equiv \frac{1}{(p-1)} \sum_{\xi \in \mu_{p-1}} \xi^{-b} (1+\xi)^{s} = \sum_{\substack{ 0 \leq j \leq s \\ j \equiv b ~\mathrm{mod}~(p-1)}} \binom{s}{j} \mod p^{t+1}.
        \end{align*}
    \end{small}%
    This completes the proof.
\end{proof}
 
 We now use the two  lemmas above to derive congruences for double binomial sums. These will be used repeatedly later.
 Following \cite{GR} for an integer $n$, let $$[n] \in \{1,2,\ldots,p-1\}$$ denote the congruence class of $n$ modulo $p-1$.
 
\begin{corollary}\label{cor: binomial sums under congruences 1}
    Let $p$ be an odd prime. Let $r,s$ be positive integers such that $r \equiv s ~\mathrm{mod}~ p^t(p-1)$ for some $t \geq 1$ and $r,s \equiv a \mod (p-1)$. For every integer $ 0 \leq m < \min\{r,s\} $, we have
    \begin{align*}
        \sum_{\substack{ 0 \leq j \leq r \\ j \equiv b ~\mathrm{mod}~(p-1)}} \binom{r}{j} \binom{j}{m} \equiv \left( \binom{r}{m} - \binom{s}{m}\right)\left(\binom{[a-m]}{[b-m]}+\delta_{[b-m],p-1} \right)+\sum_{\substack{ 0 \leq j \leq s \\ j \equiv b ~\mathrm{mod}~(p-1)}} \binom{s}{j} \binom{j}{m} \\
        ~\mathrm{mod}~ p^{t+1-v(m!)}.
    \end{align*} 
\end{corollary}
\begin{proof}
     Without loss of generality, assume that $s \leq r$.  Observe that  
    \begin{align*}
        \sum_{\substack{ 0 \leq j \leq r \\ j \equiv b ~\mathrm{mod}~(p-1)}} \binom{r}{j} \binom{j}{m} = \sum_{\substack{ m \leq j \leq r \\ j \equiv b ~\mathrm{mod}~(p-1)}} \binom{r}{j} \binom{j}{m} = \binom{r}{m} \sum_{\substack{ m \leq j \leq r \\ j \equiv b ~\mathrm{mod}~(p-1)}} \binom{r-m}{j-m}.
    \end{align*}
    By Lemma~\ref{binomial sums under congruences}, we have
    \begin{align*}
        \sum_{\substack{ m \leq j \leq r \\ j \equiv b ~\mathrm{mod}~(p-1)}} \binom{r-m}{j-m} \equiv \sum_{\substack{ m \leq j \leq s \\ j \equiv b ~\mathrm{mod}~(p-1)}} \binom{s-m}{j-m} \mod p^{t+1}.
    \end{align*}
    Thus, we obtain
    \begin{align*}
        \sum_{\substack{ 0 \leq j \leq r \\ j \equiv b ~\mathrm{mod}~(p-1)}} \binom{r}{j} \binom{j}{m}  \equiv \binom{r}{m} \sum_{\substack{ m \leq j \leq s \\ j \equiv b ~\mathrm{mod}~(p-1)}} \binom{s-m}{j-m} \mod p^{t+1}. 
    \end{align*}
    Adding and subtracting $\binom{s}{m}$ and then using the identity $\binom{s}{m} \binom{s-m}{j-m} = \binom{s}{j} \binom{j}{m}$ for $j \geq m$,  we have 
    \begin{align*}
        \binom{r}{m} \sum_{\substack{ m \leq j \leq r \\ j \equiv b ~\mathrm{mod}~(p-1)}} \binom{s-m}{j-m} =
        \left(\binom{r}{m}-\binom{s}{m}\right) \sum_{\substack{ m \leq j \leq s \\ j \equiv b ~\mathrm{mod}~(p-1)}} \binom{s-m}{j-m} + \sum_{\substack{ m \leq j \leq s \\ j \equiv b ~\mathrm{mod}~(p-1)}} \binom{s}{j} \binom{j}{m}.
    \end{align*}
    Now the corollary follows from Lemma~\ref{binomial coefficient under congruences} $(i)$ and \cite[Lemma 2.14]{GR19} (applied with $r$ there equal to $s-m$, $b$ there equal to $[b-m]$ and $m$ there equal to $0$).
\end{proof}
\begin{corollary}\label{cor: binomial sums under congruences 2}
    Let $p$ be an odd prime. Let $r,s$ be positive integers such that $r \equiv s ~\mathrm{mod}~ p^t(p-1)$ for some $t \geq 1$. For every integer $ 0 \leq m < \min\{r,s\} $, we have 
    \begin{align*}
        \sum_{\substack{ 0 \leq j \leq r \\ j \equiv b ~\mathrm{mod}~(p-1)}} \binom{r}{j} \binom{j}{m} \equiv  \sum_{\substack{ 0 \leq j \leq s \\ j \equiv b ~\mathrm{mod}~(p-1)}} \binom{s}{j} \binom{j}{m}  \mod p^{t-v(m!)}.
    \end{align*}
\end{corollary}
\begin{proof}
    By Lemma~\ref{binomial coefficient under congruences} $(i)$, we have $\binom{r}{m} \equiv \binom{s}{m} \mod p^{t-v(m!)}$. The result follows from Corollary~\ref{cor: binomial sums under congruences 1}. 
\end{proof}
Let the notation be as in Lemma~\ref{binomial sums under congruences}. We now consider the special case $j \equiv a \mod (p-1)$ and compute the double binomial sums modulo $p^{t+2}$. 
\begin{lemma}\label{lem: binomial sum Doc. Math general}
    Let $p$ be a prime and $r,s$ be positive integers with $r \equiv  s ~\mathrm{mod}~ p^t(p-1)$ for some $t\geq 1$. Let $0 \leq m \leq p-1$ be an integer. If $r \equiv a ~\mathrm{mod}~ (p-1)$ with $m+1\leq a \leq p-1+m$ and $m<\min\{r,s\}$, then 
        \begin{equation}
            \sum_{\substack{m<j<r \\ j \equiv a ~\mathrm{mod}~ (p-1)}} \binom{r}{j} \binom{j}{m} \equiv p\left\{ \binom{r}{m} - \binom{s}{m}\right\}  \frac{a-s}{a-m} + \sum_{\substack{m<j<s \\ j \equiv a ~\mathrm{mod}~ (p-1)}} \binom{s}{j} \binom{j}{m} + p \frac{s-r}{a-m} \binom{s}{m}~\mathrm{ mod }~p^{t+2}.
        \end{equation}
\end{lemma}
\begin{proof}
    We first prove the lemma in the special case $m=0$. Note that 
    \begin{equation*}
        \begin{aligned}
            (p-1)\sum_{\substack{0 \leq j \leq r \\ j \equiv a ~\mathrm{mod}~ (p-1)}}  \binom{r}{j} &= \sum_{\xi \in \mu_{p-1}} \xi^{-a}(1+\xi)^r \\
            &= \sum_{\xi \in \mu_{p-1}} \xi^{-a}(1+\xi)^s (1+\xi)^{r-s} \\
            & =  \sum_{\xi \in \mu_{p-1} \smallsetminus \{-1\}} \xi^{-a}(1+\xi)^s (1+pz_{\xi})^{\frac{r-s}{p-1}} \\
            & \equiv \sum_{\xi \in \mu_{p-1} \smallsetminus \{-1\}} \xi^{-a}(1+\xi)^s + p \frac{r-s}{p-1} \sum_{\xi \in \mu_{p-1} \smallsetminus \{-1\}} \xi^{-a}(1+\xi)^s z_{\xi} ~\mathrm{mod}~ p^{t+2} \\
            & \equiv(p-1) \sum_{\substack{0 \leq j \leq s \\ j \equiv a ~\mathrm{mod}~ (p-1)}}   \binom{s}{j} + p \frac{r-s}{p-1} \sum_{\xi \in \mu_{p-1} \smallsetminus \{-1\}} \xi^{-a}(1+\xi)^s z_{\xi} ~\mathrm{mod}~ p^{t+2}.
        \end{aligned}
    \end{equation*}
    Since $r\equiv s ~\mathrm{mod}~ p^t$ and $(1+\xi)^s \equiv (1+\xi)^a ~\mathrm{mod}~ p$, we get 
    \begin{equation*}
        \begin{aligned}
            \sum_{\substack{0 \leq j \leq r \\ j \equiv a ~\mathrm{mod}~ (p-1)}}  \binom{r}{j} \equiv \sum_{\substack{0 \leq j \leq s \\ j \equiv a ~\mathrm{mod}~ (p-1)}}  \binom{s}{j} + p \frac{r-s}{(p-1)^2} \sum_{\xi \in \mu_{p-1} \smallsetminus \{-1\}} \xi^{-a}(1+\xi)^a z_{\xi} ~\mathrm{mod}~ p^{t+2}.
        \end{aligned}
    \end{equation*}  
    Hence
    \begin{equation}\label{eq: binomial sum BG15 Lem 2.5 gen}
        \sum_{\substack{0 < j < r \\ j \equiv a ~\mathrm{mod}~ (p-1)}}  \binom{r}{j} \equiv \sum_{\substack{0 < j < s \\ j \equiv a ~\mathrm{mod}~(p-1)}}  \binom{s}{j} + p \frac{r-s}{(p-1)^2} \sum_{\xi \in \mu_{p-1} \smallsetminus \{-1\}} \xi^{-a}(1+\xi)^a z_{\xi} ~\mathrm{mod}~ p^{t+2}.
    \end{equation}  
    Taking $s=a$ and $r=a+p-1$ in the above identity we get
    \begin{equation*}
        \begin{aligned}
            \frac{p}{(p-1)}  \sum_{\xi \in \mu_{p-1} \smallsetminus \{-1\}} \xi^{-a}(1+\xi)^a z_{\xi} \equiv \binom{a+p-1}{a} \equiv \frac{p}{a} ~\mathrm{mod}~ p^2.
        \end{aligned}
    \end{equation*}
    Using this in \eqref{eq: binomial sum BG15 Lem 2.5 gen} we obtain the lemma for $m=0$.

    Now let $m$ be a non-negative integer. From the special case $m=0$, we have 
    \begin{align*}
        \sum_{\substack{m<j<r \\ j \equiv a ~\mathrm{mod}~ (p-1)}} \binom{r}{j} \binom{j}{m} &= \binom{r}{m} \sum_{\substack{m<j<r \\ j \equiv a ~\mathrm{mod}~ (p-1)}} \binom{r-m}{j-m} \\
        & \equiv \binom{r}{m} \sum_{\substack{m<j<s \\ j \equiv a ~\mathrm{mod}~ (p-1)}} \binom{s-m}{j-m} + p \frac{s-r}{a-m} \binom{r}{m} ~\mathrm{mod}~ p^{t+2}\\
        &\equiv \left\{ \binom{r}{m} - \binom{s}{m}\right\}  \sum_{\substack{m<j<s \\ j \equiv a ~\mathrm{mod}~ (p-1)}} \binom{s-m}{j-m} + \sum_{\substack{m<j<s \\ j \equiv a ~\mathrm{mod}~ (p-1)}} \binom{s}{j} \binom{j}{m} \\
        & \qquad \qquad \qquad + p \frac{s-r}{a-m} \binom{s}{m}~\mathrm{ mod }~p^{t+2},
    \end{align*}
    where we used Lemma~\ref{binomial coefficient under congruences} $(i)$ in the last step. By \cite[Lemma 2.5]{BG15}, we have
    \[
        \sum_{\substack{m<j<s \\ j \equiv a ~\mathrm{mod}~ (p-1)}} \binom{s-m}{j-m} \equiv p \frac{a-s}{a-m} ~\mathrm{mod}~ p^2.
    \]
    Substituting this above we obtain the lemma for arbitrary $m \geq 0$.
\end{proof}

In the following lemma, we prove a collection of identities involving sums of binomial coefficients, which will be useful for carrying out row operations arising in the computation of certain determinants.
    \begin{lemma}\label{lem: combinatorial id for det(B)}
        Let $N\geq0$ be an  integer. Then
        \begin{enumerate}
            \item[$(i)$] for every non-negative integer $k' \leq N$, we have
                \begin{equation*}
	                   \sum_{l'=0}^{N} (-1)^{l'}\binom{N}{l'}\binom{N-l'}{k'} =
	                \begin{cases}
		            1, & \text{if } k'=N,\\
		              0, & \text{otherwise}.
	                \end{cases}
                \end{equation*}
            \item[$(ii)$] for every real number $\alpha \neq 0,1,\ldots, N$, we have 
                \begin{equation*}
	                    \sum_{l'=0}^{N} (-1)^{l'} \binom{N}{l'} \frac{1}{\alpha-l'} = (-1)^N \frac{N!}{\alpha(\alpha-1)\cdots(\alpha-N)}.
                 \end{equation*}
            \item[$(iii)$]  for every real number $\alpha \neq 0,-1,\ldots, -N$, we have 
                \begin{equation}\label{col identity 3}
	                   \frac{1}{\alpha} - \sum_{l'=1}^{N}  \frac{(l'-1)!}{(\alpha+1)\cdots(\alpha+l')} =  \frac{N!}{\alpha(\alpha+1)\cdots(\alpha+N)}.
                \end{equation}
            \item[$(iv)$] for every  integer $k' $, we have
                \begin{equation}
                	\sum_{l'=0}^{N} (-1)^{l'}\binom{N}{l'}\binom{N+1-l'}{k'} =
                	\begin{cases}
                		1, & \text{if } k'=N,N+1,\\
                		0, & \text{otherwise}.
                	\end{cases}
               \end{equation}	
            
            \item[$(v)$] for  integers $k $ and $M \geq N$, we have
                \begin{equation}\label{col identity 5}
                	\sum_{l'=0}^{N} (-1)^{l'}\binom{N}{l'}\binom{M-l'}{k-l'} = \binom{M-N}{k}.
               \end{equation}	
            \end{enumerate}
        \end{lemma}
        \begin{proof} 
            \begin{enumerate}
                \item[$(i)$] Note that 
	            \begin{align*}
		                 \sum_{l'=0}^{N} (-1)^{l'}\binom{N}{l'}\binom{N-l'}{k'} 
		               &= \sum_{l'=0}^{N} (-1)^{l'}\binom{N}{l'}\binom{N-l'}{N-l'-k'} \\
		               &= \text{ coefficient of } x^{N-k'}  \text{ in } \sum_{l'=0}^{N} \binom{N}{l'} (-x)^{l'}(1+x)^{N-l'} \\
		               &= \text{ coefficient of } x^{N-k'}  \text{ in } 1.
	            \end{align*}
	              Now the identity follows.
                \item[$(ii)$] 	We claim that in the field of fractions of $\mathbb{Q}[X]$, we have 
                \begin{align*}
		              \sum_{l'=0}^{N} (-1)^{l'} \binom{N}{l'} \frac{1}{X-l'} = (-1)^N \frac{N!}{X(X-1)\cdots(X-N)}.
	              \end{align*}
	            Clearly the identity in $(ii)$ follows from the claim by taking $X=\alpha$. Multiplying both sides by $(-1)^N X\cdots(X-N)$, it is enough to show
	            \begin{align*}
		            F(X):=\sum_{l'=0}^{N} (-1)^{N-l'} \binom{N}{l'} \frac{{X\cdots(X-N)}}{X-l'} =  N!.
	              \end{align*}
	            Since $F(X)$ is a polynomial of degree $N$, it further suffices to show $F(X)-N!$ has more than $N$ distinct roots. Note that for every $0 \leq l' \leq N$, we have
                 \[
                    F(l') = (-1)^{N-l'} \binom{N}{l'}\times (X \cdots \widehat{(X-l')}\cdots (X-N))\big\vert_{X=l'} = \binom{N}{l'} l'! (N-l')! = N!,
                 \]	
                where $\widehat{(X-l')}$ means the term is omitted from the product. Now the proof follows.
                \item[$(iii)$] We proceed by induction on $N$. If $N=0$, then the sum is empty and we are done. Assume that \eqref{col identity 3} is true for $N=k'$. Then for $N=k'+1$, we have
	            \begin{align*}
		               \frac{1}{\alpha} - \sum_{l'=1}^{k'+1}  \frac{(l'-1)!}{(\alpha+1)\cdots(\alpha+l')} &=  \frac{1}{\alpha} - \sum_{l'=1}^{k'}  \frac{(l'-1)!}{(\alpha+1)\cdots(\alpha+l')}  - \frac{k'!}{(\alpha+1)\cdots(\alpha+k'+1)} \\
		                &= \frac{k'!}{\alpha(\alpha+1)\cdots(\alpha+k')} - \frac{k'!}{(\alpha+1)\cdots(\alpha+k'+1)} \\
		                & = \frac{(k'+1)!}{\alpha(\alpha+1)\cdots(\alpha+k'+1)}.
	              \end{align*}
	            This completes the proof by induction. 
                \item[$(iv)$] The proof is similar to $(i)$. Note that 
	            \begin{align*}
	        	    \sum_{l'=0}^{N} (-1)^{l'}\binom{N}{l'}\binom{N+1-l'}{k'} 
	        	    &= \sum_{l'=0}^{N} (-1)^{l'}\binom{N}{l'}\binom{N+1-l'}{N+1-l'-k'} \\
	        	    &= \text{ coefficient of } x^{N+1-k'}  \text{ in } \sum_{l'=0}^{N} \binom{N}{l'} (-x)^{l'}(1+x)^{N+1-l'} \\
	        	    &= \text{ coefficient of } x^{N+1-k'}  \text{ in } (1+x).
	            \end{align*}
	              Now the identity follows.
                  \item[$(v)$] The proof is similar to $(i)$ and $(iv)$. Taking $M=N,N+1$ we obtain $(i)$ and $(iv)$ respectively. If $k <0 $, then we are done. So assume $k\geq 0$. Note that 
	            \begin{align*}
	        	    \sum_{l'=0}^{N} (-1)^{l'}\binom{N}{l'}\binom{M-l'}{k-l'} 
	        	    &= \text{ coefficient of } x^{k}  \text{ in } \sum_{l'=0}^{N} \binom{N}{l'} (-x)^{l'}(1+x)^{M-l'} \\
                    &= \text{ coefficient of } x^{k}  \text{ in } (1+x)^{M-N}\sum_{l'=0}^{N} \binom{N}{l'} (-x)^{l'}(1+x)^{N-l'} \\
	        	    &= \text{ coefficient of } x^{k}  \text{ in } (1+x)^{M-N}.
	            \end{align*}
	              Now the identity follows.\qedhere
            \end{enumerate}
        \end{proof}
\subsection{Combinatorial lemmas for determinants}
In this article, we  need to solve certain linear equations or linear  congruences. We solve these  using Cramer's rule. This involves showing a matrix is invertible, which in turn involves computing its determinant. In some cases, we need to determine the $p$-adic valuation of the determinant.

We begin by recalling a result of Gessel-Viennot computing the determinant of a matrix whose entries are given by  binomial coefficients.
\begin{lemma}\cite[p. 308 \& l. -7]{Viennot}\label{GV det}
    Let $k$ be a positive integer and  $ a_1 < a_2 < \cdots < a_k $ be a sequence of positive integers. For every non-negative integer $n$, we have
    \begin{align}
        \det\limits_{1 \leq i,j \leq k}\left( \binom{a_i}{n+(j-1)}\right) = \frac{(a_1)_n\cdots(a_k)_n}{n!\cdots(n+k-1)!} \times \prod_{1\leq i< j \leq k} (a_j-a_i),
    \end{align}
    where $(a)_n = a (a-1)\cdots(a-n+1)$ is the Pochhammer symbol.
\end{lemma} 
We now derive some consequences of the above lemma that will be needed later.
\begin{corollary}\label{cor: GV det}
    Let $p$ be a prime. Let $d,k$ be positive integers and $m,n$ be non-negative integers. Then 
    \begin{equation*}
        \det\limits_{1 \leq i,j \leq k}\left( \binom{m+d(i-1)}{n+(j-1)}\right)=  \frac{(m)_n (m+d)_n \cdots(m+d(k-1))_n}{n! (n+1)!\cdots(n+k-1)!} \times d^{k(k-1)/2}\times 1!\cdots(k-1)!.
    \end{equation*}
    As a consequence, the above matrix is invertible in $M_{k}(\mathbb{Z}_p)$ if one of the following conditions holds
    \begin{enumerate}
        \item[$(i)$] $n=0$ and $p\nmid d$
        \item[$(ii)$]  $d=1$ and $n+k-1\leq m+k-1 \leq p-1$.
    \end{enumerate}
\end{corollary}
\begin{proof}
    The first assertion follows from the observation $\prod_{1\leq i< j \leq k} (a_j-a_i) = d^{k(k-1)/2}\times 1!\cdots(k-1)!$ and Lemma~\ref{GV det}. 
    \begin{enumerate}
    	\item[$(i)$] If $n=0$, then
    	\begin{align*}
    		\frac{(m)_n\cdots(m+d(k-1))_n}{n!\cdots(n+k-1)!} \times d^{k(k-1)/2}\times 1!\cdots(k-1)!= 
    		d^{k(k-1)/2}.              
    	\end{align*}
    	\item[$(ii)$]  If $d=1$, then 
    	\begin{align*}
    		\frac{(m)_n\cdots(m+d     (k-1))_n}{n!\cdots(n+k-1)!} \times d^{k(k-1)/2}\times 1!\cdots(k-1)!= 
    		\frac{(m)_n\cdots(m+(k-1))_n}{n!\cdots(n+k-1)!}  1!\cdots(k-1)! .              
    	\end{align*}
    	If $n+k \leq p$ and $n \geq 0$, then $p \nmid n!\cdots(n+k-1)!$. As $1 \leq k \leq n+k \leq p$, we get $p\nmid 1!\cdots(k-1)!$. Note that $(m+(i-1))_n = n! \binom{m+i-1}{n}$ for every $i \geq 1$. From $n+k\leq m+k \leq p$, it follows that $0 \leq n \leq m \leq m+k-1  \leq p-1$. Thus by Lucas' theorem we have $\binom{m+i-1}{n} \not\equiv 0 \mod p$. Combining this with $0 \leq n\leq p-k \leq p-1$, we get $p \nmid n! \binom{m+i-1}{n}$ for $1 \leq i \leq k$. \qedhere
    \end{enumerate}
\end{proof}

	In our later computations, given a set of integers, we need to choose another set of integers satisfying certain linear congruences. The following lemma gives a sufficient condition for when this can be done.
    
\begin{lemma}\label{lem:choice of beta}
	Let $p$ be an odd prime. Let $m'$, $n$ and $t$ be integers such that $ 0 \leq t \leq n $ and $ 0 \leq m' < p $. Let $c$, $k \geq 0$ and let $ \gamma_{c},\gamma_{c+p-1}, \ldots, \gamma_{c+k(p-1)} $  be $p$-adic integers with		
	      $$ \sum\limits_{j=0}^{k} \binom{c+j(p-1)}{m} \gamma_{c+j(p-1)} \equiv \nu_m \mod  p^{t} ~ \text{for all} ~ m=0,1, \ldots, m' $$
    for some $\nu_1,\ldots, \nu_m \in \Z_p $. If $ k \geq m' $, then there exists $  \alpha_{c}, \ldots, \alpha_{c+k(p-1)} \in \mathbb{Z}_{p} $ such that
	\begin{enumerate}
		\item[$(i)$] $ \alpha_{c+j(p-1)} \equiv \gamma_{c+j(p-1)} \mod  p^{t} $  
		\item[$(ii)$] $ \sum\limits_{j=0}^{k} \binom{c+(p-1)j}{m} \alpha_{c+j(p-1)} \equiv \nu_m \mod p^{n} $ for all $ m =0,1, \ldots, m' $.
	\end{enumerate}
\end{lemma}
\begin{proof}
	If $ n =t $, then we may take $ \alpha_{c+j(p-1)} = \gamma_{c+j(p-1)} $. So assume $ n>t $. Put $ S_{m} = \sum_{j=0}^{k} \binom{c+j(p-1)}{m} \gamma_{c+j(p-1)} $  for $ m \geq 0 $. Also set
	\begin{align*}
		\alpha_{c+j(p-1)} = \gamma_{c+j(p-1)} + p^{t} \delta_{j}
	\end{align*}
    for $j=0,\ldots,k$. To prove the lemma, it is enough to show that that the system of congruences
	\begin{align}\label{eq:choice of beta and indices reduction1}
		\sum_{j=0}^{k} \binom{c+j(p-1)}{m} \delta_{j} \equiv p^{-t} (\nu_m-S_{m}) \mod p^{n-t} \text{ for all } m=0,1,\ldots,m'.
	\end{align}
   has a solution. Set $ \delta_j=0$ for $ j \neq 0,1, \ldots, m' $.  Thus to solve \eqref{eq:choice of beta and indices reduction1}, it is enough to show that the system of congruences 
    \begin{align}\label{eq:choice of beta and indices reduction2}
		\sum_{j=0}^{m'} \binom{c+j(p-1)}{m} \delta_{j} \equiv - p^{-t} (S_{m}-\nu_m) \mod p^{n-t} \text{ for all } m=0,1,\ldots,m'.
	\end{align}
    has a solution. By Corollary~\ref{cor: GV det} $(i)$, we see that
	\begin{align*} 
		\det_{0 \leq j,m \leq m'} \Big( \binom{c+j(p-1)}{m} \Big) 
	\end{align*}
	is invertible modulo $p$. Hence we can choose $ \delta_{j} $ for $ k= 0,1, \ldots, m' $ such that \eqref{eq:choice of beta and indices reduction2} is solvable. This completes the proof of the lemma.
\end{proof}
In our applications, the $\gamma_j$ above will be linear combinations of binomial coefficients. The lemma will help us smoothen certain non-integral terms in our computations to a $p$-integral expression.  The lemma is mostly used with $\nu_m = 0$ and once (cf. Theorem~\ref{thm: above super diagonal a large}) with $\nu_m$ possibly non-zero.

     The following lemma shows that certain congruences involving linear combinations of binomial coefficients imply additional congruences.  
\begin{lemma}\label{dmm' trick}
	Let $p$ be a prime. Let $c,N$ and $t$ be non-negative integers. Let $\alpha_{0}, \ldots, \alpha_{N}$ be $p$-adic integers. Suppose 
	\begin{align*}
		\sum_{k=0}^{N} \alpha_k \binom{c+k(p-1)}{m} \equiv 0 \mod p^t
	\end{align*}
	for $m=0,\ldots, N$. Then $\alpha_{k} \equiv 0 \mod p^t$ for all $k=0,\ldots, N$. In particular, we have
	\begin{align*}
		\sum_{k=0}^{N} \alpha_k \binom{c+k(p-1)}{m} \equiv 0 \mod p^t \quad \text{for all } m.
	\end{align*}
\end{lemma}
\begin{proof}
	The given system of congruences can be expressed as
	\begin{align*}
		\begin{bmatrix} \binom{c}{0} & \cdots &  \binom{c+N(p-1)}{0} \\
			\vdots & \ddots & \vdots \\
			\binom{c}{N} & \cdots &  \binom{c+N(p-1)}{N}
		\end{bmatrix} \begin{bmatrix} \alpha_{0} \\ \vdots \\\alpha_{N} \\ \end{bmatrix} \equiv \begin{bmatrix} 0 \\ \vdots \\0 \\ \end{bmatrix}. 
	\end{align*}
	By Corollary~\ref{cor: GV det} $(i)$, the above matrix is invertible modulo $p$. Hence by Cramer's rule, we get $\alpha_k \equiv 0 \mod p^t$ for all $k$. This proves the lemma.
\end{proof}
In our computations, we often come across special matrices whose entries can be expressed as linear combinations of fixed binomial coefficients, up to an error term divisible by a power of $p$. The next two lemmas determine the exponent of $p$ dividing the determinant of these special matrices and some of its minors. 

\begin{lemma}\label{minor trick 1}
	Let $p$ be a prime. Let $c,M,N,s$ and $t$ be non-negative integers with $N<M<s$. Let
	\begin{align*}
		A=\left(\sum_{k=0}^{N} \binom{s-l}{c+k(p-1)} \binom{c+k(p-1)}{m} + O(p^t) \right)_{\substack{m=0,\ldots,M \\ l=0,\ldots,M}}.
	\end{align*}
	Then $p^{t(M-N)} \mid \det(A)$. Also, $p^{t(M-N-1)}$ divides the determinant of the minor of every entry in the last row of $A$, that is, $p^{t(M-N-1)} \mid \det(A_{M,l}) $ for all $l$.
\end{lemma}
\begin{proof}
	We first prove $p^{t(M-N)} \mid \det(A)$.  We show that the rows $m=N+1,\ldots,M$ can be expressed as a linear combination of   rows  $m =0, \ldots, N$ modulo $p^{t}$. Consider the following matrix (which appeared in the proof of Lemma~\ref{dmm' trick})
	\begin{align*}
		\begin{bmatrix} \binom{c}{0} & \cdots &  \binom{c+N(p-1)}{0} \\
			\vdots & \ddots & \vdots \\
			\binom{c}{N} & \cdots &  \binom{c+N(p-1)}{N}
		\end{bmatrix} . 
	\end{align*}
	By Corollary~\ref{cor: GV det} $(i)$,  the  matrix above is invertible modulo $p$. Thus for every $m=N+1,\ldots,M$  there exists $d_{m,m'} \in \Z_p$ such that
	\[
	\binom{c+k(p-1)}{m} = \sum_{m'=0}^{N} d_{m,m'} \binom{c+k(p-1)}{m'}.
	\]
    for all $k=0,\ldots,N$. Multiplying both sides by $\binom{s-l}{c+k(p-1)}$ and then taking the sum over $k$, we obtain 
	\[
	\sum_{k=0}^{N} \binom{s-l}{c+k(p-1)} \binom{c+k(p-1)}{m} = \sum_{m'=0}^{N} d_{m,m'} \sum_{k=0}^{N} \binom{s-l}{c+k(p-1)}\binom{c+k(p-1)}{m'}  
	\]
	for all  $m=N+1,\ldots,M$ and $l=0,\ldots,M$. This shows that the rows $m=N+1,\ldots,M$ of $A$ can be expressed as a linear combination of the first $N+1$ rows in $A$ modulo $p^{t}$. Thus after applying appropriate row operations we see that rows $m=N+1,\ldots,M$ of $A$ are multiples of $p^t$. Now the first statement follows.
	
	The second statement, namely that $p^{t(M-N-1)} \mid \det(A_{M,l}) $  for $l =0, \ldots, M$, follows by ignoring the last row and the $l$-th column at the end of the argument above.
\end{proof}
\begin{lemma}\label{minor trick 2}
	Let $p$ be a prime. Let $c,M,N,s$ and $t$ be non-negative integers with $1\leq N<M<s$. Let
	\begin{align*}
		A=\left(\begin{array}{c|c}
			\left(\sum\limits_{k=1}^{N} \binom{s-l}{c+k(p-1)} \binom{c+k(p-1)}{m} + O(p^t) \right)_{\substack{m=0,\ldots,M \\l=0,\ldots,M-1}} & \left((1-\delta_{M,m})\binom{c}{m}\right)_{m=0,\ldots,M}
		\end{array} \right)
	\end{align*}
	be a block matrix of size $(M+1)\times (M+1)$. Then  $p^{t(M-N-1)} \mid \det(A_{M,l}) $ for all $l \neq M$ and $p^{t(M-N)} \mid \det(A_{M,M})$.
\end{lemma}
\begin{proof} 
     As shown in Lemma~\ref{minor trick 1}, for every $m=N+1,\ldots,M$  there exists $d_{m,m'} \in \Z_p$ such that
	\[
	\binom{c+k(p-1)}{m} = \sum_{m'=0}^{N} d_{m,m'} \binom{c+k(p-1)}{m'}.
	\]
    for all $k=0,\ldots,N$. Using this relation, it can be checked that the rows $m=N+1,\ldots,M-1$ of $A_{M,l}$ can be expressed as a linear combination of the first $N+1$ rows in $A_{M,l}$ modulo $p^{t}$ for $l\neq M$. Now the first statement follows by the same argument as in proof of Lemma~\ref{minor trick 1}.  

    For the last statement, note that 
	\begin{align*}
		A_{M,M}&=  \left( \sum\limits_{k=1}^{N} \binom{s-l}{c+k(p-1)} \binom{c+k(p-1)}{m} + O(p^t) \right)_{\substack{m=0,\ldots,M-1 \\l=0,\ldots,M-1}} \\
		&=  \left( \sum\limits_{k=0}^{N-1} \binom{s-l}{c+p-1+k(p-1)} \binom{c+p-1+k(p-1)}{m} + O(p^t) \right)_{\substack{m=0,\ldots,M-1 \\l=0,\ldots,M-1}}.    
	\end{align*}
	By Lemma~\ref{minor trick 1}, we get $p^{t(M-N)}\mid\det(A_{M,M})$. 
\end{proof}
\subsection{The \texorpdfstring{$\theta$}{}-filtration of \texorpdfstring{$V_{r}$}{}}

 Let $V_r= \mathrm{Sym}^r\bar{\mathbb{F}}_p^2$ be the $r$-th symmetric power representation of $\gl_2(\f)$. Note that $V_r$ has a model consisting of homogeneous polynomials of degree $r$ in the variables $X,Y$ with coefficients in $\bar{\mathbb{F}}_p$. Let  $\theta = X^{p}Y- XY^{p}$ be the Dickson polynomial and $V_{r}^{(m)} = \lbrace F(X,Y) \in \bar{\mathbb{F}}_p[X,Y] : \theta^{m} \mid F(X,Y) ~ \text{in} ~ \bar{\mathbb{F}}_p[X,Y]  \rbrace$. We have the following filtration of $V_{r}$
\begin{equation}\label{theta filtration}
    V_r \supset V_r^{(1)} \supset V_r^{(2)} \supset \cdots \supset V_r^{(m)} \supset V_r^{(m+1)} \supset \cdots.
\end{equation}
We refer to the above filtration as the $\theta$-filtration of $V_r$.

In this section, we will study the  quotients $V_{r}^{(m)}/V_{r}^{(m+1)}$. The following lemma gives the structure of $V_{r}^{(m)}/V_{r}^{(m+1)}$ and describes the projection of certain polynomials onto the cosocle of $V_{r}^{(m)}/V_{r}^{(m+1)}$. 
\begin{lemma}\label{Glover-Brueil map image}
    Let $p \geq 2$, $m \geq 0$, $r \geq m(p+1)+p$ and $r \equiv b ~\mathrm{mod}~(p-1)$ with $2m+1 \leq b \leq p-1+2m$. Then we have a short exact sequence of $\Gamma$-modules
    \begin{equation}\label{exact seq. singular quotients}
    	0 \rightarrow V_{b-2m} \otimes D^m \rightarrow V_{r}^{(m)}/V_{r}^{(m+1)} \rightarrow V_{p-1-b+2m} \otimes D^{b-m} \rightarrow 0.
    \end{equation}
    and this sequence splits if and only if $ b = p-1+2m$. Furthermore, we have 
    \begin{enumerate}
    	\item[$(i)$] The monomials $X^{b-2m}, Y^{b-2m} \in V_{b-2m} \otimes D^m$ map to $\theta^m X^{r-m(p+1)}, \theta^m Y^{r-m(p+1)}$ in $V_{r}^{(m)}/V_{r}^{(m+1)}$ respectively.
    	\item[$(ii)$] The polynomials $\theta^m X^{r-m(p+1)-b+2m}Y^{b-2m}, \theta^m X^{r-m(p+1)-p+1}Y^{p-1}\in V_{r}^{(m)}/V_{r}^{(m+1)}$ map to $X^{p-1-b+2m}, (-1)^b Y^{p-1-b+2m} \in V_{p-1-b+2m} \otimes D^{b-m}$ respectively.
    \end{enumerate}
\end{lemma}
\begin{proof}
    The exact sequence and $(ii)$ follow from \cite[Lemma 2.11]{GR19}.  The proof of assertion~$(i)$ is similar to \cite[Lemma 8.5]{BG15}.
\end{proof}
The following lemma describes a $\gl_2(\f)$-generator of $V_{r}^{(m)}/V_r^{(m+1)}$. The generator described in this lemma is typical of those we will encounter in later sections.
\begin{lemma}\label{generating polynomial quotient}
    Let $p \geq 2$, $m \geq 0$ and $r \geq m(p+1)+p$. Then $\theta^m (Y^{r-m(p+1)} - X^{p-1}Y^{r-m(p+1)-(p-1)})$ generates $V_{r}^{(m)}/V_r^{(m+1)}$. 
\end{lemma}
\begin{proof}
     It is enough that the Weyl involution of the given polynomial, namely  $$F(X,Y):=\theta^m (X^{r-m(p+1)} - X^{r-m(p+1)-(p-1)}Y^{p-1}),$$  generates $V_{r}^{(m)}/V_r^{(m+1)}$. Let $r \equiv b \mod (p-1)$ with $2m+1 \leq b \leq p-1+2m$. If $b \neq p-1+2m$, then by Lemma~\ref{Glover-Brueil map image}  we obtain that $F(X,Y)$ generates $V_{r}^{(m)}/V_r^{(m+1)}$. 
    	  
    If $b = p-1+2m$, then $V_r^{(m)}/V_r^{(m+1)} \cong V_{2p-2}/V_{2p-2}^{(1)} \otimes D^{m} \cong (V_{0} \oplus V_{p-1}) \otimes D^m$ and under this isomorphism $F(X,Y)$ maps to $X^{2p-2}-X^{p-1}Y^{p-1}$. By \cite[Lemma 3.12]{AC23}, we have that $X^{2p-2}-X^{p-1}Y^{p-1} +Y^{2p-2}$ generates $ V_{0} $ and  $Y^{2p-2}$  generates $V_{p-1} $. Thus, $X^{2p-2}-X^{p-1}Y^{p-1}  = (X^{2p-2}-X^{p-1}Y^{p-1} +Y^{2p-2}) - Y^{2p-2}$ generates $V_{2p-2}/V_{2p-2}^{(1)}$. So $F(X,Y)$ generates  $V_{r}^{(m)}/V_r^{(m+1)}$ even in this case. This completes the proof of the lemma.
\end{proof}
\subsection{Hecke Operators}
In this subsection, we recall the action of the Hecke operator $T= T^++T^-$ and derive a congruence for the action of $T^+$.

 Let $R$ be a $\Z_p$-algebra and $ v = \sum\limits_{i=0}^{r} c_{i} X^{r-i} Y^{i} \in  R[X,Y] $ be a homogeneous polynomial of degree $r$.  Following the notation in \cite{Breuil1}, for $\mu \in I_n$, we have 
\begin{equation}\label{T plus formula}
	T^{+}\left( [g_{n,\mu}^{0}, v]\right)  
    = \sum_{\lambda \in I_{1}} \left[ g_{n+1,\mu+p^{n+1} \lambda }^{0}, \sum_{j=0}^{r} \left( p^{j} \sum_{i=j}^{r} c_{i} \binom{i}{j}  (- \lambda)^{i-j} \right) X^{r-j}Y^{j}  \right], 
\end{equation}
\begin{equation}\label{T minus formula}\resizebox{0.92\hsize}{!}{$
    \begin{aligned}
	      T^{-}\left( [g_{n,\mu}^{0}, v]\right) 
        & = \left[ g_{n-1,[\mu]_{n-1}}^{0}, \sum_{j=0}^{r} \left( \sum_{i=j}^{r} p^{r-i} c_{i} \binom{i}{j} \left( \frac{\mu - [\mu]_{n-1}}{p^{n-1}}\right)^{i-j} \right) X^{r-j}Y^{j} \right] 
		~(n > 0), \\
		\qquad \qquad \qquad T^{-}\left( [g_{n,\mu}^{0}, v]\right) & = \left[ g_{n,\mu}^{0} \begin{pmatrix} 1 & 0 \\ 0 & p \end{pmatrix}, \sum_{j=0}^{r}  p^{r-j} c_{j} X^{r-j} Y^{j}\right]
			\quad (n = 0). 
    \end{aligned}    $}
\end{equation}	    

We show that if $R = {\mathbb Z}_p$, then the $ T^{+} $ operator vanishes modulo a large power of $ p $ when a certain divisibility condition holds.  
\begin{lemma}\label{sum and derivative vansihing}
	Let $ F(X,Y) = \sum\limits_{j=0}^{r} a_{j} X^{r-j} Y^{j} \in \mathbb{Z}_{p}[X,Y]$ be a homogeneous polynomial of degree $ r $ in $ X, Y$ with coefficients in $\mathbb{Z}_{p}$. Then for any  $ m \geq 0$, the following are equivalent
	\begin{enumerate}
		\item[$(i)$] $ \sum\limits_{j=0}^{r} a_{j} \binom{j}{m} = 0$,
		\item[$(ii)$] $\frac{\partial^{m}}{\partial Y^{m}} F(X,Y)$ vanishes  at $ X = Y = 1 $.
	\end{enumerate}
\end{lemma}
\begin{proof}
	Note that 
	\begin{equation*} 
	   \frac{\partial^{m} }{\partial Y^{m}} F(X,Y)= \sum\limits_{j=m}^{r} a_{j} \binom{j}{m} m! X^{r-i} Y^{i-m} = \sum\limits_{i=0}^{r} a_{i} \binom{i}{m} m! X^{r-i} Y^{i-m}.
	\end{equation*}
	Now the lemma follows from the above by substituting $X = Y = 1 $ and noting that $ m! \neq 0 $ for $ m \geq 0 $.
\end{proof}
\begin{corollary}\label{sum vanishing T plus}
	Let $ F(X,Y) = \sum\limits_{j=0}^{r} a_{j} X^{r-j} Y^{j} \in \mathbb{Z}_{p}[X,Y]$ be a homogeneous polynomial of degree $ r $ in $ X, Y$ with coefficients in $\mathbb{Z}_{p}$.  If $ (X-Y)^{n}$ divides $F(X,Y) $ in $  \mathbb{Z}_{p}[X,Y]$,  then we have
	\begin{align*}
		\sum\limits_{j=0}^{r} a_{j} \binom{j}{m} = 0, \quad \text{for all }  m < n.
	\end{align*}
\end{corollary}
\begin{proof}
	Clearly, the condition $ (X-Y)^{n} \mid F(X,Y) $ implies that $\frac{\partial^{m}}{\partial Y^{m}} F(X,Y)$  vanishes  at $ X = Y =1 $, for all $ m < n$. The corollary now follows immediately from Lemma~\ref{sum and derivative vansihing}.
\end{proof}
\begin{lemma}\label{theta and T plus}
	Let $ F(X,Y) = \sum\limits_{j=0}^{r} a_{j} X^{r-j}Y^{j} $ be a homogeneous polynomial in $ \mathbb{Z}_{p} [X,Y]$ such that there exists a positive integer $b$ satisfying
	\[
		a_{j} \neq 0 ~~ \Longrightarrow ~ ~ j \equiv b \mod (p-1).
	\]
	If $ (X-Y)^{m} $ divides $F(X,Y) $ and  $a_0,\ldots, a_{m-1}=0$, then 
	\[ 
		T^{+} \left[ g_{n,\lambda}^{0}, F(X,Y) \right] \equiv 0 \mod p^{m}.
	\]
\end{lemma}
\begin{proof}
	By assumptions $a_0,\ldots, a_{m-1}=0$ and $ a_{j} \neq 0  \Rightarrow j \equiv b \mod (p-1)$, and Corollary~\ref{sum vanishing T plus} it follows that
	\begin{align}\label{T plus sum f2}
    \begin{split}
		\sum\limits_{i=j}^{r} a_{i} \binom{i}{j} [\mu]^{i-j}&=[\mu]^{b-j} \sum\limits_{i=0}^{r} a_{i} \binom{i}{j} = 0, \quad\text{for all }  j <  m \text{ and } \mu \in \mathbb{F}_{p}^\times \\
        \sum\limits_{i=j}^{r} a_i \binom{i}{j} [0]^{i-j} & = a_j = 0, \quad \text{ for all } j<m.
    \end{split}
	\end{align}  
	Applying the formula for $ T^{+} $, we see that 
	\begin{align*}
		T^{+} \left[ g_{n,\lambda}^{0}, F(X,Y) \right] 
        & = \sum_{  \mu \in \mathbb{F}_{p}}\left[ g_{n+1, \lambda + p^{n+1} [\mu] }^{0}, \sum_{j=0}^{r} p^{j} \left( \sum_{i=j}^{r} a_{i} \binom{i}{j} (-[\mu])^{i-j}  \right) X^{r-j} Y^{j} \right]  \\
		& \stackrel{\mathclap{\eqref{T plus sum f2}}}{=} \sum_{  \mu \in \mathbb{F}_{p}}\left[ g_{n+1, \lambda + p^{n+1} [\mu] }^{0}, \sum_{j=m}^{r} p^{j}\left( \sum_{i=j}^{r} a_{i} \binom{i}{j} (-[\mu])^{i-j}  \right) X^{r-j} Y^{j} \right] \\
		& \equiv 0 \mod p^{m}\mathbb{Z}_{p}.
	\end{align*}
	The last step follows because $ a_{i} \in \mathbb{Z}_{p} $. 
\end{proof}

For $ p \geq 3 $, $ r \geq m(p+1) $ and $ 1 \leq l \leq m \leq p-1 $, we define the following polynomial:
\begin{equation}\label{polynomial Fr,m,l}
    F_{r,m,l}(X,Y) :=   X^{l} Y^{r-m(p-1)-l} (X^{p-1} - Y^{p-1})^{m} = \theta^{m} X^{-(m-l)}Y^{r-m(p+1)+m-l}.
\end{equation}
Note that the exponents of $Y$ appearing in the above polynomial lie in the same congruence class modulo $(p-1)$. 
The following lemma describes the action of the $T^+$ operator on the function supported on the coset $KZ(g_{n,\lambda}^{0})^{-1}$ and taking the value $F_{r,m,l}$ 
at
$(g_{n,\lambda}^{0})^{-1}$. 
\begin{lemma}
    Let $ p \geq 3 $, $ r \geq m(p+1)+p $ and $ 1 \leq l \leq m \leq p-1 $. Then for all  $\lambda \in I_{n}$, we have 
    \begin{align}\label{vanishing of T plus on F r,m}
     	T^{+} \left[ g_{n,\lambda}^{0},  \theta^m X^{-(m-l)}Y^{r-m(p+1)+m-l} \right] \equiv 0 \mod p^{m}.
    \end{align}
\end{lemma}

\begin{proof}
    Since $ (X - Y) \mid  (X^{p-1} - Y^{p-1}) $ in $ \mathbb{Z}_{p}[X,Y] $, it follows that $ (X - Y)^m \mid \theta^m X^{-(m-l)}Y^{r-m(p+1)+m-l}$ in  $ \mathbb{Z}_{p}[X,Y] $. Also, note that $Y^{r-m(p+1)+2m-l} \mid \theta^m X^{-(m-l)}Y^{r-m(p+1)+m-l}$ and $r-m(p+1)+2m-l > m$. Thus the coefficients of $X^r,\ldots,X^{r-m}Y^m$ are zero. Now the result follows from Lemma~\ref{theta and T plus}.
\end{proof}

\subsection{The mod \texorpdfstring{$p$}{} LLC} \label{sec: mod p LLC}

In this section, we briefly recall the statement of the mod $p$ local Langlands Correspondence (mod $p$ LLC). For details, we refer the reader to  \cite{Breuil1}.  Recall that $G=\mathrm{GL}_2(\Q_p)$, $K=\mathrm{GL}_2(\Z_p)$ and $Z=\Q_p^\times \subset G$ and $\mathrm{ind}_{KZ}^{G}(V)$ denote  the compact induction of $V$. Also, recall that $V_r = \text{Sym}^r\bar{\mathbb{F}}_p^2$ is a representation of $KZ$ (where $p\in Z$ acts trivially) for $r\geq 0$. For $0\leq r \leq p-1$, $\lambda \in \bar{\mathbb{F}}_p$ and $\eta:\Q_p^\times \rightarrow \bar{\mathbb{F}}_p^\times$ a  smooth character, let $\pi(r,\lambda,\eta)$ be the  mod representation of $G$, given by
\begin{equation*}
    \pi(r,\lambda,\eta) :=\frac{\mathrm{ind}_{KZ}^{G}(V_r)}{T-\lambda}\otimes (\eta\circ\det)
\end{equation*}
where $T$ is the (mod $p$ version of the) Hecke operator from the previous section. Then under the mod $p$ LLC we have 
\begin{itemize}
    \item if $\lambda =0$, then $\text{ind}(\omega_2^{r+1}) \mapsto \pi(r,0,\eta)$
    \item if $\lambda \neq 0$, then $(\mu_\lambda \cdot \omega^{r+1} \oplus \mu_{\lambda^{-1}})\otimes \eta \mapsto \pi(r,\lambda,\eta)^{\mathrm{ss}} \oplus \pi(\llbracket p-3-r\rrbracket,\lambda^{-1},\eta\omega^{r+1})^{\mathrm{ss}} $, where $\mathrm{ss}$ in the superscript denotes the semi-simplification and $\llbracket p-3-r\rrbracket$ denotes the congruence class of $p-3-r$ mod $(p-1)$ in $\{0,\ldots,p-2\}$.
\end{itemize}

Let $V_{k,a_p}$ be the two dimensional crystalline representation of the Galois group of $\mathbb{Q}_p$ with Hodge-Tate weights $(0,k-1)$ and slope $v(a_p)$ from the introduction. The completion of the locally 
algebraic representation $\Pi_{k,a_p}$ corresponds to 
$V_{k,a_p}$ under the $p$-adic local Langlands correspondence. Moreover, by
the work of Berger \cite{Ber10} (or the functorial construction of Colmez \cite{Col10}), there is a compatibility between the $p$-adic and mod $p$ LLC with respect to the process of mod $p$ reduction. In particular,
$\bar{\Theta}_{k,a_p}$ corresponds to $\bar{V}_{k,a_p}$ under the mod $p$ LLC. This allows us to transfer the problem of 
computing the reduction $\bar{V}_{k,a_p}$ on the Galois side to the computation of $\bar{\Theta}_{k,a_p}$ on the automorphic side.

The following useful consequence of the mod $p$ LLC allows us to write down $\bar{V}_{k,a_p}$ explicitly in the introduction using the results of \S\ref{sec: good} and \S\ref{sec: bad}.

\begin{lemma}
    Let $k\geq 2$ and $r:=k-2$. If there is a surjection $\mathrm{ind}_{KZ}^{G}(V_r^{(m)}/V_r^{(m+1)}) \twoheadrightarrow \bar{\Theta}_{k,a_p}$ for some $ m \geq 0$ and $r \not \equiv 2m\pm 1 \mod (p-1)$, then $\bar{V}_{k,a_p} \simeq \mathrm{ind}(\omega_2^{[r-2m]+1+m(p+1)})$. 
\end{lemma}
\begin{proof}
    Under the assumption $r\not\equiv 2m\pm1\mod(p-1)$, it follows from \eqref{exact seq. singular quotients} that $V_r^{(m)}/V_r^{(m+1)}$ does not have a JH factor which is a twist of $V_{p-2}$. It follows from the 
    explicit description of the mod $p$ LLC that 
    one of the two JH factors of  $V_r^{(m)}/V_r^{(m+1)}$ 
    dies in $\bar{\Theta}_{k,a_p}$ and the other contributes irreducibly. This yields the lemma (amazingly, the answer is independent of the JH factor that survives). 
\end{proof}
\section{Jordan-H\"older factors of \texorpdfstring{$Q(i)$}{} } \label{sec: JH factors}
    
Recall that $\theta = X^{p}Y- XY^{p}$ and $V_{r}^{(m)} = \lbrace F(X,Y) \in \bar{\mathbb{F}}_p[X,Y] : \theta^{m} \mid F(X,Y) ~ \text{in} ~ \bar{\mathbb{F}}_p[X,Y]  \rbrace$. As in \cite{GR19}, for $0 \leq i \leq r$, let  $X_{r-i} $ be the $\bar{\mathbb{F}}_p[\gl_2(\f)]$-module generated by $X^{r-i}Y^{i}$ in $V_{r}$ and  $X_{r-i}^{(m)} = V_{r}^{(m)} \cap X_{r-i}$  for all $m \geq 0$. Furthermore, for $0 \leq i \leq p-1$, define
    \begin{align*}
        Q(i) := \frac{V_{r}}{X_{r-i}+V_{r}^{(i+1)}}.
    \end{align*}
    The importance of this module stems from the fact that if $v(a_p) \in (i,i+1)$ for $i \geq 0$,
    then there is a surjective
    map $$\mathrm{ind}_{KZ}^G Q(i) \twoheadrightarrow \bar{\Theta}_{k,a_p}$$
    which gives some handle on
    the structure of $\bar{V}_{k,a_p}$ by the mod $p$ local Langlands correspondence. For instance
    when $Q(i)$ is irreducible
    the structure of $\bar{V}_{k,a_p}$ can be written down instantly (except if it's dimension is $p-1$) as in \cite[Corollary 1.12]{GR19}. 
    In this section, we determine the Jordan-H\"older (JH) factors of $Q(i)$ for all $0 \leq i \leq p-1$ using the results obtained in \cite{GR19}.

We need the definitions of the sets $\mathcal{I}(a,i)$ and $\mathcal{J}(a,i)$ introduced in \cite[(4.11), (4.12)]{GR19}. For $1 \leq a$, $i \leq p-1$ and $i \neq a,$ $p-1$, the set $\mathcal{I}(a,i) \subseteq \lbrace 0,1, \ldots, p-1 \rbrace$ is a subset of the congruence classes  modulo $p$, given by
\begin{align}\label{interval I}
    \mathcal{I}(a,i) = 
    \begin{cases}
        \lbrace a-i+1, a-i+2,  \ldots , a-1, a \rbrace, & \mathrm{if}~ i < a-i < a, \\
        \lbrace a-i, a-i+1, \ldots , a-1,  a\rbrace, & \mathrm{if}~ a-i \leq i < a,  \\
        \lbrace a, a+1, \ldots, [a-i]-1, [a-i]\rbrace^{c}, & \mathrm{if}~ a < i < [a-i], \\
        \lbrace a, a+1, \ldots, [a-i]-2, [a-i]-1 \rbrace^{c}, & \mathrm{if}~ a < [a-i] \leq i,
    \end{cases}
\end{align}
where $c$ in the superscript  denotes the complement in  $\lbrace 0 , 1, \ldots, p-1 \rbrace$.  Since any $p-1$  consecutive numbers define a congruence classes modulo $p$, we may view $ \mathcal{I}(a,i) $ as an interval. For example, if $a<i<[a-i]$, then  $\lbrace a, a+1, \ldots, [a-i]-1, [a-i]\rbrace^{c} = \{p-1+a-i+1,p-1+a-i+2, \ldots, p-1+a-1, p-1+a\}$. Also, if $i\neq a$, $a+1$, then $\mathcal{I}(a,i-1) \subseteq \mathcal{I}(a,i)$, for $i \geq 2$.

Further, the subset $\mathcal{J}(a,i)  \subseteq \lbrace 0,1, \ldots, p-1 \rbrace$  of the congruence classes  modulo $p$ is defined as follows
\begin{align}\label{interval J}
    \mathcal{J}(a,i) =
    \begin{cases}
        \lbrace a-i, a-i+1 , \ldots,  a-2, a-1\rbrace, & \mathrm{if}~ i < a-i <a ,\\
        \lbrace  a-i-1, a-i, \ldots,  a-2, a-1\rbrace, & \mathrm{if}~a-i \leq i <a , \\ 
        \lbrace  a-1,a, \ldots, [a-i]-2, [a-i]-1 \rbrace^{c}, &\mathrm{if}~a<i < [a-i], \\
        \lbrace a-1, a, \ldots,  [a-i]-3, [a-i]-2 \rbrace^{c}, &\mathrm{if}~a<[a-i] \leq i,
    \end{cases}
\end{align}
where $c$ in the superscript again denotes the complement in  $\lbrace 0, 1, 2, \ldots, p-1 \rbrace$.  Again $ \mathcal{J}(a,i) $ can be viewed as an interval and we have $\mathcal{J}(a,i-1) \subseteq \mathcal{J}(a,i)$, for all $ i \neq a$, $a+1$ and $i \geq 2$.

First we determine the JH factors of $Q(i)$ when  $2i<a$.
\begin{lemma}\label{JH factor Q 2i < a}
	Let $p \geq 3$ and $r \equiv a~\mathrm{mod}~(p-1)$ with $1 \leq a \leq p-1$. Suppose  $ 0 \leq i < a-i \leq a$. For   $r \geq  i(p+1)+p$, we have
	\begin{enumerate}
		\item[$(i)$] If $r \not \equiv a-i+1$, $a-i+2, \ldots, a \mod p$, then $\text{JH factors of } Q(i) $ are
		\[ \{ V_{p-1-a+2l} \otimes D^{a-l} : 0 \leq l \leq i \}.\]
		\item[$(ii)$] If $r  \equiv a-i+1$, $a-i+2, \ldots, a \mod p$ or equivalently $r \equiv a-j+1 \mod p$ for some $1     \leq j \leq i$, then $\text{JH factors of } Q(i) $ are
		\[ \{ V_{p-1-a+2l} \otimes D^{a-l}: 0 \leq l <j \} \cup \text{JH factors of } \{ V_{r}^{(l)}/V_{r}^{(l+1)}: j     \leq l  \leq i \}.\]
	\end{enumerate}
As a consequence $Q(i)$ is not irreducible if $1 \leq i < a-i $.
\end{lemma}
\begin{proof}
	The case $i=0$ follows from the fact $X_r/X_r^{(1)} \equiv V_a$ and \eqref{exact seq. singular quotients}, as the hypothesis in $(i)$ is  vacuously true and $(ii)$ doesn't occur. Assume $i \geq  1 $. Note that if $1 \leq l \leq i$, then  $l \leq i < a- i \leq a-l < a$. Thus, by the first part of the definition of the interval \eqref{interval I}, we have
	\[ \mathcal{I}(a,l) = \lbrace a-l+1, a-l+2, \ldots, a \rbrace, ~ \forall ~ 1 \leq l \leq i.\]
	Also observe that $\mathcal{I}(a,1) \subset \cdots \subset \mathcal{I}(a,l) \subset  \mathcal{I}(a,i)$.
	
    Let $r \equiv r_{0}$ mod $p$ with $0 \leq r_{0} \leq p-1$. Observe that if $r \not \equiv a-i+1$, $a-i+2, \ldots, a \mod p$, then $r_{0} \not \in \mathcal{I}(a,i)$, so  $r_{0} \not \in \mathcal{I}(a,l)$, for $1 \leq l \leq i$. Thus, by \cite[Theorem 4.19]{GR19}, we have
	\[ 
        0 \rightarrow V_{p-1-a+2l} \otimes D^{a-l} \rightarrow Q(l) \rightarrow Q(l-1) \rightarrow 0, ~ \forall ~ 1 \leq l \leq i.
    \]
	Therefore, JH factors of $Q(l) = $ JH factors of $Q(l-1) \cup \lbrace V_{p-1-a+2l} \otimes D^{a-l} \rbrace$, for $1 \leq l \leq i $. Iterating this, we get 
	\[ 
        \text{JH factors of } Q(i) = \{ V_{p-1-a+2l} \otimes D^{a-l} : 1 \leq l \leq i \} \cup \text{JH factors of } Q(0).
	\]
	Since $Q(0) \cong V_{p-1-a}  \otimes D^{a}$, we obtain $(i)$. 
	
	If $ r \equiv a-j+1$ mod $p$, for some $1 \leq j \leq i$, then from the above description of the interval, we see that $a-j+1 =r_{0} \in \mathcal{I}(a,j)$. So $r_{0} \in \mathcal{I}(a,l)$, for $j \leq l \leq i$. Thus, by \cite[Theorem 4.19]{GR19}, we have
	\[
	    0 \rightarrow V_{r}^{(l)}/V_{r}^{(l+1)} \rightarrow Q(l) \rightarrow Q(l-1) \rightarrow 0, ~ \forall ~ j \leq l    \leq i.      
	\]
	Therefore, JH factors of $Q(l) = \text{JH factors of } Q(l-1) \cup  \text{JH factors of }  V_{r}^{(l)}/V_{r}^{(l+1)}$, for $j \leq l \leq i $. Iterating this, we get 
	\[
	    \text{JH factors of } Q(i) = \text{JH factors of } \{ V_{r}^{(l)}/V_{r}^{(l+1)} : j \leq l \leq i \} \cup \text{JH factors of } Q(j-1).
	  \]
	Now $(ii)$ follows from $(i)$ applied with $i$ equal to $j-1$.
\end{proof}
For emphasis, we record that if $r\equiv a$ mod $(p-1)$, then 
\begin{align}\label{Structure Q(0)}
    Q(0) \cong V_{p-1-a} \otimes D^a.
\end{align}

We now consider the case $r \equiv a$ mod $p$ with  $1 \leq 2i=a \leq p-1 $.

\begin{lemma}\label{JH factor Q 2i = a}
	Let $p \geq 3$ and $r \equiv a~\mathrm{mod}~(p-1)$ with $1 \leq a \leq p-1$. Suppose  $ a =2i$. For $r \geq i(p+1)+p$, we have
	\begin{enumerate}
		\item [$(i)$] If $r \not \equiv i-1$, $i, \ldots, a-1 \mod p$, then JH factors of $Q(i) =$ JH factors of $Q(i-1)$. More precisely,
		\begin{enumerate}
			\item [$(a)$] If $r \equiv a \mod p$, then  the $ \text{JH factors of } Q(i)$ are 
			\[ 
                \{ V_{p-1-a} \otimes D^{a} \} \cup \text{JH factors of } \{  V_{r}^{(l)}/V_{r}^{(l+1)}: 0 < l < i \}.
			\]
			\item [$(b)$] If $r \not \equiv i-1$, $i, \ldots$, $a-1$, $a \mod p$, then the $\text{JH factors of } Q(i) $ are
			\[
			    \{ V_{p-1-a+2l} \otimes D^{a-l} : 0 \leq l < i \}.
			\]
		\end{enumerate}
		\item [$(ii)$] If $r \equiv i-1 \mod p$, then  the JH factors of $Q(i) = \{ V_{p-1-[a-2l]} \otimes D^{a-l} : 0          \leq l \leq i \}$.
		\item [$(iii)$] If $r \equiv i \mod p$, then the JH factors of $Q(i) = \{ V_{p-1-a+2l} \otimes D^{a-l} : 0 \leq l         \leq i \}$.
		\item [$(iv)$]If $r  \equiv i+1$, $i+2, \ldots, a-1 \mod p$ or equivalently $r \equiv a-j+1 \mod p$ for some $2        \leq j \leq i$, then  the JH factors of  $Q(i) $ are
		\[ 
           \{ V_{p-1-a+2l}  \otimes D^{a-l}: 0 \leq l <j \}   \cup \text{JH factors of } \{ V_{r}^{(l)}/V_{r}^{(l+1)}: j          \leq l  \leq i \}.
		\] 
	\end{enumerate}
	As a consequence, $Q(i)$ is not irreducible if $a \neq 2$ and $a=2i$. Also, if $a=2$, then $Q(i)$ is irreducible if and only      if  $r  \not \equiv 0$, $1 \mod p$.
\end{lemma}
\begin{proof}
	Note that  $ a-(i-1) =a-i+1 =i+1 > i-1$. Thus, by Lemma~\ref{JH factor Q 2i < a} $(i)$, we have
	\[
	     \text{JH factors of } Q(i-1) = \lbrace V_{p-1-a+2l}\otimes D^{a-l} : 0 \leq l \leq i-1 \rbrace, \text{ if } r \not \equiv i+2, \ldots , a \mod p. 
	\]
	
	  By the second part of the definition of the interval \eqref{interval J}, we have
	\[
	    \mathcal{J}(a,i) = \lbrace i-1, i, \ldots, a-1 \rbrace.
	  \]
	\begin{enumerate}
		\item[$(i)$] By \cite[Theorem 4.23]{GR19}, we have $Q(i) \cong Q(i-1)$.  If $i=1$, then parts $(a)$ and $(b)$ hold since         $Q(0) \cong V_{p-1-a} \otimes D^{a}$. For $i \geq 2$, parts $(a)$ and $(b)$  follow from Lemma~\ref{JH factor Q 2i < a} $(ii)$ (with $i$ there equal to $i-1$ and $j=1$) and the above description of $Q(i-1)$ respectively.		
		\item[$(ii)$] Note that $[a-2i]= p-1$ and $[a-2l]=a-2l$ for $ 0 \leq l < i$. Now $(ii)$ follows from \cite[Theorem 1.8]{GR19}  and the above description of $Q(i-1)$.
		\item[$(iii)$] This follows from \cite[Theorem 4.23]{GR19}  and the above description of $Q(i-1)$.
		\item[$(iv)$] By \cite[Theorem 4.23]{GR19}, we see that JH factors of $Q(i) =$ JH factors of $ \lbrace V_{r}^{(i)}/V_{r}^{(i+1)} \rbrace$ $\cup $ JH factors of $Q(i-1)$. Now part $(iv)$ follows from the above description of 
		$Q(i-1)$ if $r \equiv i+1  = a-i+1$ mod $p$. If $r \equiv i+2, \ldots, a-1$ mod $p$, then $r \equiv a-j+1$ mod $p$ for      some $2 \leq j \leq i-1$. Now  part $(iv)$ follows from Lemma~\ref{JH factor Q 2i < a} $(ii)$ applied for $Q(i-1)$.
		\qedhere
	\end{enumerate}
\end{proof}
We next consider the case $r \equiv a$ mod $(p-1)$ with $1 \leq a \leq p-1$ and $1 \leq a-i < i < p-1$. Let $r \equiv r_0$ mod $p$ with $0 \leq r_0 \leq p-1$. Unlike in the previous two lemmas, the JH factors in this case further depend on the comparison between $a$ and $r_0$. By the second part of the definition of interval \eqref{interval J}, we have 
    \[
        \mathcal{J}(a,i) = \lbrace a-i-1, a-i, \ldots, a-1 \rbrace.
    \]
Further, by the first part of the definition of the interval \eqref{interval I}, we have
    \[
        \mathcal{I}(a,a-i) = \{ i+1, i+2, \ldots, a \},
    \]
    hence
    \[
       \mathcal{J}(a,i) \smallsetminus \mathcal{I}(a,a-i) =\lbrace a-i-1, a-i, \ldots, i \rbrace .
    \]
\begin{lemma}\label{JH factor Q 2i > a}
	  Let $p \geq 3$ and $r \equiv a~\mathrm{mod}~(p-1)$ with $1 \leq a \leq p-1$. Let $r \equiv r_{0} \mod p$ with $0 \leq r_{0} \leq p-1$. Suppose $1 \leq a-i < i < p-1$. For $r \geq i(p+1)+p$, we have
	        \begin{enumerate}
		    \item [$(i)$] If $r \not \equiv a-i-1, a-i, \ldots, a-1 \mod p$, then  $Q(i) \cong Q(a-i-1)$. More precisely, we         have
		      \begin{enumerate}
			\item [$(a)$] If $r \equiv a \mod p$, then the JH factors of $Q(i)$ are 
			\[ 
			   \{ V_{p-1-a}\otimes D^{a} \}  \cup \text{JH factors of } \{V_{r}^{(l)}/V_{r}^{(l+1)} : 0 <l \leq a-i-1 \}.
			\]
			\item[$(b)$] If $r \not \equiv a-i-1, a-i, \ldots, a-1,a \mod p$, then JH factors of $Q(i)$ are 
			\[ \{ V_{p-1-a+2l} \otimes D^{a-l} : 0 \leq l \leq a-i-1 \} .\]
		    \end{enumerate}
		      \item[$(ii)$] Assume $r \equiv  a-i-1, a-i, \ldots, i \mod p$. Then
		    \begin{enumerate}
			\item [$(a)$] If $a < 2 r_{0} +1$, then the JH factors of $Q(i)$ are 
			\[\{ V_{[a-2r_{0}]} \otimes D^{r_{0}}\} \cup 
			    \text{JH factors of } \{ V_{r}^{(l)}/V_{r}^{(l+1)} : a-r_{0} <l \leq i \},
		      \]
			$\qquad \qquad  \qquad \cup ~  \{ V_{p-1-a+2l} \otimes D^{a-l} : 0 \leq l \leq a-i-1 \}$.
			\item[$(b)$] If $a=2r_{0}+1$, then the JH factors of $Q(i)$ are
			\[ \{ V_{r}^{(l)}/V_{r}^{(l+1)} : a-r_{0} \leq l \leq i \} \cup 
			\{ V_{p-1-a+2l}\otimes D^{a-l} : 0 \leq l \leq a-i-1 \}.\]
			\item[$(c)$] If $a > 2r_{0}+1$, then the JH factors of $Q(i)$ are
			\[  
			\{ V_{p-1-[2r_{0}+2-a]} \otimes D^{r_{0}+1} \} \cup
			\text{JH factors of } \{ V_{r}^{(l)}/V_{r}^{(l+1)} : a-r_{0}  \leq l \leq  i \},
			\]
			$ \qquad \qquad \qquad \cup ~ \{ V_{p-1-a+2l}\otimes D^{a-l} : 0 \leq l \leq a-i-1 \}$ .
		    \end{enumerate}
		    \item[$(iii)$] If $r \equiv i+1$, $i+2, \ldots, a-1 \mod p$ or equivalently $r \equiv a-j+1 \mod p$ for some $ 2         \leq j \leq a-i$, then the JH factors of $Q(i)$  are
		      \[ \{ V_{r}^{(l)}/V_{r}^{(l+1)} : j  \leq l \leq i \} \cup  \{ V_{p-1-a+2l}\otimes D^{a-l} : 0 \leq l <j \}. \]
	        \end{enumerate}
	        As a consequence $Q(i)$ is not irreducible if $i \neq a-1$. Also, $Q(a-1)$ is irreducible if and only if $r \not        \equiv 0$, $1, \ldots, a-1 \mod p$.
        \end{lemma}      
        \begin{proof}
	         From \cite[(4.18), (4.21)]{GR19}, it follows that the hypothesis of $(i)$, $(ii)$ and $(iii)$ are equivalent to $r_{0} \not \in \mathcal{J}(a,i)$, $r_{0}    \in \mathcal{J}(a,i) \smallsetminus \mathcal{I}(a,a-i)$, and $r_{0} \in \mathcal{I}(a,a-i)$ and $r \not \equiv a $   mod $p$ respectively. 
	        Also note that $0 \leq a-i-1 < i+1 = a-(a-i-1)$. Thus by Lemma~\ref{JH factor Q 2i < a} $(i)$,we have 
	        \[
	           \text{ JH factors of } Q(a-i-1) = \{ V_{p-1-a+2l} \otimes D^{a-l} : 0 \leq l \leq a-i-1 \},
	        \]
	        $\text{ if } r \not \equiv i+2, \ldots, a-1 \mod p$ (so for all $r$  if $i+2 > a-1$).	 
	        \begin{enumerate}
		    \item[$(i)$] By  \cite[Theorem 4.31 (i)]{GR19}, we have $Q(i) \cong Q(a-i-1)$.
		       Now part $(a)$ follows from Lemma~\ref{JH factor Q 2i < a} $(ii)$ (with $i$ there equal to $a-i-1$ and $j=1$).  
		       Part $(b)$ follows from the above description of $Q(a-i-1)$.
		      \item[$(ii)$] Follows from \cite[Theorem 4.31 (ii)]{GR19} and  the above description of $Q(a-i-1)$.
		    \item[$(iii)$] The case $r \equiv i+1$ mod $p$ follows from  \cite[Theorem 1.9 (iii)]{GR19} and 
		      the above description of $Q(a-i-1)$. The case $r \equiv i+2, \ldots , a-1 $ mod $p$ follows from              \cite[Theorem 1.9 (iii)]{GR19} and Lemma~\ref{JH factor Q 2i < a} $(ii)$ (with $i$ there equal to $a-i-1$).
	        \end{enumerate} 
               It is easy to see that if $i \neq a-1$, then $Q(i)$ has at least two JH factors so is not irreducible. The last statement also follows by a similar argument. 
        \end{proof}
         We next determine the JH factors of $Q(a)$, where $r \equiv a$  mod $(p-1)$ and $ 1 \leq a \leq p-1$.
        \begin{lemma}\label{JH factor Q(a)}
	        Let $p \geq 3$ and $r \equiv a~\mathrm{mod}~(p-1)$ with $1 \leq a \leq p-1$.  Let $r \equiv r_{0} ~\mathrm{mod}~p$      with $0 \leq r_{0} \leq p-1$. For $r \geq a(p+1)+p$, we have
	        \begin{enumerate}
		      \item[$(i)$] If $r \equiv 0$, $1, \ldots, a-1 \mod p $ and $r \not \equiv p-2 \mod p$, then
		      \begin{enumerate}
			\item[$(a)$]  If $a< 2r_{0}+1$, then the JH factors of $Q(a)$ are 
			\[
			\text{JH factors of } \{  V_{r}^{(l)}/V_{r}^{(l+1)} : a-r_{0}< l \leq a \} \cup \{ V_{[a-2r_{0}]} \otimes                D^{r_{0}}\}. 
			\]
			\item[$(b)$]  If $a = 2r_{0}+1$, then the JH factors of $Q(a)$ are 
			\[
			   \text{JH factors of } \{  V_{r}^{(l)}/V_{r}^{(l+1)} : a-r_{0} \leq l \leq a \}.
			\]
			\item[$(c)$]  If $a > 2r_{0}+1$, then the JH factors of $Q(a)$ are 
			\[
			   \text{JH factors of } \{  V_{r}^{(l)}/V_{r}^{(l+1)} : a-r_{0} \leq l \leq a \}\cup \{ V_{p-1-[2r_{0}+2-a]} \otimes D^{r_{0}+1} \}.
			\]
		    \end{enumerate}
		    \item[$(ii)$] If $a=p-1$ and $r \equiv p-2 \mod p$, then the JH factors of  $Q(a)$ are 
	           $  \{ V_{2} \otimes D^{p-2} \}$   $\cup $ JH factors of $\{ V_{r}^{(l)}/V_{r}^{(l+1)}: 1 < l \leq p-2 \}$ $\cup$ $\{       V_{0} \}$.
		      \item[$(iii)$] If $r  \equiv a$, $a+1, \ldots, p-1 \mod p$,  then $Q(a) \cong V_{a}$.	
	        \end{enumerate}
             As a consequence, $Q(a)$ is irreducible if and only if  $r \equiv a$, $a+1, \ldots, p-1 \mod p$.
        \end{lemma}
        \begin{proof}
	        By the third part of \cite[Lemma 4.4]{GR19} (with  $i=a-1$ and $j=a$), we have $X_{r-(a-1)}^{(a)}/X_{r-(a-1)}^{(a+1)} =     X_{r}^{(a)}/X_{r}^{(a+1)}$. Thus, by \cite[Proposition 4.8]{GR19}, we have
	            \[
	             \frac{X_{r-(a-1)}^{(a)}}{X_{r-(a-1)}^{(a+1)}} =
	             \begin{cases}
	              V_{p-1-a} \otimes D^{a} & \text{if } r \equiv a, a+1, \ldots, p-1 \mod p, \\
	              0  & \text{if } r \equiv 0, 1, \ldots, a-1 \mod p.
	           \end{cases}
	          \]
	            By the  exact sequence \eqref{exact seq. singular quotients}, we have the cokernel of the inclusion $X_{r-(a-1)}^{(a)}/X_{r-(a-1)}^{(a+1)} \hookrightarrow V_{r}^{(a)}/V_{r}^{(a+1)}$ is equal to $V_{a}$ if $r \equiv a, a+1, \ldots, p-1 \mod p$ and is equal to $V_{r}^{(a)}/V_{r}^{(a+1)}$ otherwise. Thus, by  the exact sequence  \cite[(4.1)]{GR19}, we have 
	          \begin{align}\label{exact sequence Q(a-1) and P(a)}
	               0 \rightarrow W'' \rightarrow  P(a) \rightarrow Q(a-1) \rightarrow 0,
	            \end{align}
	          where $W''$ equals $V_{a} $ if $r \equiv a, a+1, \ldots, p-1 \mod p$ and is equal to $V_{r}^{(a)}/V_{r}^{(a+1)}$ otherwise.  
	
	        Let  $a=1$. Then the above sequence  determines $P(a)$ since $Q(a-1) = Q(0) \cong V_{p-1-a} \otimes D^{a}$, by               \eqref{Structure Q(0)}. The  subcase $ r\equiv 0 \mod p$ is  included in part $(i)$ $(b)$ and  follows from 
	        \cite[Theorem 4.35 (i)]{GR19}. The subcase $r \equiv 1,2, \ldots, p-1 \mod p$ is included in part $(ii)$
    	    and again follows from \cite[Theorem 4.35 (i)]{GR19}. So assume $a \geq 2$.
	
	        \begin{enumerate}
		    \item[$(i)$] By the exact sequence \eqref{exact sequence Q(a-1) and P(a)}, we have $P(a)$ is an extension of                  $V_{r}^{(a)}/V_{r}^{(a+1)}$ by $Q(a-1)$. Thus, by \cite[Theorem 4.35 (i)]{GR19}, we have JH factors of $Q(a)$             are
		        \[
		         \text{JH factors of } Q(a-1)  \cup \text{JH factors of } \{ V_{r}^{(a)}/V_{r}^{(a+1)} \} \smallsetminus \lbrace  V_{p-1-a} \otimes D^{a} \rbrace  
		    \]
	        If $a=2$, then the JH factors of $Q(a-1)$ are given by  Lemma~\ref{JH factor Q 2i = a} $(ii)$ and $(iii)$ (applied with     $i$ there equal to $1$) for $r \equiv 0 \mod p$ and $r \equiv 1 \mod p$, respectively.  These subcases are included in      part $(i)$ $(c)$ and part $(i)$ $(a)$ and we see the lemma follows in these subcases.
                If $ a-(a-1) = 1 < a-1$, then  the JH factors of $Q(a-1)$  are given by  Lemma~\ref{JH factor Q 2i > a} $(ii)$ (with    $i$ there equal to $a-1$) and again the lemma follows.	
		    \item[$(ii)$] By the exact sequence \eqref{exact seq. singular quotients} (with $m=a=p-1$), the JH factors of                            $V_{r}^{(a)}/V_{r}^{(a+1)}$ are $V_{0}$ and $V_{p-1}$. Thus, by the exact sequence  \eqref{exact sequence Q(a-1) and         P(a)}, the JH factors of $P(a)$ are $V_{0}$, $V_{p-1}$ and  the JH factors of $Q(a-1)$. It follows from \cite[Theorem  4.35 (ii)]{GR19} that the JH factors of $Q(a)$ are the same as the JH factors of $Q(a-1)$. If $p=3$, then $a-1 = p-2       = 1$ and the lemma follows from Lemma~\ref{JH factor Q 2i = a} $(iii)$  (applied with $i$ there equal to $1$). If $p       >3$, then $ a-1 = p-2  > 1= a-(a-1)$ and the lemma follows from Lemma~\ref{JH factor Q 2i > a} $(ii)$ $(a)$ (applied       with $i$ there equals $a-1$).
	        \item[$(iii)$] If $a=2$, then by Lemma~\ref{JH factor Q 2i = a} $(i)$ (with $i$ there equal to $1$), $Q(a-1) \cong Q(0)      \cong V_{p-1-a} \otimes D^{a}$, by \eqref{Structure Q(0)}. If $a>2$, then by Lemma~\ref{JH factor Q 2i > a} $(i)$            (applied with $i$ there equals $a-1$), we see that again $Q(a-1) \cong Q(0) \cong V_{p-1-a} \otimes D^{a}$. Thus, by        the exact sequence \eqref{exact sequence Q(a-1) and P(a)}, we have $P(a)$ is an extension of $V_{p-1-a} \otimes             D^{a}$ by $V_{a}$. Now the lemma  follows from \cite[Theorem 4.35 (i)]{GR19}.	
	          \end{enumerate}
	          The last assertion immediately follows from parts  $(i)$-$(iii)$.
        \end{proof}

        The four lemmas above cover all the cases $i \leq a$. In the  next four lemmas, we determine the JH factors of $Q(i)$ for $a< i \leq p-1$. Note that $[a-i] = p-1+a-i \geq a $ in this case.

        \begin{lemma}\label{JH factor Q i < [a-i]}
	        Let $p \geq 3$ and $r \equiv a~\mathrm{mod}~(p-1)$ with $1 \leq a \leq p-1$ and suppose  $a<i < [a-i] =p-1+a-i$.  Let       $r \equiv r_{0} ~\mathrm{mod}~p$ with $0 \leq r_{0} \leq p-1$. For $r \geq i(p+1)+p$, we have
	        \begin{enumerate}
				\item[$(i)$] If $r \equiv 0$, $1, \ldots, a-1 \mod p $, then
		    \begin{enumerate}
			\item[$(a)$]  If $a< 2r_{0}+1$, then the JH factors of $Q(i)$ are 
			\[
			\text{JH factors of } \{  V_{r}^{(l)}/V_{r}^{(l+1)} : a-r_{0}< l \leq i \} \cup \{ V_{[a-2r_{0}]} \otimes D^{r_{0}} \}. 
			\]
			\item[$(b)$]  If $a = 2r_{0}+1$, then the JH factors of $Q(i)$ are 
			\[
			   \text{JH factors of } \{  V_{r}^{(l)}/V_{r}^{(l+1)} : a-r_{0} \leq l \leq i \}.
			\]
			\item[$(b)$]  If $a > 2r_{0}+1$, then the JH factors of
			$Q(i)$ are 
			\[
			\text{JH factors of } \{  V_{r}^{(l)}/V_{r}^{(l+1)} : a-r_{0} \leq l \leq i \} \cup \{ V_{p-1-[2r_{0}+2-a]} \otimes          D^{r_{0}+1} \}.
			\]
		    \end{enumerate}
		    \item[$(ii)$] If $r  \equiv a$, $a+1, \ldots, [a-i] \mod p$, then the JH factors of $Q(i)$ are
		    \[
		       \{ V_{2l-a} \otimes D^{a-l} : a \leq l \leq i \} .
		      \]
		      \item[$(iii)$] If  $r \equiv [a-i]+1, \ldots, p-1$ or equivalently  $r \equiv [a-j]+1 \mod p$ for some $a+1 \leq j \leq    i$, then the JH factors of   $Q(i)$ are 
				\[
		          \{ V_{2l-a} \otimes D^{a-l} : a \leq l < j \} \cup  \text{JH factors of }  \{ V_{r}^{(l)}/V_{r}^{(l+1)} : j \leq l \leq i \}.
				\]
	        \end{enumerate}
                As a consequence, $Q(i)$ is not irreducible. 
        \end{lemma}
       \begin{proof}
	        The proof is similar to Lemma~\ref{JH factor Q 2i < a}.  Note that if $a+1 \leq l \leq i$, then $[a-l] = p-1+a-l$ and
	        $l \leq i < [a-i] \leq [a-l] <  p-1$. Thus, by the third part of \eqref{interval I}, we have
	        \begin{align}\label{interval I i < p-1+a-i}
	              \mathcal{I}(a,l) =  \{ 0, 1, \ldots, a-1 \} \cup \lbrace [a-l]+1, [a-l]+2, \ldots, p-1 \rbrace, ~ \forall ~ a+1 \leq l \leq i.
	          \end{align}    	    
	        Assume $r \not \equiv  a, a+1, \ldots, [a-i] \mod p$. Then $r_{0} \in \mathcal{I}(a,i)$. Let $a+1 \leq l \leq i$ be the     smallest integer such that $r_{0} \in \mathcal{I}(a,l)$. Since  $ \mathcal{I}(a,l) \subseteq  \mathcal{I}(a,l+1)            \subseteq \cdots \subseteq  \mathcal{I}(a,i)$, by \cite[Theorem 4.19]{GR19}, we have JH factors of $Q(l') =$ JH factors       of $Q(l'-1)$ $\cup$ JH factors  $V_{r}^{(l')}/V_{r}^{(l'+1)} $, for all $l \leq l' \leq i $. Thus,  for $r_{0} \in         \mathcal{I}(a,l)$, we have
	      \begin{align}\label{recursion JH  i< [a-i]}
	        \text{JH factors of } Q(i) = \text{ JH factors of } Q(l-1) \cup 
                 \text{ JH factors of } \{ V_{r}^{(l')}/V_{r}^{(l'+1)} : l \leq l' \leq i \}
	      \end{align}
	    \begin{enumerate}
	    	\item [$(i)$] Observe that 	$r_{0} \in \mathcal{I}(a,a+1)$. So $l=a+1$
	    	and the lemma follows from \eqref{recursion JH  i< [a-i]} and Lemma~\ref{JH factor Q(a)} $(i)$.
	    	\item [$(ii)$]  By \eqref{interval I i < p-1+a-i}, we see that $r_{0} \not \in \mathcal{I}(a,i)$. 
	    	Since $\mathcal{I}(a,i) \supseteq \mathcal{I}(a,i-1)
	    	 \supseteq \cdots \supseteq \mathcal{I}(a,a+1)$,  
	    	 we see that $r_{0} \not \in \mathcal{I}(a,l)$ for all 
	    	 $a+1 \leq l \leq i$. 
	    	 By \cite[Theorem 4.19]{GR19}, we have 
	    	\[
	    	    0 \rightarrow V_{p-1-[a-2l]} \otimes D^{a-l} \rightarrow Q(l)
	    	    \rightarrow Q(l-1) \rightarrow0, ~ \forall ~ a+1 \leq l \leq i.
	    	\]
	    	It can be checked that $[a-2l] = [a-l]-l = p-1+a-2l$, for  all
	    	$a+1 \leq l \leq i$. 
	    	So the JH factors of $Q(l-1)$ = JH factors of $Q(l) \cup 
	    	\{ V_{2l-a} \otimes D^{a-l} \}$, for all $a+1 \leq l \leq i$. 
	    	Iterating this, we get
	    	\begin{align}\label{JH of Q(i) interms of Q(a)}
	    	\text{ JH factors of } Q(i) = \{ V_{2l-a} \otimes  D^{a-l} : a+1 \leq l \leq i \} \cup  \text{JH factors of } Q(a).
	    	\end{align}
	    	Now $(ii)$ follows from Lemma~\ref{JH factor Q(a)} $(iii)$.  	
          	\item[$(iii)$]  By \eqref{interval I i < p-1+a-i}, we have $r_{0} \in \mathcal{I}(a,j)$.  If $j-1>a$, then
                 $[a-j]+1 = p-1+a-j > a-1$ and $[a-j]+1  <  [a-j]+2 = [a-(j-1)]+1 $. So  $r_{0} \not \in \mathcal{I}(a,j-1)$,
          	 by \eqref{interval I i < p-1+a-i}. So $l=j$.
          	 Thus, by 	\eqref{recursion JH  i< [a-i]}, we have 
          	\[
          	\text{JH factors of } Q(i) = \text{ JH factors of } Q(j-1) \cup \text{ JH factors }
          	\{ V_{r}^{(l')}/V_{r}^{(l'+1)} : j \leq l' \leq i \}
          	\]
          	If $j-1=a$, then $(iii)$ follows from Lemma~\ref{JH factor Q(a)} $(iii)$ as $r \equiv [a-j]+1 = p-1 \mod p$. Assume $j-1 > a$. Then  $j-1 < [a- (j-1)]$ as $a< j-1 < i$, so $j-1$ satisfies the hypothesis of the lemma.  Also note that $[a-(j-1)] = p-1+a-j+1 = [a-j]+1$. Now $(iii)$ follows from part $(ii)$ of the lemma applied with $i$ equal to $j-1$. 
               \end{enumerate}
                The last assertion is an immediate consequence of the parts $(i)$-$(iii)$.
        \end{proof}

        \begin{remark}
	\label{compatibility of JH Q(a) and i<[a-i]}
	     Observe that if $a< p-1$, then part $(ii)$ of \Cref{JH factor Q(a)} doesn't occur. Taking $i=a$  in the assertions $(i)$   and $(ii)$ of   \Cref{JH factor Q i < [a-i]}, we obtain $(i)$ and $(iii)$ of Lemma~\ref{JH factor Q(a)}, respectively.     Thus the  lemma above holds for  $a \leq i < [a-i] $.
        \end{remark}

        \begin{lemma}
	\label{JH factor Q i = [a-i]}
	Let $p \geq 3$ and $r \equiv a~\mathrm{mod}~(p-1)$ with $1 \leq a \leq p-1$ and suppose $a<i = [a-i] =p-1+a-i$.  Let $r         \equiv r_{0} ~\mathrm{mod}~p$ with $0 \leq r_{0} \leq p-1$. For $r \geq i(p+1)+p$, we have
        \begin{enumerate}
	\item[$(i)$] If  $r \equiv 0, \ldots ,a-2 \mod p$, then the 
	JH factors of $Q(i)$  are
	\[  
	\text{JH factors   of } Q(a) \cup  \text{JH factors of }
		\{ V_{r}^{(l)}/V_{r}^{(l+1)} :  a+1 \leq l \leq i \}.
	\]	
	\item[$(ii)$] If $r \equiv a-1$, $a, \ldots, i-2  \mod p$, then
	$Q(i) \cong Q(i-1)$. More precisely,
	\begin{enumerate}
		\item [$(a)$] If $r \equiv a-1 \mod p$, then the JH factors of $Q(i)$
		are 
		\[
		   \{ V_{p+1-a} \otimes D^{a-1} \} \cup  
		   \text{ JH factors of } \{  V_{r}^{(l)}/V_{r}^{(l+1)} : 2 \leq l \leq i-1 \} .		    
		\]
		\item[$(b)$]  If $r \equiv a$, $ a+1, \ldots, i-2 \mod p$, then the JH factors of $Q(i)$
		are 
		\[
		     \{ V_{2l-a} \otimes D^{a-l} : a \leq l \leq i-1 \}.
		\]
	\end{enumerate}
	 \item[$(iii)$] 	If  $r \equiv i-1 \mod p $, 	then the JH factors of $Q(i) =
	   \{ V_{p-1-[a-2l]} \otimes D^{a-l} : a \leq l \leq i\}$. 
	   \item[$(iv)$] 	If  $r \equiv i \mod p $, 	then the JH factors of $Q(i) =
	   \{ V_{2l-a} \otimes D^{a-l} : a \leq l \leq i\}$.
	\item[$(v)$] If $r \equiv i+1, \ldots, p-1 \mod p$ or equivalently
	$r \equiv [a-j]+1 \mod p$ for some $a+1 \leq j \leq i$,  then the 
	JH factors of $Q(i) $ are 
	\[ \text{JH factors  of } Q(i-1) \cup  \text{JH factors of }
      V_{r}^{(i)}/V_{r}^{(i+1)}.
      \]	
      More precisely, the JH factors of $Q(i)$ are 
	\[
	   \text{JH factors of } \{ V_{r}^{(l)}/V_{r}^{(l+1)} : j \leq l \leq i  \} \cup 
	   \{ V_{p-1-[a-2l]} \otimes D^{a-l} : a \leq l < j \}.
	\]																																					
\end{enumerate}	
As a consequence, $Q(i)$ is not irreducible.
\end{lemma}	
\begin{proof}
	Since $ a< i =[a-i] $, by the fourth part of \eqref{interval J}, we have 
	\[
	    \mathcal{J}(a,i) = \{ a-1, a, \ldots, [a-i]-2 \}^{c}
	    =\{ 0,1, \ldots, a-2\} \cup \{ i-1, i, \ldots, p-1 \}.
	\]
	Note that $[a-(i-1)] = p+a-i = [a-i]+1 > [a-i] = i > i-1 \geq a$. So, either $i-1 =a$ or $i-1$ satisfies the hypotheses of Lemma~\ref{JH factor Q i < [a-i]}. Thus, by \Cref{JH factor Q i < [a-i]}  $(ii)$  and  Remark~\ref{compatibility of JH Q(a) and i<[a-i]}, we have
	\[
	      \text{ JH factors of } Q(i-1) = \{ V_{2l-a} \otimes D^{a-l} : a \leq l \leq i-1 \},
	      \text{ if } r \equiv a, a+1, \ldots, i+1.
	\]
    \begin{enumerate}
    	    \item [$(i)$] By \cite[Theorem 1.8]{GR19}, we have JH factors of $Q(i)$ = JH factors of $Q(i-1)$ $\cup$  JH factors      of  $V_{r}^{(i)}/V_{r}^{(i+1)} $. If $i-1 =a$, then we are done. Else note that $r_{0} \in \{ 0,1, \ldots, a-2 \}       \subset \mathcal{I}(a,a+1)$  by \eqref{interval I i < p-1+a-i}. So $l =a+1$ in the notation just above                  \eqref{recursion JH  i< [a-i]}. Now  $(i)$ follows  from \eqref{recursion JH  i< [a-i]} (with  $i$ equal to $i-1$        and $l$ equal to $a+1$).
  
    	    \item[$(ii)$]  As $r_0 \not \in \mathcal{J}(a,i)$, by \cite[Theorem 1.8]{GR19}, we have $Q(i) \cong Q(i-1)$.  Since        $p-1+a-i =i$,  it forces $a$ to be even. Therefore $a \geq 2$ or equivalently $a < 2(a-1)+1$.Now part (a)  follows      from Lemma~\ref{JH factor Q i < [a-i]} $(i)$ $(a)$ (with $i$ there equal to $i-1$) and Remark~\ref{compatibility of JH Q(a) and i<[a-i]}. Part  $(b)$ follows from  the above description of $Q(i-1)$.  
    	    \item[($iii)$] Follows from \cite[Theorem 1.8]{GR19} and from the above description of $Q(i-1)$, noting that $p-1-       [a-2l] =2l-a$  for $a \leq  l \leq  i-1$ and $p-1-[a-2i] =0$.
    	    \item[$(iv)$] Follows from \cite[Theorem 1.8]{GR19}  and  from  the above description of $Q(i-1)$. 
    	    \item[$(v)$] By \cite[Theorem 1.8]{GR19}, we have  JH factors of $Q(i)$ = JH factors of $Q(i-1)$ $\cup $ JH factors        of $V_{r}^{(i)}/V_{r}^{(i+1)}$. The case $r_{0}=i+1$, follows from the above description of $Q(i-1)$. The other         cases follow from Lemma~\ref{JH factor Q i < [a-i]} $(iii)$ as $r \equiv [a-j]+1$ mod $p$ for some $a+1 \leq j <i$.
              \end{enumerate} 
	          The non-irreducibility  of $Q(i)$  is obvious  for parts $(i)$, $(iii)$-$(v)$  since it can be checked that in each of  the above cases $Q(i)$ has at least two JH factors as $a<i$. Observe that part $(ii)$ is vacuous if $i<a+1$. If $i > a+1$, then $i \geq 3$ so the non-irreducibility  also holds in part $(ii)$ $(a)$ and $(ii)$ $(b)$. If $i =a+1$, then  $(ii)$ $(b)$ is vacuous and we have  $2 \leq a+1 = i = [a-i] = p-2$, so $p \geq 5 $ and $i \geq 3$,  and  again non-irreducibility holds in part $(ii)$ $(a)$.
        \end{proof}  

        \begin{lemma}\label{JH factor Q i > [a-i]}
	        Let $p \geq 3$ and $r \equiv a~\mathrm{mod}~(p-1)$ with $1 \leq a \leq p-1$ and suppose $a<  i < p-1$ and $i>[a-i] =p-1+a-i$.  Let $r \equiv r_{0} ~\mathrm{mod}~p$  with $0 \leq r_{0} \leq p-1$. For $r \geq i(p+1)+p$, we have
	     \begin{enumerate}
	     	\item[$(i)$] If $r \equiv 0$, $1, \ldots, a-2 \mod p$, then the JH factors of $Q(i)$ are 
	     	\[
	     	\text{ JH factors of } \{ V_{r}^{(l)}/V_{r}^{(l+1)} : a+1 \leq l \leq i \} 
	     	\cup  \text{JH factors of } Q(a).
	     	\] 
	     	\item[$(ii)$] If $r \equiv a-1$, $a, \ldots, [a-i]-2 \mod p$, then
	     	$Q(i) \cong Q([a-i]-1)$. More precisely, 
	     	\begin{enumerate}
	     		\item[$(a)$] If $p \nmid r \equiv a-1 \mod p$, 
	     		then  the JH factors of $Q(i)$	are 
	     		\[
	     		\{ V_{p+1-a} \otimes D^{a-1} \} \cup  
	     		\text{ JH factors of } \{  V_{r}^{(l)}/V_{r}^{(l+1)} : 2 \leq l \leq [a-i]-1 \} .	
	     		\]
	     		If $p \mid r \equiv a-1 \mod p$, then the 
	     		JH factors of $Q(i)$	are 
	     		\[
	     		\text{JH factors of  } \{  V_{r}^{(l)}/V_{r}^{(l+1)} : 1 \leq l \leq [a-i]-1 \} .
	     		\]
	     		\item[$(b)$] If $r \equiv a$, $a+1, \ldots, [a-i]-2 \mod p$, then 
	     		the JH factors of $Q(i)$	are  
	     		\[
	     		\{ V_{p-1-[a-2l]} \otimes D^{a-l} : a \leq l \leq [a-i]-1 \}.
	     		\]
	     	\end{enumerate}
     	    \item[$(iii)$] Assume $r \equiv [a-i]-1, \ldots, i \mod p$.
     	    \begin{enumerate}
     	    	\item[$(a)$] If $[a-r_{0}] < r_{0}+1$, then  the JH factors of $Q(i)$ 
     	    	are 
     	    	\[ 
     	    	      \text{JH factors of } \{ V_{r}^{(l)}/V_{r}^{(l+1)} : [a-r_{0}] < l \leq i \}
     	    	      \cup     \{ V_{[a-2r_{0}]} \otimes D^{r_{0}}  \}   	 
     	    	  \]
     	    	   \qquad \qquad \qquad  $ \cup$  $\{ V_{p-1-[a-2l]} \otimes D^{a-l} : a \leq l \leq [a-i]-1 \}   $.
     	    	  \item[$(b)$] If $[a-r_{0}] = r_{0}+1$, then  the JH factors of $Q(i)$  are 
     	    	  \[ 
     	    	   \text{JH factors of } \{ V_{r}^{(l)}/V_{r}^{(l+1)} : [a-r_{0}] \leq l \leq i \}
     	    	  \, \cup \,
     	    	  \{ V_{p-1-[a-2l]} \otimes D^{a-l} : a \leq l \leq [a-i]-1 \} .
     	    	  \]
     	    	   \item[$(c)$]  If $[a-r_{0}] > r_{0}+1$, then  the JH factors of $Q(i)$ 
     	    	   are
     	    	   \[ 
     	    	   \text{JH factors of } \{ V_{r}^{(l)}/V_{r}^{(l+1)} : [a-r_{0}] \leq l \leq i \}
     	    	   \, \cup \,
     	    	   \{ V_{p-1 -[2r_{0}+2-a]} \otimes D^{r_{0}+1} \}      	 
     	    	   \]
     	    	  \qquad \qquad \qquad  $ \cup$  $ \{ V_{p-1-[a-2l]} \otimes D^{a-l} : a \leq l \leq [a-i]-1 \}   $.     	    
     	    \end{enumerate} 
     	     \item[$(iv)$] If $r \equiv i+1$, $i+2, \ldots, p-1 \mod p$  or equivalently
     	     $r \equiv [a-j]+1 \mod p$ for some $a+1 \leq j \leq [a-i]$, then the
     	     JH factors of $Q(i)$ are 
     	     \[
     	         \text{ JH factors of } \{ V_{r}^{(l)}/V_{r}^{(l+1)} :  j \leq l \leq i \}
     	         \cup \{ V_{p-1-[a-2l]} \otimes D^{a-l} : a \leq l < j \}.
     	     \]
	           \end{enumerate}
            \noindent As a consequence, $Q(i)$ is not irreducible.
        \end{lemma}
        \begin{proof}
	         By the fourth part of \eqref{interval J}, we have 
	        \[
	           \mathcal{J}(a,i) = \{ a-1, a, \ldots, [a-i]-2 \}^{c} =\{ 0,1, \ldots, a-2 \} \cup \{ [a-i]-1, \ldots, p-1\}.
	        \]
	        By  hypothesis, we have $a< i  < p-1$  and $a \leq p-1+a-i-1=[a-i]-1$. Also, $[a-([a-i]-1)] = i+1 > [a-i]-1$. So      either $[a-i]-1 =a$ or $[a-i]-1$ satisfies the hypothesis of \Cref{JH factor Q i < [a-i]}. Thus, by                     Remark~\ref{compatibility of JH Q(a) and i<[a-i]}, the JH factors of $Q([a-i]-1)$ are determined by Lemma~\ref{JH       factor Q i < [a-i]}.  By Lemma~\ref{JH factor Q i < [a-i]} $(ii)$, we have
	        \[
	        \text{JH factors of } Q([a-i]-1) = \{ V_{p-1-[a-2l]} \otimes D^{a-l} : a \leq l \leq [a-i]-1 \}  \text{ if }  r         \equiv  a,  \ldots, i+1 \mod p.
	        \]
	        \begin{enumerate}
		        \item[$(i)$] By the third part of \eqref{interval I}, we have  $r_{0} \in \mathcal{I}(a,[a-i])$. Thus by \cite[Theorem 4.31 (iii)]{GR19}, we have $\text{ JH factors of } Q(i) $ are 
		        \[
		        \text{JH factors of } \{V_{r}^{(l)}/V_{r}^{(l+1)}:[a-i] \leq l \leq i\} \cup \text{JH factors of } Q([a-i]-1).
		        \]
		     If $[a-i]-1 =a$, then $(i)$ follows. Assume $[a-i]-1 >a$.
		     By \eqref{interval I i < p-1+a-i}, we see that
		     $r_{0} \in \{ 0,1, \ldots, a-2 \} \subset \mathcal{I}(a,a+1)$.
		     So $l =a+1$ in the notation 
		     just above \eqref{recursion JH  i< [a-i]}. Now 
		      $(i)$ follows from \eqref{recursion JH  i< [a-i]} (with $i$ there equal to $[a-i]-1$
		      and $l =a+1$).
		      \item[$(ii)$] The first statement follows from \cite[Theorem 4.31 (i)]{GR19}
		      as $r_0 \not \in \mathcal{J}(a, i)$.Assertion $(a)$ follow from parts $(i)$ $(a)$ (if $p \nmid r$) and $(i)$ $(b)$ (if $p \mid r$) of Lemma~\ref{JH factor Q i < [a-i]} (with $i$ there equal to $[a-i]-1$) and Remark~\ref{compatibility of JH Q(a) and i<[a-i]}.  Assertion $(b)$ follows from the above description of $Q([a-i]-1)$ as $i>[a-i]$.
		   \item[$(iii)$] By \cite[(4.18),(4.20)]{GR19}, we have $r_{0} \in \mathcal{J}(a,i)\smallsetminus \mathcal{I}(a,[a-i])$. 
		   Now  $(iii)$ follows from \cite[Theorem 4.31 (ii)]{GR19} and the above description of  $Q([a-i]-1)$.
		   \item[$(iv)$]  By the third part of \eqref{interval I}, we have  $r_{0} \in \mathcal{I}(a,[a-i])$.  Thus, by               \cite[Theorem 1.9 (iii)]{GR19}, we have  the JH factors of $Q(i)$ equals the
		       \[
		       \text{JH factors of } \{ V_{r}^{(l)}/V_{r}^{(l+1)} :  [a-i] \leq  l \leq i \} \cup
		       \text{JH factors of } Q([a-i]-1).
		       \]
		   If $r \equiv i+1 \mod p$, then $(iv)$ follows from above description of $Q([a-i]-1)$. Else $(iv)$ follows from Lemma~\ref{JH factor Q i < [a-i]} $(iii)$ (applied with $i$ there equal to $[a-i]-1$) and Remark~\ref{compatibility of JH Q(a) and i<[a-i]} as $r \equiv [a-j]+1 \mod p$ for $a+1 \leq j \leq [a-i]-1$.
	\end{enumerate}
         The non-irreducibility of $Q(i)$  follows, noting that $[a-i]-1 \geq a \geq 2$ in part (ii) (a) if $p \nmid r$,  and 
         $[a-r_{0}] = r_{0} +1 \leq i$ in part (iii) (b) since $[a-i]<i$.
\end{proof}
 
 \begin{remark}
 \label{remark compatibility  in equal and greater case}	
 	       In view of  parts $(i)$ and $(ii)$ $(b)$  of  \Cref{JH factor Q i = [a-i]}
 	       we see that the corresponding parts 
 	       of  \Cref{JH factor Q i > [a-i]} are also true when  $a<i=[a-i]$.
 	       Similarly for  parts $(ii)$ $(a)$ 
 	       of  \Cref{JH factor Q i = [a-i]} 
 	       and  \Cref{JH factor Q i > [a-i]}, noting that necessarily  $p \nmid r$
 	       in the former lemma 
 	       as the condition $a < i =[a-i]$  implies 
 	       $a  $ is even, so at least $2$.
 	       Similarly for parts $(iii)$, $(iv)$ and $(v)$ of \Cref{JH factor Q i = [a-i]}
 	       and parts $(iii)$ $(c)$, $(iii)$ $(a)$ and $(iv)$ of  \Cref{JH factor Q i > [a-i]}, respectively. 
 	       Part $(iii)$ $(b)$ of the latter lemma does not occur when $a< i =[a-i]$
 	       as this condition forces $a$ to be even whereas  the
 	        condition $[a-r_0] = r_{0} + 1$ forces $a$ to be odd. 
          Summarizing,   \Cref{JH factor Q i > [a-i]}   is valid even if $a< i =[a-i]$.
 \end{remark}	
 Finally we determine the JH factors of $Q(p-1)$. Recall that the case $p-1 = a$ was
 treated in \Cref{JH factor Q(a)}. 
 
\begin{lemma}
\label{JH factors of Q(p-1)}	
	      Let $p \geq 3$ and $r \equiv a~\mathrm{mod}~(p-1)$
	      with $1 \leq a < p-1$. Let $r \equiv r_{0} \mod p$, with
	      $0 \leq r_{0} \leq p-1$. For $r \geq (p-1)(p+1)+p$, we have
	      \begin{enumerate}
	      	       \item[$(i)$] If $r \equiv 0, 1, \ldots, a-2 \mod p$, then 
	      	       the JH factors of $Q(p-1)$ are 
	      	       \[
	      	            \{ V_{r}^{(l)}/V_{r}^{(l+1)} : a+1 \leq l \leq p-1 \} \cup
	      	            \text{JH factors of } Q(a).
	      	       \]
	      	       \item[$(ii)$] If $r \equiv a-1 \mod p$, then 
	      	       JH factors of $Q(p-1)$ = JH factors of $Q(p-2)$ 
	      	       $\smallsetminus$ JH factors of $V_{r}^{(a)}/V_{r}^{(a+1)}$.
	      	       More precisely,
	      	       \begin{enumerate}
	      	       	       \item[$(a)$] If $p \mid r \equiv a-1 \mod p$, then the
	      	       	            JH factors of $Q(p-1)$ are
	      	       	            \[
	      	       	                 \text{JH factors of } \{  V_{r}^{(l)}/V_{r}^{(l+1)} : 1 \leq l < a \}
	      	       	                 \cup \{ V_{p-1-a} \otimes D^{a} \}.
	      	       	            \]
	      	       	       \item[$(b)$]  if $p \nmid r \equiv a-1 \mod p$, then the
	      	       	       JH factors of $Q(p-1)$ are 
	      	       	       \[
	      	       	           \text{JH factors of } \{  V_{r}^{(l)}/V_{r}^{(l+1)} : 1 < l  < a \}
	      	       	           \cup \{ V_{p-1-a} \otimes D^{a} , V_{p+1-a} \otimes D^{a-1} \}.
	      	       	       \]
	      	       \end{enumerate}
      	         \item[$(iii)$] If $r \equiv a, a+1, \ldots, p-2 \mod p$, then  we have
      	         \begin{enumerate}
      	         	\item[$(a)$] If $[a-r_{0}] < r_{0}+1$, then the JH factors of $Q(p-1)$ 
      	         	 are 
      	         	 \[
      	         	      \text{ JH factors of }
      	         	      \{ V_{r}^{(l)}/V_{r}^{(l+1)} : [a-r_{0}] < l \leq p-1 \} \cup
      	         	      \{ V_{[a-2r_{0}]} \otimes D^{r_{0}} \}.
      	         	 \]
      	         	\item[$(b)$] If $[a-r_{0}] = r_{0}+1$, then the JH factors of $Q(p-1)$ 
      	         	are 
      	         	\[
      	         	\text{ JH factors of } \{ V_{r}^{(l)}/V_{r}^{(l+1)} : [a-r_{0}] \leq l \leq p-1 \} .
      	         	\]
      	         	\item[$(c)$]  If $[a-r_{0}] > r_{0}+1$, then the JH factors of $Q(p-1)$ 
      	         	are 
      	         	\[
      	         	\text{ JH factors of } \{ V_{r}^{(l)}/V_{r}^{(l+1)} : [a-r_{0}] \leq l \leq p-1 \} 
      	         	\cup  \{ V_{p-1-[2r_{0}+2-a]} \otimes D^{r_{0}+1} \}.
      	         	\]      	      
      	         \end{enumerate}
               \item[$(iv)$] If $r \equiv p-1 \mod p$, then the JH factors of $Q(p-1)$
               are
               \[
                    	\text{ JH factors of } \{ V_{r}^{(l)}/V_{r}^{(l+1)} : a+1 \leq l \leq p-1 \} 
                    	\cup \{ V_{a} \}.
               \]
	      \end{enumerate} 
      As a consequence $Q(p-1)$ is irreducible if and only if 
      $r \equiv 1 \mod (p-1)$ and $p \mid r$.
\end{lemma} 
\begin{proof}
	     By the second part of \cite[Lemma 4.4]{GR19}, we have 
	     $X_{r-(p-2)}^{(p-1)}/X_{r-(p-2)}^{(p)} = X_{r-a}^{(p-1)}/X_{r-a}^{(p)} $.
	     Thus, by \cite[Proposition 4.34]{GR19}, we have 
	     $X_{r-(p-2)}^{(p-1)}/X_{r-(p-2)}^{(p)} = V_{a}$ if $r \equiv a-1 \mod p$
	     and 0 otherwise. By  the exact sequence \eqref{exact seq. singular quotients} (with $m =p-1$),
	     the cokernel of $X_{r-(p-2)}^{(p-1)}/X_{r-(p-2)}^{(p)} \hookrightarrow
	     V_{r}^{(p-1)}/V_{r}^{(p)}$ is equal to $V_{p-1-a} \otimes D^{a}$ 
	     if $r \equiv a-1 \mod p$ and
	     $V_{r}^{(p-1)}/V_{r}^{(p)}$ otherwise.
	     Thus, by  the exact sequence \cite[(3.2)]{GR19} and the discussion below it, we have 
	     \begin{align}
	     \label{exact sequence P and Q for i=p-2}
	     	     0 \rightarrow W' \rightarrow P(p-1) \rightarrow Q(p-2) \rightarrow 0.
	     \end{align}
	     where $W' = V_{p-1-a} \otimes D^{a}$ if $r \equiv a-1 \mod p$ and
	     $V_{r}^{(p-1)}/V_{r}^{(p)}$ otherwise. If $a=p-2$, then $Q(p-2)$ is determined by 
	     Lemma~\ref{JH factor Q(a)}. If $a< p-2$, then $[a-(p-2)] = a+1 \leq p-2$
	     so by Remark~\ref{remark compatibility  in equal and greater case},
	      the JH factors of $Q(p-2)$ can be obtained 
	     from Lemma~\ref{JH factor Q i > [a-i]}. 
	     \begin{enumerate}
	     	\item[$(i)$] By \cite[Theorem 4.39 (i)]{GR19}, we have 
	     	$Q(p-1) \cong P(p-1)$. By the exact sequence
	     	\eqref{exact sequence P and Q for i=p-2}, we have JH factors of 
	     	$Q(p-1)  = $ JH factors of $Q(p-2)$ $\cup $ JH factors of 
	     	$V_{r}^{(p-1)}/V_{r}^{(p)}$. Now $(i)$ follows from 
	     	Lemma~\ref{JH factor Q i > [a-i]} $(i)$ (applied with $i$ there equal to $p-2$)
	     	in the case 	$a< p-2$
	     	 and 
	     	 is obvious in the case  $a=p-2$.
	     	\item[$(ii)$]  By \cite[Theorem 4.39 (ii)]{GR19}, we have 
	     	JH factors of $Q(p-1) =$ JH factors of $P(p-1) \smallsetminus 
	     	 \text{ JH factors of } V_{r}^{(a)}/V_{r}^{(a+1)}$. Further, by  the exact
	     	 sequence \eqref{exact sequence P and Q for i=p-2}, we have
	     	 $P(p-1)$ is an extension of $Q(p-2)$ by $V_{p-1-a} \otimes D^{a-1}$.
	     	 If $a =p-2$, then $Q(p-2) = Q(a)$.
	     	 If $a< p-2$, then by Lemma~\ref{JH factor Q i > [a-i]} $(ii)$, we have 
	     	 $Q(p-2) \cong Q(a)$. Thus,  the JH factors of $Q(p-1)$ are 
	     	 \[
	     	     \{ V_{p-1-a} \otimes D^{a-1} \} \cup \text{JH factors of } Q(a) \, \smallsetminus 
	     	      \text{JH factors of } V_{r}^{(a)}/V_{r}^{(a+1)}.
	     	 \] 
	     	 Now $(ii)$ $(a)$ and $(ii)$ $(b)$ follow 
	     	 from parts $(ii)$ $(b)$ and $(ii)$ $(a)$ of Lemma~\ref{JH factor Q(a)}, respectively.
	     	\item[$(iii)$] By \cite[Theorem 4.39 (iii)]{GR19} and  the exact sequence \eqref{exact sequence P and Q for i=p-2}, we have the JH factors of $Q(p-1)$ equals
	     	\[
	     	    \text{JH factors of } Q(p-2) \cup \text{JH factors of } 
	     	    V_{r}^{(p-1)}/ V_{r}^{(p)} \smallsetminus \{ V_{a} \}.
	     	\]
	     	Part $(iii)$ follows from Lemma~\ref{JH factor Q i > [a-i]} $(iii)$ (applied with $i=p-2$) if $a< p-2$ 
                noting that $p-1-[a-2a] = a$. If $a=p-2$, then $r \equiv p-2 \mod p$, so $[a-r_{0}] = p-1 = p-2+1 = r_{0}+1$, and we are in  part $(iii)$ $(b)$. This case follows from Lemma~\ref{JH factor Q(a)} $(iii)$ (applied with $a=p-2$).
	     	
	     	\item[$(iv)$] As in $(i)$ the JH factors of $Q(p-1)=$ JH factors of $Q(p-2)$ $\cup $ JH factors of  
                $V_{r}^{(p-1)}/V_{r}^{(p)}$.  If $a=p-2$, then $(iv)$ follows from Lemma~\ref{JH factor Q(a)} $(iii)$. If $a< p-2$, then $(iv)$ follows from Lemma~\ref{JH factor Q i > [a-i]} $(iv)$ (with $i=p-2$) noting that $r \equiv p-1 = [a-(a+1)] +1 \mod p$.
	     \end{enumerate}
     If $ r \equiv a-1 \mod p$ and  $p \mid r $, then $Q(p-1)$ is 
     irreducible by the assertion $(ii)$ $(a)$. The converse follows  from
     the above assertions noting that $Q(p-1)$ contains at least two 
     JH factors if $p \nmid r$ or $ r \not \equiv a-1 \mod p$.
    \end{proof}

    \section{Good congruence classes}\label{sec: good}

Let $k\geq 2$ be an integer and 
        $v(a_p) \in (i,i+1)$ for some integer $0 \leq i < p-1$. In this section, we determine the shape of $\bar{V}_{k,a_p}$ when the congruence class of $r=k-2$ lies in certain `good' congruence classes modulo $p$. Recall that there is a surjection 
	\begin{align}\label{surjection from Q(i)}
		\ind_{KZ}^{G}  Q(i) \twoheadrightarrow \bar{\Theta}_{k,a_{p}},
	\end{align}
	where
	\[
	Q(i) = \frac{V_{r}}{X_{r-i} + V_{r}^{(i+1)}}.
	\] 
    We have chosen to work with these congruence classes of $r$ mod $p$ first, since the structure of $Q(i)$ in these classes is a bit simpler. 
        
        Let us say that $r$ is a good congruence class mod $p$ if for each of the sub-quotients $V_r^{(m)}/ V_r^{(m+1)}$ appearing in $V_r/V_r^{(i+1)}$, the socle of $V_r^{(m)}/ V_r^{(m+1)}$ does not contribute to the final quotient $Q(i)$. That is, if $V_r^{(m)}/ V_r^{(m+1)}$ is non-split, then the socle of $V_r^{(m)}/ V_r^{(m+1)}$ does not contribute to $Q(i)$, and if it is split, then neither of the JH factors of $V_r^{(m)}/ V_r^{(m+1)}$ contribute to $Q(i)$.

        A case by case inspection of the results in the previous chapter shows that these congruence classes can be described in the following elementary fashion, which for all practical purposes, can then be taken as the definition of the good classes.

        \begin{definition}
        Let $r \equiv a $ mod $(p-1)$ with $1\leq a \leq p-1$.
        Let 
        \begin{eqnarray*}
           b = \begin{cases}
               a & \text{ if } 0 \leq i < a \\
               p-1+a & \text{ if }a \leq i \leq p-2.
           \end{cases} 
        \end{eqnarray*}
        Then we say $r$ is in a {\it good congruence class mod $p$} if
        \begin{eqnarray*}
            r \not\equiv \begin{cases}
                b-i+1, b-i, \ldots, b & \text{ if } b > 2i \\
                b-i-1, b-i, \cdots, b  & \text{ if } b \leq  2i 
            \end{cases} \quad \mod p
        \end{eqnarray*} 
        \end{definition}


        We next prove a technical lemma, which guarantees the existence of certain constants $\beta_{l}$ needed in the proof of the next theorem.
    	\begin{lemma}\label{conj}
		Let $ p \geq 3 $,  $ m \geq 0 $  and $  2p \leq r \equiv a \mod (p-1)$ with $ m+1 \leq a-m \leq p-1 $. If $ r \not \equiv  a-m,  \ldots, a-1, a  \mod p$, then one can choose $ p $-adic integers $ \beta_{0},\beta_{1},
		\ldots, \beta_{m} \in \mathbb{Z}_{p}$, such that 
		\begin{enumerate}
			\item[$(i)$] $   
			\sum\limits_{\substack{0 < j < r-m\\ j \equiv r-m ~\mathrm{mod}~ (p-1)}}  \left( \beta_m \binom{r-m}{j} + p \sum\limits_{l=0}^{m-1} \beta_l \binom{r-l}{j} \right) \binom{j}{n}  \equiv 0  \mod p^2$, for all $ 0 \leq n < m $,
			\item[$(ii)$] $ 
			\sum\limits_{\substack{0 < j < r-m \\ j \equiv r-m ~\mathrm{mod}~ (p-1)}} \left( \beta_m \binom{r-m}{j} + p \sum\limits_{l=0}^{m-1} \beta_l \binom{r-l}{j} \right) \binom{j}{m}  \equiv p \mod p^2$.
		\end{enumerate}
	\end{lemma}
	\begin{proof}
	For $ 0 \leq n \leq m $, we have $r-m -n \equiv a-m -n ~\mathrm{mod}~ (p-1)$ and $ 1 \leq m+1 -n  \leq a-m-n\leq p-1-n \leq p-1 $. By \cite[Lemma 2.5]{BG15}, we have
		\begin{align}\label{beta first term sum}
		\begin{split}
			\sum_{\substack{0 < j < r-m \\ j \equiv r-m ~\mathrm{mod}~ (p-1)}}\binom{r-m}{j} \binom{j}{n} 
            & = \binom{r-m}{n} \sum_{\substack{0 < j < r-m \\ j \equiv r-m ~\mathrm{mod}~ (p-1)}} \binom{r-m-n}{j-n} \\
			& = \binom{r-m}{n} \sum_{\substack{0 < j < r-m-n \\ j \equiv r-m-n ~\mathrm{mod}~ (p-1)}}\binom{r-m-n}{j} \\
			& \equiv p \binom{r-m}{n}  \frac{a-r}{a-m-n} \mod p^{2},
		\end{split}
		\end{align}
		where the second step is obvious if $ n= 0 $ and for $ n > 1 $ it follows from the observation that the smallest positive integer congruent to $ r-m $ modulo $ (p-1)  $ is $ a-m > m \geq n $ so $a - m -n > 0$.
		
        Next we observe that for $ 1 \leq l \leq  m $ and $ 0 \leq n \leq m $, since $a > 2m $, by \cite[Lemma 2.15]{GR19}, we have
		\begin{align}\label{beta second part sum1}
		\begin{split}
			     \sum_{\substack{0 < j < r-m \\ j \equiv r-m ~\mathrm{mod}~ (p-1)}} \binom{r-l}{j} \binom{j}{n} 
                & = \sum_{\substack{0 < j \leq r-m \\ j \equiv r-m ~\mathrm{mod}~ (p-1)}} \binom{r-l}{j} \binom{j}{n} - \binom{r-l}{r-m} \binom{r-m}{n} \\
				& = \sum_{\substack{0 \leq j \leq r-m \\ j \equiv r-m ~\mathrm{mod}~ (p-1)}} \binom{r-l}{j} \binom{j}{n} - \binom{r-l}{r-m} \binom{r-m}{n} \\
				& \qquad \qquad \qquad \qquad  \qquad \qquad \qquad - \binom{0}{n}  \delta_{[r-m],p-1}  \\
				& \equiv \binom{r-l}{m} \binom{[a-l-n]}{a-m-n}  + \binom{r-l}{n} \delta_{[r-m-n],p-1}  \\
				& \qquad \quad -\binom{r-l}{r-m} \binom{r-m}{n}
				- \delta_{0,n}  ~ \delta_{[r-m],p-1} \mod p.
			\end{split}
		\end{align}
		Since  $ 0 \leq n \leq m < a-m \leq p-1 $, it follows that $ [r-m-n] = [[r-m]-n] = p-1 $ if and only if  $ [r-m] = p-1 $ and $ n = 0 $. Thus the $\binom{r-l}{n} \delta_{[r-m-n],p-1}$ and $-\delta_{0,n}  ~ \delta_{[r-m],p-1} $ in \eqref{beta second part sum1} cancel out. Also, if $ a-l-n \leq p-1 $, then $ \binom{[a-l-n]}{a-m-n} = \binom{a-l-n}{a-m-n} $. If $ a-l-n > p-1 $, then since $ l \leq m < p-1 $ and $ a-l-n\leq a \leq p-1+m \leq 2p-2$. Thus  $ \binom{[a-l-n]}{a-m-n}  = \binom{a-l-n-p+1}{a-m-n}  = 0$, and by Lucas' theorem (Lemma~\ref{Lucas}),   $ \binom{a-l-n}{a-m-n} \equiv  \binom{a-l-n-p}{a-m-n} \equiv 0 \mod p$. Hence \eqref{beta second part sum1} reduces to 
		\begin{align}
			\label{beta second part sum2}
			\sum_{\substack{0 < j < r-m \\ j \equiv r-m ~\mathrm{mod}~  (p-1)}}
			\binom{r-l}{j} \binom{j}{n}  \equiv 
			\binom{r-l}{n} \binom{a-l-n}{a-m-n} -
			\binom{r-l}{r-m} \binom{r-m}{n} \mod p.
		\end{align}
		It follows from  \eqref{beta first term sum} and 
		\eqref{beta second part sum2} that for $ 0 \leq n \leq m$ we have
        \begin{multline*}
				\sum\limits_{\substack{0 < j < r-m \\ j \equiv r-m ~\mathrm{mod}~ (p-1)}}   \left( \beta_m \binom{r-m}{j} + p \sum\limits_{l=0}^{m-1} \beta_l \binom{r-l}{j} \right) \binom{j}{n}  \\
               \quad \equiv p\beta_{m}
				\binom{r-m}{n}  \frac{a-r}{a-m-n}  + p\sum_{l=0}^{m-1} 
				\beta_{l} \Bigg\lbrace \binom{r-l}{n}  \binom{a-l-n}{a-m-n} 
				 - \binom{r-l}{r-m} \binom{r-m}{n} \Bigg\rbrace 
			\mod p^2.
        \end{multline*}
		Observe that to prove the lemma, it is enough to solve the following  
		equations  
		\begin{align}
			\label{linear equations alpha's}
			\beta_{m}
			\binom{r-m}{n}  \frac{a-r}{a-m-n}  + \sum_{l=0}^{m-1} 
			\beta_{l} \Bigg\lbrace \binom{r-l}{n}  \binom{a-l-n}{a-m-n} 
			- \binom{r-l}{r-m} \binom{r-m}{n} \Bigg\rbrace = \delta_{m,n}.
		\end{align}
        for $ 0 \leq n \leq m $.
		Putting this in matrix form, we must solve
		\[
		A ( \beta_{0}, \beta_{1}, \ldots, \beta_{m-1},\beta_{m})^{t} =
		(0,0, \ldots, 0, 1)^{t},
		\]
		where  $ A = (B \vert v^{t}) $ and
		\begin{align*}		
			B &= \left(\binom{r-l}{n}  \binom{a-l-n}{a-m-n} 
			- \binom{r-l}{r-m} \binom{r-m}{n}\right)_{\substack{0 \leq n \leq m \\ 0 \leq l \leq m-1}} \\
            v &= \left(  \binom{r-m}{0} \frac{a-r}{a-m}, \ldots, \binom{r-m}{n} \frac{a-r}{a-m-n},
			\ldots, \binom{r-m}{m} \frac{a-r}{a-m-m} \right).
		\end{align*}
		To show that the system of linear equations \eqref{linear equations alpha's} has a solution, it is enough to show that $ A $ is invertible if  $r \not \equiv a-m,  \ldots, a-1, a \mod p $. To do this, we
		compute the determinant of $ A $.
		
		Pulling out the factor of $\frac{1}{(r-m-n)! n!}$ from the $ n^{\rm{th}} $-row and 
		$ \frac{(r-l) !}{(m-l)!} $ from the $ l^{\rm{th}}$-column,  we see that
		\begin{align*}
			\mathrm{det}(A) = \left( \prod_{l=0}^{m}   \frac{(r-l) !}{(m-l) !}
			\right) \left( \prod_{n=0}^{m}  \frac{1}{(r-m-n)! n!} \right)
			\det (C),
		\end{align*}
		where $ C = (D \vert w^{t}) $ and
		\begin{align*}
			D & = \left(  \frac{(r-m-n)! (a-l-n)! }{ (r-l-n)! (a-m-n)! } -1  \right)_{\substack{0  \leq n \leq m  
					\\ 0 \leq l \leq m-1 }}, \\
			w & = \left( \frac{a-r}{a-m}, \ldots,  \frac{a-r}{a-m-n}, \ldots, 
			\frac{a-r}{a-2m}\right).
		\end{align*}
		Performing the column operations subtracting  $l^{\rm{th}}$-column from
		$ (l+1)^{\rm{th}}$-column for $ l \leq m-2 $, that is, subtracting the first  column from 
		the second   column,  the second column from the third column, 
		and so forth, and the third-to-last column from the second-to-last column, and then pulling out a sign from each of the last two columns, we see that
		\begin{align*}
			\mathrm{det}(A) = \left( \prod_{l=0}^{m}   \frac{(r-l) !}{(m-l) !}
			\right) \left( \prod_{n=0}^{m}  \frac{1}{(r-m-n)! n!} \right) \cdot \det(E),
		\end{align*}
		where   
		\begin{align*}
			E =  \left( (r-a) \frac{(r-m-n)!}{(r-l-n)!} 
			\frac{(a-l-1-n)!}{(a-m-n)!}\right)_{0	\leq n,l \leq m}.
		\end{align*}
		Multiplying the $n^{\mathrm{th}}$-row of $E$ by  $\frac{1}{(r-m-n)! n!}$ and the
		$ l^{\mathrm{th}} $-column of $ E $ by $ (r-l) ! $, we see that
		\begin{align*}
		    \mathrm{det}(A) 
            &=  \left( \prod_{l=0}^{m} \frac{1}{{(m-l) !}} \right) \times \det \left( \frac{(r-a)}{n!} \frac{(r-l)!}{(r-l-n)!} \frac{(a-l-n-1)!}{(a-m-n)!} \right)_{0 \leq n,l \leq  m} \\
            &= \left( \prod_{l=0}^{m} \frac{1}{{(m-l) !}} \right) \times \det \left( (r-a) \frac{(m-n)!(a-1-l-m)!}{(a-m-n)!}\binom{r-l}{n} \binom{a-1-l-n}{m-n} \right)_{0 \leq n,l \leq  m} \\
            &= (r-a)^{m+1} \left( \prod_{l=0}^{m} \frac{(a-1-l-m)!}{{(m-l) !}} \right) \times \left( \prod_{n=0}^{m} \frac{(m-n)!}{(a-m-n)!} \right) \\
            & \qquad \qquad \qquad \times \det \left( \binom{r-l}{n} \binom{a-1-l-n}{n-m} \right)_{0 \leq n,l \leq  m}.
		\end{align*}
        The last quantity is non-zero modulo $p$. Indeed, by \cite[Proposition 2.16 (ii)]{GR19}, the above determinant is non-zero modulo $p$ if $r \not \equiv a-m, \ldots, a-2,a-1 \mod p$ and the fudge factor 
        in front doesn't vanish modulo $p$
        since $r \not \equiv a \mod p$.
	\end{proof}

    The following important theorem gives conditions  which will allow us to kill the JH factors 
    in $V_r^{(m)}/V_r^{(m+1)}$
    for $m \geq 0$.
    It will be used both in 
    the `good' and `bad' cases
    for $m$ small.

    \begin{theorem}\label{Elimination i < a and not in interval}
		Let $p\geq 3$ and $ 0 \leq m \leq p-1 $ be an integer. Let  $v(a_{p}) > 0 $ and $v(a_p)\not\in \Z$. Let $r \geq  (\lfloor v(a_{p}) \rfloor)(p+1)+p$ and $r \equiv a \mod (p-1) $  with $ m+1 \leq a-m \leq p-1  $.
		If the following conditions hold 
		\begin{itemize}
			\item[$(i)$] $ m +1 <  v(a_{p})$ and $ \lfloor v(a_{p}) \rfloor  \leq a-m-1 $ 
            \item[$(ii)$] if $ \lfloor v(a_p) \rfloor  = a-m-1 $, then $p \mid \binom{r-m}{\lfloor v(a_{p}) \rfloor +1 }$
			\item[$(iii)$]  
			$ r \not \equiv  a-m, a-m+1, \ldots, a-1,a \mod p$, 
		\end{itemize}
		then the image of  $\ind_{KZ}^{G}  V_{r}^{(m)} $ in $ \bar{\Theta}_{k,a_{p}} $ 
		is the same as the image of  $\ind_{KZ}^{G} V_{r}^{(m+1)} $  in $ \bar{\Theta}_{k,a_{p}} $.
	\end{theorem}
	
	\begin{proof}
       Assume first that $r \geq (\lfloor v(a_{p}) \rfloor +2)(p-1)+a $ or equivalently $r > (\lfloor v(a_{p}) \rfloor +1)(p-1)+a $.
    Recall from \eqref{polynomial Fr,m,l} that
        \begin{align*}
            F_{r,m+2,m} =  X^m Y^{r-(m+2)(p-1)-m} (X^{p-1}-Y^{p-1})^{m+2}=\theta^{m+2} X^{-2}Y^{r-(m+2)(p+1)+2}.
        \end{align*}
        Let $F_{r,m}:= (-1)^m F_{r,m+2,m}$. As $\lfloor v(a_p) \rfloor \geq m+1 \geq 1$, we have $r > (\lfloor v(a_{p}) \rfloor +1)(p-1)+a \geq 2p-1$. Moreover, hypothesis $(iii)$ holds. Thus there exist $ \beta_{0}$, $\beta_{1}, \ldots$, $\beta_{m}  \in \Z_{p}$ as in Lemma~\ref{conj}. Consider the following function 
        \begin{align}\label{f2r}
	            f_{2} =\sum_{\lambda \in \mathbb{F}_{p}}^{} \left[ g_{2, p[\lambda]}^{0} , \frac{\beta_{m}}{p^{m+1}} F_{r,m}  + \left( \sum_{l=0}^{m-1} \frac{\beta_{l}}{p^{l}} [\lambda]^{p-1-(m-l)} (X^{-1}Y)^{m-l} F_{r,m} \right) \right].        
        \end{align}
        We now compute the action of the Hecke operator $T$ on the above function. Since  $ (X-Y)^{l+2} \mid (X^{-1}Y)^{m-l} F_{r,m} $ for $ 0 \leq l \leq m $, it follows from Lemma~\ref{theta and T plus} that $ T^{+} f_{2} \equiv 0  \mod p$. Here we used $m+1 \leq \lfloor v(a_p) \rfloor$ and  $r\geq \lfloor v(a_{p}) \rfloor (p+1)+p \geq (m+1) (p+1)+p $ to conclude that the coefficients of $X^r,\ldots,X^{r-m-1}Y^{m+1}$ are zero. Clearly $a_pf_2 \equiv 0 \mod p$ since $m+1 < v(a_p)$ by hypothesis $(i)$. 
        Also, it follows from formula \eqref{T minus formula} for $ T^{-} $ that
        \begin{align*}
	            T^{-} f_{2} & \equiv \sum_{\lambda \in \mathbb{F}_{p}}^{} T^{-}  \left[ g_{2,p[\lambda]}^{0} , \frac{\beta_{m}}{p^{m+1}}  X^{m} Y^{r-m} + \sum_{l=0}^{m-1} \frac{\beta_{l}}{p^{l}} [\lambda]^{p-1-(m-l)} X^{l} Y^{r-l}  \right] \mod p\\
	            & \equiv \sum_{\lambda \in \mathbb{F}_{p}}^{}  \Bigg[ g_{1,0}^{0}, \frac{\beta_{m}}{p} \sum_{j = 0}^{r-m} \binom{r-m}{j} [\lambda]^{r-m-j} X^{r-j}Y^{j} + \sum_{l=0}^{m-1} \beta_{l} \sum_{j = 0}^{r-l} \binom{r-l}{j} [\lambda]^{p-1+r-m-j}  X^{r-j}Y^{j}  \Bigg] \\
	            & \equiv \left[ g_{1,0}^{0}, \frac{p-1}{p} \sum\limits_{\substack{0 \leq  j < r - m \\  j \equiv r - m ~\mathrm{mod}~  (p-1)}}^{} \left(\beta_m \binom{r-m}{j} + \sum_{l=0}^{m-1} p \beta_l \binom{r-l}{j}  \right)X^{r-j} Y^{j} \right]  \\
                 & \qquad \qquad \qquad + \left[ g_{1,0}^{0},\left( \beta_{m} + \sum\limits_{l=0}^{m-1} (p-1) \beta_{l} \binom{r-l}{r-m} \right)X^{m} Y^{r-m} \right]  \mod p.
        \end{align*}
        The last term dies as it is integral and $X^mY^{r-m} \in X_{r-i}$. Thus we obtain
        \begin{equation}\label{T f2r}
             \begin{aligned}
		        T^{-} f_{2}  \equiv  & \left[ g_{1,0}^{0}, \frac{p-1}{p} \sum\limits_{\substack{0 \leq  j < r - m \\  j \equiv r - m ~\mathrm{mod}~  (p-1)}}^{} \left(\beta_m \binom{r-m}{j} + \sum_{l=0}^{m-1} p \beta_l \binom{r-l}{j}  \right)X^{r-j} Y^{j} \right] \mod p. 
	        \end{aligned}
        \end{equation}
        
        Note that the $ X^{r} $ term occurs in \eqref{T f2r} only if  
		$  r - m \equiv p-1 $ mod $(p-1)$. Let
		\begin{align}
			f_{0} =
			\begin{cases}
				\left[ \mathrm{Id}, \beta_{m}\left(\frac{1-p}{p}\right)
				(X^{r} - X^{r-(p-1)}Y^{p-1})\right], & \mathrm{if} ~ 
				a - m= p-1, \vspace*{2mm} \\
				0, &\mathrm{otherwise}.                
			\end{cases}
		\end{align}
		Since $ r > p $ and $ v(a_{p})> 1$, $ T^{-}f_{0}  $ and $ a_{p} f_{0} $ 
		vanish modulo $ p $. Also, it can be checked that
		\begin{align*}
			T^{+} f_{0} = \begin{cases}
				\left[   g_{1,0}^{0}, \beta_{m} \left( \frac{1-p}{p}  \right)X^{r} \right] 
				\mod \ind_{KZ}^{G} \langle X^{r-1} Y , p \rangle, & \mathrm{ if }
				~ a - m =p-1,  \\
				0, & \mathrm{otherwise}. 
			\end{cases}
		\end{align*}
        Thus $T^{-}f_2+T^{+}f_0$ is as in \eqref{T f2r} but  without the   $j=0$ term.
        
        We now wish to introduce a function $f_1$ in radius $1$ such that $T^{-} f_2- a_p f_1 + T^{+}f_0$ is integral  and $T^{+}f_1$ vanishes modulo $p$. 
             Note that for all $n$, we have
       \begin{align*}
            \sum_{\substack{0 < j < r-m \\ j \equiv r-m ~\mathrm{mod}~  (p-1)}}^{} \left( \beta_m \binom{r-m}{j} + \sum_{l=0}^{m-1} p \beta_l \binom{r-l}{j} \right) \binom{j}{n} 
             & \equiv \beta_{m} \sum_{\substack{0 < j < r-m \\ j \equiv r-m ~\mathrm{mod}~  (p-1)}}^{} \binom{r-m}{j} \binom{j}{n} \mod p.
       \end{align*}
       If  $0 \leq n < [r-m]$, then by \cite[Lemma 2.5]{BG15} we have 
        \begin{align}\label{beta sum 0}
        \begin{split}
            \sum_{\substack{0 < j < r-m \\ j \equiv r-m ~\mathrm{mod}~  (p-1)}}^{} \binom{r-m}{j} \binom{j}{n} &\equiv   \binom{r-m}{n} \sum_{\substack{0 <  j < r-m \\ j \equiv r-m ~\mathrm{mod}~  (p-1)}}^{}  \binom{r-m-n}{j-n} \\
            &= \binom{r-m}{n}\sum_{\substack{0 < j' < r-m-n \\ j' \equiv r-m-n ~\mathrm{mod}~ (p-1)}}^{}  \binom{r-m-n}{j'} \equiv 0 \mod p.
        \end{split}    
        \end{align}
        This also holds if $p \mid \binom{r-m}{n}$. From $i+1 \leq a-m$ and $(ii)$, it follows that
        \begin{align}
		        \sum_{\substack{0 \leq j < r-m \\ j \equiv r-m ~\mathrm{mod}~ (p-1)}}^{} \left( \beta_m \binom{r-m}{j} + \sum_{l=0}^{m-1} p \beta_l \binom{r-l}{j} \right) \binom{j}{n}  \equiv  0 \mod  p~\text{for } n=0,\ldots,i+1.
        \end{align}
           Take $c=a-m$, $t=1$, $m'=\lfloor v(a_p) \rfloor+1$, $\nu_0,\ldots, \nu_{m'}=0$, $n=\lfloor v(a_p)\rfloor+2$, $k=\frac{r-a}{p-1}-1$ and 
         \begin{equation*}
             \gamma_{c+j(p-1)} = \beta_m \binom{r-m}{a-m+j(p-1)} + \sum_{l=0}^{m-1} p \beta_l \binom{r-l}{a-m+j(p-1)}  \quad \text{for } j=0,\ldots,\frac{r-a}{p-1}-1
         \end{equation*}
         in Lemma~\ref{lem:choice of beta}. Note that the condition $r\geq \lfloor (v(a_p) \rfloor+2) (p-1)+a$ is equivalent to the condition $k\geq m'$ in Lemma~\ref{lem:choice of beta}. By the lemma, 
         there exist $\alpha_j$  for $0 <  j < r - m \text{ with } j \equiv r - m ~\mathrm{mod}~(p-1)$ such that 
         \begin{enumerate}
             \item $\alpha_j \equiv  \beta_m \binom{r-m}{j} + \sum\limits_{l=0}^{m-1} p \beta_l \binom{r-l}{j}   \mod p \text{ for all } 0 <  j < r - m \text{ with } j \equiv r - m ~\mathrm{mod}~(p-1)$
             \item $\sum_j \binom{j}{n}\alpha_j \equiv 0 ~\mathrm{mod}~p^{\lfloor v(a_p)\rfloor+2}$ for $n=0,\ldots,\lfloor v(a_p) \rfloor+1$.
         \end{enumerate}

        Define
		\begin{align*}
			f_{1} = \left[ g_{1,0}^{0},   \frac{p-1}{p a_{p}} 
			\sum\limits_{\substack{0 <  j < r - m \\  j \equiv r - m \mod (p-1)}}^{}
			\alpha_{j} X^{r-j} Y^{j} \right]. 
		\end{align*}
        Note that the valuation of the denominator of $f_{1}$ is  
        $v(a_p)+1$ and the smallest power of $X$ appearing in $f_{1}$ is at least $m+p-1$. Since $v(a_p)  + 1 < p-1 + m $ in all cases except $a=p-1$, $\lfloor v(a_p) \rfloor =p-2$, $m=0$, it follows that $T^{-} f_{1}$ vanishes generically. 
        We now check that $T^{-}f_{1}$ vanishes modulo $p$ even in the case $a=p-1, \lfloor v(a_p) \rfloor =p-2, m=0$. Since $v(a_p) < p-1$, it suffices to show 
        $v(\alpha_{r-m-(p-1)}) +m+p-1  > v(a_p)+1$. This follows if $v(\alpha_{r-(p-1)})=v(\alpha_{r-m-(p-1)}) \geq 1$. 
        But $\alpha_{r-(p-1)} \equiv  \beta_0 \binom{r}{r-(p-1)} \mod p$. But, from hypothesis $(iii)$, we have $r \not \equiv p-1 \mod p$ as $a=p-1$. 
        Thus, by Lucas' theorem,
        $ p \mid \binom{r}{p-1}$, as desired.
        
    By formula \eqref{T plus formula} and property (2) of the $\alpha_j$ above, we see that the $\lambda \neq 0$ terms in $ T^{+}f_{1} $ vanish modulo $ p $. For the $\lambda=0$ term, note that the smallest power of $Y$ appearing in $f_1$ is $a-m$. 
    If $\lfloor v(a_p) 
    \rfloor + 1 <a-m$, then  this term dies. If   $\lfloor v(a_p) 
    \rfloor + 1 =a-m$, by property (1) of the $\alpha_j$ above and hypothesis $(ii)$, we have $\alpha_{a-m} \equiv \beta_m \binom{r-m}{a-m} \equiv 0 \mod p$, so again this term dies.

		 Therefore
		\begin{align*}
			(T-a_{p})(f_{2}+f_{1}+f_{0}) & \equiv  T^{-} f_{2} 
			- a_{p} f_{1} + T^{+} f_{0}  \equiv  [ g_{1,0}^{0},   
			F(X,Y)  ] \mod \ind_{KZ}^{G} \langle X^{r-1} Y , p \rangle, 
		\end{align*}   
		where 
        \begin{align*}
           F(X,Y)= \frac{p-1}{p } 
			\sum\limits_{\substack{0 <  j < r - m \\  j \equiv r - m \mod (p-1)}}^{}
			\left( \beta_m \binom{r-m}{j} + \sum\limits_{l=0}^{m-1} p \beta_l \binom{r-l}{j} - \alpha_{j}\right) X^{r-j} Y^{j}.
        \end{align*}
        By property (1) of the $\alpha_j$ above, $F(X,Y)$ is integral.
        The smallest exponent of $ Y $ in the above expression is $ a-m > m $, so the coefficients
		of $ X^{r}, X^{r-1}Y, \ldots, X^{r-m}Y^{m} $ are zero. The smallest exponent of 
		$ X $ in the above expression is $ r-m-(p-1) > m$
		so the coefficient of $ Y^{r}, XY^{r-1}, \ldots, X^{m}Y^{r-m} $ are zero. 
		So the polynomial satisfies the first condition of \cite[Lemma 2.8]{GR19}. By Lemma~\ref{conj} $(i)$ and property (2) of the $\alpha_j$ above, $\theta^m  \mid \overline{F(X,Y)}$. Noting that the coefficient of $ X^{{m}}Y^{r-{m}} $ and $ X^{r-{m}} Y^{{m}} $
		are zero, by  \cite[Lemma 2.12]{GR19}, we have 
		\begin{multline*}
			F(X,Y) \equiv \frac{p-1}{p} \left(\sum\limits_{\substack{0 <  j < r - m \\  j \equiv r - m \mod (p-1)}}^{}
			\left( \beta_m \binom{r-m}{j} + \sum\limits_{l=0}^{m-1} p \beta_l \binom{r-l}{j} - \alpha_{j}\right) \binom{j}{{m}} \right)
			\\ \times \theta^{m} X^{r -m(p+1)-a+2{m}} Y^{a-2{m}}   \mod V_{r}^{({m}+1)}.
		\end{multline*}
        By Lemma~\ref{conj} $(ii)$ and property (2) of the $\alpha_j$ above, $\overline{F(X,Y)} \equiv - \theta^{m} X^{r -m(p+1)-a+2{m}} Y^{a-2{m}}   \mod V_{r}^{({m}+1)} $.
		Thus by \cite[Lemma 2.12]{GR19}, Lemma~\ref{Glover-Brueil map image} $(ii)$ and the assumption  $ a - {m} >{m} $
		we see that $(T-a_{p}) (f_{2}+f_{1}+f_{0})$ mod $p$
		generates  $\ind_{KZ}^{G} \left( \frac{V_{r}^{({m})}}{V_{r}^{({m}+1)}} \right)$.
        But any such function dies
        in $\bar\Theta_{k,a_p}$.
        This finishes the proof if $r> (\lfloor v(a_p) \rfloor+1)(p-1)+a$.
    
     We now show that the bound can be improved to $r\geq(\lfloor v(a_p)\rfloor )(p+1)+p$. Note that 
    \begin{align}\label{eq: bound inequality good}
        \frac{r-a}{p-1} \geq \frac{(\lfloor v(a_p)\rfloor )(p+1)+p-a}{p-1} = \lfloor v(a_p)\rfloor +1 + \frac{2\lfloor v(a_p)\rfloor+1 -a}{p-1}.
    \end{align}
    If $a\leq2\lfloor v(a_p)\rfloor$, then $\frac{r-a}{p-1}>\lfloor v(a_p)\rfloor +1$ and the previous bound on $r$ holds. Assume $a> 2\lfloor v(a_p)\rfloor$. By hypothesis $(i)$, we have $a\leq m+p-1\leq \lfloor v(a_p)\rfloor+p-2  $. Thus $-(p-2)\leq \lfloor v(a_p)\rfloor+1-(p-2) \leq  2\lfloor v(a_p)\rfloor+1-a    $.  Thus by \eqref{eq: bound inequality good}, we have
    \(
         \frac{r-a}{p-1} > \lfloor v(a_p) \rfloor
    \).
    Thus the only case which is not covered above when $a>2\lfloor v(a_p)\rfloor$ is $r= (\lfloor v(a_p)\rfloor+1)(p-1)+a$. In this case, the proof still works.  Indeed, $T^{-}f_2+T^+f_0$  is already integral since every $0 <  j < r - m  $ with $ j \equiv r - m \mod (p-1)$ can be expressed as $a-m+k'(p-1)$ for some $0 \leq k' \leq \lfloor v(a_p) \rfloor$,  so by hypothesis $(i)$ and Lucas' theorem, for $0 \leq k' \leq \lfloor v(a_p) \rfloor$, we have 
    \begin{align*}
 	   \binom{r-m}{a-m+k'(p-1)} = \binom{(\lfloor v(a_p) \rfloor+1)p+a-\lfloor v(a_p) \rfloor-1-m}{k'p+a-m-k'} \equiv 0 \mod p.
    \end{align*}
     So we may take $f_1=0$ and proceed as before.
    \end{proof}
    
The next two corollaries allow us to show that all but one of the JH factors in $Q(i)$ die in $\bar{\Theta}_{k,a_p}$ for $0 \leq i < p-1$, for 
    good congruence classes of $r$ mod $p$. The first treats
    the case $i<a$.
    
	\begin{corollary}
		\label{Shape of Theta i < a not in interval}		
        Let $ p \geq 3 $ and
		$  v(a_{p}) \in (i,i+1)$ for some $ 0 \leq i < p-1 $.	Let 
		$r \geq i(p+1)+p$ and $   k-2 =:r\equiv a \mod (p-1) $  with $ 1 \leq a \leq p-1  $.
		Further assume that  $ i < a $.
		\begin{enumerate}
			\item[$(i)$] If $  a > 2i $ and $ r \not \equiv a-i+1, a-i+2, \ldots,a-1, a  
			\mod p $, then
			$
			\ind_{KZ}^{G}(V_{p-1-a+2i} \otimes D^{a-i}) \twoheadrightarrow 
			\bar{\Theta}_{k,a_{p}}.
			$
			\item[$(ii)$] If $ a = 2i $ and $ r \not \equiv a-i-1, a-i, \ldots, a-1,a
			\mod p $, then
			$
			\ind_{KZ}^{G}(V_{p-3-a+2i} \otimes D^{i+1}) \twoheadrightarrow 
			\bar{\Theta}_{k,a_{p}}.
			$
			\item[$(iii)$] If $ a < 2i $ and $ r \not \equiv a-i-1, a-i, \ldots, a-1, a 
			\mod p $, then
			$
			\ind_{KZ}^{G}(V_{p-3+a-2i} \otimes D^{i+1}) \twoheadrightarrow 
			\bar{\Theta}_{k,a_{p}}.
			$
		\end{enumerate}
	\end{corollary}
	\begin{proof}
        Note that $ r-m \equiv  a-m \mod (p-1)$ where $ 1 \leq a-m \leq p-1 $ for all $ 0 \leq m < a $. 
		\begin{enumerate}
			\item[$(i)$] By Lemma~\ref{JH factor Q 2i < a} $(i)$, the JH factors of $ Q(i) $ are
		     \[ \{ V_{p-1-a+2l} \otimes D^{a-l} : 0 \leq l \leq i \},\]
            which are the cosocles of $V_r^{(l)}/V_r^{(l+1)}$
            for $0 \leq l \leq i$.  If $ i=0 $, then $ Q(0) \cong V_{p-1-a} \otimes D^{a} $ and $(i)$ follows from \eqref{surjection from Q(i)}.
			So assume $ i > 1 $.
			Since $ a > 2i $ it can be checked that $ 0 \leq m \leq i-1 $
			satisfy the hypothesis of Theorem~\ref{Elimination i < a and not in interval}.
			Hence,  the images of $\ind_{KZ}^{G} V_{r}$,  $\ind_{KZ}^{G} V_{r}^{(1)},
			\ldots,  \ind_{KZ}^{G} V_{r}^{(i)} $
			are the same in $ \bar{\Theta}_{k,a_{p}} $, that is,  the JH factors above for $l=0,\ldots,i-1$ die in $ \bar{\Theta}_{k,a_{p}} $.
			Thus, $(i)$ again follows from \eqref{surjection from Q(i)}.  
			\item[$(ii)$] Since $2i=a\geq 1$ we have $i\geq 1$. By Lemma~\ref{JH factor Q 2i = a} $(i)$ part $(b)$, we have $ Q(i) \cong Q(i-1) $ and the JH factors of $ Q(i) $ are
		     \[ \{ V_{p-1-a+2l} \otimes D^{a-l} : 0 \leq l \leq i-1 \},\]
            which are the cosocles of $V_r^{(l)}/V_r^{(l+1)}$ for $0 \leq l \leq i-1$.                
			If $ i=1 $, then  $a=2$ and $ Q(i-1) = Q(0) \cong V_{p-1-a} \otimes D^{a} 
			= V_{p-3} \otimes D^{2}$
			and $(ii)$ follows immediately from \eqref{surjection from Q(i)}. So assume $ i > 1  $. 
			Since $ a = 2i $  it can be checked that $ 0 \leq m \leq i-2 $
			satisfy the hypothesis of Theorem~\ref{Elimination i < a and not in interval}.
			Hence, the images of $\ind_{KZ}^{G} V_{r}$, $\ind_{KZ}^{G} V_{r}^{(1)},
			\ldots$,   $\ind_{KZ}^{G} V_{r}^{(i-1)} $
			are the same in $ \bar{\Theta}_{k,a_{p}} $, that is,  the JH factors above for $l=0,\ldots,i-2$ die in $ \bar{\Theta}_{k,a_{p}} $. Thus, $(ii)$ again follows from \eqref{surjection from Q(i)}.  
			\item[$(iii)$] Since $i<a$ we have $a-i\geq 1$. By Lemma~\ref{JH factor Q 2i > a} $(i)$ part $(b)$, we have $ Q(i) \cong Q(a-i-1) $ and the JH factors of $ Q(i) $ are
		     \[ \{ V_{p-1-a+2l} \otimes D^{a-l} : 0 \leq l \leq a-i-1 \},\]
            which are the cosocles of $V_r^{(l)}/V_r^{(l+1)}$ for $0 \leq l \leq a-i-1$. 
			If $ a-i =1 $, then $ Q(a-i-1) = Q(0) \cong V_{p-1-a} \otimes D^{a}
			= V_{p-3+a-2i} \otimes D^{i+1}$ and $(iii)$ follows immediately from \eqref{surjection from Q(i)}. 
			So assume $ a-i \geq 2  $. As $ a<2i $, we have
			$ a- i+1 < i+1 = a-(a-i-1) $. Then it can be checked 
			that $ 0 \leq m \leq a-i-2 $ satisfy the hypothesis of 
			Theorem~\ref{Elimination i < a and not in interval}.	Hence, the 
			images of $ \ind_{KZ}^{G} V_{r},  \ind_{KZ}^{G} V_{r}^{(1)},
			\ldots, \ind_{KZ}^{G} V_{r}^{(a-i-2)},  \ind_{KZ}^{G} V_{r}^{(a-i-1)} $
			are the same in $ \bar{\Theta}_{k,a_{p}} $, that is,  the JH factors above for $l=0,\ldots,a-i-2$ die in $ \bar{\Theta}_{k,a_{p}} $.  Thus, $(iii)$ again follows from \eqref{surjection from Q(i)}. \qedhere
		\end{enumerate}
	\end{proof}

	\begin{remark}
		\label{remark on cor shape of Theta i < a not in interval}
		In the case   $ i=1 $ and $a = 2$ (resp. $ i =a-1 $ and $i\geq 2$), the statement of part $(ii)$ 
		(resp. part $(iii)$) of Corollary~\ref{Shape of Theta i < a not in interval} 
		is valid even when $ r \equiv a $ mod $ p $. This follows from  
		Lemma~\ref{JH factor Q 2i = a} $(i)$ part $(a)$ and Lemma~\ref{JH factor Q 2i > a} $(i)$ part $(a)$.
	\end{remark}
	
    In the next two theorems, we take a momentary detour and provide some further information about the case $r \equiv 2i+1 \mod (p-1)$ in Corollary~\ref{Shape of Theta i < a not in interval}. This is precisely the case where the JH factor surjecting onto $\bar{\Theta}_{k,a_p}$has dimension $p-1$. The argument given below is based on the proof of \cite[Theorem 9.1]{BG15}.
	\begin{theorem}\label{thm: a=2i+1 good cases and non-half integer slopes}
		Let $ p > 3 $ and $ v = v(a_{p}) \in (i,i+1)$ for some $ 0 \leq i \leq (p-2)/2$.	Let 
		$k-2=r \geq i(p+1)+p$ $ r \equiv a \mod (p-1) $ with $1\leq a\leq p-1 $. Assume $a=2i+1$, $v(a_p) \neq i+\frac{1}{2}$ and $r \not \equiv i+2,\ldots, 2i+1 \mod p$. Then the surjection $\ind_{KZ}^{G}(V_{p-2} \otimes D^{i+1}) \twoheadrightarrow 
		\bar{\Theta}_{k,a_{p}}$ factors as $$ \frac{\ind_{KZ}^{G}(V_{p-2} \otimes D^{i+1})}{T} \twoheadrightarrow 
		\bar{\Theta}_{k,a_{p}}.$$
	\end{theorem}
	\begin{proof}
		  Let 
		\begin{align*}
			f_2 = \sum_{  \lambda \in \mathbb{F}_{p}} \left[ g^0_{2,p[\lambda]}, \frac{1}{a_p} (-\theta)^{i}(Y^{r-i(p+1)}-X^{r-i(p+1) -1}Y)\right].
		\end{align*}
		By Lemma~\ref{theta and T plus}, we have $T^+f_2$ vanishes modulo $p$. Furthermore,
		\begin{align*}
			T^{-} f_2 \equiv \left[ g^0_{1,0}, \frac{(p-1)p^i}{a_p} \sum_{ \substack{0 \leq j < r-i \\ j \equiv r-i \mod (p-1)} } \binom{r-i}{j} X^{r-j}Y^j \right] \mod p.
		\end{align*}
		By \cite[Lemma 2.5]{BG15} and the usual argument, we have 
		\[
		\sum_{ \substack{0 < j < r-i \\ j \equiv r-i \mod (p-1)} } \binom{r-i}{j}\binom{j}{m} \equiv 0 \mod p \text{ for } m=0,\ldots,i.
		\]
		By Lemma~\ref{lem:choice of beta}, there exists $\alpha_j$ for $0<j<r-i$ and $j \equiv r-i \mod (p-1)$ such that 
		\begin{enumerate}
			\item[$(i)$] $\alpha_j \equiv \binom{r-i}{j} \mod p$
			\item[$(ii)$] $\sum \alpha_j \binom{j}{m} \equiv 0 \mod p^{i+2}$ for  $m=0,\ldots,i$.
		\end{enumerate}
		Let 
		\begin{align*}
			f_1 = \left[ g^0_{1,0}, \frac{(p-1)p^i}{a_p^2} \sum_{ \substack{0 \leq j < r-i \\ j \equiv r-i \mod (p-1)} } \alpha_{j} X^{r-j}Y^j \right].
		\end{align*}
		It can be checked that $T^{-} f_2-a_p f_1 $ and $T^{-}f_1$ vanish modulo $p$. Furthermore, using $i\leq p-3$ we have
		\begin{align*}
			T^{+}f_1 & \equiv \left[g^0_{2,0}, \frac{(p-1)p^{2i+1}}{a_p^2}\alpha_{i+1} X^{r-i-1}Y^{i+1} \right] \\
			& \qquad + \sum_{  \lambda \in \mathbb{F}_{p}^{\times}} \left[g^0_{2,p[\lambda]}, \frac{(p-1)p^{2i+1}}{a_p^2} \sum_{ \substack{0 < j < r-i \\ j \equiv r-i \mod (p-1)} } \alpha_{j}\binom{j}{i+1} X^{r-i-1}Y^{i+1} \right] \mod p.
		\end{align*}
		Note that $\alpha_{i+1} \equiv \binom{r-i}{i+1} \mod p$ and 
		\begin{align*}
			\sum_{\substack{0 < j < r-i \\ j \equiv r-i \mod (p-1)} } \alpha_{j}\binom{j}{i+1} \equiv \sum_{\substack{0 < j < r-i \\ j \equiv r-i \mod (p-1)} } \binom{r-i}{j}\binom{j}{i+1} \equiv \binom{r-i}{i+1}  \mod p,
		\end{align*}
		where the last congruence follows from \cite[Lemma 2.14]{GR19}. Thus 
		\begin{align}\label{eq: good congruence class a=2i+1 hecke operator final}
        \begin{split}
			(T-a_p)(f_2+f_1) & \equiv -a_p f_2 +T^{+} f_1  \\
			& \equiv - \sum_{  \lambda \in \mathbb{F}_{p}} \left[ g^0_{2,p[\lambda]},  (-\theta)^{i}(Y^{r-i(p+1)}-X^{r-i(p+1) -1}Y) + \frac{p^{2i+1}}{a_p^2}\binom{r-i}{i+1}X^{r-i-1}Y^{i+1}\right] \mod p. 
        \end{split}    
		\end{align}
        Observe that the first term $(-\theta)^{i}(Y^{r-i(p+1)}-X^{r-i(p+1) -1}Y)$ has projection $(-1)^{i+1}X^{p-2}$ under the map $V_r^{(i)}/V_r^{(i+1)}  \twoheadrightarrow V_{p-2} \otimes D^{i+1}$ by Lemma~\ref{Glover-Brueil map image} $(ii)$.
        
        If $v(a_p^2) <2i+1$ or $p\mid \binom{r-i}{i+1}$, then the second term dies and the map $\ind_{KZ}^{G}(V_{p-2} \otimes D^{i+1}) \twoheadrightarrow 
		\bar{\Theta}_{k,a_{p}}$ factors through $T$ and we are done.

	    So now assume that  $v(a_p) > i+\frac{1}{2}$ and $p\nmid \binom{r-i}{i+1}$. Thus $r\not \equiv i,i+1 \mod p$ at least if $i\geq 1$. In this case $(T-a_p)(f_2+f_1)$ is not integral. However, we can use the following modified function $f= f_2' +f_1'$:
	   \begin{align*}
	   	     f_2' &= \frac{a_p^2}{p^{2i+1}} f_2 =  \sum_{  \lambda \in \mathbb{F}_{p}} \left[ g^0_{2,p[\lambda]}, \frac{a_p}{p^{2i+1}}(-\theta)^{i}(Y^{r-i(p+1)}-X^{r-i(p+1)-1}Y)\right] \\
	   	     f_1' &=  \frac{a_p^2}{p^{2i+1}} f_1 =\left[ g^0_{1,0}, \frac{(p-1)}{p^{i+1}} \sum_{ \substack{0 \leq j < r-i \\ j \equiv r-i \mod (p-1)} } \alpha_{j} X^{r-j}Y^j \right].
	   \end{align*} 
	   Multiplying \eqref{eq: good congruence class a=2i+1 hecke operator final} by $a_p^2/p^{2i+1}$, we get
	   \begin{align*}
	   	(T-a_p)(f_2'+f_1') & \equiv T^{+} f_1' 
	   	\equiv -\sum_{  \lambda \in \mathbb{F}_{p}} \left[ g^{0}_{2,p[\lambda]},  \binom{r-i}{i+1}X^{r-i-1}Y^{i+1}\right] \mod p. 
	   \end{align*}
       If $i=0$, then by Lemma~\ref{Glover-Brueil map image} $(ii)$, $X^{r-1}Y$ maps to  $X^{p-2}$ under the surjection $V_r^{(i)}/V_r^{(i+1)}  \twoheadrightarrow V_{p-2} \otimes D^{i+1}$. As $p\nmid \binom{r-i}{i+1}$, we are done if $i=0$. So assume $i\geq 1 $. Note that the above function doesn't take values in $V_r^{(i)}$. To compute the projection, we need the following: \\

	   \textbf{Claim:} There exists $F(X,Y) \in X_{r-i}$ such that 
	  \begin{enumerate}
	  	\item[$(i)$] $X^{i+1}  Y^{r-i-1}+ F(X,Y) \in V_r^{(i)}$
	  	\item[$(ii)$] The image of $X^{i+1}  Y^{r-i-1} + F(X,Y) \in V_r^{(i)}$ under the map $V_r^{(i)}/V_r^{(i+1)}   \twoheadrightarrow V_{p-2} \otimes D^{i+1}$ equals $-\binom{r-i-1}{i}Y^{p-2}$. Hence $\begin{psmallmatrix}
	  	0 & 1 \\ 1 & 0
	  \end{psmallmatrix} (X^{i+1}  Y^{r-i-1}+ F(X,Y))$ maps to $-(-1)^{i+1}\binom{r-i-1}{i} X^{p-2}$.
	  \end{enumerate}
\underline{Proof of claim:} Since $i \leq p-3$, by \cite[Proposition 2.16 (ii)]{GR19} we have the following matrix 
	  \begin{align*}
	  	\left( \binom{r-l}{m} \binom{2i+1-l-m}{i-m}\right)_{0\leq m,l \leq i-1}
	  \end{align*}
	  is invertible modulo $p$ if $r \not\equiv i+2, \ldots, 2i+1 \mod p$. Thus there exists $\beta_0, \ldots, \beta_{i-1}$ such that 
	  \begin{align*}
	  	\sum_{l=0}^{i-1} \beta_l \binom{r-l}{m} \binom{2i+1-l-m}{i-m} \equiv \binom{r-i-1}{m} \mod p \quad \text{ for } m=0,\ldots,i-1.
	  \end{align*}
	  Consider the polynomial
	 \begin{align*}
	 	 F(X,Y) &:= \sum_{l=0}^{i-1} \beta_{l} \sum_{  \lambda \in \mathbb{F}_{p}^{\times}} \lambda^{-i} X^{l}(X+\lambda Y)^{r-l} = (p-1) \sum_{l=0}^{i-1} \beta_{l}  \sum_{\substack{i \leq j \leq  r-i-1 \\ j \equiv i \mod (p-1)} }  \binom{r-l}{j} X^{r-j}Y^j.
	 \end{align*}
     We show that $F(X,Y)$ satisfies the claim. 
     Clearly $F(X,Y) \in X_{r-i}$. We check $X^{i+1}Y^{r-i-1}+F(X,Y)$ satisfies $(i)$. Clearly the coefficients of $X^r,\ldots,X^{r-i+1}Y^{i-1}$ and $Y^r,\ldots, X^iY^{r-i}$ are zero in $F(X,Y)$. Furthermore, for $m=0,\ldots,i-1$ we have 
	  \begin{multline*}
	  	\binom{r-i-1}{m} + (p-1) \sum_{l=0}^{i-1} \beta_l \sum_{ \substack{i\leq  j \leq  r-i-1 \\ j \equiv i ~\mathrm{mod}~ (p-1)} } \binom{r-l}{j} \binom{j}{m} \\ \equiv \binom{r-i-1}{m} + (p-1) \sum_{l=0}^{i-1} \beta_l \binom{r-l}{m} \binom{2i+1-l-m}{i-m} \equiv 0 \mod p,
	  \end{multline*}
	  where the first congruence follows from \cite[Lemma 2.14]{GR19} and the second congruence follows from the choice of $\beta_l$. This proves $(i)$.
	  
	  By \cite[Lemma 2.12]{GR19}, we have the projection of $X^{r-i}Y^{i}+F(X,Y)$ under $V_r^{(i)}/V_r^{(i+1)} \twoheadrightarrow V_{p-2} \otimes D^{i+1}$ is the same as the image of 
	  \begin{align*}
	  	\theta^i \left( \binom{r-i-1}{i} + (p-1)\sum_{l=0}^{i-1} \beta_l \sum_{ \substack{0 < j \leq r-i-1 \\ j \equiv i ~\mathrm{mod}~ (p-1)} } \binom{r-l}{j} \binom{j}{i} - (p-1)\sum_{l=0}^{i-1} \beta_l \binom{r-l}{i}\right)X^{r-i(p+1)-(p-1)}Y^{p-1}.
	  \end{align*}
	  By \cite[Lemma 2.14]{GR19}, we have
	  \begin{align*}
	  	   \sum_{ \substack{0 < j \leq r-i-1 \\ j \equiv i ~\mathrm{mod}~(p-1)} } \binom{r-l}{j} \binom{j}{i} \equiv \binom{r-l}{i} \mod p.
	  \end{align*}
	   Thus by Lemma~\ref{Glover-Brueil map image} $(ii)$, the projection of $X^{r-i}Y^{i}+F(X,Y)$ under $V_r^{(i)}/V_r^{(i+1)} \twoheadrightarrow V_{p-2} \otimes D^{i+1}$ equals $-\binom{r-i-1}{i} Y^{p-2}$. This proves $(ii)$ and the claim. 

    From the claim above, it follows that $\overline{(T-a_p)(f_2+f_1)}$ maps to
	   \begin{align*}
	   (-1)^{i+1}	\binom{r-i}{i+1}\binom{r-i-1}{i}\sum_{  \lambda \in \mathbb{F}_{p}} \left[g^{0}_{2,p[\lambda]}, X^{p-2} \right] \in \mathrm{ind}_{KZ}^{G}(V_{p-2} \otimes D^{i+1}). 
	   \end{align*}
       As $r\not\equiv i,i+1, i+2,\ldots,2i+1 \mod p$, we have $p \nmid \binom{r-i}{i+1}\binom{r-i-1}{i}$ the surjection $\ind_{KZ}^{G}(V_{p-2} \otimes D^{i+1}) \twoheadrightarrow 
		\bar{\Theta}_{k,a_{p}}$ factors through $T$. This completes the proof of the theorem.
	\end{proof}
    We now consider the case $v(a_p) = i+\frac{1}{2}$. Following \cite{Gha20}, we set
	\begin{align*}
		c = \frac{a_p^2 -  \binom{r-i-1}{i} \binom{r-i}{i+1} p^{2i+1}}{pa_p} \quad \mathrm{and}\quad  \tau = v(c),
\end{align*}
	and  $t=v(r-2i-1)$. By  \cite[Conjecture 1.1]{Gha20}, for $t\gg 0$, we have 
	\begin{align*} 
		\bar{V}_{k,a_p} \simeq
		\begin{cases}
			\mathrm{ind}_{KZ}^{G}(\omega_2^{i+2+ip}) &\text{ if } \tau \in (t+i-1,t+i), \\
			\mu_\lambda \omega^{i+1} \oplus	\mu_{\lambda^{-1}} \omega^{i+1}	&\text{ if } \tau \geq t+i,
		\end{cases}
	\end{align*}
for some $\lambda \in \mathbb{F}_p$.  This was proved in \cite[Theorem 1.2]{Gha22}.  We now show 
that this also holds when $r \not \equiv i+2, \ldots, 2i+1 \mod p$ and $\tau  =i-\frac{1}{2}$.

\begin{theorem}\label{a=2i+1 good cases and half integer slopes}
	Let $ p > 3 $, $ a_{p} \in \bar{\mathbb{Q}}_{p} $ be such that 
	$ v = v(a_{p}) = i+\frac{1}{2}$ for some $ 0 \leq i \leq (p-2)/2$.	Let 
	$i(p+1) + p \leq r  \equiv a \mod (p-1) $ with $1 \leq a \leq p-1$. Assume $a=2i+1$ and $r \not \equiv i+2,\ldots, 2i+1 \mod p$. If $\tau  =i-\frac{1}{2}$, then we have
	\begin{align*}
		\bar{V}_{k,a_p} \simeq \mathrm{ind}(\omega_2^{i+2+ip}).
    \end{align*}
\end{theorem}
\begin{proof}
	Let $f_2$ and $f_1$ be as defined in Theorem~\ref{thm: a=2i+1 good cases and non-half integer slopes}. By \eqref{eq: good congruence class a=2i+1 hecke operator final}, we have
	\begin{align*}
		(T-a_p) f' & \equiv -a_p f'_2 + T^+f'_1 \\
		& \equiv
		-  \sum_{  \lambda \in \mathbb{F}_{p}}  \left[ g^0_{2,p[\lambda]},  (-\theta)^{i}(Y^{r-i(p+1)}-X^{r-i(p+1) -1}Y)- \frac{(p-1)p^{2i+1}}{a_p^2}  \binom{r-i}{i+1}X^{r-i-1}Y^{i+1}\right] \mod p.
	\end{align*}
	As in Theorem~\ref{thm: a=2i+1 good cases and non-half integer slopes}, it can be checked that the first summand $(T-a_p)f'$ maps to  $(-1)^{i+1} T([g^0_{1,0},X^{p-2}])$ under  $ \mathrm{ind}_{KZ}^{G}(V_r^{(i)}/V_r^{(i+1)}) \rightarrow \mathrm{ind}_{KZ}^{G}(V_{p-2} \otimes D^{i+1})$. By the claim proved in Theorem~\ref{thm: a=2i+1 good cases and non-half integer slopes}  the second summand (after a suitable modification) maps to $(-1)^{i+1} \frac{(p-1)p^{2i+1}}{a_p^2}  \binom{r-i}{i+1} \binom{r-i-1}{i}T([g^0_{1,0},X^{p-2}])$. 
    This completes the proof as $v(a_p^2 -  \binom{r-i-1}{i} \binom{r-i}{i+1} p^{2i+1}) = 2i+1$.
\end{proof}

    Returning to the main development,
    we now treat the good congruence classes of $r$ mod $p$ for $i \geq a$.
    	\begin{corollary}
		\label{Shape of Theta i > a not in interval} 
		Let $ p \geq 3 $, $ a_{p} \in \bar{\mathbb{Q}}_{p} $ be such that $ v = v(a_{p}) \in (i,i+1)$
		for some $ 0 \leq i < p-1 $.	Let $k-2=:r\geq i(p+1)+p$ and $ r \equiv a \mod (p-1) $  
		with $ 1 \leq a \leq p-1  $ and $ b = p-1+a $. Further assume $ i \geq a $.
		\begin{enumerate}
			\item[$(i)$] If $ b > 2i  $ and $ r \not \equiv  b-i+1, b-i+2, \ldots, 
			b-1, b \mod  p $, then
			$
			\ind_{KZ}^{G}(V_{p-1-b+2i} \otimes D^{b-i}) \twoheadrightarrow 
			\bar{\Theta}_{k,a_{p}}.
			$
			\item[$(ii)$] If $ b = 2i $ and $ r \not \equiv b-i-1, b-i, \ldots, 
			b-1, b \mod  p $, then
			$
			\ind_{KZ}^{G}(V_{p-3-b+2i} \otimes D^{i+1}) \twoheadrightarrow 
			\bar{\Theta}_{k,a_{p}}.
			$
			\item[$(iii)$] If $ b < 2i $ and $ r \not \equiv b-i-1, b-i, \ldots, 
			b-1, b \mod  p $, then
			$
			\ind_{KZ}^{G}(V_{p-3+b-2i} \otimes D^{i+1}) \twoheadrightarrow 
			\bar{\Theta}_{k,a_{p}}.
			$
		\end{enumerate}
	\end{corollary}
	\begin{proof}
		The proof is similar to Corollary~\ref{Shape of Theta i < a not in interval}. Note that $p-1+a>i+a\geq2a$ and $ r-m \equiv  b-m \mod (p-1)$ where $ 1 \leq b-m \leq p-1 $ for all $ a \leq m \leq p-1 $. 
		\begin{enumerate}
			\item[$(i)$] Note that  $ r \not \equiv  b-i+1, b-i+2, \ldots, 
			b-1, b \mod  p $ is equivalent to $ r \equiv  b+1, b+2, \ldots, 
			 p+b-i \mod  p $, that is, $ r \equiv a,\ldots,[a-i] \mod p$. 
            If $ i =a $, then  by Lemma~\ref{JH factor Q(a)} $(iii)$, we have $Q(i)=V_a$ which is the cosocle of $V_{r}^{(a)}/V_{r}^{(a+1)}$. Now $(i)$ follows from \eqref{surjection from Q(i)}. Assume $i>a$. By Lemma~\ref{JH factor Q i < [a-i]} $(ii)$, the JH factors of $ Q(i) $ are
		     \[ \{ V_{2l-a} \otimes D^{a-l} : a \leq l \leq i \},\]
            which are the cosocles of $V_r^{(l)}/V_r^{(l+1)}$
            for $a \leq l \leq i$.
			If $ i > a $, then it can be checked that $ a \leq m \leq i-1  $
			satisfy the hypothesis of Theorem~\ref{Elimination i < a and not in interval}
			with $ a $ there equal to $ b = p-1+a $. Thus the images of 
			$ \ind_{KZ}^{G}V_{r}^{(a)}, \ind_{KZ}^{G}V_{r}^{(a+1)}, \ldots, 
			\ind_{KZ}^{G}V_{r}^{(i)}$ are the same in
			$ \bar{\Theta}_{k,a_{p}} $, that is,  the JH factors above for $l=a,\ldots,i-1$ die in $ \bar{\Theta}_{k,a_{p}} $.
			Thus, $(i)$ again follows from \eqref{surjection from Q(i)}.     
			\item[$(ii)$] Note that  $ r \not \equiv  b-i-1, b-i, \ldots, 
			b-1, b \mod  p $ is equivalent to $ r  \equiv  b+1, b+2, \ldots, p+b-i-2 \mod  p $, that is, $ r \equiv a,\ldots,i-2 \mod p$.  Since $ i < p-1  $, we have
			$ a = b-(p-1) =2i -(p-1)   < i$.  Thus, by Lemma~\ref{JH factor Q i = [a-i]} $(ii)$ part $(b)$,  we have $ Q(i) \cong Q(i-1) $ and the JH factors of $ Q(i) $ are
		     \[ \{ V_{2l-a} \otimes D^{a-l} : a \leq l \leq i-1 \},\]
            which are the cosocles of $V_r^{(l)}/V_r^{(l+1)}$ for $a \leq l \leq i-1$.   
			If $ i= a+ 1 $, then  $ a = p-3 $ and $ Q(i)  \cong V_{a}  
			= V_{p-3}  $  and $(ii)$ follows immediately from  \eqref{surjection from Q(i)}.
			So assume $ i > a+1  $.
			Since $ b = 2i $ it can be checked that $ a \leq m \leq i-2 $
			satisfy the hypothesis of Theorem~\ref{Elimination i < a and not in interval}
			with $ a $ there equal to $ b = p-1+a $.
			Hence, the images of $\ind_{KZ}^{G} V_{r}^{(a)}$, $\ind_{KZ}^{G} V_{r}^{(a+1)},
			\ldots$,   $\ind_{KZ}^{G} V_{r}^{(i-1)} $
			are the same in $ \bar{\Theta}_{k,a_{p}} $, that is,  the JH factors above for $l=a,\ldots,i-2$ die in $ \bar{\Theta}_{k,a_{p}} $. Thus, $(ii)$ again follows from \eqref{surjection from Q(i)}.
			\item[$(iii)$] Note that  $ r \not \equiv  b-i-1, b-i, \ldots, 
			b-1, b \mod  p $ is equivalent to $ r  \equiv  b+1, b+2, \ldots, p+b-i-2 \mod  p $, that is, $ r \equiv a,\ldots,[a-i]-2 \mod p$.  
        Thus, by Lemma~\ref{JH factor Q i > [a-i]} $(ii)$ part $(b)$,  we have $ Q(i) \cong Q(b-i-1) $ and the JH factors of $ Q(i) $ are
		     \[ \{ V_{2l-a} \otimes D^{a-l} : a \leq l \leq [a-i]-1 \},\]
             which are the cosocles of $V_r^{(l)}/V_r^{(l+1)}$ for $a \leq l \leq b-i-1$.
			If $ b-i =a+1 $, then  $i=p-2$ and $ Q(i) =  V_{a} 
			= V_{p-3+b-2i} \otimes D^{i+1}$ and $(iii)$ follows immediately
			from \eqref{surjection from Q(i)}. 
			So assume $ b-i \geq a+2  $. 
            Since $ b < 2i $ it can be checked 
			that $ 0 \leq m \leq b-i-2 $ satisfy the hypothesis of 
			Theorem~\ref{Elimination i < a and not in interval}
			with $ a $ there equal to $ b = p-1+a $.	Hence the 
			images of $ \ind_{KZ}^{G} V_{r}^{(a)},  \ind_{KZ}^{G} V_{r}^{(a+1)},
			\ldots,   \ind_{KZ}^{G} V_{r}^{(b-i-1)} $
			are the same in $ \bar{\Theta}_{k,a_{p}} $, that is,  the JH factors above for $l=a,\ldots,[a-i]-2$ die in $ \bar{\Theta}_{k,a_{p}} $. Thus, $(iii)$ again follows from \eqref{surjection from Q(i)}.
			\qedhere 	     
		\end{enumerate}
	\end{proof}
    Finally, we note that in 
    part $(i)$ above, if $b = 2i+1$,
    then a JH factor of dimension $p-1$ occurs, and more information is required to 
    predict the structure of 
    $\bar{\Theta}_{k,a_p}$. 
    Unlike the case 
    $i < a$ treated in the two theorems before the corollary, we do not explore
    what happens when $i \geq a$.
    In any case, the two corollaries
    above give a complete treatment of the
    structure of $\bar{\Theta}_{k,a_p}$ for
    good congruence classes of $r$ mod $p$ away from such cases.

      	\section{Bad congruence classes when \texorpdfstring{$i<a$}{}}\label{sec: bad}

    Recall that $r \equiv a$ mod $(p-1)$ with $1  \leq a \leq p-1$ and $v(a_p) \in (i,i+1)$ for some $i \in \Z_{\geq 0}$.  To keep this paper at a reasonable length we shall assume throughout this section that
    $$ i < a.$$  
    In the previous chapter, we treated the good cases 
    $$ r \not\equiv \begin{cases}
          a-i+1, a-i+2, \ldots,a \mod p & \text{if } a > 2i, \\
         a-i-1, a-i, \ldots,a \mod p & \text{if } a \leq 2i.
    \end{cases}$$
    We now turn our attention to the remaining congruence classes of $r \mod p$, namely 
     $$r \equiv \begin{cases}
         a-i+1, a-i+2, \ldots,a \mod p & \text{if } a > 2i, \\
         a-i-1, a-i, \ldots,a \mod p & \text{if } a \leq 2i,
    \end{cases}$$
    which we refer to as the
    {\it bad} congruence classes of $r$ mod $p$.
    Thus, we shall assume that 
	\begin{align*}
		r \equiv a-i+n \mod p \text{ with } -1 \leq n \leq i.
	\end{align*}
     
    
	
	 Let $s =a-i+n+(i-n)p$. Thus $r \equiv s$ mod $p(p-1)$. We can use results on local constancy to describe the shape of $\bar{V}_{r+2, a_p}$
     when $v(r-s)$ is large, at least when $1 \leq n \leq i$. This is a good first approximation to what happens in general.

     \begin{lemma}
         Let  $p \geq 2$ and $v(a_p) \in (i,i+1)$ with $0 \leq i < p-1$ be an integer. Let $ r \geq i(p+1)+p$ and $r\equiv a-i+n+(i-n)p$ mod $p(p-1)$ with $ 1 \leq a \leq p-1$ and $1\leq n < i$. Assume $i<a$. Let $s =a-i+n+(i-n)p$. If $s \geq 4i+2$ and $v(r-s) \gg 0$, then 
    \begin{equation*}
        \bar{V}_{k, a_p}  \sim 
        \begin{cases}
           \ind(\omega_2^{s+1}) \quad \qquad \qquad &\text{if } (p+1)\nmid (s+1),\\
           (\mu_{\sqrt{-1}} \oplus \mu_{-\sqrt{-1}}) \otimes \omega^{\frac{s+1}{p+1}}     &\text{if } (p+1)\mid (s+1).   \end{cases}
    \end{equation*}
     \end{lemma}
     
     \begin{proof}
     By  \cite[Theorem B]{Ber12}  (see also \cite{Ber}), it can be checked that if $s \geq 4i+2$ and $v(r-s) \gg 0$, then 
    \begin{equation*}
        \bar{V}_{r+2, a_p}  \cong \bar{V}_{s+2, a_p}. 
    \end{equation*}   
     Indeed, under these conditions 
     $$\alpha(s+1) :=  \sum_{n\geq 1} \Big\lfloor \frac{s+1}{p^n(p-1)}\Big\rfloor = i-n+ \Big\lfloor \frac{a+1}{p-1} \Big\rfloor \leq i-n+ 1 \leq i,$$
     where we have used $s+1 \leq  p^2-3p+3$, $a\leq p-1$ and  $n\geq 1$. Thus $s+2 > 3v+\alpha(s+1)+1$.    Also, by  \cite[Theorem 1.1.1]{BL22},  we have $\bar{V}_{s+2, a_p} \cong \bar{V}_{s+2, 0}$ since $v(a_p)>i\geq  i-n+1 \geq\lfloor \frac{s+1}{p-1}\rfloor$.
    The lemma follows from \cite[Proposition 6.2]{Breuil}.
    \end{proof}

    However, in the lemma, one does not have an estimate on  how large $v(r-s)$ should be. See also \cite{Bha20}, \cite{GS24}. In this chapter, we make precise this estimate. In fact, we determine the shape of $\bar{V}_{k,a_p}$ as $v(r-s)$ varies through all non-negative integers.
	
	\subsection{The case \texorpdfstring{$a \leq 2i-2n-1$}{}}\label{sec: bad small}
     In this section, we determine $\bar{V}_{k,a_p}$ when $a \leq 2i-2n-1$ and $i<a$. As  $a\leq 2i-2n-1$ and $n \geq -1$, we get $a-i+n \leq i-n-1 \leq i$. In the setting, we are dealing with the congruence classes $r \equiv a-i-1,\ldots, i \mod p$. Also note that we must have $n<i$ since $1\leq a \leq 2i-2n-1$.
     
     We first recall the JH factors of $Q(i)$ in the setting of this subsection.
    \begin{lemma}\label{JH factors bad congruences a<2i case}
		Let  $p \geq 3$ and $0 \leq i \leq p-1$ be an integer. Let $ r \geq i(p+1)+p$ and $r\equiv a-i+n+(i-n)p$ mod $p(p-1)$ with $ 1 \leq a \leq p-1$ and $-1\leq n < i$. Assume that $1\leq a-i<i $. 
		\begin{enumerate}
			\item[$(i)$] If $a=2i-2n-1$, then the JH factors of $Q(i)$ are 
			\begin{center}
				JH factors of $\{V_r^{(l)}/V_r^{(l+1)}: i-n \leq l \leq i \} ~\cup ~\{ V_{p-1-a+2l} \otimes D^{a-l} :  0\leq l \leq a-i-1\}$.
			\end{center}
			\item[$(ii)$] If $a < 2i-2n-1$, then the JH factors of $Q(i)$ are 
			\begin{center}
				JH factors of $\{V_r^{(l)}/V_r^{(l+1)}: i-n \leq l \leq i \} ~\cup ~ \{ V_{2i-2n-2-a} \otimes D^{a-i+n+1} \} $ \newline $ \cup~ \{ V_{p-1-a+2l} \otimes D^{a-l}:0 \leq l \leq a-i-1 \}$.
			\end{center}
		\end{enumerate}
	\end{lemma}
	\begin{proof}
		Follows from parts $(b)$ and $(c)$ of Lemma~\ref{JH factor Q 2i > a} $(ii)$ by taking $r_0$ there equal to $a-i+n$.
	\end{proof}

    Next, we show that the shallow JH factors of $Q(i)$ vanish in $\bar{\Theta}_{k,a_p}$.
    \begin{theorem}\label{shallow JH die a-i<i in bad cases}
		Let  $p \geq 3$ and $v(a_p) \in (i,i+1)$ for some $0 \leq i \leq p-1$. Let $ k-2=:r \geq i(p+1)+p$ and $r\equiv a~\mathrm{mod}~ (p-1)$ with $ 1 \leq a \leq p-1$. If $1\leq a-i<i $ and $r \equiv a-i-1, \ldots, i \mod p$, then the image of $\mathrm{ind}_{KZ}^{G}(V_r)$ is the same as the image of $\mathrm{ind}_{KZ}^{G}(V_r^{(a-i)})$ in $\bar{\Theta}_{k,a_p}$.
	\end{theorem}	
    \begin{proof}
       By Theorem~\ref{Elimination i < a and not in interval}, we see that image of the JH factors coming from $V_{r}^{(m)}/V_{r}^{(m+1)}$ vanish for $0 \leq m \leq a-i-1$. Note that if $m = a-i-1$, then $ \binom{r-m}{i+1} \equiv \binom{n+1}{i+1} \equiv 0 \mod p$ by Lucas' theorem since $n < i$, so hypothesis $(ii)$ of that theorem holds. 
    \end{proof}
    As an immediate consequence, we have the following result when $n=-1$:
    \begin{theorem}\label{a<2i bad n=-1}
        Let  $p \geq 3$ and $v(a_p) \in (i,i+1)$ for some $0 \leq i < p-1 $. Let $ k-2=:r \geq i(p+1)+p$ and $r \equiv a \mod (p-1)$ with $ 1 \leq a \leq p-1$. Assume $1\leq a-i<i$ and $r\equiv a-i-1 \mod p$. Then $\mathrm{ind}_{KZ}^{G}(V_{2i-a} \otimes D^{a- i}) \twoheadrightarrow \bar{\Theta}_{k,a_p}$.
    \end{theorem}
    \begin{proof}
        This follows from Lemma~\ref{JH factors bad congruences a<2i case} $(ii)$ and Theorem~\ref{shallow JH die a-i<i in bad cases}.
    \end{proof}
	The remaining congruence classes $r\equiv a-i, \ldots, a$ mod $p$ are trickier. We prove the following theorem with details provided in the next two subsections. 
    
	\begin{theorem}[Diagonal conjecture] \label{diagonal conj a-i<i}
		Let  $p \geq  3$ and $v(a_p) \in (i,i+1)$ for some $0 \leq i < p-1$. Let $ k-2=: r \geq i(p+1)+p$ and $r \equiv a ~\mathrm{mod}~(p-1)$ with $ 1 \leq a \leq p-1$. Assume  $i<a$ and $r\equiv a-i+n \mod p$ for some $0 \leq n <i$. Set $s = (a-i+n)+(i-n)p $ and $t = v(r-s)$.  
		\begin{enumerate}
			\item[$(i)$] If $a = 2i-2n-1 $, then 
			\begin{itemize}
				\item[$(a)$]  $\mathrm{ind}_{KZ}^{G}(V_r^{(i-t+1)}/V_r^{(i-t+2)}) \twoheadrightarrow \bar{\Theta}_{k,a_p}$ for $  t \leq n$ 
				\item[$(b)$]  $\mathrm{ind}_{KZ}^{G}(V_r^{(i-n)}/V_r^{(i-n+1)}) \twoheadrightarrow \bar{\Theta}_{k,a_p}$ for $ t \geq n+1$.
			\end{itemize}
			\item[$(ii)$] If $a < 2i-2n-1 $, then 
			\begin{itemize}
				\item[$(a)$] $\mathrm{ind}_{KZ}^{G}(V_r^{(i-t+1)}/V_r^{(i-t+2)}) \twoheadrightarrow \bar{\Theta}_{k,a_p}$ for $  t \leq n+1$ 
				\item[$(b)$]  $\mathrm{ind}_{KZ}^{G}(V_r^{(i-n-1)}/V_r^{(i-n)}) \twoheadrightarrow \bar{\Theta}_{k,a_p}$ for $ t \geq n+2$.
			\end{itemize}  
		\end{enumerate}
	\end{theorem} 
    
        The above result explicitly describes which sub-quotient of $V_{r}/V_r^{(i+1)}$ survives in $\bar{\Theta}_{k,a_p}$, namely, given $t=v(r-s)$, it gives the unique value of $T$ such that we have a surjection $\mathrm{ind}_{KZ}^{G}(V_r^{(i-T)}/V_r^{(i-T+1)}) \twoheadrightarrow \bar{\Theta}_{k,a_p}$.  By Theorem~\ref{shallow JH die a-i<i in bad cases}, we know that the JH factors of the sub-quotient $V_{r}/V_r^{(a-i)}$ vanish in $\bar{\Theta}_{k,a_p}$. 
        
        The following beautiful picture explains which of the remaining JH factors $V_r^{(i-T)}/V_r^{(i-T+1)}$ survive when $a=2i-2n-1$.  We list these JH factors along the horizontal axis (according to  increasing $T$) and the possibilities for $t=v(r-s)$ along the vertical axis.
     \begin{figure}[hbt!]
     \centering
    \resizebox{.6\textwidth}{!}{
        \begin{tikzpicture}[font=\small]
  \def\rows{12}
  \def\cols{12}
  \def\cwA{2}   
  \def\cw{1}    
  \def\rhA{1.2} 
  \def\rh{1}    

  \draw (0,0) -- (0,-\rhA-\rows*\rh);
  \draw (\cwA,0) -- (\cwA,-\rhA-\rows*\rh);
  \foreach \j in {1,...,\cols}{ \draw (\cwA+\j*\cw,0) -- (\cwA+\j*\cw,-\rhA-\rows*\rh); }
  \draw (0,0) -- (\cwA+\cols*\cw,0);
  \draw (0,-\rhA) -- (\cwA+\cols*\cw,-\rhA);
  \foreach \i in {1,...,\rows}{ \draw (0,-\rhA-\i*\rh) -- (\cwA+\cols*\cw,-\rhA-\i*\rh); }

  \draw (0,0) rectangle (\cwA,-\rhA);
  \draw (0,0) -- (\cwA,-\rhA); 

  \pgfmathsetmacro{\tx}{\cwA/3}
  \pgfmathsetmacro{\ty}{-2*\rhA/3}
  \pgfmathsetmacro{\Tx}{2*\cwA/3}
  \pgfmathsetmacro{\Ty}{-1*\rhA/3}

  \node at (\tx,\ty) {\(t\)};  
  \node at (\Tx,\Ty) {\(T\)};  

  \node at (\cwA+0.5*\cw,-0.5*\rhA) {$0$};
  \node at (\cwA+1.5*\cw,-0.5*\rhA) {$1$};
  \node at (\cwA+2.5*\cw,-0.5*\rhA) {$2$};
  \node at (\cwA+3.5*\cw,-0.5*\rhA) {$3$};
  \node at (\cwA+4.5*\cw,-0.5*\rhA) {$\cdot$};
  \node at (\cwA+5.5*\cw,-0.5*\rhA) {$\cdot$};
  \node at (\cwA+6.5*\cw,-0.5*\rhA) {$\cdot$};
  \node at (\cwA+7.5*\cw,-0.5*\rhA) {$\cdot$};
  \node at (\cwA+8.5*\cw,-0.5*\rhA) {$n-3$};
  \node at (\cwA+9.5*\cw,-0.5*\rhA) {$n-2$};
  \node at (\cwA+10.5*\cw,-0.5*\rhA) {$n-1$};
  \node at (\cwA+11.5*\cw,-0.5*\rhA) {$n$};

  \node at (0.5*\cwA,-\rhA-0.5*\rh) {$t=1$};
  \node at (0.5*\cwA,-\rhA-1.5*\rh) {$t=2$};
  \node at (0.5*\cwA,-\rhA-2.5*\rh) {$t=3$};
  \node at (0.5*\cwA,-\rhA-3.5*\rh) {$t=4$};
  \node at (0.5*\cwA,-\rhA-4.5*\rh) {$\cdot$};
  \node at (0.5*\cwA,-\rhA-5.5*\rh) {$\cdot$};
  \node at (0.5*\cwA,-\rhA-6.5*\rh) {$\cdot$};
  \node at (0.5*\cwA,-\rhA-7.5*\rh) {$\cdot$};
  \node at (0.5*\cwA,-\rhA-8.5*\rh) {$t=n-2$};
  \node at (0.5*\cwA,-\rhA-9.5*\rh) {$t=n-1$};
  \node at (0.5*\cwA,-\rhA-10.5*\rh) {$t=n$};
  \node at (0.5*\cwA,-\rhA-11.5*\rh) {$t\ge n+1$};

  \foreach \i in {1,...,\cols}{
    \foreach \j in {1,...,\rows}{
      \pgfmathsetmacro{\cx}{\cwA+(\i-0.5)*\cw}
      \pgfmathsetmacro{\cy}{-\rhA-(\j-0.5)*\rh}

      \ifnum\i>4
        \ifnum\i<9
          \ifnum\j>4
            \ifnum\j<9
            \else
              \ifnum\i=\j \node at (\cx,\cy){\cmark}; \else \node at (\cx,\cy){\xmark}; \fi
            \fi
          \else
            \ifnum\i=\j \node at (\cx,\cy){\cmark}; \else \node at (\cx,\cy){\xmark}; \fi
          \fi
        \else
          \ifnum\i=\j \node at (\cx,\cy){\cmark}; \else \node at (\cx,\cy){\xmark}; \fi
        \fi
      \else
        \ifnum\i=\j \node at (\cx,\cy){\cmark}; \else \node at (\cx,\cy){\xmark}; \fi
      \fi

    }
  }

  \foreach \d in {-3,...,3}{ 
    \pgfmathtruncatemacro{\imin}{max(5,5+\d)}
    \pgfmathtruncatemacro{\imax}{min(8,8+\d)}
    \ifnum\imin<\imax
      \pgfmathsetmacro{\xstart}{\cwA+(\imin-1)*\cw}
      \pgfmathsetmacro{\ystart}{-\rhA-(\imin-\d-1)*\rh}
      \pgfmathsetmacro{\xend}{\cwA+(\imax)*\cw}
      \pgfmathsetmacro{\yend}{-\rhA-(\imax-\d)*\rh}
      \ifnum\d=0
        \draw[red,dotted,very thick] (\xstart,\ystart) -- (\xend,\yend);
      \else
        \draw[dotted,thick] (\xstart,\ystart) -- (\xend,\yend);
      \fi
    \fi
  }
\end{tikzpicture}
}
\caption{\label{fig:diagonal conj. a=2i-2n-1} Contribution of JH factors when $a=2i-2n-1$ and $r\equiv a-i+n \mod p$.}
\end{figure}

\noindent 
The symbol $\times$ at position $(t,T)$ indicates that, for this $t$, the image of $\mathrm{ind}_{KZ}^{G}(V_r^{(i-T)}/V_r^{(i-T+1)})$ vanishes  in $\bar{\Theta}_{k,a_p}$. Thus for a given value of $t$, the sub-quotient that survives in $\bar{\Theta}_{k,a_p}$ is marked by \textcolor{red}{$\checkmark$}. In the above picture, we notice that \textcolor{red}{$\checkmark$} always appears along the diagonal ($t=T+1$).

\newpage
Similarly,  for $a<2i-2n-1$, we have
    \begin{figure}[hbt!]
     \centering
    \resizebox{.6\textwidth}{!}{
        \begin{tikzpicture}[font=\small]
  \def\rows{12}
  \def\cols{12}
  \def\cwA{2}   
  \def\cw{1}    
  \def\rhA{1.2} 
  \def\rh{1}    

  \draw (0,0) -- (0,-\rhA-\rows*\rh);
  \draw (\cwA,0) -- (\cwA,-\rhA-\rows*\rh);
  \foreach \j in {1,...,\cols}{ \draw (\cwA+\j*\cw,0) -- (\cwA+\j*\cw,-\rhA-\rows*\rh); }
  \draw (0,0) -- (\cwA+\cols*\cw,0);
  \draw (0,-\rhA) -- (\cwA+\cols*\cw,-\rhA);
  \foreach \i in {1,...,\rows}{ \draw (0,-\rhA-\i*\rh) -- (\cwA+\cols*\cw,-\rhA-\i*\rh); }

  \draw (0,0) rectangle (\cwA,-\rhA);
  \draw (0,0) -- (\cwA,-\rhA); 

  \pgfmathsetmacro{\tx}{\cwA/3}
  \pgfmathsetmacro{\ty}{-2*\rhA/3}
  \pgfmathsetmacro{\Tx}{2*\cwA/3}
  \pgfmathsetmacro{\Ty}{-1*\rhA/3}

  \node at (\tx,\ty) {\(t\)};  
  \node at (\Tx,\Ty) {\(T\)};  

  \node at (\cwA+0.5*\cw,-0.5*\rhA) {$0$};
  \node at (\cwA+1.5*\cw,-0.5*\rhA) {$1$};
  \node at (\cwA+2.5*\cw,-0.5*\rhA) {$2$};
  \node at (\cwA+3.5*\cw,-0.5*\rhA) {$3$};
  \node at (\cwA+4.5*\cw,-0.5*\rhA) {$\cdot$};
  \node at (\cwA+5.5*\cw,-0.5*\rhA) {$\cdot$};
  \node at (\cwA+6.5*\cw,-0.5*\rhA) {$\cdot$};
  \node at (\cwA+7.5*\cw,-0.5*\rhA) {$\cdot$};
  \node at (\cwA+8.5*\cw,-0.5*\rhA) {$n-2$};
  \node at (\cwA+9.5*\cw,-0.5*\rhA) {$n-1$};
  \node at (\cwA+10.5*\cw,-0.5*\rhA) {$n$};
  \node at (\cwA+11.5*\cw,-0.5*\rhA) {$n+1$};

  \node at (0.5*\cwA,-\rhA-0.5*\rh) {$t=1$};
  \node at (0.5*\cwA,-\rhA-1.5*\rh) {$t=2$};
  \node at (0.5*\cwA,-\rhA-2.5*\rh) {$t=3$};
  \node at (0.5*\cwA,-\rhA-3.5*\rh) {$t=4$};
  \node at (0.5*\cwA,-\rhA-4.5*\rh) {$\cdot$};
  \node at (0.5*\cwA,-\rhA-5.5*\rh) {$\cdot$};
  \node at (0.5*\cwA,-\rhA-6.5*\rh) {$\cdot$};
  \node at (0.5*\cwA,-\rhA-7.5*\rh) {$\cdot$};
  \node at (0.5*\cwA,-\rhA-8.5*\rh) {$t=n-1$};
  \node at (0.5*\cwA,-\rhA-9.5*\rh) {$t=n$};
  \node at (0.5*\cwA,-\rhA-10.5*\rh) {$t=n+1$};
  \node at (0.5*\cwA,-\rhA-11.5*\rh) {$t\ge n+2$};

  \foreach \i in {1,...,\cols}{
    \foreach \j in {1,...,\rows}{
      \pgfmathsetmacro{\cx}{\cwA+(\i-0.5)*\cw}
      \pgfmathsetmacro{\cy}{-\rhA-(\j-0.5)*\rh}

      \ifnum\i>4
        \ifnum\i<9
          \ifnum\j>4
            \ifnum\j<9
            \else
              \ifnum\i=\j \node at (\cx,\cy){\cmark}; \else \node at (\cx,\cy){\xmark}; \fi
            \fi
          \else
            \ifnum\i=\j \node at (\cx,\cy){\cmark}; \else \node at (\cx,\cy){\xmark}; \fi
          \fi
        \else
          \ifnum\i=\j \node at (\cx,\cy){\cmark}; \else \node at (\cx,\cy){\xmark}; \fi
        \fi
      \else
        \ifnum\i=\j \node at (\cx,\cy){\cmark}; \else \node at (\cx,\cy){\xmark}; \fi
      \fi

    }
  }

  \foreach \d in {-3,...,3}{ 
    \pgfmathtruncatemacro{\imin}{max(5,5+\d)}
    \pgfmathtruncatemacro{\imax}{min(8,8+\d)}
    \ifnum\imin<\imax
      \pgfmathsetmacro{\xstart}{\cwA+(\imin-1)*\cw}
      \pgfmathsetmacro{\ystart}{-\rhA-(\imin-\d-1)*\rh}
      \pgfmathsetmacro{\xend}{\cwA+(\imax)*\cw}
      \pgfmathsetmacro{\yend}{-\rhA-(\imax-\d)*\rh}
      \ifnum\d=0
        \draw[red,dotted,very thick] (\xstart,\ystart) -- (\xend,\yend);
      \else
        \draw[dotted,thick] (\xstart,\ystart) -- (\xend,\yend);
      \fi
    \fi
  }
\end{tikzpicture}
}
    \caption{\label{fig:diagonal conj. a<2i-2n-1} Contribution of JH factors when $a<2i-2n-1$ and $r\equiv a-i+n \mod p$.}
\end{figure}

The two pictures above explain why Theorem~\ref{diagonal conj a-i<i} was in the run up to our proof of it called the \textquote{Diagonal Conjecture}. 
 \begin{remark} We give some numerical examples.
 \begin{itemize}
     \item  When $p = 7$, $r = 65$ and $i = 4$,
        so that $a = 5$, $n = 1$, $s = 23$ and $t = 1$, then 
    $V_r^{(4)}/V_r^{(5)}$ 
        contributes,
        corroborating part $(i)$ $(a)$ of  Theorem~\ref{diagonal conj a-i<i}.
        \item 
        When $p = 5$, $r = 116$ and $i = 3$,
        so that $a = 4$, $n = 0$, $s = 16$ and $t = 2$, then 
        \cite{Roz} shows that  $V_r^{(2)}/V_r^{(3)}$ contributes,
        corroborating part $(ii)$ $(b)$ of  Theorem~\ref{diagonal conj a-i<i}. 
        \end{itemize}
        \end{remark}
        
The proof of the above theorem is lengthy. It will be based on the following two steps.  In Theorem~\ref{eliminating JH dcad} in \S\ref{sec: above diagonal} (resp. Theorem~\ref{eliminating JH dcbd} in \S\ref{sec: below diagonal}), we eliminate the JH factors lying strictly Above (resp. strictly Below) the diagonal, i.e., we prove: 
\begin{itemize}
	\item[(A)] If $1 \leq t \leq T$, $T \leq n+1$ and $n+T \leq 2i-a-1$, then  $V_{r}^{(i-T)}/V_{r}^{(i-T+1)}$ dies in $\bar{\Theta}_{k,a_p}$.
	\item[(B)] If $0  \leq T \leq t-2$, $T\leq n$ and $n+T < 2i-a-1$, then  $V_{r}^{(i-T)}/V_{r}^{(i-T+1)}$ dies in $\bar{\Theta}_{k,a_p}$.
\end{itemize}
  Assuming these results, we prove Theorem~\ref{diagonal conj a-i<i}. 

\begin{proof}[Proof of Theorem~\ref{diagonal conj a-i<i}] By Theorem~\ref{shallow JH die a-i<i in bad cases}, we see that the JH factors coming from $V_{r}^{(m)}/V_{r}^{(m+1)}$ vanish for $0 \leq m \leq a-i-1$.  
	\begin{itemize}
		\item[$(i)$] By Lemma~\ref{JH factors bad congruences a<2i case} $(i)$, we have 
		JH factors of $Q(i)$ are 
		\begin{center}
			JH factors of $\{V_r^{(l)}/V_r^{(l+1)}: i-n \leq l \leq i \} ~\cup  ~\{ V_{p-1-a+2l} \otimes D^{a-l} :  0\leq l \leq a-i-1\}$.
		\end{center}
		Since $V_{p-1-a+2l} \otimes D^{a-l}$ is the cosocle of $V_r^{(l)}/V_r^{(l+1)}$ for $0\leq l \leq a-i-1$, from above it follows that these JH factors  vanish in $\bar{\Theta}_{k,a_p}$. To prove the theorem in this case, we show all but one of the following quotients
		\begin{align*}
			\{V_r^{(l)}/V_r^{(l+1)}: i-n \leq l \leq i \}
		\end{align*}
		vanish in $\bar{\Theta}_{k,a_p}$.
		\begin{itemize}
			\item[$(a)$] If $T \leq n$, then we have $n+T \leq 2n  =2i-a-1$.  Thus it follows from (A) that the JH factors of $V_{r}^{(i-T)}/V_{r}^{(i-T+1)}$ vanish in $\bar{\Theta}_{k,a_p}$ for all $t \leq T \leq n$. Also, if $0 \leq T \leq t-2$, then $n+T \leq n+t-2 \leq 2n -2 <  2i-a-1$.  Hence it follows from (B) that the JH factors of $V_{r}^{(i-T)}/V_{r}^{(i-T+1)}$ vanish in $\bar{\Theta}_{k,a_p}$ for all $0 \leq T \leq t-2$. This shows that all the $V_{r}^{(i-T)}/V_{r}^{(i-T+1)}$ vanish in $\bar{\Theta}_{k,a_p}$ for $0 \leq T \leq n$ except for $T=t-1$.
			\item[$(b)$] If $0 \leq T \leq n-1$, then $n+T \leq 2n-1 =  2i-a-2 $ and $T \leq n-1 \leq t-2$. Thus by (B) we have  the JH factors of $V_{r}^{(i-T)}/V_{r}^{(i-T+1)}$ vanish in $\bar{\Theta}_{k,a_p}$ for all $0 \leq T \leq n-1$. The theorem follows in this case also. 
		\end{itemize}
		\item[$(ii)$] A similar argument as in part $(i)$ using  Lemma~\ref{JH factors bad congruences a<2i case}  $(ii)$, we see that it is enough to show all but one of the following 
		\begin{align*}
			\{V_r^{(l)}/V_r^{(l+1)}: i-n \leq l \leq i \} \cup  \{ V_{2i-2n-2-a} \otimes D^{a-i+n+1} \} 
		\end{align*}
		vanish in $\bar{\Theta}_{k,a_p}$. 
		\begin{itemize}
			\item[$(a)$] The proof is very similar to part $(a)$ of $(i)$. If $T \leq n+1$, then we have $n+T \leq 2n+1  \leq 2i-a-1$.  Thus it follows from (A) that the JH factors of $V_{r}^{(i-T)}/V_{r}^{(i-T+1)}$ vanish in $\bar{\Theta}_{k,a_p}$ for all $t \leq T \leq n+1$. Also, if $0 \leq T \leq t-2$, then $n+T \leq n+t-2 \leq 2n -1 <  2i-a-1$.  Hence it follows from (B) that the JH factors of $V_{r}^{(i-T)}/V_{r}^{(i-T+1)}$ vanish in $\bar{\Theta}_{k,a_p}$ for all $0 \leq T \leq t-2$. This shows that all the $V_{r}^{(i-T)}/V_{r}^{(i-T+1)}$ vanish in $\bar{\Theta}_{k,a_p}$ for $t \leq T \leq n+1$ except for $T=t-1$.
			\item[$(b)$] The proof is very similar to part $(b)$ of $(i)$. If $0 \leq T \leq n$, then $n+T \leq 2n <  2i-a-1 $ and $T \leq n \leq t-2$.  Thus by (B) we have  the JH factors of $V_{r}^{(i-T)}/V_{r}^{(i-T+1)}$ vanish in $\bar{\Theta}_{k,a_p}$ for all $0 \leq T \leq n$. The theorem follows in this case also. \qedhere
		\end{itemize}
	\end{itemize}
\end{proof}  

\subsubsection{Above the diagonal.}\label{sec: above diagonal}

In this subsection, we eliminate the Jordan-H\"{o}lder factors coming from  $V_{r}^{(i-T)}/V_{r}^{(i-T+1)}$ for $t \leq T$ (above the diagonal in Figures~\ref{fig:diagonal  conj. a=2i-2n-1} and \ref{fig:diagonal conj. a<2i-2n-1}).

We introduce some notation. 
For a set $A$, the indicator function of $A$ denoted by $\mathbf{1}_{A}(\cdot)$ is  given by
\begin{align*}
	\mathbf{1}_{A}(x) = \begin{cases}
		1 \quad &\mathrm{if}~ x \in A,\\
		0 &\mathrm{if}~ x\not\in A.
	\end{cases}
\end{align*}
Let $\mathbb{N}$ denote the set of natural numbers. In the next few lemmas we show that certain systems of linear equations/congruences are solvable. 

\begin{lemma}\label{adhoc eqn solvable dcad}
	Let  $r \equiv a \mod{(p-1)}$ with $ 1 \leq a \leq p-1 $ and $r \geq i(p+1)+p $ with $ v(a_{p})  \in (i,i+1) $. Let $ s =  a-i+n+(i-n)p$ and $v(r-s) =t $ for some $t\geq 1$. Assume that 
	\begin{enumerate}
		\item[$(i)$] $t \leq T \leq n+1$
		\item[$(ii)$]  $  0 \leq n < i < a$ and $n+T \leq 2i-a-1$.
	\end{enumerate}
	Then there exists $ \beta_0, \ldots, \beta_{i-T} \in \mathbb{Z}_{p}$  with $p \mid \beta_{l}$ for $a-i+n+1 \leq l \leq i-T-1$ satisfying
	\begin{small}
		\begin{align*}
			\sum_{l=0}^{i-T} &\beta_l  \left( \sum_{ \substack{a-i+T <  j < s-i+T \\  j \equiv a-i+T ~\mathrm{mod}~ (p-1)} }   \binom{s-l}{j} \binom{j}{m} + \Delta_{m,l} \right) = \delta_{i-T,m} p^t \quad \mathrm{for} ~ m=0,\ldots, i-T  
		\end{align*}
	\end{small}%
	where
	\begin{small}
		\begin{align*}
			\Delta_{m,l} = &- \frac{(r-s)(a-i+n-l)!(p-1)!}{(p+a-2i+n+T)!(i-T-l)!}  \binom{s-i+T}{m} \\ &\qquad + \frac{(r-s)(a-i+n-l)!(p-1)!}{m!(p+a-m+n-i-l)!}
			\mathbf{1}_{\N}(l+m-a+1) \binom{p+a-l-m-1}{p+a-i+T-m-1} .
		\end{align*}
	\end{small}%
\end{lemma}
\begin{proof}
	Note that $a-i+n+(i-n)p = a +(i-n)(p-1) \leq a+i(p-1) \leq (i+1) p < r $. Hence $r>s$. 
	Let 
	\begin{align*}
		a_{m,l} &= \sum\limits_{ \substack{a-i+T <  j < s-i+T \\  j \equiv a-i+T ~\mathrm{mod}~ (p-1)} }   \binom{s-l}{j} \binom{j}{m}  + \Delta_{m,l} \\
		&= 
		\sum_{k=1}^{i-n-1} \binom{s-l}{a-i+T+k(p-1)} \binom{a-i+T+k(p-1)}{m}  + \Delta_{m,l}. 
	\end{align*}
	Set $A = (a_{m,l})_{0 \leq m,l \leq i-T}$. Using the hypotheses of the lemma it  can be checked that $a_{m,l} \in \Z_p$. 
	Then the system of equations in the lemma can be written as 
	\begin{align*}
		A \begin{bmatrix}
			\beta_0 \\ \vdots \\ \beta_{i-T}
		\end{bmatrix} = \begin{bmatrix}
			0 \\ \vdots \\ 0 \\ p^t
		\end{bmatrix}.
	\end{align*}
	Let $A_{m,l}$ be the minor of $A$ 
	corresponding to the entry $a_{m,l}$. Then, by Cramer's rule, we have $$\beta_{l} = \dfrac{(-1)^{i-T+l} p^t \det(A_{i-T,l})}{\det(A)}$$ provided $\det(A) \neq 0$. Thus, to prove the lemma, it suffices to show that
	\begin{enumerate}[label=(\arabic*)]
		\item $\det(A) \neq 0$. This ensures that the equations have solutions in $\Q_p$.
		\item  $v(\det(A)) \leq  t + v(\det(A_{i-T,l})) $ for all $l$. This shows that all the solutions belong to $\mathbb{Z}_p$.
		\item   $p \mid \beta_{l}$ for $a-i+n+1 \leq l \leq i-T-1$. 
	\end{enumerate}
	\textbf{Computing the determinant of $A$:} Observe that $A$ factors as $BC$, where  $B =(b_{m,k})$ and $C =(c_{k,l})$ with 
	\begin{align}\label{entry of B dcad}
		b_{m,k} = \begin{cases}
			\binom{a-i+T+k(p-1)}{m} ~ &\mathrm{if}~ 0 \leq m \leq  i-T~\mathrm{and}~ 1 \leq k \leq i-n, \\
			\frac{r-s}{m!} \delta_{k-1,m}~ &\mathrm{if}~ 0 \leq m \leq  i-T~\mathrm{and}~ i-n+1 \leq k \leq i-T+1, 
		\end{cases}
	\end{align}
	and 
	\begin{small}
		\begin{align}\label{entry of C dcad}
			c_{k,l} = \begin{cases}
				\binom{s-l}{a-i+T+k(p-1)} ~ &\mathrm{if}~   1\leq k \leq i-n-1 ~\& ~0 \leq l \leq i-T,    \\
				\frac{(a-i+n-l)!(s-r)(p-1)\cdots(p+a-2i+n+T)}{(i-T-l)!}~ &\mathrm{if}~   k = i-n ~\mathrm{and} ~ 0 \leq l \leq  i-T, 
				\vspace{2mm} \\
				(a-i+n-l)!(p-1)\cdots(p+a-i+n+2-l-k) &\mathrm{if}~  i-n+1 \leq k \leq i-T+1 
				\\ ~~~~ \times \mathbf{1}_{\mathbb{N}}(l+k-a) \binom{p+a-l-k}{p+a-i+T-k} &\qquad \qquad \quad  \mathrm{and} ~ 0 \leq l \leq  i-T.\\
			\end{cases}
		\end{align}
	\end{small}%
	Hence, $\det(A) = \det(B)\det(C)$.
	
	\noindent \textbf{Computing the determinant of $B$:}
	Let 
	\begin{alignat*}{2}
		B_1 &= (b_{m,k})_{\substack{m=0, \ldots, i-n-1 \\ k=1, \ldots, i-n}}, \qquad &&B_{2} =(b_{m,k})_{\substack{m=0, \ldots, i-n-1 \\ k=i-n+1 \ldots i-T+1}}, \\
		B_{3} &= (b_{m,k})_{\substack{m=i-n \ldots  i-T \\ k=1, \ldots, i-n}}, \qquad  &&B_4= (b_{m,k})_{\substack{m=i-n,\ldots, i-T \\ k=i-n+1 \ldots i-T+1} }
	\end{alignat*}
	be the submatrices of $B$.  Note that $\delta_{k-1,m} = 0$ whenever $0 \leq m \leq i-n-1$ and $ i-n+1 \leq k \leq i-T+1$. Thus, $B_2$ is the zero matrix, and 
	\begin{align*}
		B = \left(\begin{array}{c|c} B_1 & B_2 \\ \hline B_3 & B_4 \end{array}\right) = \left(\begin{array}{c|c} B_1 & \mathbf{0}\\ \hline B_3 & B_4 \end{array}\right). 
	\end{align*}
	Using the determinant formula for block matrices, we have $\det(B) = \det(B_1)  \det(B_4)$. Since $\delta_{k-1,m} \neq 0$ if and only if $k=m+1$, it follows from \eqref{entry of B dcad} that $B_4$ is a diagonal matrix and $$\det(B_4) = (r-s)^{n+1-T} \prod\limits_{i-n \leq m \leq i-T} \frac{1}{m!}.$$
	Applying Corollary~\ref{cor: GV det}
	with $d$ there equal to $p-1$, $k=i-n$, $m=a-i+T+(p-1)$ and $n=0$, we see that 
	$\det(B_1) = (p-1)^{(i-n)(i-n-1)/2}$. Thus, 
	\begin{align}\label{det B dcad final}
		\det(B) = (p-1)^{(i-n)(i-n-1)/2} (r-s)^{n+1-T} \prod_{i-n \leq m \leq i-T} \frac{1}{m!}.
	\end{align}
	\textbf{Computing the determinant of $C$:} Consider the following submatrices of $C$
	\begin{alignat*}{3}
		C_1 &= (c_{k,l})_{\substack{ k=1, \ldots, i-n \\ l =0, \ldots,a-i+T-1 }}, \quad && C_{2} =(c_{k,l})_{\substack{ k=1, \ldots, i-n \\ l=a-i+T, \ldots, a-i+n}}, \quad 
		C_{3} &&= (c_{k,l})_{\substack{ k=1, \ldots, i-n \\ l= a-i+n+1, \ldots ,i-T }},  \\ 
		C_4 &= (c_{k,l})_{\substack{ k= i-n+1, \ldots, i-T+1 \\ l =0, \ldots,a-i+T-1 }}, \quad &&C_5= (c_{k,l})_{\substack{k= i-n+1, \ldots, i-T+1 \\ l= a-i+T,\ldots, a-i+n}} \quad C_{6} &&= (c_{k,l})_{\substack{ k=i-n+1, \ldots, i-T+1  \\ l= a-i+n+1, \ldots ,i-T }}.
	\end{alignat*}
	Since  $\mathbf{1}_{\mathbb{N}}(l+k-a)=0$ if $l+k \leq a$, we obtain $c_{k,l} = 0$ if  $0 \leq l \leq a-i+T$  and $i-n+1 \leq k \leq i-T$. Also,  $(a-i+n-l)! =0$ if $l\geq a-i+n+1$, so we get $c_{k,l} = 0$ if  $a-i+n+1 \leq l \leq i-T$  and $i-n+1 \leq k \leq i-T$. Hence $C_4$ and $C_6$ are zero matrices. Thus  we have 
	\begin{align*}
		C = \left(\begin{array}{c|c|c} C'_{1} & C'_{2} & C'_{3} \\ \hline \mathbf{0}_{(n-T+1)\times(a-i+T)} &  C'_5 & \mathbf{0}_{(n-T+1)\times(2i-a-n-T)}  \end{array} \right). 
	\end{align*} 
	By interchanging the columns and then applying the  determinant formula for block matrices, we get 
	\begin{align*}
		\det(C) = (-1)^{(2i-a-n-T)(n-T+1)} \det\left( \begin{array}{c|c}
			C_1 & C_3
		\end{array}\right) \det(C_5).
	\end{align*}
	Again using the fact $\mathbf{1}_{\mathbb{N}}(l+k-a)=0$ if $l+k\leq a$ we obtain $C_5$ is an anti-lower triangular matrix, i.e. entries above anti-diagonal ($l+k = a+1$) are zero and 
	\begin{align}\label{det C dcad 1}
		\begin{split}	
			\det(C_5) &= (-1)^{\binom{n+1-T}{2}} \prod_{k=i-n+1}^{i-T+1} (k-i+n-1)! (p-1) \cdots (p-i+n+1)\binom{p-1}{p+a-i+T-k} \\ 
			&= (-1)^{\binom{n+1-T}{2}} \prod_{k=0}^{n-T} k! (p-1) \cdots (p-i+n+1)\binom{p-1}{2i-a-n-T+k}.
		\end{split}	
	\end{align}
	We now determine $\det\left( \begin{array}{c|c}
		C_1 & C_3
	\end{array}\right) $. Pulling out the constant $(s-r) (p-1)\cdots (p+a-2i+n+T)$ from the last row of $\left( \begin{array}{c|c}
		C_1 & C_3
	\end{array}\right) $ and using \eqref{entry of C dcad}, we get 
	\begin{align}\label{det C dcad 2}
		\det\left( \begin{array}{c|c}
			C_1 & C_3
		\end{array}\right) = (-1)^{a-2i+n+T} (2i-a-n-T)! (r-s)  \det\left( \begin{array}{c|c}
			C'_1 & C_3 \\
		\end{array}\right)~\mathrm{mod}~p^{t+1},
	\end{align}
	where $C'_1=(c'_{k,l})_{\substack{ k=1 \ldots  i-n \\ l = 0, \ldots, a-i+T-2 }}$ with 
	\begin{align*}
		c'_{k,l} = \begin{cases}
			\binom{s-l}{a-i+T+k(p-1)} &~\mathrm{if}~ k=1, \ldots, i-n-1  ~ \mathrm{and}~ l=0, \ldots, a-i+T-1 \\
			\frac{(a-i+n-l)!}{(i-T-l)!} & ~\mathrm{if}~ k=i-n  ~\mathrm{and}~ l=0, \ldots, a-i+T-1.
		\end{cases} 
	\end{align*} 
	We now compute $\det\left( \begin{array}{c|c}
		C'_1 & C_3 \\
	\end{array}\right)$ modulo $p$ and show that it doesn't vanish. By Lucas' theorem, modulo $p$ we have 
	\begin{align*}
		\binom{s-l}{a-i+T+k(p-1)} \equiv 
		\begin{cases}
			\binom{i-n}{k} \binom{a-i+n-l}{a-i+T-k} &\mathrm{if}~ 1 \leq k \leq a-i+T \\
			& \qquad \mathrm{and}~  0\leq l \leq a-i+T-1,\\
			\binom{i-n-1}{k}  \binom{p+a-i+n-l}{a-i+T-k} & \mathrm{if}~  1 \leq k \leq a-i+T \\
			&\qquad \mathrm{and}~ a-i+n+1\leq  l \leq i-T,\\
			\binom{i-n}{k-1} \binom{a-i+n-l}{p+a-i+T-k} &\mathrm{if}~ a-i+T+1 \leq k \leq i-n-1 \\ & \qquad \mathrm{and}~ 0\leq l \leq a-i+T-1,\\
			\binom{i-n-1}{k-1} \binom{p+a-i+n-l}{p+a-i+T-k} & \mathrm{if}~  a-i+T+1 \leq k \leq i-n-1 \\
			&\qquad \mathrm{and}~  a-i+n+1\leq l \leq i-T.
		\end{cases}
	\end{align*}
	If $k \leq i-n-1$, then $a-i+n-l < a-k-l < p+a-i+T-k$. Thus $\binom{a-i+n-l}{p+a-i+T-k} =0$ if $a-i+T+1 \leq k \leq i-n-1$ and $0 \leq l \leq a-i+T-1$.  Thus modulo $p$ we have
	\begin{align*}
		\left(\begin{array}{c|c}
			C'_1 & C_3 \\
		\end{array}\right) \equiv  \left( \begin{array}{c|c} \begin{pmatrix}  \binom{i-n}{k}\binom{a-i+n-l}{a-i+T-k}\end{pmatrix}_{\substack{k=1,\ldots, a-i+T \\ l=0,\ldots, a-i+T-1}}  &  \begin{pmatrix} \binom{i-n-1}{k}\binom{p+a-i+n-l}{a-i+T-k} \end{pmatrix}_{\substack{k=1,\ldots, a-i+T \\ l=a-i+n+1,\ldots, i-T}}  \\ \hline \mathbf{0}_{(2i-a-n-T-1) \times (a-i+T)} & \begin{pmatrix}  \binom{i-n-1}{k-1}\binom{p+a-i+n-l}{p+a-i+T-k}\end{pmatrix}_{\substack{ k = a-i+T+1,\ldots, i-n-1 \\ l=a-i+n+1, \ldots, i-T}} \\ \hline 
			\left(\frac{(a-i+n-l)!}{(i-T-l)!}\right)_{l=0,\ldots, a-i+T-1}  & \mathbf{0}_{1 \times(2i-a-n-T) }
		\end{array}\right).
	\end{align*}
	We wish to make all but one entry in the last row of the matrix above zero by applying some column operations. We need the following identity below while 
	\begin{small}
		\begin{multline*}
			\binom{i-n}{k} \binom{a-i+n-l}{a-i+T-k} -  \frac{a-i+n-l}{i-T-l}\binom{i-n}{k}\binom{a-i+n-l-1}{a-i+T-k} \\ = \frac{i-n}{i-T-l} \binom{i-n-1}{k} \binom{a-i+n-l}{a-i+T-k}.
		\end{multline*}
	\end{small}%
	Successively applying the column operations \begin{align*}
		C_0 &\rightarrow C_0 - \frac{a-i+n}{i-T} C_1 \\
		& \vdots \\
		C_{l} &\rightarrow C_{l} - \frac{a-i+n-l}{i-T-l} C_{l+1} \\
		& \vdots \\
		C_{a-i+T-2} &\rightarrow C_{a-i+T-2} - \frac{n+2-T}{2i-a-2T+2} C_{a-i+T-1}
	\end{align*} 
	and expanding out the determinant using $(a-i+T-1)^{\mathrm{th}}$ entry of  the last row we see that 
	\begin{align*}
		\det\left( \begin{array}{c|c}	C'_1 & C_3 \end{array}\right) & \equiv  (-1)^{a-n+T}  \frac{(n+1-T)!}{(2i-a-2T+1)!} \prod_{k=1}^{a-i+T}  \binom{i-n-1}{k} \\ & \quad \times \prod_{l=0}^{a-i+T-2} \frac{i-n}{i-T-l} \prod_{k=a-i+T+1}^{i-n-1} \binom{i-n-1}{k-1} \\ &  \quad \times \det\left( \begin{array}{c|c} \begin{pmatrix}  \binom{a-i+n-l}{a-i+T-k}\end{pmatrix}_{\substack{k=1,\ldots, a-i+T \\ l=0,\ldots, a-i+T-2}}  &  \begin{pmatrix} \binom{p+a-i+n-l}{a-i+T-k} \end{pmatrix}_{\substack{k=1,\ldots, a-i+T \\ l=a-i+n+1,\ldots, i-T}}  \\ \hline \mathbf{0}_{(2i-a-n-T-1) \times (a-i+T-1)} & \begin{pmatrix}  \binom{p+a-i+n-l}{p+a-i+T-k}\end{pmatrix}_{\substack{ k = a-i+T+1,\ldots, i-n-1 \\ l=a-i+n+1, \ldots, i-T}}
		\end{array}\right)  ~\mathrm{mod}~ p.
	\end{align*} 
	Since $\binom{p+a-2i+n+T}{p+a-i+T-k} =0 $ for $a-i+T+1 \leq k \leq i-n-1$, we get that the last column of the matrix in the  lower right corner block is zero. Thus by moving the last column of $C_3$ to the last column of $C_1'$ to obtain a block upper triangular matrix  and then using  the determinant formula for a block matrix, we see that 
	\begin{align*}
		\det\left( \begin{array}{c|c}	C'_1 & C_3 \end{array}\right) & = (-1)^{a-n+T} (-1)^{2i-a-T-n-1} \frac{(n+1-T)!}{(2i-a-2T+1)!} \\  & \quad \prod_{k=1}^{a-i+T}  \binom{i-n-1}{k} \prod_{l=0}^{a-i+T-2} \frac{i-n}{i-T-l}  \prod_{k=a-i+T+1}^{i-n-1} \binom{i-n-1}{k-1} \\
		& \qquad \times  \det\left( \begin{array}{c|c} \begin{pmatrix}  \binom{a-i+n-l}{a-i+T-k}\end{pmatrix}_{\substack{k=1,\ldots,a-i+T \\ l=0,\ldots, a-i+T-2}}  &  \quad \begin{pmatrix} \binom{p+a-2i+n+T}{a-i+T-k}  \end{pmatrix}_{k=1,\ldots,a-i+T}
		\end{array} \right) \\
		& \quad \times \det_{\substack{ k= a-i+T+1 , \ldots, i-n-1 \\ l=a-i+n+1,\ldots, i-T-1}} \begin{pmatrix}  \binom{p+a-i+n-l}{p+a-i+T-k}\end{pmatrix}  \mod p.
	\end{align*}  
	Since $\binom{p+a-i+n-l}{p+a-i+T-k} = 0$ if $l>n-T+k$, we see that $\begin{pmatrix}  \binom{p+a-i+n-l}{p+a-i+T-k}\end{pmatrix}_{\substack{ k=a-i+T+1,\ldots, i-n-1 \\ l=a-i+n+1,\ldots, i-T-1}} $ is a lower triangular matrix with all diagonal entries equal to $1$. Hence 
	\begin{align*}
		\det_{\substack{ k=a-i+T+1,\ldots, i-n-1 \\ l=a-i+n+1,\ldots, i-T-1}} \begin{pmatrix}  \binom{p+a-i+n-l}{p+a-i+T-k}\end{pmatrix} = 1.
	\end{align*}
	Reversing the rows and reversing all columns except the last,  we see that
	\begin{multline*}
		\det \left(\begin{array}{c|c}
			\begin{pmatrix}  \binom{a-i+n-l}{a-i+T-k}\end{pmatrix}_{\substack{k=1,\ldots,a-i+T \\ l=0,\ldots, a-i+T-2}}  &  \begin{pmatrix} \binom{p+a-2i+n+T}{a-i+T-k}  \end{pmatrix}_{k=1,\ldots,a-i+T}
		\end{array} \right) \\ 
		= (-1)^{a-i+T-1} 
		\det \left(\begin{array}{c|c}
			\begin{pmatrix}  \binom{n-T+2+l}{k}\end{pmatrix}_{\substack{k=0,\ldots,a-i+T-1 \\ l=0,\ldots, a-i+T-2}}  &  \begin{pmatrix} \binom{p+a-2i+n+T}{k}  \end{pmatrix}_{k=0,\ldots,a-i+T-1}
		\end{array} \right),
	\end{multline*}
	which by Lemma~\ref{GV det} applied with $b$ there equal to $0$, $k$ there equal to $a-i+T$ and $ a_1= n+2-T, \ldots, a_{k-1} = a-i+n, a_k = p+a-2i+n+T$, equals 
	\begin{align*}
		(-1)^{a-i+T-1}  \frac{\prod\limits_{l=1}^{a-i+T-1}(p-1-i+T+l)}{(a-i+T-1)!} \equiv \frac{\prod\limits_{l=0}^{a-i+T-2}(i-T-l)}{(a-i+T-1)!} \mod p.
	\end{align*}
	Combining all of the above,  we have  modulo $p$
	\begin{align*}
		\det\left( \begin{array}{c|c}	C'_1 & C_3 \end{array}\right) & \equiv -\frac{(n+1-T)! (i-n)^{a-i+T-1} }{(2i-a-2T+1)! (a-i+T-1)!}  \binom{i-n-1}{a-i+T}  \prod_{k=1}^{i-n-2} \binom{i-n-1}{k}. 
	\end{align*}
	Thus it follows from \eqref{det C dcad 1} and \eqref{det C dcad 2}
	\begin{align}
		\det(C) &\equiv (-1)^{n+1-T} (s-r)  \binom{2i-a-2T+1}{n+1-T}^{-1} \frac{ (i-n)^{a-i+T} }{(a-i+T-1)!}   \prod_{k=1}^{i-n-2} \binom{i-n-1}{k}  \nonumber \\
		& \qquad \quad \times \binom{i-n-1}{a-i+T}   \prod_{k=0}^{n-T} \frac{k! (p-1)!} {(p-i+n)!} \binom{p-1}{2i-a-n-T+k} \mod p^{t+1}.
	\end{align}
	Since $n+1<i<a\leq p-1$, we see that all the terms except $(r-s)$ in the above expression are non-zero modulo $p$. Thus $v(\det(C)) = v(r-s) = t$. Combining this with \eqref{det B dcad final} we get $$v(\det(A)) = v(\det(B))+ t = (n+2-T)v(r-s)=(n+2-T)t$$
	which is finite. Therefore $\det(A) \neq 0$ and this proves (1).
	
	We now prove (2) holds, i.e., $v(\det(A_{i-T,l}))+t \geq v(\det(A))$. If $T=n+1$, we see that $v(\det(A)) = t$ and $v(\det(A_{i-T,l})) \geq 0$ so (2) holds. If $T<n+1$, then by Lemma~\ref{minor trick 1} (applied with $c = a-i+T+p-1$, $M=i-T$ and $N=i-n-2$), we see that $v(\det(A_{i-T,l})) + t \geq (n+2-T)t = v(\det(A))$. Now (2) follows even in the case $T<n+1$. 
	
	We show recursively that $p \mid \beta_l$ for $a-i+T+1 \leq l \leq i-T-1$. Since $p \mid (r-s)$, we have $\Delta_{m,l} \equiv 0 \mod p$. Thus from above we have 
	$$ \sum_{l=0}^{i-T} \beta_l \sum_{ \substack{a-i+T <  j < s-i+T \\  j \equiv a-i+T ~\mathrm{mod}~ (p-1)} }   \binom{s-l}{j} \binom{j}{m} \equiv 0 \mod p.$$
	By \cite[Lemma 2.15]{GR19}, for $a-i+T \leq m \leq i-n-1$ we have
	\begin{multline*}
		\sum_{ \substack{a-i+T <  j < s-i+T \\  j \equiv a-i+T ~\mathrm{mod}~ (p-1)} }   \binom{s-l}{j} \binom{j}{m}   \\
		=\sum_{ \substack{0 \leq   j \leq s-i+T \\  j \equiv a-i+T ~\mathrm{mod}~ (p-1)} }   \binom{s-l}{j} \binom{j}{m}- \binom{s-l}{a-i+T} \binom{a-i+T}{m}  - \binom{s-l}{s-i+T} \binom{s-i+T}{m} \\
		\equiv \binom{s-l}{m} \binom{[a-l-m]}{p+a-i+T-m-1} + \delta_{p-1,[a-i+T-m]}\binom{s-l}{m} - \binom{s-l}{a-i+T} \binom{a-i+T}{m} \\
		- \binom{s-l}{s-i+T} \binom{s-i+T}{m} \mod p. 
	\end{multline*}
	Assume that $a-i+T \leq m \leq i-n-1$. Then $\binom{a-i+T}{m} = \delta_{a-i+T,m} = \delta_{p-1,[a-i+T-m]}$. Thus the second and third terms in the above expression cancel each other.  If $l=i-T$, then the first and fourth term in above expression cancel each other. If $ 0 \leq  l \leq a-m-1$, then $[a-l-m] = a-l-m < p+a-i+T-m$. So the first term vanishes in the above expression for $0 \leq l \leq a-m-1$. Hence 
	\begin{small}
		\begin{align*}
			0 & \equiv \sum_{l=0}^{i-T} \beta_l \sum_{ \substack{a-i+T <  j < s-i+T \\  j \equiv a-i+T ~\mathrm{mod}~ (p-1)} }   \binom{s-l}{j} \binom{j}{m}  \\
			& \equiv \sum_{l=a-m}^{i-T-1} \beta_l \binom{s-l}{m} \binom{p+a-l-m-1}{p+a-i+T-m-1} - \sum_{l=0}^{i-T-1} \beta_l \binom{s-l}{s-i+T} \binom{s-i+T}{m} \mod p. 
		\end{align*}
	\end{small}%
	Taking $m=a-i+T$, we see that $\sum_{l=0}^{i-T-1} \beta_l \binom{s-l}{s-i+T} \binom{s-i+T}{a-i+T} \equiv 0 $ mod $p$. By Lucas' theorem, $\binom{s-i+T}{a-i+T} \not \equiv 0 $ mod $p$. Thus $\sum_{l=0}^{i-T-1} \beta_l \binom{s-l}{s-i+T} \equiv 0 $ mod $p$. Hence
	\begin{small}
		\begin{align*}
			\sum_{l=a-m}^{i-T-1} \beta_l \binom{s-l}{m} \binom{p+a-l-m-1}{p+a-i+T-m-1} \equiv 0 \mod p \quad \mathrm{for}~ m=a-i+T+1, \ldots, i-n-1.
		\end{align*}
	\end{small}%
	Taking $m=a-i+T+1$ and noting that  $\binom{s-i+T+1}{a-i+T+1}$ and $\binom{p-1}{p+a-i+T-m-1}$ are non-zero, we obtain $\beta_{i-T-1} \equiv 0$ mod $p$. Assume that $\beta_{i-T-1}, \ldots, \beta_{l} \equiv 0$ mod $p$ for some $a-i+n+1<l\leq i-T-1$. Taking $ m =a-l+1$, we see that $\beta_{l-1} \binom{s-l+1}{a-l+1}\binom{p-1}{p-i+T+l}$ is a linear combination of  $\beta_{i-T-1}, \ldots, \beta_{l} $. Hence $\beta_{l-1} \binom{s-l+1}{a-l+1}\binom{p-1}{p-i+T+l} \equiv 0$ mod $p$. Since $\binom{p-1}{p-i+T+l} \equiv (-1)^{i-T+l-1}$ mod $p$  and $\binom{s-l+1}{a-l+1} \equiv \binom{p+a-i+n-l+1}{a-l+1} \not \equiv 0$ mod $p$ we get $\beta_{l-1} \equiv 0$ mod $p$. This proves part (3) by recursion. 
\end{proof}
The following lemma shows that $\beta_0, \ldots, \beta_{i-T}$ given by Lemma~\ref{adhoc eqn solvable dcad} satisfy certain additional congruences. 

\begin{lemma}\label{choice beta for dcad}
	Let  $r \equiv a \mod{(p-1)}$ with $ 1 \leq a \leq p-1 $ and $r \geq i(p+1)+p $ with $ v(a_{p})  \in (i,i+1) $. Let $ s =  a-i+n+(i-n)p$ and $v(r-s) =t $ for some $t\geq 1$. Assume that 
	\begin{enumerate}
		\item[$(i)$] $t \leq T \leq n+1$
		\item[$(ii)$] $  0 \leq n < i < a $  and $n+T \leq 2i-a-1$.
	\end{enumerate}
	Then there exist $ \beta_0, \ldots, \beta_{i-T} \in \mathbb{Z}_{p}$  with $p \mid \beta_{l}$ for $a-i+n+1 \leq l \leq i-T-1$, satisfying
	\begin{small}
		\begin{align*}
			\sum_{l=0}^{i-T} &\beta_l  \sum_{ \substack{a-i+T <  j < r-i+T \\  j \equiv a-i+T ~\mathrm{mod}~ (p-1)} }   \binom{r-l}{j} \binom{j}{m} \equiv \begin{cases}
				0 \mod p^t ~ & \mathrm{for}~ m=0,\ldots, a-i+T-1\\
				\delta_{i-T,m} p^t  \mod p^{t+1} ~ &\mathrm{for}~ m=a-i+T,\ldots, i-T \\
				0 \mod p^{t-v(m!)} ~ & \mathrm{for}~ m=i-T+1,\ldots, i+t.
			\end{cases}
		\end{align*}
	\end{small}%
\end{lemma}
\begin{proof}
	Let $\beta_l$ be defined as in Lemma~\ref{adhoc eqn solvable dcad}. We show that these $\beta_l$ solve the linear congruences mentioned in the lemma. 
	
	By Corollary~\ref{cor: binomial sums under congruences 2}, we see that for all $m\geq 0$
	\begin{align*}
		\sum_{ \substack{a-i+T \leq  j \leq r-i+T \\  j \equiv a-i+T ~\mathrm{mod}~ (p-1)} }   \binom{r-l}{j} \binom{j}{m} \equiv \sum_{ \substack{a-i+T \leq   j \leq s-i+T \\  j \equiv a-i+T ~\mathrm{mod}~ (p-1)} }   \binom{s-l}{j} \binom{j}{m} ~\mathrm{mod}~p^{t-v(m!)}.
	\end{align*}
	Since $r \equiv s$ mod $p^t$, it follows from Lemma~\ref{binomial coefficient under congruences} $(i)$ that $\binom{r-l}{a-i+T} \equiv \binom{s-l}{a-i+T}$ mod $p^t$, $\binom{r-i+T}{m} \equiv\binom{s-i+T}{m}$ mod $p^{t-v(m!)}$ and $\binom{r-l}{r-i+T} = \binom{r-l}{i-T-l} \equiv \binom{s-l}{i-T-l} \equiv \binom{s-l}{s-i+T}$ mod $p^t$. Thus, $\binom{r-l}{a-i+T} \equiv \binom{s-l}{a-i+T}$ mod $p^t$ and $\binom{r-l}{r-i+T}\binom{r-i+T}{m} \equiv \binom{s-l}{s-i+T}\binom{s-i+T}{m} \mod p^{t-v(m!)}$. This shows that for all $m \geq 0$, we have 
	\begin{align}\label{eq: red beta dcad 1}
		\sum_{ \substack{a-i+T <  j < r-i+T \\  j \equiv a-i+T ~\mathrm{mod}~ (p-1)} }   \binom{r-l}{j} \binom{j}{m} \equiv \sum_{ \substack{a-i+T <  j < s-i+T \\  j \equiv a-i+T ~\mathrm{mod}~ (p-1)} }   \binom{s-l}{j} \binom{j}{m} ~\mathrm{mod}~p^{t-v(m!)} .
	\end{align}		
	Noting that $p^t\mid(r-s)$, it follows from Lemma~\ref{adhoc eqn solvable dcad} that 
	\begin{align}\label{vanishing mod p dcad}
		\sum_{l=0}^{i-T} &\beta_l  \sum_{ \substack{a-i+T <  j < s-i+T \\  j \equiv a-i+T ~\mathrm{mod}~ (p-1)} }   \binom{s-l}{j} \binom{j}{m} \equiv 0 ~\mathrm{mod}~p^t ~~~~ \mathrm{for}~ m=0,1, \ldots, i-T. 
	\end{align}		 		
	The congruence in the lemma for $m=0,\ldots, a-i+T-1$  is immediate from \eqref{eq: red beta dcad 1} and \eqref{vanishing mod p dcad} noting that $v(m!) =0 $ 
	
	We now consider the cases $m=a-i+T,\ldots, i-T$. We claim that for each $0 \leq l \leq i-T$, we have 
	\begin{align}\label{eq: red dcad 2}
		\beta_l \sum_{ \substack{a-i+T  <  j < r-i+T \\  j \equiv a-i+T ~\mathrm{mod}~ (p-1)} }   \binom{r-l}{j} \binom{j}{m} \equiv \beta_l \sum_{ \substack{a-i+T <  j < s-i+T \\  j \equiv a-i+T ~\mathrm{mod}~ (p-1)} }   \binom{s-l}{j} \binom{j}{m} + \beta_{l} \Delta_{m,l}~\mathrm{mod}~p^{t+1},
	\end{align}
	where $\Delta_{m,l}$ is as defined in Lemma~\ref{adhoc eqn solvable dcad}. Clearly, the  congruences in the lemma for $m=a-i+T,\ldots, i-T$ follow from this claim. Noting that $a-i+T \leq m \leq i-T \leq p-2$ and applying Corollary~\ref{cor: binomial sums under congruences 1}, we get
	\begin{small}
		\begin{multline}\label{eq: binomial sum cong dcad 1}
			\beta_l \sum_{ \substack{0  \leq   j \leq r \\  j \equiv a-i+T ~\mathrm{mod}~ (p-1)} }   \binom{r-l}{j} \binom{j}{m} \equiv   \beta_l  \left(\binom{r-l}{m}-\binom{s-l}{m} \right) \left(\binom{[a-l-m]}{p+a-i+T-m-1}+\delta_{a-i+T,m}\right) \\
			+ \beta_l \sum_{ \substack{0 \leq  j \leq s \\  j \equiv a-i+T ~\mathrm{mod}~ (p-1)} }  \binom{s-l}{j} \binom{j}{m}  \mod p^{t+1}.
		\end{multline}
	\end{small}%
	Subtracting $\beta_l \binom{r-l}{a-i+T} \binom{a-i+T}{m}+\beta_l \binom{r-l}{r-i+T} \binom{r-i+T}{m}$ on both sides of congruence  
	\begin{small}
		\begin{align*}
			\beta_l &\sum_{ \substack{a-i+T <   j < r-i+T \\  j \equiv a-i+T ~\mathrm{mod}~ (p-1)} }  \binom{r-l}{j} \binom{j}{m} \\ & \equiv    \beta_l  \left(\binom{r-l}{m}-\binom{s-l}{m} \right) \binom{[a-l-m]}{p+a-i+T-m-1} + \beta_l \sum_{ \substack{a-i+T <   j < s-i+T \\  j \equiv a-i+T ~\mathrm{mod}~ (p-1)} }  \binom{s-l}{j} \binom{j}{m} \\ \nonumber
			& \qquad  + \beta_l  \left(\binom{s-l}{s-i+T} \binom{s-i+T}{m} -\binom{r-l}{r-i+T} \binom{r-i+T}{m} \right)  ~\mathrm{mod}~p^{t+1}  ,
		\end{align*}
	\end{small}%
	where we have separated out the $j=a-i+T$ term (which cancels away since $\binom{a-i+T}{m} =\delta_{a-i+T,m}$) and the $j=s-i+T$ term in the last sum.
	Thus to  show \eqref{eq: red dcad 2}, it remains to check that
	\begin{multline}\label{eq: red dcad 3}
		\beta_l \Delta_{m,l} \equiv \beta_l  \left(\binom{s-l}{s-i+T} \binom{s-i+T}{m} -\binom{r-l}{r-i+T} \binom{r-i+T}{m} \right) \\ 
		+ \beta_l  \left(\binom{r-l}{m}-\binom{s-l}{m} \right) \binom{[a-l-m]}{p+a-i+T-m-1} ~\mathrm{mod}~p^{t+1}  .
	\end{multline}
	If $l=i-T$, then the two terms on right side of \eqref{eq: red dcad 3} cancel. Since $\Delta_{m,i-T}=0$, we get \eqref{eq: red dcad 3} holds for $l=i-T$. If $l=a-i+n+1, \ldots, i-T-1$, then $p\mid \beta_l$. Using Lemma~\ref{binomial coefficient under congruences}$(i)$ we see that each term on the right hand side of \eqref{eq: red dcad 3} vanish modulo $p^{t+1}$. As $\Delta_{m,a-i+n+1} = \cdots = \Delta_{m,i-T-1} =0$, we get \eqref{eq: red dcad 3} holds for $l=a-i+n+1, \ldots, i-T-1$. Thus it remains to show \eqref{eq: red dcad 3} holds for $l=0, \ldots, a-i+n$. 
	Note that $(r-s) \mid \binom{r-i+T}{m} -\binom{s-i+T}{m} $ and by Lucas' theorem we have $\binom{r-l}{r-i+T}  \equiv 0 $ mod $p$. Thus    $ \binom{r-l}{r-i+T}\binom{r-i+T}{m}  \equiv \binom{r-l}{r-i+T} \binom{s-i+T}{m} $ mod $p^{t+1}$. Hence
	\begin{small}
		\begin{multline*}
			\binom{s-l}{s-i+T} \binom{s-i+T}{m} -\binom{r-l}{r-i+T} \binom{r-i+T}{m}   \\
			\equiv  \left(\binom{s-l}{s-i+T}-\binom{r-l}{r-i+T}\right) \binom{s-i+T}{m} 	~\mathrm{mod}~ p^{t+1}.
		\end{multline*}
	\end{small}%
	From the inequality $a<2i-n-T$, it follows that $s-i+T+1 \leq  s-(a-i+n) \leq s-l$ for all $0 \leq l \leq a-i+n$. As $p\mid (s-(a-i+n))$ it follows from Lemma~\ref{binomial coefficient under congruences} $(iii)$, that
	\begin{small} 
		\begin{align*}
			\binom{s-l}{s-i+T}-\binom{r-l}{r-i+T} &= \binom{s-l}{i-T-l}-\binom{r-l}{i-T-l} \\
			&\equiv \frac{s-r}{s-(a-i+n)}  \binom{s-l}{i-T-l}  \\
			& =  \frac{s-r}{(i-T-l)!} \times \frac{(s-l) \cdots  (s-(a-i+n))  \cdots (s-i+T+1)}{(s-(a-i+n)) }  \\
			& \equiv \frac{(s-r) (a-i+n-l)! (p-1) \cdots (p+a-2i+n+T+1)}{(i-T-l)!} ~\mathrm{mod}~ p^{t+1}, 
		\end{align*}
	\end{small}%
	where we used $s \equiv a-i+n$ mod $p$ and $p^t \mid (r-s)$. Thus, we have
	\begin{small}
		\begin{multline}\label{eq: congruence dcad 1}
			\binom{s-l}{s-i+T} \binom{s-i+T}{m} -\binom{r-l}{r-i+T} \binom{s-i+T}{m}   \\
			\equiv  
			(s-r) \frac{ (a-i+n-l)! (p-1) \cdots (p+a-2i+n+T+1)}{(i-T-l)!} \binom{s-i+T}{m} \mod p^{t+1}.
		\end{multline}
	\end{small}%
	Having simplified the first term on the right-hand side of \eqref{eq: red dcad 3}, we next simplify the second term. Observe that 
	\begin{align*}
		\binom{[a-l-m]}{p+a-i+T-m-1} = \begin{cases}
			\binom{a-l-m}{p+a-i+T-m-1} \quad &\mathrm{if}~ l+m<a, \\
			\binom{p+a-l-m-1}{p+a-i+T-m-1} \quad & \mathrm{if}~ l+m\geq a.
		\end{cases}
	\end{align*}
	Since $0< i-T-l \leq i < a \leq p-1$, we get $\binom{[a-l-m]}{p+a-i+T-m-1} =0$ if $l+m<a$. Thus 
	\begin{small}
		\begin{align}\label{eq: congruence dcad 2}
			\binom{[a-l-m]}{p+a-i+T-m-1} = \mathbf{1}_\mathbb{N}(l+m+1-a) \binom{p+a-l-m-1}{p+a-i+T-m-1} .
		\end{align} 
	\end{small}%
	If $l+m \geq  a$ and $0 \leq l \leq a-i+n$, then $s-l-m+1 \leq s-(a-i+n) \leq s-l$. Thus, by Lemma~\ref{binomial coefficient under congruences} $(iii)$, for $0 \leq l \leq a-i+n$ and $l+m \geq a$ we have 
	\begin{small}
		\begin{align}\label{eq: congruence dcad 3}
			\binom{r-l}{m} - \binom{s-l}{m} &\equiv \frac{r-s}{s-(a-i+n)} \binom{s-l}{m}  \nonumber \\
			& \equiv  \frac{(r-s)(s-l) \cdots (s-(a-i+n)+1)(s-(a-i+n)-1)\cdots(s-l-m+1)}{m!} \nonumber \\
			& \equiv (r-s)  \frac{ (a-i+n-l)! (p-1)!}{m!(p+a-i+n-l-m)!}   \mod p^{t+1}.
		\end{align}%
	\end{small}%
	Substituting  \eqref{eq: congruence dcad 1}, \eqref{eq: congruence dcad 2} and \eqref{eq: congruence dcad 3} into \eqref{eq: red dcad 3}, and comparing with the expression for $\Delta_{m,l}$ completes the proof of \eqref{eq: red dcad 3} in the case $l=0, \ldots, a-i+n$. Hence, the congruences in the lemma follow for $m=a-i+T,\ldots, i-T$. 
	
	We next treat the remaining case $m=i-T+1,\ldots, i+t$. By Corollary~\ref{cor: GV det} $(i)$, we have 
	\begin{align*}
		\det_{\substack{m=0,\ldots, i-n-1 \\ k=1,\ldots, i-n}} \left( \binom{a-i+T+k(p-1)}{m}\right) 
	\end{align*}
    is invertible modulo $p$. Hence the rows of the above matrix span $\mathbb{Z}_{p}^{(i-n)}$. Thus for every $m \geq 0$  there exist constants $d_{m,m'} \in \mathbb{Z}_{p}$ such that 
	\begin{align*}
		\binom{a-i+T+k(p-1)}{m} = \sum_{m'=0}^{i-n-1} d_{m,m'}  \binom{a-i+T+k(p-1)}{m'} \quad \text{for all } k=1, \ldots, i-n.
	\end{align*}
	Noting that every $a-i+T < j < s-i+T$ with $j \equiv a-i+T$ mod $(p-1)$ can be expressed as $a-i+T+k(p-1)$ for some $1 \leq k \leq i-n-1$.  Thus 
	\begin{align}\label{eq: dcad 2}
		\sum_{ \substack{a-i+T <  j < s-i+T \\  j \equiv a-i+T ~\mathrm{mod}~ (p-1)} }   \binom{s-l}{j} \binom{j}{m} = \sum_{m'=0}^{i-n-1} d_{m,m'} \sum_{ \substack{a-i+T <  j < s-i+T  \\  j \equiv a-i+T ~\mathrm{mod}~ (p-1)} }   \binom{s-l}{j} \binom{j}{m'}.
	\end{align}
	By \eqref{vanishing mod p dcad} and noting that $i-n-1 \leq i-T$, we get $$\sum_{ \substack{a-i+T <  j < s-i+T  \\  j \equiv a-i+T ~\mathrm{mod}~ (p-1)} }   \binom{s-l}{j} \binom{j}{m'} \equiv 0  ~\mathrm{mod} ~p^t ~ \mathrm{for}~ 0 \leq m' \leq i-n-1. $$ 
	Thus the left hand side of \eqref{eq: dcad 2} vanishes modulo $p^t$. Hence the congruence in the lemma follows also in the case $m=i-T+1,\ldots, i+t$ by \eqref{eq: red beta dcad 1}. 
\end{proof}

\begin{lemma}\label{choice beta and gamma dcad}
	Let  $r \equiv a \mod{(p-1)}$ with $ 1 \leq a \leq p-1 $ and $r \geq i(p+1)+p $ with $ v(a_{p})  \in (i,i+1) $. Let $ s =  a-i+n+(i-n)p$ and $v(r-s) =t $ for some $t\geq 1$. Assume that 
	\begin{enumerate}
		\item[$(i)$] $t \leq T \leq n+1$
		\item[$(ii)$] $  0 \leq n < i < a $  and $n+T \leq 2i-a-1$.
	\end{enumerate}
	Then there exist $ \beta_0, \ldots, \beta_{i-T} \in \mathbb{Z}_{p}$  with $p \mid \beta_{l}$ for $a-i+n+1 \leq l \leq i-T-1$ and $\gamma_0, \ldots, \gamma_{i-T}$ satisfying
	\begin{small}
		\begin{itemize}
			\item[$(i)$] $\sum\limits_{l=0}^{i-T} (\beta_l+p^t \gamma_l)  \sum\limits_{ \substack{a-i+T <  j < r-i+T \\  j \equiv a-i+T ~\mathrm{mod}~ (p-1)} }   \binom{r-l}{j} \binom{j}{m}  \equiv 
			\delta_{i-T,m} p^t  \mod p^{t+1} ~ \mathrm{for}~ m=0,\ldots, i-T$ 
			
			\item[$(ii)$] $\sum\limits_{l=0}^{i-T} (\beta_l+p^t \gamma_l)  \sum\limits_{ \substack{a-i+T <  j < r-i+T \\  j \equiv a-i+T ~\mathrm{mod}~ (p-1)} }   \binom{r-l}{j} \binom{j}{m} \equiv  0 \mod p^{t-v(m!)} ~  \mathrm{for}~ m=i-T+1,\ldots, i+t$.
		\end{itemize}		
	\end{small}%
\end{lemma}
\begin{proof}
	Let $ \beta_0, \ldots, \beta_{i-T} \in \mathbb{Z}_{p}$ be as defined in Lemma~\ref{choice beta for dcad}. Set $\gamma_0 = \cdots = \gamma_{n-T}=0$ and $\gamma_{a-i+n+1} = \cdots = \gamma_{i-T}=0$.
	Note that  by Lemma~\ref{choice beta for dcad}, we have 
	$$ \sum_{l=0}^{i-T} \beta_l\sum_{ \substack{a-i+T <  j < r-i+T \\  j \equiv a-i+T ~\mathrm{mod}~ (p-1)} }   \binom{r-l}{j} \binom{j}{m} \equiv 0 ~\mathrm{mod}~ p^t~ \text{for all } 0 \leq m \leq i-T.$$
	
	For $(i)$, we need to show that 
	\begin{small}
		\begin{multline}\label{eq: gamma red dcad}
			\sum_{l=n+1-T}^{a-i+n}  \gamma_l \sum_{ \substack{a-i+T <  j < r-i+T \\  j \equiv a-i+T ~\mathrm{mod}~ (p-1)} }   \binom{r-l}{j} \binom{j}{m}  \\
			\equiv \delta_{i-T,m}- p^{-t} \sum_{l=0}^{i-T} \beta_l\sum_{ \substack{a-i+T <  j < r-i+T \\  j \equiv a-i+T ~\mathrm{mod}~ (p-1)} }   \binom{r-l}{j}\binom{j}{m}~\mathrm{mod}~p
		\end{multline}
	\end{small}%
	for $m=0, \ldots,i-T$.
	
	Now assume $m=0, \ldots, a-i+T-1$. By Lucas' theorem we have  $\binom{r-l}{a-i+T} \equiv \binom{a-i+n-l}{a-i+T} \equiv 0$ mod $p$ as $l \geq n+1-T$. Again by Lucas' theorem $\binom{r-l}{r-i+T} = \binom{r-l}{i-T-l} \equiv \binom{a-i+n-l}{i-T-l} \equiv 0$ mod $p$ since $a<2i-n-T$. Thus for $m=0, \ldots, a-i+T-1$ and $l=n+1-T, \ldots, a-i+n$ it follows from \cite[Lemma 2.15]{GR19} that
	\begin{align*}
		\sum_{ \substack{a-i+T <  j < r-i+T \\  j \equiv a-i+T ~\mathrm{mod}~ (p-1)} }   \binom{r-l}{j} \binom{j}{m} & \equiv \sum_{ \substack{a-i+T \leq  j \leq r-i+T \\  j \equiv a-i+T ~\mathrm{mod}~ (p-1)} }   \binom{r-l}{j} \binom{j}{m} \\
		& \equiv  \binom{r-l}{m} \binom{a-l-m}{a-i+T-m} ~\mathrm{mod}~ p \\
		& \equiv \binom{a-i+n-l}{m} \binom{a-l-m}{a-i+T-m} ~\mathrm{mod}~ p .
	\end{align*}
	Thus to show  \eqref{eq: gamma red dcad} holds in the case $m=0,\ldots,a-i+T-1$, it suffices to find $\gamma_{n-T+1},\ldots,  \gamma_{a-i+n} \in \mathbb{Z}_{p}$ satisfying
	\begin{small}
		\begin{align}\label{eq: choice gamma dcad}
			\sum_{l=n-T+1}^{a-i+n} \gamma_{l} \binom{a-i+n-l}{m} \binom{a-l-m}{a-i+T-m} \equiv - \frac{1}{p^{t}} \sum_{l=0}^{i-T} \beta_l\sum_{ \substack{a-i+T <  j < r-i+T \\  j \equiv a-i+T ~\mathrm{mod}~ (p-1)} }   \binom{r-l}{j} \binom{j}{m}~\mathrm{mod}~p
		\end{align} 
	\end{small}%
	for $m=0, \ldots, a-i+T-1$. Note that the matrix $$\left(\binom{a-i+n-l}{m} \binom{a-l-m}{a-i+T-m} \right)_{\substack{m=0,\ldots, a-i+T-1 \\ l=n-T+1,\ldots, a-i+n}}$$ 
	is upper anti-triangular, i.e., entries below the anti-diagonal ($l+m=a-i+n$) vanish. Furthermore, the anti-diagonal entries are non-zero modulo $p$ as $n+T < 2i-a-1$. Hence, \eqref{eq: choice gamma dcad} admits a solution in $\Z_p$. This proves $(i)$ in the case $m=0,\ldots, a-i+T-1$.		  
	
	If $m=a-i+T, \ldots, i-T$, then by Lemma~\ref{choice beta for dcad} we see that the right side of \eqref{eq: gamma red dcad} vanishes modulo $p$. Observe that by Lucas' theorem, we have $\binom{r-l}{m} \equiv \binom{a-i+n-l}{m} \equiv 0$ mod $p$ for $m=a-i+T, \ldots, i-T$ and $l=n+1-T, \ldots, a-i+n$. Hence
	\begin{align*}
		\sum_{ \substack{a-i+T <  j < r-i+T \\  j \equiv a-i+T ~\mathrm{mod}~ (p-1)} }   \binom{r-l}{j} \binom{j}{m} = \binom{r-l}{m}  \sum_{ \substack{a-i+T< j < r-i+T \\  j \equiv a-i+T ~\mathrm{mod}~ (p-1)} }   \binom{r-l-m}{j-m} \equiv 0 ~\mathrm{mod }~p 
	\end{align*}
	for $m=a-i+T, \ldots, i-T$ and $l=n+1-T, \ldots, a-i+n$. So the left side of  \eqref{eq: gamma red dcad}  vanishes modulo $p$. Hence \eqref{eq: gamma red dcad} holds for $m=a-i+T, \ldots, i-T$. This proves $(i)$.
	
	Clearly $(ii)$ follows from Lemma~\ref{choice beta for dcad} noting that $ t \geq t-v(m!)$.
	This finishes the proof of the lemma.
\end{proof}

\begin{theorem}\label{eliminating JH dcad}
	Let  $r \equiv a \mod{(p-1)}$ with $ 1 \leq a \leq p-1 $ and $r \geq i(p+1)+p $ with $ v(a_{p})  \in (i,i+1) $. Let $ s =  a-i+n+(i-n)p$ and $v(r-s) =t $ for some $t\geq 1$. Assume that
	\begin{enumerate}
		\item[$(i)$] $t \leq T \leq n+1$
		\item[$(ii)$] $ 1 \leq i < a $  and $n+T \leq 2i-a-1$.
	\end{enumerate}
	Then the image of $ V_{r}^{(i-T)} $ is the same as the image of $ V_{r}^{(i-T+1)} $ in $ \bar{\Theta}_{k,a_{p}} $.
\end{theorem}
\begin{proof}
	Let $\beta_l$ and $\gamma_l$ be the $p$-adic integers as chosen in Lemma~\ref{choice beta and gamma dcad}. We wish to apply Lemma~\ref{lem:choice of beta}. To do this we need to check the hypotheses of the lemma. We do this check here and leave all such future checks  in this chapter to the reader. Note that $s-a= (i-n)(p-1) < i(p+1) <r-a $. Thus $r>s$. So $r\geq s+p(p-1) \geq a+(p+1)(p-1)$ and $\frac{r-a}{p-1}-2 \geq p-1$. In Lemma~\ref{lem:choice of beta}, take $c=a-i+T+(p-1)$, $k=\frac{r-a}{p-1}-2$, $m'=\min\{i+t,p-1\}$, $\gamma_{c}$ there equal to $\sum_{l=0}^{i-T} \left(\beta_{l}+p^t \gamma_l \right) \binom{r-l}{c}$, $\ldots$, $\gamma_{c+k(p-1)}$ there equal to $\sum_{l=0}^{i-T} \left(\beta_{l}+p^t \gamma_l \right) \binom{r-l}{c+k(p-1)}$, $\nu_0,\ldots,\nu_{m'}=0$  and $n=i+t+1$. So $k\geq p-1\geq m'$.  By Lemma~\ref{choice beta and gamma dcad} $(i)$ and $(ii)$, we also have 
    \begin{align*}
        \sum_{k'=0}^k \gamma_{c+k'(p-1)} \binom{c+k'(p-1)}{m} = 
        \sum_{l=0}^{i-T} \sum_{k'=0}^k \left(\beta_{l}+p^t \gamma_l \right) \binom{r-l}{c+k'(p-1)}\binom{c+k'(p-1)}{m} \equiv 0 \mod p^t
    \end{align*}
    for $m=0,\ldots,m'$. This checks all the hypotheses of Lemma~\ref{lem:choice of beta}.
    Hence, there exist $\alpha_{j} \in \mathbb{Z}_{p}$, for all $j \equiv a-i+T $ mod $(p-1)$ with $ a-i+T< j < r-i+T$ satisfying
	\begin{enumerate}
		\item[$(1)$]  $\alpha_j \equiv \sum_{l=0}^{i-T} \left(\beta_{l}+p^t \gamma_l \right) \binom{r-l}{j} ~\mathrm{mod}~p^{t} $, for all $j$ as above, 
		\item[$(2)$]   $\sum_{j}^{} \alpha_j \binom{j}{m} \equiv 0 $ mod $p^{i+t+1-m}$ for $m=0,\ldots, \min\{i+t, p-1\}$.
	\end{enumerate}
	From the congruence condition $(1)$ and Lemma~\ref{choice beta and gamma dcad}, we also have 
	\begin{enumerate}
		\item[$(2')$] $\sum_{j}^{} \alpha_j \binom{j}{m} \equiv 0 $ mod $p^{t-1}$ for $m=p,\ldots, i+t$. 
	\end{enumerate}
	Let
	\begin{align*}
		f_2 &=  
		\sum_{  \lambda \in \mathbb{F}_{p}^{\times}} 
		\Bigg[ g_{2,p[\lambda]}^{0}, \sum_{l=0}^{i-T} \frac{[\lambda]^{l-(i-T)}}{p^{l+t}} 
		(\beta_{l}+p^t \gamma_l) (-\theta)^{l+t+2} X^{-t-2}Y^{r-(l+t+2)(p+1)+t+2} \Bigg]  \\
		& \quad + \left[ g_{2,0}^{0}, \frac{(1-p)}{p^{i-T+t}}  \sum_{l=0}^{i-T} (\beta_{l} +p^t \gamma_l) \binom{r-l}{r-i+T}  (-\theta)^{i+t+1} X^{-T-t-1}Y^{r-(i+t+1)(p+1)+T+t+1} \right]  \\
		f_1 &=  \left[g_{1,0}^0, \frac{p-1}{p^t a_p} \sum_{ \substack{a-i+T <  j < r-i+T \\  j \equiv a-i+T ~\mathrm{mod}~ (p-1)} } \alpha_{j} X^{r-j} Y^{j} \right] \\
		f_{0} &= \left[ \mathrm{id}, \frac{1-p}{p^{a-i+T+t}}  \sum_{l=0}^{i-T} \left(\beta_{l}+p^t \gamma_l \right) \binom{r-l}{a-i+T} \theta^{a-i+T+t+1} X^{r-(a-i+T+t+1)(p+1)+t+1} Y^{-t-1}\right].
	\end{align*}
	From Lemma~\ref{theta and T plus} it follows that $T^+ f_2$ vanishes modulo $p$. It easy to see that $-a_p f_2$, $T^{+}f_1$, $T^{-}f_1$, $-a_p f_0 $ and $T^{-} f_0$ all vanish modulo $p$ using the $(2)$ and $(2')$ above and the hypotheses $a \leq 2i-T-n-1$, $t \leq T  \leq  n+1$. It can be checked that 
	\begin{align*}
		T^{-} f_2 -a_p f_1 + T^+ f_0 \equiv   \left[g_{1,0}^0, F(X,Y) \right] ~\mathrm{mod}~ p. 
	\end{align*}
	where 
	\begin{align*}
		F(X,Y) = \frac{p-1}{p^t} \sum_{ \substack{a-i+T <  j < r-i+T \\  j \equiv a-i+T ~\mathrm{mod}~ (p-1)} } \left( \sum_{l=0}^{i-T}  \left(\beta_{l}+p^t \gamma_l \right) \binom{r-l}{j}- \alpha_{j} \right) X^{r-j} Y^{j}
	\end{align*}
	By $(1)$ above, we have $F(X,Y) \in \mathbb{Z}_p[X,Y]$.  To prove the proposition it is enough to show that $\bar{F}(X,Y)$ generates $V_{r}^{(i-T)}/V_r^{(i-T+1)}$. Using the properties of $\beta_l$ and $\gamma_l$ from Lemma~\ref{choice beta and gamma dcad} and the choice of $\alpha_{j}$, for $m=0,1,\ldots, i-T$, we have 		     
	\begin{align*}
		\sum_{ \substack{a-i+T <  j < r-i+T \\  j \equiv a-i+T ~\mathrm{mod}~ (p-1)} } &\binom{j}{m} \left( \sum_{l=0}^{i-T}  \left(\beta_{l}+p^t \gamma_l \right) \binom{r-l}{j}- \alpha_{j} \right) \equiv \delta_{i-T,m}p^t ~\mathrm{mod}~p^{t+1}.
	\end{align*}		     
	Since $i-T< a-i+T+ p-1 $, it follows from \cite[Lemma 2.8]{GR19} that $\theta^{i-T} \mid \overline{F(X,Y)}$. From hypothesis $(ii)$ of the theorem, it follows that $0 < a-i+T \leq i-T <  p-1$. Thus either $a-i+T = i-T$ or $ a-i+T \not \equiv i-T ~\mathrm{mod}~(p-1)$. In both cases, we obtain that the coefficient of $X^{r-i+T}Y^{i-T}$ in $F(X,Y)$ is zero. Applying \cite[Lemma 2.12]{GR19}, with $m$ there equal $i-T$ and $l$ there equal to $a-i+T$ we obtain 
	\begin{align*}
		\overline{F(X,Y)} & \equiv (p-1) \theta^{i-T} X^{r-(i-T)(p+1)-[a-2i+2T]}Y^{[a-2i+2T]} ~\mathrm{mod}~V_r^{(i-T+1)} \\
		& = (p-1) \theta^{i-T} X^{r-(i-T)(p+1)-(p-1+a-2i+2T)}Y^{p-1+a-2i+2T}  ~\mathrm{mod}~V_r^{(i-T+1)},
	\end{align*}  
	since $a \leq 2i-2T $ by the hypotheses $T \leq n+1$ and $n+T \leq 2i-a-1$.
	Applying Lemma~\ref{Glover-Brueil map image} with $m$ there equal to $i-T$ and $b$ there equal to $ p-1+a$, it follows that $\overline{F(X,Y)} $ generates $V_r^{(i-T)}/V_r^{(i-T+1)}$. This finishes the proof of the theorem.
\end{proof}

\subsubsection{Below the diagonal.}\label{sec: below diagonal}

Let the setting be as in Theorem~\ref{diagonal conj a-i<i}. In this subsection, we eliminate the Jordan H\"{o}lder factors coming from $V_{r}^{(i-T)}/V_{r}^{(i-T+1)}$ whenever $0 \leq T \leq t-2$ (below the diagonal in Figures~\ref{fig:diagonal  conj. a=2i-2n-1} and \ref{fig:diagonal conj. a<2i-2n-1}). In order to do this, we start with a preparatory lemma which guarantees a certain choice of $p$-adic integers.  
\begin{lemma}\label{choice beta dcbd}
	Let  $r \equiv a \mod{(p-1)}$ with $ 1 \leq a \leq p-1 $ and $r \geq i(p+1)+p $ with $ v(a_{p})  \in (i,i+1) $. Let $ s =  a-i+n+(i-n)p$ and $v(r-s) =t $ for some $t\geq 1$. Assume that 
	\begin{enumerate}
		\item[$(i)$]  $t \geq T+2$ 
		\item[$(ii)$]  $ 1 \leq i < a $, $0 \leq T \leq n < i-1$ and $n+T \leq 2i-a-1$.
	\end{enumerate}
	Then there exist $\beta_{n-T}, \ldots, \beta_{i-T-1} \in \mathbb{Z}_{p}$ with $p\mid \beta_{i-T-1}$ satisfying
	\begin{enumerate}
		\item[$(1)$] $\sum\limits_{l=n-T}^{i-T-1} \beta_{l} \sum\limits_{ \substack{a-i+T <   j < r-i+T \\  j \equiv a-i+T ~\mathrm{mod}~ (p-1)} }\binom{r-l}{j} \binom{j}{m} ~\equiv 0~ \mathrm{mod}~ p^{T+2-v(m!)} $ for $m=0, \ldots, i+T+3$
		\item[$(2)$] $\sum\limits_{l=n-T}^{i-T-1} \beta_{l} \binom{r-l}{r-i+T} \equiv -p~\mathrm{mod}~p^{T+2}$.
	\end{enumerate}
\end{lemma}
\begin{proof}
	Let $s=a-i+n+(i-n)p$.  As $r \equiv s$ mod $p^{T+2}$, by Corollary~\ref{cor: binomial sums under congruences 2} we have 
	\begin{align*}
		\sum\limits_{ \substack{a-i+T <   j < r-i+T \\  j \equiv a-i+T ~\mathrm{mod}~ (p-1)} }\binom{r-l}{j} \binom{j}{m} \equiv \sum\limits_{ \substack{a-i+T <   j < s-i+T \\  j \equiv a-i+T ~\mathrm{mod}~ (p-1)} }\binom{s-l}{j} \binom{j}{m}~\mathrm{mod}~p^{T+2-v(m!)}
	\end{align*}
	for $l=n-T,\ldots, i-T-1$. Also by Lemma~\ref{binomial coefficient under congruences} we have $\binom{r-l}{r-i+T} = \binom{r-l}{i-T-l} \equiv \binom{s-l}{i-T-l}  = \binom{s-l}{s-i+T} $ mod $p^{T+2}$. Thus it is enough to show that the lemma holds for the special value  $s$.  We first show the existence of $\beta_{n-T}, \ldots, \beta_{i-T-1} \in \mathbb{Z}_{p}$ with $p \mid \beta_{i-T-1}$ satisfying
	\begin{align}
		\sum\limits_{l=n-T}^{i-T-1} \beta_{l} \sum\limits_{ \substack{a-i+T <   j < s-i+T \\  j \equiv a-i+T ~\mathrm{mod}~ (p-1)} }\binom{s-l}{j} \binom{j}{m} & = 0~ \quad  \mathrm{for}~ m=0, \ldots, i-n-2 \label{red eq 1 dcbd}\\
		\sum\limits_{l=n-T}^{i-T-1} \beta_{l} \binom{s-l}{s-i+T} = -p.\label{red eq 2 dcbd}
	\end{align}
	Then above system of equations  can be written as 
	\begin{align*}
		A \begin{bmatrix}
			\beta_{n-T} \\ \vdots \\ \beta_{i-T-1}
		\end{bmatrix} = \begin{bmatrix}
			0 \\ \vdots \\ 0 \\ -p
		\end{bmatrix}.
	\end{align*}
	where 
	\begin{align*}
		A = \left[ \begin{array}{c}
			\left(\sum\limits_{ \substack{a-i+T <   j < s-i+T \\  j \equiv a-i+T ~\mathrm{mod}~ (p-1)} }\binom{s-l}{j} \binom{j}{m}\right)_{\substack{m=0,\ldots, i-n-1\\ l=n-T, \ldots, i-T-1}} \\ \hline
			\left( \binom{s-l}{s-i+T} \right)_{l=n-T, \ldots, i-T-1}
		\end{array}\right]. 
	\end{align*}
	Noting that every $a-i+T < j < s-i+T$ with $j \equiv a-i+T$ mod $(p-1)$ can be expressed as $a-i+T+k(p-1)$ for some $1 \leq k \leq i-n-1$, we have $A = BC$ with 
	\begin{small}
		\begin{align*}	
			B &= \left[\begin{array}{c|c}
				\left(\binom{a-i+T+k(p-1)}{m}\right)_{\substack{m=0, \ldots, i-n-2 \\ k=1, \ldots, i-n-1}} & \mathbf{0}_{(i-n-1)\times 1} \\ \hline
				\mathbf{0}_{1 \times (i-n-1)} & 1 
			\end{array}\right] \\
			C &= \left[ \binom{s-l}{a-i+T+k(p-1)}_{\substack{k=1, \ldots, i-n \\ l=n-T, \ldots, i-T-1}}\right]. 
		\end{align*}
	\end{small}%
	Applying Corollary~\ref{cor: GV det}
	with $d$ there equal to $p-1$, $k=i-n-1$, $m=a-i+T+p-1$ and $n=0$ we see that $\det(B) = (p-1)^{(i-n-1)(i-n-2)/2}$. Hence $\det(B)$ is a unit in $\Z_p$.  Applying \cite[Lemma 8]{Viennot} $(a-i+T+1)$ times we see that
	\begin{align}\label{det C dcbd}
		\det(C) 
		& =  \prod_{u=0}^{a-i+T} \left( \prod_{l=n-T}^{i-T-1}(s-l-u) \times  \prod_{k=1}^{i-n} \frac{1}{(a-i+T+k(p-1)-u)} \right) \nonumber \\
		& \qquad \qquad \times \det_{\substack{k=1, \ldots, i-n \\ l=n-T, \ldots, i-T-1}} \left( \binom{s-l-(a-i+T+1)}{k(p-1)-1}\right)\nonumber  \\
		& = \prod_{u=0}^{a-i+T} \left( \prod_{l=1}^{i-n}(s-n+T+1-l-u) \times \prod_{k=1}^{i-n} \frac{1}{(a-i+T+k(p-1)-u)} \right) \nonumber \\
		& \qquad \qquad \times \det_{\substack{k=1, \ldots, i-n \\ l=1, \ldots, i-n}} \left( \binom{s-l-(a-i+n)}{k(p-1)-1}\right) \\
		& = \prod_{k=1}^{i-n} \frac{(s-n+T+1-k)_{a-i+T+1}}{(a-i+T+k(p-1))_{a-i+T+1}}  \times \det_{\substack{k=1, \ldots, i-n \\ l=1, \ldots, i-n}} \left( \binom{s-l-(a-i+n)}{k(p-1)-1}\right)  \nonumber.
	\end{align} 
	By Lucas' theorem and noting $s=(a-i+n)+(i-n)p$ we have
	\begin{align*}
		\binom{s-l-(a-i+n)}{k(p-1)-1} = \binom{(i-n)p-l}{k(p-1)-1} \equiv \binom{i-n-1}{k-1} \binom{p-l}{p-1-k} ~\mathrm{mod}~p.
	\end{align*}
	Thus 
	\begin{small}
		\begin{align}\label{det C' dcbd}
			\det_{\substack{k=1, \ldots, i-n \\ l=1, \ldots, i-n}} \left( \binom{s-l-(a-i+n)-1}{k(p-1)-1}\right) \equiv \prod_{k=1}^{i-n} \binom{i-n-1}{k-1} \times  \det_{\substack{k=1, \ldots, i-n \\ l=1, \ldots, i-n}} \left( \binom{p-l}{p-1-k}\right) ~\mathrm{mod}~p.
		\end{align}
	\end{small}%
	Reversing  the order of rows and columns and then applying Lemma~\ref{GV det} on the transpose with $b$ there equal to $p-1-i+n$, $k=i-n$ and $a_1 = p-i+n, a_2 = p-i+n+1, \ldots, a_k = p-1$ we see that
	\begin{align*}
		\det_{\substack{k=1, \ldots, i-n \\ l=1, \ldots, i-n}} \left( \binom{p-l}{p-1-k}\right) 
		&=\det_{\substack{l=1, \ldots, i-n  \\ k=1, \ldots, i-n}} \left( \binom{p-1-i+n+l}{p-2-i+n+k}\right)\\
		&= \frac{(p-i+n)_{p-1-i+n} \cdots (p-1)_{p-1-i+n}}{(p-1-i+n)!\cdots(p-2)!} \prod_{1 \leq l<l'\leq i-n} (l'-l).
	\end{align*}
	Clearly the above expression is a $p$-adic unit. So the right side of \eqref{det C' dcbd} is a $p$-adic unit.  Thus by \eqref{det C dcbd} we have
	\begin{align*}
		v(\det(C)) & = \sum_{k=1}^{i-n} v\left(\frac{(s-n+T+1-k)!}{(s-(a-i+n)-k)!}\right) - v\left(\frac{(a-i+T+k(p-1))!}{(k(p-1)-1)!}\right) 
	\end{align*}
	Note $ s= a-i+n +(i-n)p$ and $v((s-n+T+1-k)!) = \lfloor \frac{(s-n+T+1-k)}{p} \rfloor + \lfloor \frac{(s-n+T+1-k)}{p^2} \rfloor + \cdots$. Thus 
	\begin{align*}
		v((s-n+T+1-k)!) = 
		\begin{cases}
			(i-n) & \mathrm{if}~ 1\leq k\leq a-i+T+1, \\
			(i-n-1) & \mathrm{if}~ a-i+T+1 < k \leq i-n.
		\end{cases}
	\end{align*} 
	Similarly $v((s-(a-i+n)-k)!) = v(((i-n)p-k)!) = i-n-1$ for $1 \leq k \leq i-n$. Thus  
	\begin{align*}
		v((s-n+T+1-k)!) - v((s-(a-i+n)-k)!)= 
		\begin{cases}
			1 & \mathrm{if}~ 1\leq k\leq a-i+T+1, \\
			0 & \mathrm{if}~ a-i+T+2 \leq k \leq i-n.
		\end{cases}
	\end{align*} 
	A similar calculation shows that
	\begin{align*}
		v((a-i+T+k(p-1))!)-v((k(p-1)-1)!) = \begin{cases}
			1 &\mathrm{if}~ 1\leq k\leq a-i+T, \\ 
			0 & \mathrm{if}~a-i+T +1 \leq  k \leq i-n.
		\end{cases} 		  	     
	\end{align*} 
	Thus
	\begin{align*} v\left(\frac{(s-n+T+1-k)!}{(s-(a-i+n)-k)!}\right) - v\left(\frac{(a-i+T+k(p-1))!}{(k(p-1)-1)!}\right) 
		= \begin{cases}
			1 &\mathrm{if}~  k =  a-i+T+1, \\ 
			0 & \mathrm{otherwise}.
		\end{cases} 		
	\end{align*} 
	Since $p\nmid \det(B)$, we obtain $v( \det(A)) = v( \det(B))+v( \det(C)) = 1$. By Cramer's rule, we have $$\beta_{l} = \pm \frac{(-p) \cdot \det(\text{minor of } a_{i-n,l} \text{ entry of } A)}{\det(A)} .$$ 
	Thus $\beta_{l} \in \mathbb{Z}_{p}$.  Hence \eqref{red eq 1 dcbd} and \eqref{red eq 2 dcbd}  have solutions in $\mathbb{Z}_{p}$. 
	
	We next show that $\beta_{i-T-1} \in p\mathbb{Z}_{p}$. By $(ii)$, we have $0 < i-1 < a-1 \leq p-2 $. 
	Thus by \eqref{red eq 1 dcbd} and \eqref{red eq 2 dcbd} we have
	\begin{align*}
		0 &\equiv \sum\limits_{l=n-T}^{i-T-1} \beta_{l} \sum\limits_{ \substack{a-i+T <   j  < s-i+T \\  j \equiv a-i+T ~\mathrm{mod}~ (p-1)}  } \binom{s-l}{j} \binom{j}{a-i+T+1}  \\
		&= \sum\limits_{l=n-T}^{i-T-1} \beta_{l} \sum\limits_{ \substack{a-i+T \leq    j  \leq  s-i+T \\  j \equiv a-i+T ~\mathrm{mod}~ (p-1)}  } \binom{s-l}{j} \binom{j}{a-i+T+1} - \sum\limits_{l=n-T}^{i-T-1} \beta_{l} \binom{s-l}{s-i+T} \binom{s-i+T}{a-i+T+1} \\
		& \equiv \sum\limits_{l=n-T}^{i-T-1} \beta_{l} \sum\limits_{ \substack{a-i+T \leq    j  \leq  s-i+T \\  j \equiv a-i+T ~\mathrm{mod}~ (p-1)}  } \binom{s-l}{j} \binom{j}{a-i+T+1} .
	\end{align*}	    
	By \cite[Lemma 2.15]{GR19}, for $n-T \leq  l \leq i-T-1$ we have
	\begin{align*}
		\sum\limits_{ \substack{a-i+T \leq   j \leq s-i+T \\  j \equiv a-i+T ~\mathrm{mod}~ (p-1)} }\binom{s-l}{j} \binom{j}{a-i+T+1} 
		& \equiv  \binom{s-l}{a-i+T+1} \binom{[i-T-1-l]}{p-2} \\
		& \equiv   \binom{s-l}{a-i+T+1} \cdot (p-1) \delta_{l,i-T-1} ~\mathrm{mod}~p.
	\end{align*} 
	Substituting this above 
	\begin{align*}
		(p-1) \beta_{i-T-1} \binom{s-i+T+1}{a-i+T+1}  \equiv 0~\mathrm{mod}~p.
	\end{align*}
	By Lucas' theorem $\binom{s-i+T+1}{a-i+T+1} \equiv \binom{p+a-2i+n+T+1}{a-i+T+1} \not \equiv 0$ mod $p$. This shows that $\beta_{i-T-1} \in p\mathbb{Z}_{p}$.
	
	We next treat the remaining cases, i.e.,  $m=i-n-1,\ldots, i+T+3$. 
	As we saw, the top left block of $B$ has unit determinant in $\Z_p$. Therefore its rows  span $\mathbb{Z}_{p}^{(i-n-1)}$. Thus for every $m\geq 0$,  there exist constants $d_{m,m'} \in \mathbb{Z}_{p}$ such that 
	\begin{align*}
		\binom{a-i+T+k(p-1)}{m} = \sum_{m'=0}^{i-n-2} d_{m,m'}  \binom{a-i+T+k(p-1)}{m'} \quad \text{for all } k=1, \ldots, i-n-1.
	\end{align*}
	Thus 
	\begin{align*}
		\sum_{ \substack{a-i+T <  j < s-i+T \\  j \equiv a-i+T ~\mathrm{mod}~ (p-1)} }   \binom{s-l}{j} \binom{j}{m} = \sum_{m'=0}^{i-n-2} d_{m,m'} \sum_{ \substack{a-i+T <  j < s-i+T  \\  j \equiv a-i+T ~\mathrm{mod}~ (p-1)} }   \binom{s-l}{j} \binom{j}{m'}.
	\end{align*}
	Hence the congruence in the lemma follows for $m=i-n-1,\ldots, i+T+3$ from  \eqref{red eq 1 dcbd}. 
\end{proof}
In the following theorem, we show that the JH factors coming from $V_{r}^{(i-T)}/V_{r}^{(i-T+1)}$ die in $\ind_{KZ}^{G}(\bar{\Theta}_{k,a_{p}})$ whenever $T<t-1$ and the assumptions of Theorem~\ref{diagonal conj a-i<i} hold.
\begin{theorem}\label{eliminating JH dcbd}
	Let  $r \equiv a \mod{(p-1)}$ with $ 1 \leq a \leq p-1 $ and $r \geq i(p+1)+p $ with $ v(a_{p})  \in (i,i+1) $. Let $ s =  a-i+n+(i-n)p$ and $v(r-s) =t $ for some $t\geq 1$. Assume that
	\begin{enumerate}
		\item[$(i)$] $ t \geq T+2$  
		\item[$(ii)$] $ 1 \leq  i < a $,  $0 \leq T \leq n < i-1$ and $n+T < 2i-a-1$.
	\end{enumerate}
	Then the image of $ V_{r}^{(i-T)} $ is the same as the image of $ V_{r}^{(i-T+1)} $ in $ \bar{\Theta}_{k,a_{p}} $.
\end{theorem}
\begin{proof}
	Let $\beta_{n-T}, \ldots, \beta_{i-T-1}$ be the integers as chosen in Lemma~\ref{choice beta dcbd}. Then by Lemma~\ref{lem:choice of beta}  there exists $\alpha_{j} \in \mathbb{Z}_{p}$, for all $j \equiv a-i+T $ mod $(p-1)$, with $ a-i+T< j < r-i+T$ satisfying
	\begin{enumerate}
		\item[(1)] $\alpha_{j} \equiv \sum\limits_{l=n-T}^{i-T-1} \beta_{l} \binom{r-l}{j}$ mod $p^{T+2}$, for all $j$ as above, 
		\item[(2)] $\sum\limits_{\substack{ a-i+T < j < r-i+T \\ j \equiv a-i+T ~\mathrm{mod}~(p-1)}}^{} \alpha_{j} \binom{j}{m} \equiv 0$ mod $p^{i+T+3-m}$ for $m=0, \ldots, \min\{p-1, i+T+3\}$.
	\end{enumerate} 
	From the congruence condition (1) and Lemma~\ref{choice beta dcbd} it follows that 
	\begin{enumerate}
		\item[$(2')$] $\sum\limits_{j}^{} \alpha_{j} \binom{j}{m} \equiv 0$ mod $p^{T+1}$ for $m=p, \ldots, i+T+3$.
	\end{enumerate}
	Let		  
	\begin{align*}
		f_2 &= \left[ g_{2,0}^{0}, \frac{(p-1)}{a_p}   (-\theta)^{i+1} X^{-T-1}Y^{r-(i+1)(p+1)+T+1} \right] \\ 
		& \quad + 
		\sum_{  \lambda \in \mathbb{F}_{p}^{\times}} 
		\Bigg[ g_{2,p[\lambda]}^{0},\sum_{l=n-T}^{i-T-1} \beta_l \frac{[\lambda]^{l-(i-T)}}{a_p} p^{i-T-l-1} 
		(-\theta)^{l+T+2} X^{-T-2} Y^{r-(l+T+2)(p+1)+T+2} \Bigg] \\
		f_1 &=  \left[g_{1,0}^0, \frac{(p-1)p^{i-T-1}}{a_p^2} \sum_{ \substack{a-i+T <  j < r-i+T \\  j \equiv a-i+T ~\mathrm{mod}~ (p-1)} } \alpha_{j} X^{r-j} Y^{j} \right] \\
		f_{0} &= \left[ \mathrm{id}, \frac{p^{2i-2T-a-1}(1-p)}{a_p}   \sum_{l=n-T}^{i-T-1} \beta_{l} \binom{r-l}{a-i+T}  \theta^{a-i+2T+2} X^{r-(a-i+2T+2)(p+1)+T+2} Y^{-T-2}\right].
	\end{align*}
	From Lemma~\ref{theta and T plus} it follows that $T^+ f_2$ vanishes modulo $p$. Using  property $(ii)$ of Lemma~\ref{choice beta dcbd} we see that the coefficient of $X^{i-T}Y^{r-i+T}$ coming from the $\lambda =0$ and the $\lambda \neq 0$ terms in $T^{-}f_2$  cancel. Hence
	\begin{align*}
		T^{-}f_2 \equiv \left[ g_{1,0}^0, \frac{p^{i-T-1}(p-1)}{a_p} \sum_{l=n-T}^{i-T-1} \beta_l \sum_{ \substack{a-i+T \leq   j < r-i+T \\  j \equiv a-i+T ~\mathrm{mod}~ (p-1)} } \binom{r-l}{j} X^{r-j} Y^{j} \right]~\mathrm{mod}~p.
	\end{align*}
	On the other hand $T^+ f_0$ kills the $j=a-i+T$ term in $T^- f_2$ above. Thus 
	\begin{align*}
		T^{-}&f_2 -a_p f_1 +T^{+}f_0 \\ 
		& \equiv \left[ g_{1,0}^0, \frac{p^{i-T-1}(p-1)}{a_p} \sum_{l=n-T}^{i-T-1} \beta_l \sum_{ \substack{a-i+T <   j < r-i+T \\  j \equiv a-i+T ~\mathrm{mod}~ (p-1)} } \left( \binom{r-l}{j}-\alpha_{j} \right) X^{r-j} Y^{j} \right]~\mathrm{mod}~p.
	\end{align*}
	which vanishes by (1) above. It can  also be checked  that $ T^{-}f_1, -a_p f_0 $ and $T^{-} f_0$ all vanish modulo $p$ using  the hypotheses $a <  2i-T-n-1$ and  $T \leq n$.  Also $T^{+}f_1$ vanishes by $(2)$ and $(2')$. Hence $(T-a_p)(f_2 + f_1 + f_0) \equiv - a_p f_2$ mod $p$. Using $p\mid \beta_{i-T-1}$ it can be shown that $[\lambda] \neq 0$ vanish modulo $p$, so we have 
	\begin{align*}
		-a_p f_2 \equiv \left[ g_{2,0}^{0}, (-\theta)^{i-T} \left( \sum_{k=0}^{T+1} (-1)^k \binom{T+1}{k} X^{k(p-1)}Y^{r-(i-T)(p+1)-k(p-1)}\right) \right] ~\mathrm{mod}~p.
	\end{align*}
	Since $\sum_{k=0}^{T+1} (-1)^k \binom{T+1}{k} =0$, it follows from \cite[Lemma 2.12]{GR19}  that 
	\begin{multline*}
		\sum_{k=0}^{T+1} (-1)^k \binom{T+1}{k} X^{k(p-1)}Y^{r-(i-T)(p+1)-k(p-1)} \\
		\equiv - X^{p-1}Y^{r-(i-T)(p+1)-(p-1)} 
		+Y^{r-(i-T)(p+1)} ~\mathrm{mod}~V_{r-(i-T)(p+1)}^{(1)}.
	\end{multline*} 
	By Lemma~\ref{Glover-Brueil map image} $(ii)$, the  first monomial generates the cosocle of $V_{r-(i-T)(p+1)}/V_{r-(i-T)(p+1)}^{(1)}$ and the second monomial dies in it. Hence $-a_p f_2 $ generates $V_r^{(i-T)}/V_r^{(i-T+1)}$. This finishes the proof.
\end{proof}

	\subsection{The case \texorpdfstring{$2i-2n \leq a < 2i$}{}}\label{sec: bad medium}

In this section, we prove a more complicated conjecture than 
the diagonal conjecture of the previous section. We do use some results from the previous section. For instance to eliminate JH factors above the diagonal, we   use Theorem~\ref{eliminating JH dcad}. The final results and pictures appear at the end of the section.  
\subsubsection{Below the diagonal } 
We begin by eliminating JH factors below the diagonal.

\begin{lemma}\label{choice beta hcbd}
	Let  $r \equiv a \mod{(p-1)}$ with $ 1 \leq a \leq p-1 $ and $r \geq i(p+1)+p $ with $ v(a_{p})  \in (i,i+1) $. Let $ s =  a-i+n+(i-n)p$ and $v(r-s) =t $ for some $t\geq 1$ and $0 \leq T < n < i-1 <a-1$. If   $a<2i-n-T$ and $t \geq T+2$, then there exist $\beta_0, \ldots, \beta_{i-n-1} \in \mathbb{Z}_p$ such that 
	\begin{enumerate}
		\item[$(i)$] $\sum\limits_{l=0}^{i-n-1} \beta_l  \sum\limits_{ \substack{a-i+T <  j < r-i+T \\  j \equiv a-i+T ~\mathrm{mod}~ (p-1)} }   \binom{r-l}{j} \binom{j}{m}   \equiv 0 \mod p^{T+2-v(m!)} ~ \mathrm{for}~ m=0,\ldots, i+T+2$ 
		
		\item[$(ii)$] $\sum\limits_{l=0}^{i-n-1} \beta_l  \binom{r-l}{r-i+T} \equiv  -p \mod p^{T+2}$.
	\end{enumerate}
\end{lemma}
\begin{proof}
	Let $s=a-i+n+(i-n)p$.  As $r \equiv s$ mod $p^{T+2}(p-1)$, by Corollary~\ref{cor: binomial sums under congruences 2} we have 
	\begin{align*}
		\sum\limits_{ \substack{a-i+T <   j < r-i+T \\  j \equiv a-i+T ~\mathrm{mod}~ (p-1)} }\binom{r-l}{j} \binom{j}{m} \equiv \sum\limits_{ \substack{a-i+T <   j < s-i+T \\  j \equiv a-i+T ~\mathrm{mod}~ (p-1)} }\binom{s-l}{j} \binom{j}{m}~\mathrm{mod}~p^{T+2-v(m!)}
	\end{align*}
	for $l=0,\ldots, i-n-1$. Also by Lemma~\ref{binomial coefficient under congruences} we have $\binom{r-l}{r-i+T} = \binom{r-l}{i-T-l} \equiv \binom{s-l}{i-T-l}  = \binom{s-l}{s-i+T} $ mod $p^{T+2}$. Thus it is enough to show that lemma hold for the special value  $s$.  We first show the existence of $\beta_{0}, \ldots, \beta_{i-n-1} \in \mathbb{Z}_{p}$ satisfying
	\begin{align}
		\sum\limits_{l=0}^{i-n-1} \beta_{l} \sum\limits_{ \substack{a-i+T <   j < s-i+T \\  j \equiv a-i+T ~\mathrm{mod}~ (p-1)} }\binom{s-l}{j} \binom{j}{m} & = 0~ \quad  \mathrm{for}~ m=0, \ldots, i-n-2 \label{red eq 1 hcbd}\\
		\sum\limits_{l=0}^{i-n-1} \beta_{l} \binom{s-l}{s-i+T} = -p.\label{red eq 2 hcbd}
	\end{align}
	Then the above system of equations  can be written as 
	\begin{align*}
		A \begin{bmatrix}
			\beta_{0} \\ \vdots \\ \beta_{i-n-1}
		\end{bmatrix} = \begin{bmatrix}
			0 \\ \vdots \\ 0 \\ -p
		\end{bmatrix}.
	\end{align*}
	where 
	\begin{align*}
		A = \left[ \begin{array}{c}
			\left(\sum\limits_{ \substack{a-i+T <   j < s-i+T \\  j \equiv a-i+T ~\mathrm{mod}~ (p-1)} }\binom{s-l}{j} \binom{j}{m}\right)_{\substack{m=0,\ldots, i-n-2\\ l=0, \ldots, i-n-1}} \\ \hline
			\left( \binom{s-l}{s-i+T} \right)_{l=0, \ldots, i-n-1}
		\end{array}\right]. 
	\end{align*}
	By Cramer's rule, we have $\beta_l = (-1)^{i-n+l+1} (-p) \det(A_{i-n,l})/\det(A) $, where $A_{i-n,l}$ corresponds to the minor of the $(i-n-1,l)$-entry namely  $\binom{s-l}{s-i+T}$. To show $\beta_l \in \mathbb{Z}_p$, it is enough to show that $v(\det(A)) \leq 1 + v( \det(A_{i-n,l}))$. 
	
	We now compute $v(\det(A))$. Noting that every $a-i+T < j \leq s-i+T$ with $j \equiv a-i+T$ mod $(p-1)$ can be expressed as $a-i+T+k(p-1)$ for some $1 \leq k \leq i-n$, we have $A = BC$ with 
	\begin{small}
		\begin{align*}	
			B &= \left[\begin{array}{c|c}
				\left(\binom{a-i+T+k(p-1)}{m}\right)_{\substack{m=0, \ldots, i-n-2 \\ k=1, \ldots, i-n-1}} & \mathbf{0}_{(i-n-1)\times 1} \\ \hline
				\mathbf{0}_{1 \times (i-n-1)} & 1 
			\end{array}\right] \\
			C &= \left[ \binom{s-l}{a-i+T+k(p-1)}_{\substack{k=1, \ldots, i-n \\ l=0, \ldots, i-n-1}}\right]. 
		\end{align*}
	\end{small}%
	Applying Corollary~\ref{cor: GV det} with $d$ there equal to $p-1$, $k=i-n-1$, $m=a-i+T+p-1$ and $n=0$ we see that
	\begin{align}\label{det B hcbd}
		\det(B) = \det\left(\binom{a-i+T+k(p-1)}{m}\right)_{\substack{m=0, \ldots, i-n-2 \\ k=1, \ldots, i-n-1}}    = (p-1)^{(i-n-1)(i-n-2)/2}
	\end{align}
	Hence $\det(B)$ is a unit in $\Z_p$. Applying \cite[Lemma 8]{Viennot} $a-i+n+1$ times we see that
	
	\begin{align}\label{det C hcbd}
		\begin{split}
			\det(C) 
			& =  \prod_{u=0}^{a-i+n} \left( \prod_{l=0}^{i-n-1}(s-l-u) \times  \prod_{k=1}^{i-n} \frac{1}{(a-i+T+k(p-1)-u)} \right)  \\
			& \qquad \qquad \times \det \left( \binom{s-l-(a-i+n+1)}{a-i+T+k(p-1)-(a-i+n+1)}\right)_{\substack{k=1, \ldots, i-n \\ l=0, \ldots, i-n-1}}
		\end{split}
	\end{align} 
	
	By Lucas' theorem and noting $s=(a-i+n)+(i-n)p$ we have
	\begin{align*}
		\binom{s-l-(a-i+n+1)}{a-i+T+k(p-1)-(a-i+n+1)} &\equiv \binom{(i-n-1)p+p-1-l}{(k-1)p+p+T-n-1-k} \\ &\equiv \binom{i-n-1}{k-1} \binom{p-1-l}{p-1+T-n-k} ~\mathrm{mod}~p
	\end{align*}
	where we used $l \leq i-n-1 < i \leq p-2$ and $1 \leq T+1 \leq p+T-i-1 \leq p+T-n-1-k< p-2-k \leq p-1$. Thus by \eqref{det C hcbd} we obtain 
	\begin{align}\label{det C hcbd 1}
		\begin{split}
			\det(C) 
			& \equiv  \prod_{u=0}^{a-i+n} \left( \prod_{l=0}^{i-n-1}(s-l-u) \times  \prod_{k=1}^{i-n} \frac{1}{(a-i+T+k(p-1)-u)} \right)  \\
			& \qquad \qquad \times  \prod_{k=1}^{i-n} \binom{i-n-1}{k-1} \times  \det \left( \binom{p-1-l}{p-1+T-n-k}\right)_{\substack{k=1, \ldots, i-n \\ l=0, \ldots, i-n-1}}  ~\mathrm{mod}~p.
		\end{split}
	\end{align}
	Reversing  the order of rows and column
	\begin{align*}
		\det\left( \binom{p-1-l}{p-1+T-n-k}\right)_{\substack{k=1, \ldots, i-n \\ l=0, \ldots, i-n-1}}  
		&=\det \left( \binom{p-i+n+l}{p-1-i+T+k}\right)_{\substack{k=0, \ldots, i-n-1 \\ l=0, \ldots, i-n-1 }}\\
	\end{align*}
	which is a $p$-adic unit by Corollary~\ref{cor: GV det} $(ii)$. Noting that $\det(B) \in \Z_p^\times$, it follows from \eqref{det C hcbd 1} that
	\begin{equation}\label{valutaion det(A) hcbd}
		\begin{split}
			v(\det(A)) &=v(\det(B))+v(\det(C)) \\ &=  v\left(\prod_{u=0}^{a-i+n} \left( \prod_{l=0}^{i-n-1}(s-l-u) \times  \prod_{k=1}^{i-n} \frac{1}{(a-i+T+k(p-1)-u)} \right)  \right).
		\end{split}
	\end{equation}
	We now compute the determinant of the $(i-n-1,l)$-minor of $A$ which is given by 
	\begin{align*}
		\left(\sum\limits_{ \substack{a-i+T <   j < s-i+T \\  j \equiv a-i+T ~\mathrm{mod}~ (p-1)} }\binom{s-l}{j} \binom{j}{m}\right)_{\substack{m=0,\ldots, i-n-2\\ l'=0, \ldots, i-n-1, l' \neq l }} 
	\end{align*}
	It can be checked $A_{i-n-1,l} = B' C_{i-n-1,l}$ with 
	\begin{small}
		\begin{align*}	
			B' &= 
			\left(\binom{a-i+T+k(p-1)}{m}\right)_{\substack{m=0, \ldots, i-n-2 \\ k=1, \ldots, i-n-1}}  \\
			C_{i-n-1,l} &= \left[ \binom{s-l'}{a-i+T+k(p-1)}_{\substack{k=1, \ldots, i-n-1 \\ l'=0, \ldots, i-n-1, l' \neq l}}\right]. 
		\end{align*}
	\end{small}%
	By \eqref{det B hcbd}, we have $\det(B') \in \mathbb{Z}_p^\times$. We next compute $v(\det(C_{i-n-1,l}))$. This is very similar to $v(\det(C))$. Applying \cite[Lemma 8]{Viennot} $a-i+n+1$ times we see that
	
	\begin{align}\label{det minor C hcbd}
		\begin{split}
			\det(C_{i-n-1,l}) 
			& =  \prod_{u=0}^{a-i+n} \left( \prod_{\substack{l'=0\\ l' \neq l}}^{i-n-1}(s-l'-u) \times  \prod_{k=1}^{i-n} \frac{1}{(a-i+T+k(p-1)-u)} \right)  \\
			& \qquad \qquad \times \det \left( \binom{s-l'-(a-i+n+1)}{a-i+T+k(p-1)-(a-i+n+1)}\right)_{\substack{k=1, \ldots, i-n \\ l'=0, \ldots, i-n-1, l' \neq l }}
		\end{split}
	\end{align} 
	
	By Lucas' theorem and noting $s=(a-i+n)+(i-n)p$ we have
	\begin{align*}
		\binom{s-l'-(a-i+n+1)}{a-i+T+k(p-1)-(a-i+n+1)} &\equiv \binom{(i-n-1)p+p-1-l'}{(k-1)p+p+T-n-1-k} \\ &\equiv \binom{i-n-1}{k-1} \binom{p-1-l'}{p+T-n-1-k} ~\mathrm{mod}~p
	\end{align*}
	where we used $l' \leq i-n-1 < i \leq p-2$ and $1 \leq T+1 \leq p+T-i-1 \leq p+T-n-1-k < p-1-k \leq p-1$. Thus by \eqref{det minor C hcbd} we obtain 
	\begin{align}\label{det minor C hcbd 1}
		\begin{split}
			\det(C_{i-n-1,l}) 
			& \equiv  \prod_{u=0}^{a-i+n} \left( \prod_{\substack{l'=0\\l' \neq l}}^{i-n-1}(s-l'-u) \times  \prod_{k=1}^{i-n} \frac{1}{(a-i+T+k(p-1)-u)} \right)  \\
			& \qquad \qquad \times  \prod_{k=1}^{i-n} \binom{i-n-1}{k-1} \times  \det\left( \binom{p-1-l'}{p-1+T-n-k}\right)_{\substack{k=1, \ldots, i-n \\ l'=0, \ldots, i-n-1, l' \neq l}}  ~\mathrm{mod}~p.
		\end{split}
	\end{align}
	Reversing  the order of rows and columns and then applying Lemma~\ref{GV det} we get 
	\begin{align*}
		\det \left( \binom{p-1-l'}{p-1+T-n-k}\right)_{\substack{k=1, \ldots, i-n-1 \\ l'=0, \ldots, i-n-1, l' \neq l}} 
		&=\det \left( \binom{p-i+n+l'}{p-1-i+T+k}\right)_{\substack{k=0, \ldots, i-n-2 \\ l'=0, \ldots, i-n-1, l' \neq l  }}\\
		&= \frac{(p-i+n)_{p-1-i+T} \cdots \widehat{{(p-1-i+n+l)_{p-1-i+T}}} \cdots (p-1)_{p-1-i+T}}{(p-1-i+T)!\cdots(p-1+T-n-2)!} \\ & \quad \times \prod_{\substack{0 \leq l''<l'\leq i-n-1 \\ l',l'' \neq l}} (l'-l'),
	\end{align*}
	where we omit $(p-1-i+n+l)_{p-1-i+T}$ in the product.
	Since $0 \leq T \leq n-1 \leq i-3\leq a-4\leq p-5$, we obtain the above quantity is a $p$-adic unit.
	Noting that $\det(B') \in \Z_p^\times$, it follows from \eqref{det C hcbd 1} that
	\begin{equation}\label{valutaion minor det(A) hcbd}
		\begin{split}
			v(\det(A_{i-n-1,l})) &=v(\det(B'))+v(\det(C_{i-n-1,l})) \\ &=  v\left(\prod_{u=0}^{a-i+n} \left( \prod_{\substack{l'=0\\l'\neq l}}^{i-n-1}(s-l'-u) \times  \prod_{k=1}^{i-n} \frac{1}{(a-i+T+k(p-1)-u)} \right)  \right).
		\end{split}
	\end{equation}
	Hence by \eqref{valutaion det(A) hcbd} and \eqref{valutaion minor det(A) hcbd} we get
	\begin{align*}
		v \left(\frac{\det(A)}{\det(A_{i-n-1,l})}\right) = v\left(\prod_{u=0}^{a-i+n} (s-l-u)\right). 
	\end{align*}
	Noting that $s = a-i+n+(i-n)p$ and $0 \leq u+l \leq a-i+n +l < a-i+n+p$ for $u=0,\ldots,a-i+n$, we get  $p\mid (s-l-u)$ if and only if $l+u =a-i+n$. Thus
	\begin{align*}
		v\left(\prod_{u=0}^{a-i+n} (s-l-u)\right) = \begin{cases}
			1 \quad &\text{ if } 0 \leq l \leq a-i+n,\\
			0 \quad &\text{ otherwise}. 
		\end{cases}
	\end{align*}
	Thus $v(\det(A)) \leq v (\det(A_{i-n-1,l})) +1$. This shows that there exist $\beta_0, \ldots, \beta_{i-n-1} \in \mathbb{Z}_p$ satisfying \eqref{red eq 1 hcbd} and \eqref{red eq 2 hcbd}. By \eqref{det B hcbd}, we see that the rows of $B'$ span $\mathbb{Z}_p^{(i-n-1)}$. Thus
	for every $m=0, \ldots, i+T+2$, there exists $d_{m,m'} \in \mathbb{Z}_p$ satisfying \begin{align*}
		\binom{a-i+T+k(p-1)}{m} = \sum_{m'=0}^{i-n-2} d_{m,m'}  \binom{a-i+T+k(p-1)}{m'} \quad \text{for all } k=1, \ldots, i-n-1.
	\end{align*}
	Thus 
	\begin{align*}
		\sum_{ \substack{a-i+T <  j < s-i+T \\  j \equiv a-i+T ~\mathrm{mod}~ (p-1)} }   \binom{s-l}{j} \binom{j}{m} = \sum_{m'=0}^{i-n-2} d_{m,m'} \sum_{ \substack{a-i+T <  j < s-i+T  \\  j \equiv a-i+T ~\mathrm{mod}~ (p-1)} }   \binom{s-l}{j} \binom{j}{m'}.
	\end{align*}
	Hence the congruence in  part $(i)$ of the lemma  holds for all $m=0,\ldots, i+T+2$ from  \eqref{red eq 1 hcbd} when $r=s$. This proves the lemma for arbitrary $r$ as observed at the beginning of the proof.
\end{proof}

\begin{theorem}\label{eliminating JH hcbd}
	Let  $r \equiv a \mod{(p-1)}$ with $ 1 \leq a \leq p-1 $ and $r \geq i(p+1)+p $ with $ v(a_{p})  \in (i,i+1) $. Let $ s =  a-i+n+(i-n)p$ and $v(r-s) =t $ for some $t\geq 1$ and $ 0 \leq T < n < i-1 < a-1$. If   $a<2i-n-T$ and $t \geq T+2$, then the image of $ \mathrm{ind}_{KZ}^{G}(V_{r}^{(i-T)}) $ is the same as the image of $ \mathrm{ind}_{KZ}^{G}(V_{r}^{(i-T+1)}) $ in $ \bar{\Theta}_{k,a_{p}} $.
\end{theorem}
\begin{proof}
	Let $\beta_{0}, \ldots, \beta_{i-n-1} \in \mathbb{Z}_p$ be   as chosen in Lemma~\ref{choice beta hcbd}. Then by Lemma~\ref{lem:choice of beta}  there exist $\alpha_{j} \in \mathbb{Z}_{p}$ for all $j \equiv a-i+T $ mod $(p-1)$, with $ a-i+T< j < r-i+T$ satisfying
	\begin{enumerate}
		\item[(1)] $\alpha_{j} \equiv \sum\limits_{l=0}^{i-n-1} \beta_{l} \binom{r-l}{j}$ mod $p^{T+2}$, for all $j$ as above, 
		\item[(2)] $\sum\limits_{\substack{ a-i+T < j < r-i+T \\ j \equiv a-i+T ~\mathrm{mod}~(p-1)}}^{} \alpha_{j} \binom{j}{m} \equiv 0$ mod $p^{i+T+3-m}$ for $m=0, \ldots, \min\{p-1, i+T+2\}$.
	\end{enumerate} 
	From the congruence condition (1) and Lemma~\ref{choice beta hcbd}  it follows that 
	\begin{enumerate}
		\item[$(2')$] $\sum\limits_{j}^{} \alpha_{j} \binom{j}{m} \equiv 0$ mod $p^{T+2 -v(m!)}$ for $m=p, \ldots, i+T+2$.
	\end{enumerate}
	Since $i+T+2 \leq 2i-1 \leq 2p-5$, we obtain $T+2 - v(m!) = T+1$ for $m=p, \ldots,i+T+2$.  Thus we get $m+T+2 - v(m!)  \geq m+T+1  \geq p+T+1 \geq i+T+3$ for $m=p, \ldots, i+T+2$. 
	
	Let		  
	\begin{align*}
		f_2 &= \left[ g_{2,0}^{0}, \frac{(p-1)}{a_p}   (-\theta)^{i+1} X^{-T-1}Y^{r-(i+1)(p+1)+T+1} \right] \\ 
		& \quad + 
		\sum_{  \lambda \in \mathbb{F}_{p}^{\times}} 
		\Bigg[ g_{2,p[\lambda]}^{0},\sum_{l=0}^{i-n-1} \beta_l \frac{[\lambda]^{l-(i-T)}}{a_p} p^{i-T-1-l} 
		(-\theta)^{l+T+2} X^{-T-2} Y^{r-(l+T+1)(p+1)+T+2} \Bigg] \\
		f_1 &=  \left[g_{1,0}^0, \frac{(p-1)p^{i-T-1}}{a_p^2} \sum_{ \substack{a-i+T <  j < r-i+T \\  j \equiv a-i+T ~\mathrm{mod}~ (p-1)} } \alpha_{j} X^{r-j} Y^{j} \right] \\
		f_{0} &= \left[ \mathrm{id}, \frac{p^{2i-2T-a-1}(1-p)}{a_p}   \sum_{l=0}^{i-n-1} \beta_{l} \binom{r-l}{a-i+T}  \theta^{a-i+2T+1} X^{r-(a-i+2T+1)(p+1)+T+1} Y^{-T-1}\right].
	\end{align*}
	From Lemma~\ref{theta and T plus} it follows that $T^+ f_2$ vanishes modulo $p$. Using  property $(ii)$ of Lemma~\ref{choice beta hcbd} we see that the coefficient of $X^{i-T}Y^{r-i+T}$ coming from the \textquote{$\lambda =0$} and the $\lambda \neq 0$ terms in $T^{-}f_2$  cancel. Hence
	\begin{align*}
		T^{-}f_2 \equiv \left[ g_{1,0}^0, \frac{p^{i-T-1}(p-1)}{a_p} \sum_{l=0}^{i-n-1} \beta_l \sum_{ \substack{a-i+T \leq   j < r-i+T \\  j \equiv a-i+T ~\mathrm{mod}~ (p-1)} } \binom{r-l}{j} X^{r-j} Y^{j} \right]~\mathrm{mod}~p.
	\end{align*}
	On the other hand, it can be checked that $T^+ f_0$ kills the $j=a-i+T$ term in $T^- f_2$ above. Thus we have
	\begin{align*}
		T^{-}&f_2 -a_p f_1 +T^{+}f_0 \\ 
		& \equiv \left[ g_{1,0}^0, \frac{p^{i-T-1}(p-1)}{a_p} \sum_{l=0}^{i-n-1} \beta_l \sum_{ \substack{a-i+T <   j < r-i+T \\  j \equiv a-i+T ~\mathrm{mod}~ (p-1)} } \left( \binom{r-l}{j}-\alpha_{j} \right) X^{r-j} Y^{j} \right] \equiv 0 ~\mathrm{mod}~p,
	\end{align*}
	by $(1)$. It can  also be checked  that $ -a_p f_0 $ and $T^{-} f_0$ vanish modulo $p$ using  the hypotheses $a \leq  2i-T-n-1$ and  $T < n$.  We note that $T^{-}f_1$ vanishes if $(p-1)+2(i-T)-1 \geq 2(i+1) \Leftrightarrow 2T \leq p-4$. Since $2T \leq T+n-1  \leq 2i-a-2 < a -2 \leq p-3 $, we get $T^{-}f_1$ vanishes. Also $T^{+}f_1$ vanishes by $(2)$ and (the discussion below)  $(2')$. Hence $(T-a_p)(f_2 + f_1 + f_0) \equiv - a_p f_2$ mod $p$. Using $T\leq n-1$ it can be shown that $\lambda \neq 0$ terms vanish modulo $p$ in $-a_p f_2$, so we have 
	\begin{align*}
		-a_p f_2 \equiv \left[ g_{2,0}^{0}, (-\theta)^{i-T} \left( \sum_{k=0}^{T+1} (-1)^k \binom{T+1}{k} X^{k(p-1)}Y^{r-(i-T)(p+1)-k(p-1)}\right) \right] ~\mathrm{mod}~p.
	\end{align*}
	Since $\sum_{k=0}^{T+1} (-1)^k \binom{T+1}{k} =0$, it follows from \cite[Lemma 2.12]{GR19}  that 
	\begin{multline*}
		\sum_{k=0}^{T+1} (-1)^k \binom{T+1}{k} X^{k(p-1)}Y^{r-(i-T)(p+1)-k(p-1)} \\
		\equiv (- X^{p-1}Y^{r-(i-T)(p+1)-(p-1)} 
		+Y^{r-(i-T)(p+1)} ~\mathrm{mod}~V_{r-(i-T)(p+1)}^{(1)}.
	\end{multline*} 
	By Lemma~\ref{generating polynomial quotient}, $\theta^{i-T}$ times the above polynomial generates $V_{r}^{(i-T)}/V_{r}^{(i-T+1)}$. Hence $-a_p f_2 $ generates $V_r^{(i-T)}/V_r^{(i-T+1)}$. This finishes the proof.
\end{proof}

\subsubsection{Below the superdiagonal}
\begin{lemma}\label{choice beta sdcbd}
	Let  $r \equiv a \mod{(p-1)}$ with $ 1 \leq a \leq p-1 $ and $r \geq i(p+1)+p $ with $ v(a_{p})  \in (i,i+1) $. Let $ s =  a-i+n+(i-n)p$ and $v(r-s) =t $ for some $t\geq 1$ and $ 0 \leq T < n < i <  p-1$. If  $2i-n -T \leq a $,  $ i \leq a$ and $t \geq T+1$, then there exist $\beta_{n-T-1}, \ldots, \beta_{i-T-1}$ satisfying the following congruences
	\begin{enumerate}
		\item[$(i)$] $\sum\limits_{l=n-T-1}^{i-T-1} \beta_{l} \sum\limits_{ \substack{a-i+T <   j < r-i+T \\  j \equiv a-i+T ~\mathrm{mod}~ (p-1)} }\binom{r-l}{j} \binom{j}{m} ~\equiv 0~ \mathrm{mod}~ p^{T+1-v(m!)} $ for $m=0, \ldots, i+T+1$.
		\item[$(ii)$] $\sum\limits_{l=n-T-1}^{i-T-1} \beta_{l} \binom{r-l}{r-i+T} \equiv 1 ~\mathrm{mod}~p^{T+1}$.
		\item[$(iii)$] $\sum\limits_{l=n-T-1}^{i-T-1} \beta_{l} \binom{r-l}{a-i+T} \equiv 0 ~\mathrm{mod}~p^{T+1}$.
		\item[$(iv)$] $\sum\limits_{l=n-T-1}^{i-T-1} \beta_{l} \binom{r-l}{r-i+T-(p-1)} \equiv 0 ~\mathrm{mod}~p^{T}$. 
	\end{enumerate}
\end{lemma}
\begin{proof}
	As observed in the proof of Lemma~\ref{choice beta dcbd}, using $r\equiv s \mod p^{T+1}$, Lemma~\ref{binomial sums under congruences} and Corollary~\ref{cor: binomial sums under congruences 1}, it is enough to show that the lemma holds when $r=s$. We first show the existence of $\beta_{n-T-1}, \ldots, \beta_{i-T-1} \in \mathbb{Z}_{p}$  satisfying
	\begin{align}
		\sum\limits_{l=n-T-1}^{i-T-1} \beta_{l} \sum\limits_{ \substack{a-i+T <   j < s-i+T \\  j \equiv a-i+T ~\mathrm{mod}~ (p-1)} }\binom{s-l}{j} \binom{j}{m} &= 0 ~\text{ for } m=0,\ldots, i-n-2 \label{red eq 1 sdcbd}\\
		\sum\limits_{l=n-T-1}^{i-T-1} \beta_{l} \binom{s-l}{s-i+T} &= 1 \label{red eq 2 sdcbd} \\
		\sum\limits_{l=n-T-1}^{i-T-1} \beta_{l} \binom{s-l}{a-i+T} &= 0. \label{red eq 3 sdcbd}
	\end{align}
	Note that if $n=i-1$, then \eqref{red eq 1 sdcbd} is vacuously true. Then the above system of equations  can be written as 
	\begin{align*}
		A \begin{bmatrix}
			\beta_{n-T-1} \\ \vdots \\ \beta_{i-T-1}
		\end{bmatrix} = \begin{bmatrix}
			1 \\ 0 \\ \vdots \\ 0 
		\end{bmatrix}.
	\end{align*}
	where 
	\begin{align*}
		A = \left[ \begin{array}{c}
			\left( \binom{s-l}{a-i+T} \right)_{l=n-T-1, \ldots, i-T-1} \\ \hline
			\left(\sum\limits_{ \substack{a-i+T <   j < s-i+T \\  j \equiv a-i+T ~\mathrm{mod}~ (p-1)} }\binom{s-l}{j} \binom{j}{m}\right)_{\substack{m=0,\ldots, i-n-2\\ l=n-T-1, \ldots, i-T-1}} \\ \hline
			\left( \binom{s-l}{s-i+T} \right)_{l=n-T-1, \ldots, i-T-1}
		\end{array}\right]. 
	\end{align*}
	Noting that every $a-i+T \leq j \leq s-i+T$ with $j \equiv a-i+T$ mod $(p-1)$ can be expressed as $a-i+T+k(p-1)$ for some $0 \leq k \leq i-n$, we have $A = BC$ with 
	\begin{small}
		\begin{align*}	
			B &= \left[\begin{array}{c|c|c}
				1 & \mathbf{0}_{1 \times (i-n-1)} & 0 \\ \hline \mathbf{0}_{(i-n-1) \times 1} &\left(\binom{a-i+T+k(p-1)}{m}\right)_{\substack{m=0, \ldots, i-n-2 \\ k=1, \ldots, i-n-1}} & \mathbf{0}_{(i-n-1)\times 1} \\ \hline
				0 & \mathbf{0}_{1 \times (i-n-1)} & 1 
			\end{array}\right] \\
			C &= \left[ \binom{s-l}{a-i+T+k(p-1)}_{\substack{k=0, \ldots, i-n \\ l=n-T-1, \ldots, i-T-1}}\right]. 
		\end{align*}
	\end{small}%
	Applying Corollary~\ref{cor: GV det}  with $d$ there equal to $(p-1)$, $k=i-n-1$, $m=a-i+T+p-1$ and $n=0$ we see that $\det(B) = (p-1)^{(i-n-1)(i-n-2)/2}$. Hence $\det(B)$ is a unit in $\Z_p$. By Lucas' theorem, for $n-T-1 \leq l \leq i-T-1$ and $0 \leq k \leq i-n$ we have $\binom{s-l}{a-i+T+k(p-1)} \equiv \binom{i-n}{k} \binom{a-i+n-l}{a-i+T-k} \mod p$, where we have used $2i-n-T \leq a$. Hence
	\begin{align*}
		\det(C) \equiv \left( \prod_{k=0}^{i-n} \binom{i-n}{k} \right) \det\left( \binom{a-i+n-l}{a-i+T-k}\right)_{\substack{k=0, \ldots, i-n \\ l=n-T-1, \ldots, i-T-1}}. 
	\end{align*}
	Reversing the order of rows and columns we get 
	\[
	\det\left( \binom{a-i+n-l}{a-i+T-k}\right)_{\substack{k=0, \ldots, i-n \\ l=n-T-1, \ldots, i-T-1}} = \det\left( \binom{a-2i+n+T+1+l}{a-2i+n+T+k}\right)_{\substack{k=0, \ldots, i-n \\ l=0, \ldots, i-n}}.
	\]
	Applying Corollary~\ref{cor: GV det} $(ii)$, we see that $\det(C) \in \Z_p^\times$. This shows that \eqref{red eq 1 sdcbd}, \eqref{red eq 2 sdcbd} and \eqref{red eq 3 sdcbd} has a solution in $\mathbb{Z}_p$.  As  observed towards the end of the proof of Lemma~\ref{choice beta dcbd}, it can be shown that \eqref{red eq 1 sdcbd} implies $(i)$.

    We now prove $(iv)$ holds. By Lemma~\ref{binomial coefficient under congruences} $(i)$ we have 
		\begin{align*}
			\binom{r-l}{r-i+T-(p-1)} = \binom{r-l}{i-T-l+(p-1)} \equiv \binom{s-l}{i-T-l+(p-1)} = \binom{s-l}{s-i+T-(p-1)} \mod p^T.
		\end{align*}
		Thus it is enough to prove that $(iv)$ holds for the special value $s$. By \eqref{red eq 1 sdcbd}, we have 
        \begin{equation*}
            \sum_{k=1}^{i-n-1} \binom{a-i+T+k(p-1)}{m} \sum_{l=n-T-1}^{i-T-1} \beta_{l} \binom{s-l}{a-i+T+k(p-1)} =0 \quad \text{ for } m=0,\ldots,i-n-2.
        \end{equation*}
        Now by Lemma~\ref{dmm' trick}, we obtain 
        \begin{equation*}
            \sum_{l=n-T-1}^{i-T-1} \beta_{l} \binom{s-l}{a-i+T+k(p-1)} \equiv 0 \mod p^t \quad \text{ for } k=1,\ldots,i-n-1.
        \end{equation*}
        Taking $k=i-n-1$, we get $(iv)$ holds for $s$, as $s-i+T-(p-1) = a-i+T+(i-n-1)(p-1)$.  This proves $(iv)$ as noted earlier.
	This finishes the proof of the lemma.
\end{proof}

\begin{theorem}\label{eliminating JH sdcbd}
	Let  $r \equiv a \mod{(p-1)}$ with $ 1 \leq a \leq p-1 $ and $r \geq i(p+1)+p $ with $ v(a_{p})  \in (i,i+1) $. Let $ s =  a-i+n+(i-n)p$ and $v(r-s) =t $ for some $t\geq 1$ and $ 0 \leq T < n < i < a$. If  $2i-n -T \leq a $,   and $t \geq T+1$, then the image of $ \mathrm{ind}_{KZ}^{G}(V_{r}^{(i-T)}) $ is the same as the image of $\mathrm{ind}_{KZ}^{G} (V_{r}^{(i-T+1)}) $ in $\bar{\Theta}_{k,a_{p}} $.
\end{theorem}
\begin{proof}
	Consider the function
	\begin{align*}
		f_{2} = \sum_{\lambda \in \mathbb{F}_{p}^{\times}} 
		\Bigg[ g_{2,p[\lambda]}^{0}, & \frac{1}{a_{p}} \sum_{l=n-T-1}^{i-T-1} 
		\bigg( \frac{p}{[\lambda]}\bigg)^{i-T-l}  
		\beta_{l} (-\theta)^{l+T+1} X^{-T-1}Y^{r-(l+T+1)(p+1)+T+1} \Bigg] \\
		& + \frac{(1-p)}{a_{p}}  \Bigg[ g_{2,0}^{0},  (-\theta)^{i+1} X^{-T-1}Y^{r-(i+1)(p+1)+T+1} \Bigg]
	\end{align*}
	where $\beta_l \in \Z_p$ for $n-T-1 \leq l \leq i-T-1$ are as chosen in Lemma~\ref{choice beta sdcbd}. 
	Applying Lemma~\ref{lem:choice of beta} with $m$ there equal to $\min\{i+T+1, p-1\}$, we obtain $\alpha_j$ satisfying
	\begin{enumerate}
		\item[$(1)$] $\alpha_{j} \equiv \sum\limits_{l=n-T-1}^{i-T-1} \beta_{l} \binom{r-l}{j}$ mod $p^{T+1}$, for all     $a-i+T <  j < r-i+T$ and $j \equiv a-i+T \mod (p-1)$ and $\alpha_j =0$ otherwise 
		\item[$(2)$] $\sum\limits_{\substack{ a-i+T < j < r-i+T \\ j \equiv a-i+T ~\mathrm{mod}~(p-1)}}^{} \alpha_{j} \binom{j}{m} \equiv 0$ mod $p^{i+T+2-m}$ for $m=0, \ldots, \min\{p-1, i+T+1\}$.
	\end{enumerate} 
	From the congruence condition $(1)$ and Lemma~\ref{choice beta sdcbd}  it follows that 
	\begin{enumerate}
		\item[$(2')$] $\sum\limits_{j}^{} \alpha_{j} \binom{j}{m} \equiv 0$ mod $p^{T+1- v(m!)}$ for $m=p, \ldots, i+T+1$.
	\end{enumerate}	
	Note that $i+T+2-m \leq T+1-v(m!) =T$ for $m=p, \ldots, i+T+1 $ since $i+2 \leq p$.  Let 
	\begin{align*}
		f_1 =  \left[g_{1,0}^0, \frac{(p-1)p^{i-T}}{a_p^2} \sum_{ \substack{a-i+T <   j < r-i+T \\  j \equiv a-i+T ~\mathrm{mod}~ (p-1)} } \alpha_{j} X^{r-j} Y^{j} \right]
	\end{align*}
	From Lemma~\ref{theta and T plus} it follows that $T^+ f_2$ vanishes modulo $p$. Using  property $(ii)$ of Lemma~\ref{choice beta sdcbd} we see that the coefficient of $X^{i-T}Y^{r-i+T}$ coming from the \textquote{$\lambda =0$} and the $\lambda \neq 0$ terms in $T^{-}f_2$  cancel. Using  property $(iii)$ of Lemma~\ref{choice beta sdcbd} we see that the coefficient of $X^{r-(a-i+T)}Y^{a-i+T}$  in $T^{-}f_2$  vanishes modulo $p$. Hence
	\begin{align*}
		T^{-}f_2 -a_p f_1 \equiv \left[ g_{1,0}^0, \frac{p^{i-T}(p-1)}{a_p}  \sum_{ \substack{a-i+T <   j < r-i+T \\  j \equiv a-i+T ~\mathrm{mod}~ (p-1)} } \left(\sum_{l=n-T-1}^{i-T-1} \beta_l \binom{r-l}{j} - \alpha_j \right) X^{r-j} Y^{j} \right] \equiv  0 ~\mathrm{mod}~p.
	\end{align*}
	From $(2)$ and $(2')$ it follows that $T^+f_1 $ vanishes modulo $p$. Since $i \leq p-2$, it follows from Lemma~\ref{choice beta sdcbd} $(iv)$ that $T^{-}f_1$ also vanishes modulo $p$. Thus 
	$$(T-a_p)(f_2 + f_1) \equiv -a_p f_2 \equiv -\Bigg[ g_{2,0}^{0},  (-\theta)^{i+1} X^{-T-1}Y^{r-(i+1)(p+1)+T+1} \Bigg] \mod p . $$
	Note 
	\begin{align*}
		(-\theta)^{i+1} X^{-T-1}Y^{r-(i+1)(p+1)+T+1} &= (-\theta)^{i-T} (\sum\limits_{j=0}^{T+1} (-1)^{j} \binom{T+1}{j} X^{(p-1)j}Y^{r-(i-T)(p+1)-j(p-1)})\\
		& \equiv (-\theta)^{i-T}  (Y^{r-(i-T)(p+1)} - X^{p-1}Y^{r-(i-T)(p+1)-(p-1)}) \mod{V}_r^{(i-T+1)}.
	\end{align*}
	By Lemma~\ref{generating polynomial quotient},  the above polynomial generates ${V}_r^{(i-T)}/{V}_r^{(i-T+1)}$. This completes the proof of the theorem.
\end{proof}

We now apply the above theorem to the case $a= 2i-1$ and $r \equiv i \mod p$ (that is $r \equiv a-i+n \mod p$ with $n=1$). Taking $ T = 0$ in Theorem~\ref{eliminating JH sdcbd} we get $\mathrm{ind}_{KZ}^{G}(V_{r}^{(i)}/V_{r}^{(i+1)})$ vanishes $\bar{\Theta}_{k,a_{p}} $. By Theorem~\ref{Elimination i < a and not in interval} we have the image of $\mathrm{ind}_{KZ}^{G}(V_{r}^{(m)})$ is the same as the image of $\mathrm{ind}_{KZ}^{G}(V_{r}^{(m+1)})$ for $m=0, \ldots, a-i-1$. Applying Lemma~\ref{JH factor Q 2i > a} $(ii)$ part $(a)$  we obtain $\mathrm{ind}_{KZ}^{G}(V_{p-2} \otimes D^{i}) \twoheadrightarrow \bar{\Theta}_{k,a_{p}} $. We now study how the above surjection factors, that is, whether it factors through $T$ or $T^2-cT+1$ for some $c \in \bar{\mathbb{F}}_p$. 

In the next few results, the following number 
$$d = -\frac{1}{p} \binom{r-i+1}{i} - \frac{(-1)^{i+1}}{i}$$
plays an important role. Observe that when $r \equiv i \mod p$ and $i \geq 2$, then $\binom{r-i+1}{i} \equiv 0 \mod p$ by Lucas' theorem. Thus $d $ is a $p$-adic integer if $2 \leq i\leq p-1$.

We first prove two combinatorial lemmas. 
\begin{lemma}\label{beta l for a = 2i-1, n=1}
	Let  $r \equiv 2i-1 \mod{(p-1)}$ and $r \equiv i \mod p$ with $ 2 \leq i \leq (p-1)/2 $ and $r \geq p $. Let $\beta_l = \frac{(-1)^{i-l+1}}{i-l} \binom{i}{l}$ for $l=0, \ldots, i-2$ and $\beta_{i-1} = \frac{1}{p}$. Then we have the following
	\begin{enumerate}
		\item[$(i)$] $p \sum\limits_{l=0}^{i-1} \beta_l  \sum\limits_{\substack{i \leq j < r-i+1 \\ j \equiv i \mod (p-1)  }} \binom{r-l}{j} \binom{j}{m}  \equiv 0 \mod p^2 $ for $0 \leq m \leq i-1 $.
		\item[$(ii)$] $p \sum\limits_{l=0}^{i-1} \beta_l  \sum\limits_{\substack{i \leq j < r-i+1 \\ j \equiv i \mod (p-1)  }} \binom{r-l}{j} \binom{j}{i}  \equiv -pd \mod p^2 $. 
	\end{enumerate}
\end{lemma}
\begin{proof}
	We first compute the partial sum $$ p \sum\limits_{l=0}^{i-2} \beta_l  \sum\limits_{\substack{i \leq j < r-i+1 \\ j \equiv i \mod (p-1)  }} \binom{r-l}{j} \binom{j}{m} \text{ for } m=0, \ldots, i-2.$$ 
	Let $m=0, \ldots, i-1$ and $l=0,\ldots, i-2$. By \cite[Lemma 2.14]{GR19}, for $l=0, \ldots, i-2$ and  we have 
	\begin{align}\label{sum a=2i-1, n=1}
		\begin{split}
			\sum\limits_{\substack{i \leq j < r-i+1 \\ j \equiv i \mod (p-1)  }} \binom{r-l}{j} \binom{j}{m}  &\equiv  \binom{r-l}{m} \binom{2i-1-l-m}{i-m} - \binom{r-l}{r-i+1} \binom{r-i+1}{m} \\
			&  \equiv \binom{i-l}{m} \binom{2i-1-l-m}{i-m} - \binom{i-l}{i-1-l}\binom{1}{m}  \mod p,	          	  		
		\end{split}	       
	\end{align}
	where we have used $r \equiv i \mod p$ and Lucas' theorem in the last step. Note that for $0 \leq m \leq i-1$ and $l=0, \ldots,i-2$, we have  
	\begin{align*}
		\beta_l \binom{i-l}{m} \binom{2i-1-l-m}{i-m} &= (-1)^{i-l+1} \frac{i!(2i-1-l-m)!}{l!(i-l-m)! m! (i-m)!(i-l)!} \\ &= \frac{(-1)^{i-l+1}}{i-m} \binom{i}{m} \binom{i-m}{l} \binom{2i-1-l-m}{i-m-1}.
	\end{align*}
	Multiplying both sides of  \eqref{sum a=2i-1, n=1} by $p\beta_{l}$ and using the above identity we obtain
	\begin{multline*}
		p \beta_{l} \sum\limits_{\substack{i \leq j < r-i+1 \\ j \equiv i \mod (p-1)  }} \binom{r-l}{j} \binom{j}{m}  \\ \equiv  p (-1)^{i-l+1} \left( \frac{1}{i-m} \binom{i}{m} \binom{i-m}{l} \binom{2i-1-l-m}{i-m-1} -  \binom{i}{l}\binom{1}{m} \right) \mod p^2.
	\end{multline*}
	Taking sum over $l =0, \ldots, i-2$ we obtain
	\begin{equation}\label{sum a=2i-1, n=1 red1}
		\begin{multlined}
			\sum_{l=0}^{i-2} p \beta_{l} \sum\limits_{\substack{i \leq j < r-i+1 \\ j \equiv i \mod (p-1)  }} \binom{r-l}{j} \binom{j}{m}  \\
			\equiv  \sum_{l=0}^{i}  p (-1)^{i-l+1} \left( \frac{1}{i-m} \binom{i}{m} \binom{i-m}{l} \binom{2i-1-l-m}{i-m-1} -  \binom{i}{l}\binom{1}{m} \right) \qquad \qquad \quad  \\
			\qquad -  \sum_{l=i-1}^{i}  p (-1)^{i-l+1} \left( \frac{1}{i-m} \binom{i}{m} \binom{i-m}{l} \binom{2i-1-l-m}{i-m-1} -  \binom{i}{l}\binom{1}{m} \right) \\
			= \sum_{l=0}^{i}  p (-1)^{i-l+1}  \frac{1}{i-m} \binom{i}{m} \binom{i-m}{l} \binom{2i-1-l-m}{i-m-1} - \sum_{l=0}^{i}  p (-1)^{i-l+1} \binom{i}{l}\binom{1}{m}  \\
			- p\left( \binom{i}{m} \binom{i-m}{i-1} - \frac{1}{i-m} \binom{i}{m} \binom{i-m}{i} -(i-1) \binom{1}{m}\right) \mod p^2.      	  	   
		\end{multlined}
	\end{equation}
	Note that $\sum\limits_{l=0}^{i}  p (-1)^{i-l+1} \binom{i}{l}\binom{1}{m} =0$. Also $\binom{i}{m} \binom{i-m}{i-1} = \binom{1}{m} i$ and $\binom{i}{m} \binom{i-m}{i} =  \delta_{m,0} $. Thus
	\begin{multline*}
		\sum_{l=0}^{i-2} p \beta_{l} \sum\limits_{\substack{i \leq j < r-i+1 \\ j \equiv i-1 \mod (p-1)  }} \binom{r-l}{j} \binom{j}{m}  \\ \equiv 
		\frac{p}{i-m} \binom{i}{m} \sum_{l=0}^{i}   (-1)^{i-l+1}    \binom{i-m}{l} \binom{2i-1-l-m}{i-m-1} - p \left( \binom{1}{m} - \frac{1}{i} \delta_{m,0} \right) \mod p^2.
	\end{multline*}
	Observe that 
	\begin{align*}
		\sum_{l=0}^{i}   (-1)^{i-l+1}    \binom{i-m}{l} \binom{2i-1-l-m}{i-l} & = \text{coefficient of } x^{i} \text{ in } \sum_{l=0}^{i} (-1)^{i-l+1} \binom{i-m}{l} x^{l} (1+x)^{2i-l-m-1} \\
		& = \text{coefficient of } x^{i} \text{ in } (-1)^{i+1} (1+x)^{i-1} \sum_{l=0}^{i}  \binom{i-m}{l} (-x)^{l} (1+x)^{i-m-l} \\
		& = \text{coefficient of } x^{i} \text{ in } (-1)^{i+1} (1+x)^{i-1} \cdot 1 =0.
	\end{align*}
	Noting that $\binom{2i-1-l-m}{i-m-1} = \binom{2i-1-l-m}{i-l}$, we obtain 
	\begin{align}\label{sum a=2i-1, n=1 red2}
		\sum_{l=0}^{i-2} p \beta_{l} \sum\limits_{\substack{i \leq j < r-i+1 \\ j \equiv i \mod (p-1)  }} \binom{r-l}{j} \binom{j}{m}  =  - p \left( \binom{1}{m} - \frac{1}{i} \delta_{m,0} \right) \mod p^2.
	\end{align}
	By \cite[Lemma 2.5]{BG15}, for $m=0, \ldots, i-1$ we have
	\begin{align*}
		p \beta_{i-1} \sum\limits_{\substack{i \leq j < r-i+1 \\ j \equiv i \mod (p-1)  }} \binom{r-i+1}{j} \binom{j}{m}  &= \binom{r-i+1}{m}  \sum\limits_{\substack{i-m \leq j' < r-i+1-m \\ j' \equiv i-m \mod (p-1)  }} \binom{r-i+1-m}{j'} \\
		&\equiv p \binom{r-i+1}{m} \frac{(i-(r-i+1))}{i-m} \mod p^2 \\
		&\equiv p \binom{1}{m} \frac{i-1}{i-m} = p \left( \delta_{m,1} + \frac{i-1}{i} \delta_{m,0} \right).
	\end{align*}
	Now $(i)$ follows from \eqref{sum a=2i-1, n=1 red2} and the above congruence.
	
	We now prove $(ii)$.  First note that $0 < 2i-1-l-i = i-1-l < p-1$ for $l=0,\ldots, i-2$. Thus by \cite[Lemma 2.14]{GR19} with $r$ there equal to $r-l$, $m$ there equal to $i$ and $b$ there equal to $i$, we have
	\begin{align*} 
		\sum\limits_{\substack{i \leq j < r-i+1 \\ j \equiv i \mod (p-1)  }} \binom{r-l}{j} \binom{j}{i}  &\equiv  \binom{r-l}{i} \left( \binom{i-1-l}{p-1} +1 \right)- \binom{r-l}{r-i+1} \binom{r-i+1}{i} \\
		&  \equiv \binom{i-l}{i} - \binom{i-l}{i-1-l}\binom{1}{i}  \equiv  \delta_{l,0} \mod p, 		   \end{align*}
	where we have used Lucas' theorem and $r \equiv i \mod p$ in the penultimate step and $i \geq 2$ in the last step. Multiplying both sides by $p\beta_l$ and the taking sum  over $l=0, \ldots,i-2$ we get
	\begin{align}\label{sum a=2i-1, n=1 red3}
		\sum_{l=0}^{i-2} p \beta_{l} \sum\limits_{\substack{i \leq j < r-i+1 \\ j \equiv i \mod (p-1)  }} \binom{r-l}{j} \binom{j}{i} \equiv p \beta_0 = p \frac{(-1)^{i+1}}{i}\mod p^2.  
	\end{align}
	By \cite[Lemma 2.5]{BG15},  we have 
	\begin{align*}
		p \beta_{i-1} \sum\limits_{\substack{i \leq j < r-i+1 \\ j \equiv i \mod (p-1)  }} \binom{r-i+1}{j} \binom{j}{i}  &= \binom{r-i+1}{i}  \sum\limits_{\substack{0 \leq j' < r-2i+1 \\ j' \equiv 0 \mod (p-1)  }} \binom{r-2i+1}{j'} \\
		& \equiv  p \binom{r-i+1}{i} \frac{(p-1-(r-2i+1))}{p-1} + \binom{r-i+1}{i}  \\
		& \equiv p \binom{1}{i} \frac{i-2}{p-1} + \binom{r-i+1}{i}   \equiv \binom{r-i+1}{i} \mod p^2 
	\end{align*}
	as $i\geq 2$. Now $(ii)$ follows \eqref{sum a=2i-1, n=1 red3} and the above congruence. 
\end{proof}
\begin{lemma}\label{gamma l for a = 2i-1, n=1}
	Let  $r \equiv 2i-1 \mod{(p-1)}$ and $r \equiv i \mod p$ with $ 2 \leq i \leq (p-1)/2 $ and $r \geq p $. Let $\gamma_l =(-1)^{l} \binom{r}{l}$ for $l=0, \ldots, i-2$. Then we have the following
	$$ \sum\limits_{l=0}^{i-2} \gamma_l  \sum\limits_{\substack{i < j < r-i+1 \\ j \equiv i \mod (p-1)  }} \binom{r-l}{j} \binom{j}{m}  \equiv 0 \mod p  \text{ for } 0 \leq m \leq i+1 .$$                
\end{lemma}
\begin{proof}
	We first prove the congruences in the lemma for $m=0, \ldots,i$. By \cite[Lemma 2.14]{GR19}, for $l=0, \ldots, i-2$ and $m=0, \ldots,i$ and noting $1 \leq 2i-1-l-m  < p-1$ we have 
	\begin{align}\label{sum2 a=2i-1, n=1}
		\begin{split}
			\sum\limits_{\substack{i < j < r-i+1 \\ j \equiv i \mod (p-1)  }} \binom{r-l}{j} \binom{j}{m} &= \sum\limits_{\substack{i \leq  j \leq  r-i+1 \\ j \equiv i \mod (p-1)  }} \binom{r-l}{j} \binom{j}{m} - \binom{r-l}{i} \binom{i}{m} - \binom{r-l}{r-i+1} \binom{r-i+1}{m} \\  
			&\equiv  \binom{r-l}{m} \left( \binom{[2i-1-l-m]}{[i-m]} + \delta_{i,m}\right)- \binom{r-l}{i} \binom{i}{m}- \binom{r-l}{r-i+1} \binom{r-i+1}{m} \\
			&\equiv  \binom{r-l}{m} \binom{2i-1-l-m}{i-m} - \binom{r-l}{i} \binom{i}{m}- \binom{r-l}{r-i+1} \binom{r-i+1}{m} \\
			&  \equiv \binom{i-l}{m} \binom{2i-1-l-m}{i-m} - \binom{i-l}{i}\binom{i}{m}- \binom{i-l}{i-1-l}\binom{1}{m}  \mod p,	          	 	
		\end{split}	       
	\end{align}
	where we used $r \equiv i \mod p$ and Lucas' theorem in the last step. Note that 
	\begin{align*}
		\binom{i}{l} \binom{i-l}{m} \binom{2i-1-l-m}{i-m} =  \frac{i!}{l!(i-l-m)!m!} \binom{2i-1-l-m}{i-m} =  \binom{i}{m} \binom{i-m}{l} \binom{2i-1-l-m}{i-m}.
	\end{align*}
	Also note that $\gamma_l \equiv (-1)^l \binom{i}{l} \mod p $ for $l=0, \ldots, i-2$. Thus multiplying both sides of \eqref{sum2 a=2i-1, n=1} by $\gamma_l$ and then using this congruence and above identity we obtain
	\begin{multline*}
		\gamma_l \sum\limits_{\substack{i < j < r-i+1 \\ j \equiv i \mod (p-1)  }} \binom{r-l}{j} \binom{j}{m} \equiv  (-1)^{l}  \binom{i}{m} \binom{i-m}{l} \binom{2i-1-l-m}{i-m} - (-1)^{l} \binom{i}{l} \binom{i-l}{i} \binom{i}{m} \\
		\qquad  \qquad \qquad - (-1)^{l} \binom{i}{l} \binom{i-l}{i-1-l}\binom{1}{m} \\
		\qquad \qquad \qquad \qquad \qquad\equiv  (-1)^{l}  \Bigg\lbrace \binom{i}{m} \binom{i-m}{l} \binom{2i-1-l-m}{i-m} -  \delta_{l,0} \binom{i}{m}  -  i \binom{i-1}{l} \binom{1}{m} \Bigg\rbrace \mod p.
	\end{multline*}
	Taking sum over $l=0, \ldots, i-2$, we get 
	\begin{equation}\label{sum gamma a=2i-1 and n=1}
		\begin{split}
			\sum_{l=0}^{i-2} \gamma_l & \sum\limits_{\substack{i < j < r-i+1 \\ j \equiv i-1 \mod (p-1)  }} \binom{r-l}{j} \binom{j}{m} \\ 
			&\equiv  -\binom{i}{m} + \sum_{l=0}^{i-2} (-1)^{l}  \Bigg\lbrace \binom{i}{m} \binom{i-m}{l} \binom{2i-1-l-m}{i-m}   -  i \binom{i-1}{l} \binom{1}{m} \Bigg\rbrace  \\
			&\equiv  -\binom{i}{m} + \sum_{l=0}^{i} (-1)^{l}  \Bigg\lbrace \binom{i}{m} \binom{i-m}{l} \binom{2i-1-l-m}{i-m}   -  i \binom{i-1}{l} \binom{1}{m} \Bigg\rbrace  \\
			&\qquad \qquad- \sum_{l=i-1}^{i} (-1)^{l}  \Bigg\lbrace \binom{i}{m} \binom{i-m}{l} \binom{2i-1-l-m}{i-m}   -  i \binom{i-1}{l} \binom{1}{m} \Bigg\rbrace \\
			&\equiv -\binom{i}{m} + \sum_{l=0}^{i} (-1)^{l}   \binom{i}{m} \binom{i-m}{l} \binom{2i-1-l-m}{i-m} \\  & \qquad \qquad - \sum_{l=i-1}^{i} (-1)^{l}  \Bigg\lbrace \binom{i}{m} \binom{i-m}{l} \binom{2i-1-l-m}{i-m}   -  i \binom{i-1}{l} \binom{1}{m} \Bigg\rbrace \mod p,
		\end{split}
	\end{equation}
	where we used $\sum\limits_{l=0}^{i-1} (-1)^l \binom{i-1}{l} =0$ in the last step. Note that 
	\begin{align*}
		\sum_{l=i-1}^{i} (-1)^{l}  \Bigg\lbrace \binom{i}{m} \binom{i-m}{l} \binom{2i-1-l-m}{i-m}   -  i \binom{i-1}{l} \binom{1}{m} \Bigg\rbrace = 0. 
	\end{align*}
	Also note that 
	\begin{align*}
		\sum_{l=0}^{i} (-1)^{l}   \binom{i}{m} \binom{i-m}{l} \binom{2i-1-l-m}{i-1-l} &= \text{ coefficient of } x^{i-1} \text{ in } \left( \binom{i}{m} \sum_{l=0}^{i} (-1)^{l}    \binom{i-m}{l} x^l (1+x)^{2i-1-l-m} \right) \\
		&= \text{ coefficient of } x^{i-1} \text{ in } \left( (1+x)^{i-1} \binom{i}{m} \sum_{l=0}^{i}     \binom{i-m}{l} (-x)^l (1+x)^{i-m-l} \right) \\ 
		& = \binom{i}{m}.
	\end{align*}
	Noting that $\binom{2i-1-l-m}{i-1-l}=\binom{2i-1-l-m}{i-m} $ and substituting this in \eqref{sum gamma a=2i-1 and n=1}, we obtain the sum vanishes. This proves the congruences in the lemma for $m=0,\ldots,i$.
	
	Note that by Lucas' theorem, $\binom{r-l}{i+1} \equiv \binom{i-l}{i+1} \equiv 0 \mod p$ for $l=0,\ldots i-2$. Thus 
	\begin{align*}
		\sum_{l=0}^{i-2} \gamma_l \sum\limits_{\substack{i < j < r-i+1 \\ j \equiv i \mod (p-1)  }} \binom{r-l}{j} \binom{j}{i+1} &= \sum_{l=0}^{i-2} \gamma_l  \binom{r-l}{i+1}\sum\limits_{\substack{i < j < r-i+1 \\ j \equiv i \mod (p-1)  }} \binom{r-l-(i+1)}{j-(i+1)} \\
		& \equiv 0 \mod p.
	\end{align*}
	This proves the congruences in the lemma for $m=i+1$.
\end{proof}  
\begin{theorem}\label{a=2i-1, left and right half, irred}
	Let  $r \equiv 2i-1 \mod{(p-1)}$ and $r\equiv i \mod p$ with $ 2 \leq i \leq \frac{p-1}{2} $ and $r \geq i(p+1)+p $ with $ v(a_{p})  \in (i,i+1) $. Let $-d = \frac{1}{p} \binom{r-i+1}{i} + \frac{(-1)^{i+1}}{i}$.  Assume that one of the following holds
	\begin{enumerate}
		\item[$(i)$] $v(a_p^2)< 2i+1$
		\item[$(ii)$] $v(a_p^2) \geq 2i+1$ and $p\mid d$.
	\end{enumerate} 
	Then we have
	\[ \frac{\ind_{KZ}^{G}(V_{p-2} \otimes D^{i})}{T} \twoheadrightarrow \bar{\Theta}_{k,a_p}.\]
\end{theorem}
\begin{proof} 
	Note that $p \geq 5$ since $4 \leq 2i \leq p-1 $. Let $\beta_0, \ldots, \beta_{i-1}$ be as in Lemma~\ref{beta l for a = 2i-1, n=1}, that is, $\beta_l = \frac{(-1)^{i-l+1}}{(i-l)} \binom{i}{l} $ for $l=0,\ldots,i-2$ and  $\beta_{i-1}=1/p$. Then by Lemma~\ref{lem:choice of beta}, there exist $\alpha_j$ such that
	\begin{enumerate}
		\item[(1)] $ \alpha_j \equiv \sum\limits_{l=0}^{i-1} p\beta_l \binom{r-l}{j} \mod p^2$ for $i \leq j<r-i+1$ and $j \equiv i \mod (p-1)$
		\item[$(2)$] $\sum\limits_{\substack{i \leq j<r-i+1 \\ j \equiv i \mod p} } \alpha_j \binom{j}{m} \equiv 0 \mod p^{i+3-m}$ for $m=0, \ldots,i-1$.
	\end{enumerate}
	Using $(1)$ and Lemma~\ref{beta l for a = 2i-1, n=1} we obtain
	\begin{enumerate}
		\item[$(2')$] $\sum\limits_{\substack{i \leq j<r-i+1 \\ j \equiv i \mod p} } \alpha_j \binom{j}{i} \equiv -pd \mod p^{2}$.
	\end{enumerate}
	Noting that by Lucas' theorem $p \beta_l \binom{r-l}{i} \equiv 0 \mod p^2$ for $1\leq l \leq i-2$ and noting that $p\beta_{i-1} =1$, we get $\alpha_i \equiv p\beta_0 \binom{r}{i}+ p\beta_{i-1} \binom{r-i+1}{i} \equiv p\beta_0 + \binom{r-i+1}{i} \equiv  -pd \mod p^2 $.  
	
	Let $\gamma_0, \ldots, \gamma_{i-2}$ be as in Lemma~\ref{gamma l for a = 2i-1, n=1}, that is, $\gamma_l = (-1)^l \binom{r}{l}$. Then again by Lemma~\ref{lem:choice of beta}, there exist $\alpha'_j$ such that
	\begin{enumerate}
		\item[$(3)$] $ \alpha'_j \equiv \sum\limits_{l=0}^{i-2}\gamma_l \binom{r-l}{j} \mod p$ for $i<j<r-i+1$ and $j \equiv i \mod (p-1)$
		\item[$(4)$] $\sum\limits_{\substack{i<j<r-i+1 \\ j \equiv i \text{ mod } (p-1)} } \alpha'_j \binom{j}{m} \equiv 0 \mod p^{i+3-m}$ for $m = 0, \ldots,i+1$.
	\end{enumerate}
	Consider
	\begin{align*}
		f_3 &=  
		d \sum_{  \lambda \in \mathbb{F}_{p}} \Bigg\lbrace \sum_{  \mu \in \mathbb{F}_{p}^{\times}}
		\Bigg[ g_{3,p[\lambda]+p^2[\mu]}^{0}, \sum_{l=0}^{i-2} \gamma_l  \frac{p^{2i-l}[\mu]^{l-i+1}}{a_p^2(p-1)} 
		(-\theta)^{l+2} X^{-2}Y^{r-(l+2)(p+1)+2} \Bigg]   \\
		& \qquad \qquad - \left[ g_{3,p[\lambda]}^{0}, \frac{p^{i+1}}{a_p^{2}}  \sum_{l=0}^{i-2} \gamma_{l}  \binom{r-l}{r-i+1}  (-\theta)^{i+1} X^{-2}Y^{r-(i+1)(p+1)+2} \right] \Bigg\rbrace \\
		f_{21} &=  
		\sum_{  \lambda \in \mathbb{F}_{p}^{\times}} 
		\Bigg[ g_{2,p[\lambda]}^{0}, \sum_{l=0}^{i-1} \frac{[\lambda]^{l+1-i}p^{i-l}}{a_p(p-1)} 
		\beta_{l} (-\theta)^{l+2} X^{-2}Y^{r-(l+2)(p+1)+2} \Bigg]  \\
		& \qquad \qquad - \left[ g_{2,0}^{0}, \frac{p}{a_p}  \Bigg( \sum_{l=0}^{i-1} \beta_{l}  \binom{r-l}{r-i+1}\Bigg)   (-\theta)^{i+1} X^{-2}Y^{r-(i+1)(p+1)+2} \right]  \\
		f_{22} &=  \sum_{  \lambda \in \mathbb{F}_{p}}  \left[g_{2,p[\lambda]}^0, d \frac{p^{2i}}{a_p^3} \sum_{ \substack{i <   j < r-i+1 \\  j \equiv i ~\mathrm{mod}~ (p-1)} } \alpha'_{j} X^{r-j} Y^{j} \right] \\
		f_1 &=  \left[g_{1,0}^0, \frac{p^{i-1}}{a_p^2} \sum_{ \substack{i \leq   j < r-i+1 \\  j \equiv i ~\mathrm{mod}~ (p-1)} } \alpha_{j} X^{r-j} Y^{j} \right]. 
	\end{align*}
	It follows from Lemma~\ref{theta and T plus} that $T^{+}f_3, T^+f_{21}$ vanish modulo $p$ if $v(a_p) \in (i,i+1)$. Clearly $-a_p f_3$ vanishes modulo $p$ if $v(a_p) \in (i,i+1)$. Also, using $(4)$ it follows that $T^+f_{22}$ also vanishes modulo $p$ if $v(a_p) \in (i,i+1)$. Using $(2)$,  $(2')$ and the fact $v(a_p^2) < 2i+2$ it can be checked from the formula of $T^+$ \eqref{T plus formula} that
	\begin{align*}
		T^+ f_1 &\equiv \left[  g_{2,0}^0 , \frac{p^{2i-1}}{a_p^2} \alpha_i X^{r-i}Y^{i}\right] + \sum_{  \lambda \in \mathbb{F}_{p}^\times }  \left[  g_{2,p[\lambda]}^0 , -\frac{p^{2i}}{a_p^2} d X^{r-i}Y^{i}\right] \\
		& \qquad \qquad + \sum_{  \lambda \in \mathbb{F}_{p}^\times }  \left[  g_{2,p[\lambda]}^0 , (-[\lambda])^{-1} \frac{p^{2i}}{a_p^2} \sum\limits_{\substack{i \leq j<r-i+1 \\ j \equiv i \text{ mod } (p-1)} } \alpha_j \binom{j}{i+1} X^{r-i-1}Y^{i+1}\right] \\ & \qquad \qquad + \sum_{  \lambda \in \mathbb{F}_{p}^\times }  \left[  g_{2,p[\lambda]}^0 , (-[\lambda])^{-2} \frac{p^{2i+1}}{a_p^2} \sum\limits_{\substack{i \leq j<r-i+1 \\ j \equiv i \text{ mod } (p-1)} } \alpha_j \binom{j}{i+2} X^{r-i-2}Y^{i+2}\right] + O(p).
	\end{align*}
	We now show the terms involving $X^{r-i-1}Y^{i+1}$ and $X^{r-i-2}Y^{i+2}$ in $T^+ f_1$ vanish modulo $p$. It follows from $(1)$  that
	\begin{align*}
		\sum_{\substack{i\leq j < r-i+1 \\ j \equiv i \mod (p-1)}}^{} \alpha_{j} \binom{j}{i+1} & \equiv \sum_{l=0}^{i-1} p \beta_{l} \sum_{\substack{i\leq j < r-i+1 \\ j \equiv i \text{ mod } (p-1)}}^{} \binom{r-l}{j} \binom{j}{i+1} \\
		& \equiv \sum_{l=0}^{i-1} p \beta_{l} \binom{r-l}{i+1} \sum_{\substack{i\leq j < r-i+1 \\ j \equiv i \text{ mod } (p-1)}}^{}  \binom{r-l-i-1}{j-i-1} \\
		& \equiv \binom{r-(i-1)}{i+1} \sum_{\substack{i\leq j< r-i+1 \\ j \equiv i \text{ mod } (p-1)}}^{}  \binom{r-2i}{j-i-1} \equiv 0 \mod p^2,
	\end{align*}
	where the penultimate  congruence follows since $p\mid p\beta_l$ and $p\mid \binom{r-l}{i+1}$ for $l=0, \ldots, i-2$ and the last step follows from \cite[Lemma 2.5]{BG15} and $\binom{r-(i-1)}{i+1} \equiv 0 \mod p$. Thus $X^{r-i-1}Y^{i+1}$ term in $T^{+}f_1$ vanishes modulo $p$. A similar check using $\binom{r-l}{i+2} \equiv 0 \mod p$ for $l=0,\ldots,i-1$, shows that $X^{r-i-2}Y^{i+2}$ term in $T^{+}f_1$ also vanishes modulo $p$. From the discussion below $(2')$, we have  $\alpha_i \equiv -pd \mod p^2$. Thus
	\begin{align}\label{T plus f_1, a=2i-1, n=1}
		T^+ f_1 &\equiv \sum_{  \lambda \in \mathbb{F}_{p}}  \left[  g_{2,p[\lambda]}^0 , -\frac{p^{2i}}{a_p^2} d X^{r-i}Y^{i}\right] \mod p.
	\end{align}            
	Using $p \mid p\beta_l$ for $0 \leq l \leq i-2$ and $p\beta_{i-1} =1$, we get 
	\begin{align*}
		-a_p f_{21} \equiv \sum_{  \lambda \in \mathbb{F}_{p}} \left[ g_{2,p[\lambda]}^{0},      (-\theta)^{i+1} X^{-2}Y^{r-(i+1)(p+1)+2} \right] \mod p.
	\end{align*}
	
	Note that the coefficient of $X^{i-1}Y^{r-i+1}$ coming from the \textquote{$\mu =0$} and the $\mu \neq 0$ terms in $T^{-}f_3$  cancels.  Thus 
	\begin{align*}
		T^-f_3 \equiv \sum_{  \lambda \in \mathbb{F}_{p}} \left[ g_{2,p[\lambda]}^{0},  d \frac{p^{2i}}{a_p^2}\sum_{\substack{0 \leq    j < r-i+1 \\  j \equiv i ~\mathrm{mod}~ (p-1)}} \sum_{l=0}^{i-2} \gamma_l \binom{r-l}{j} X^{r-j}Y^{j} \right]  \mod p.
	\end{align*}
	By Lucas' theorem and $r \equiv i \mod p$, it follows that $p \mid \binom{r-l}{i}$ for $1 \leq l \leq i-2$. Thus $\sum\limits_{l=0}^{i-2} \gamma_l \binom{r-l}{i} \equiv \gamma_0 \binom{r}{i} \equiv \gamma_0 \equiv 1 \mod p$. Thus the $X^{r-i}Y^{i}$  term in $T^- f_3$ cancels with $T^{+}f_1$ if $v(a_p^2) < 2i+1$ or if $v(a_p^2) \geq 2i+1$ and $p\mid d$. Thus
	\begin{align*}
		T^-f_3 -a_p f_{22} +T^+f_1 \equiv \sum_{  \lambda \in \mathbb{F}_{p}} \left[ g_{2,p[\lambda]}^{0},  {d} \frac{p^{2i}}{a_p^2} \sum_{\substack{i <    j < r-i+1 \\  j \equiv i ~\mathrm{mod}~ (p-1)}} \left( \sum_{l=0}^{i-2}  \gamma_l \binom{r-l}{j} - \alpha'_j \right) X^{r-j}Y^{j} \right]  \mod p.
	\end{align*}
	If $v(a_p^2) < 2i+1$, then the above expression vanishes modulo $p$ by (3).  If $v(a_p^2) \geq 2i+1$ and $p\mid d$, then once again the above expression vanishes modulo $p$.  Thus under the hypotheses of the lemma we have $T^-f_3 -a_p f_{22} +T^+f_1  \equiv 0 \mod p$. 
	
	It is easy to see that the coefficient of $X^{i-1}Y^{r-i+1}$ coming from the \textquote{$\lambda =0$} and the $\lambda \neq 0$ terms in $T^{-}f_{21}$  cancel each other. Also, we see that $T^- f_{22}$ vanishes modulo $p$  since the smallest power of $X$ appearing in $f_1$ equals $i-1+p-1 $ which is at least $i+3$ as $p\geq5$. Thus  
	\begin{align*}
		T^-f_{21}+T^{-}f_{22} -a_p f_{1}  \equiv  \left[ g_{1,0}^{0},  \frac{p^{i-1}}{a_p} \sum_{\substack{0 \leq    j < r-i+1 \\  j \equiv i ~\mathrm{mod}~ (p-1)}} \left( \sum_{l=0}^{i-1}  p \beta_l \binom{r-l}{j} - \alpha_j \right) X^{r-j}Y^{j} \right]  \equiv 0 \mod p
	\end{align*}
	by $(2)$.
	
	Similarly, $T^{-}f_1$ vanishes modulo $p$ using $p\geq 5$. Thus we finally get
	\begin{align*}
		(T-a_p) (f_3 + f_{21} + f_{22} +f_1) & \equiv -a_p f_{21} 
		\equiv  \sum_{  \lambda \in \mathbb{F}_{p}} \left[ g_{2,p[\lambda]}^{0},   F(X,Y) \right]    \mod p,            
	\end{align*}
	where $ F(X,Y) = (-\theta)^{i+1} X^{-2}Y^{r-(i+1)(p+1)+2} $. Note that $ F(X,Y) = (-\theta)^{i-1}  (Y^{r'} - 2X^{p-1} Y^{r'-(p-1)} + X^{2p-2} Y^{r'-2(p-1)})  $ where $r'=r-(i-1)(p+1)$.  Using Lemma~\ref{Glover-Brueil map image} applied with $r$ there equal to $r'$ and $m=0$ we get that the image of $F(X,Y)$ under the composition $V_r^{(i-1)}/V_{r}^{(i)} \rightarrow V_{r'}/V_{r'}^{(1)} \otimes D^{i-1} \twoheadrightarrow V_{p-2} \otimes D^i$ equals $(-1)^{i}X^{p-2}$. Thus $(T-a_p) (f_3 + f_{21} + f_{22} +f_1)$ maps to $(-1)^{i} \sum_{  \lambda \in \mathbb{F}_{p}} \left[ g_{2,p[\lambda]}^{0}, X^{p-2} \right] $ in $\ind_{KZ}^{G}(V_{p-2} \otimes D^{i})$. Noting that $ \sum_{  \lambda \in \mathbb{F}_{p}} [   g_{2,p[\lambda]}^{0}, X^{p-2} ] \equiv T([g_{1,0}^{0}, X^{p-2} ] ) \mod p$, we obtain the theorem.             
\end{proof}                    
\begin{remark}
	We now make some remarks on the proof of the above theorem.
	\begin{enumerate}
		\item For the convenience of the reader we summarize the idea of the proof since it is a bit more complicated than previous proofs. One starts with the function  $f_{21}$. To smoothen $T^- f_{21}$  we introduce $f_1$. As $T^+f_1$ doesn't vanish, one defines $f_3$ so that $T^{-}f_3$ cancels $T^+ f_1$. To smoothen the remaining terms of   $T^{-}f_3$ one introduces $f_{22}$.
		\item If $p^2 \mid d$ and $v(a_p^2) \geq 2i+1$, then $(T-a_p)(f_{21}+f_1)$ suffices to conclude. Indeed, if $p^2 \mid d$, then by \eqref{T plus f_1, a=2i-1, n=1}, we have $T^+ f_1$ vanishes. Thus we don't need $f_3$ 
		and hence  $f_{22}$. 
		\item By Lucas' theorem,
		\begin{align}\label{d mod p}
			\begin{split}
				-d &= \frac{(r-i)(r-i+1)}{i(i-1)p} \binom{r-i-1}{i-2}  + \frac{(-1)^{i+1}}{i} \\
				& \equiv \frac{(r-i)}{i(i-1)p} (-1)^i + \frac{(-1)^{i+1}}{i} \equiv \frac{(-1)^{i}(r-i-(i-1)p)}{i(i-1)p}  \mod p.   
			\end{split}
		\end{align}
		Thus we conclude that if $r \equiv i+(i-1)p =s \mod p^2$, then $p\mid d$ and the theorem applies.
	\end{enumerate}
\end{remark}
It remains to investigate what happens if $2i+1 \leq v(a_p^2)$ and $p\mid d$. We now show that the conclusion of the previous theorem still holds when $p\nmid d$ and $2i+1 < v(a_p^2)$. As we shall see later, the conclusion of the previous theorem does not always hold if $v(a_p^2)= 2i + 1$ (and  $p \nmid d$).

\begin{theorem}\label{a=2i-1, right half, p not divides d}
	Let  $r \equiv 2i-1 \mod{(p-1)}$ and $r\equiv i \mod p$ with $ 1 \leq i \leq \frac{p-1}{2} $ and $r \geq i(p+1)+p $ with $ v(a_{p})  \in (i,i+1) $. Let $-d = \frac{1}{p} \binom{r-i+1}{i} + \frac{(-1)^{i+1}}{i}$. If $v(a_p^2) > 2i+1$ and $p \nmid d$, then 
	\[ \frac{\ind_{KZ}^{G}(V_{p-2} \otimes D^{i})}{T} \twoheadrightarrow \bar{\Theta}_{k,a_p}.\]
\end{theorem}
\begin{proof}
	To treat the \textquote{right half interval} we multiply the functions $f_3$, $f_{21}$, $f_{22}$, $f_{1}$ used in the proof of Theorem~\ref{a=2i-1, left and right half, irred} in the \textquote{left half interval} by a scalar (this trick was also used in  \cite{BG15}, \cite{GR}). Consider
	\begin{align*}
		f'_3 = \frac{a_p^2}{p^{2i+1}}f_{3} &=  
		d \sum_{  \lambda \in \mathbb{F}_{p}^{\times}} \sum_{  \mu \in \mathbb{F}_{p}^{\times}}
		\Bigg[ g_{3,p[\lambda]+p^2[\mu]}^{0}, \sum_{l=0}^{i-2} \gamma_l  \frac{[\mu]^{l-i+1}}{p^{l+1}(p-1)} 
		(-\theta)^{l+2} X^{-2}Y^{r-(l+2)(p+1)+2} \Bigg]  \\
		& \qquad \qquad + d \left[ g_{3,0}^{0}, \frac{1}{p^i}  \sum_{l=0}^{i-2} \gamma_{l}  \binom{r-l}{r-i+1}  (-\theta)^{i+1} X^{-2}Y^{r-(i+1)(p+1)+2} \right]  \\
		f'_{21} =  \frac{a_p^2}{p^{2i+1}}f_{21} &=
		\sum_{  \lambda \in \mathbb{F}_{p}^{\times}} 
		\Bigg[ g_{2,p[\lambda]}^{0}, \sum_{l=0}^{i-1} \frac{[\lambda]^{i-1-l}a_p}{p^{i+1+l}(p-1)} \beta_{l} (-\theta)^{l+2} X^{-2}Y^{r-(l+2)(p+1)+2} \Bigg]  \\
		& \qquad \qquad - \left[ g_{2,0}^{0}, \frac{a_p}{p^{2i}}  \sum_{l=0}^{i-1} \beta_{l}  \binom{r-l}{r-i+1}  (-\theta)^{i+1} X^{-2}Y^{r-(i+1)(p+1)+2} \right]  \\
		f'_{22} =  \frac{a_p^2}{p^{2i+1}}f_{22} &= \sum_{  \lambda \in \mathbb{F}_{p}}  \left[g_{2,p[\lambda]}^0, d \frac{1}{pa_p} \sum_{ \substack{i <   j < r-i+1 \\  j \equiv i ~\mathrm{mod}~ (p-1)} } \alpha'_{j} X^{r-j} Y^{j} \right] \\
		f'_1 = \frac{a_p^2}{p^{2i+1}}f_{1} &= \left[g_{1,0}^0, \frac{1}{p^{i+2}} \sum_{ \substack{i \leq   j < r-i+1 \\  j \equiv i ~\mathrm{mod}~ (p-1)} } \alpha_{j} X^{r-j} Y^{j} \right].  
	\end{align*}
	Since $v(a_p^2) > 2i+1$, it follows that $T^+f'_3, -a_p f'_3,T^+f'_{21}$ and $T^+f'_{22}$ all vanish modulo $p$ as observed in Theorem~\ref{a=2i-1, left and right half, irred}. In addition, $-a_p f'_{21}$ now vanishes modulo $p$. A check similar to that in the proof of Theorem~\ref{a=2i-1, left and right half, irred} shows that
	\begin{align*}
		T^-f'_3 -a_p f'_{22} +T^+f'_1 \equiv \sum_{  \lambda \in \mathbb{F}_{p}} \left[ g_{2,p[\lambda]}^{0}, \frac{d}{p} \sum_{\substack{i <    j < r-i+1 \\  j \equiv i ~\mathrm{mod}~ (p-1)}} \left( \sum_{l=0}^{i-2}  \gamma_l \binom{r-l}{j} - \alpha'_j \right) X^{r-j}Y^{j} \right]  \mod p.
	\end{align*}
	Again as in the proof of Theorem~\ref{a=2i-1, left and right half, irred}, it follows 
	\begin{align*}
		T^-f'_{21}+T^{-}f'_{22} -a_p f'_{1}  \equiv  \left[ g_{1,0}^{0},  \frac{a_p}{p^{i+2}} \sum_{\substack{0 \leq    j < r-i+1 \\  j \equiv i ~\mathrm{mod}~ (p-1)}} \left( \sum_{l=0}^{i-1}  p\beta_l \binom{r-l}{j} - \alpha_j \right) X^{r-j}Y^{j} \right]  \equiv 0 \mod p.
	\end{align*}
	Note that $T^{-}f'_1$ vanishes modulo $p$, since the smallest power of $X$ appearing in $f'_1$ equals $i-1+p-1 $ which is at least $i+3$  as $p\geq 5$. Thus we get
	\begin{align*}
		(T-a_p) (f'_3 + f'_{21} + f'_{22} +f'_1) & \equiv  T^-f'_3 -a_p f'_{22} +T^+f'_1 \\
		&  \equiv \sum_{  \lambda \in \mathbb{F}_{p}} \left[ g_{2,p[\lambda]}^{0}, \frac{d}{p} \sum_{\substack{i < j < r-i+1 \\  j \equiv i ~\mathrm{mod}~ (p-1)}} \left( \sum_{l=0}^{i-2}  \gamma_l \binom{r-l}{j} - \alpha'_j \right) X^{r-j}Y^{j} \right]  \mod p.               
	\end{align*}
	Since the above function is not necessarily in $\mathrm{ind}_{KZ}^{G}(V_r^{(i-1)}/V_r^{(i)})$, we need to modify it slightly. Note that for   $0 \leq l \leq i-2$ $$ \sum\limits_{\substack{i \leq  j < r-i+1 \\  j \equiv i ~\mathrm{mod}~ (p-1)}} \binom{r-l}{j} X^{r-j}Y^{j} \equiv  -\sum\limits_{\lambda \in \mathbb{F}_p^\times} \lambda^{-i} X^l(X+\lambda Y)^{r-l} - \binom{r-l}{r-i+1} X^{i-1} Y^{r-i+1} \in \langle X^{r-i}Y^{i} \rangle . $$ We now claim that there exist $A_0, \ldots, A_{i-2} \in \mathbb{Z}_p$ such that
	\begin{align}\label{polynomial modification a=2i-1, n=1}
		F(X,Y) := \frac{d}{p} \sum_{\substack{i <    j < r-i+1 \\  j \equiv i ~\mathrm{mod}~ (p-1)}} \left( \sum_{l=0}^{i-2}  \gamma_l \binom{r-l}{j} - \alpha'_j \right) X^{r-j}Y^{j} - \sum_{l=0}^{i-2} A_l \sum_{\substack{i \leq  j < r-i+1 \\  j \equiv i ~\mathrm{mod}~ (p-1)}} \binom{r-l}{j} X^{r-j}Y^{j}
	\end{align}
	lies in $V_r^{(i-1)}$ and maps to a non-zero multiple of $X^{p-2}$ under the map $V_r^{(i-1)}/V_r^{(i)} \twoheadrightarrow V_{p-2} \otimes D^i $. By \cite[Lemma 2.7]{GR19} and \cite[Lemma 2.12]{GR19} it is sufficient to show that 
	\begin{align*}           	
		&\sum_{l=0}^{i-2} A_l \sum_{\substack{i \leq j < r-i+1 \\  j \equiv i ~\mathrm{mod}~ (p-1)}} \binom{r-l}{j}  \binom{j}{m} - \frac{d}{p} \sum_{\substack{i <    j < r-i+1 \\  j \equiv i ~\mathrm{mod}~ (p-1)}} \left( \sum_{l=0}^{i-2}  \gamma_l \binom{r-l}{j} - \alpha'_j \right) \binom{j}{m}  \in p\mathbb{Z}_{p} ~\text{ for } 0 \leq m \leq i-2 \\
		&\sum_{l=0}^{i-2} A_l \sum_{\substack{i \leq j < r-i+1 \\  j \equiv i ~\mathrm{mod}~ (p-1)}} \binom{r-l}{j}  \binom{j}{i-1} - \frac{d}{p} \sum_{\substack{i <    j < r-i+1 \\  j \equiv i ~\mathrm{mod}~ (p-1)}} \left( \sum_{l=0}^{i-2}  \gamma_l \binom{r-l}{j} - \alpha'_j \right) \binom{j}{i-1} \in \mathbb{Z}_{p}^\times.           	
	\end{align*}
	Since $ \sum_{ \substack{i <   j < r-i+1 \\  j \equiv i ~\mathrm{mod}~ (p-1)} }\alpha'_{j} \binom{j}{m} \equiv 0 \mod p^2$ for $m=0,1, \ldots, i-1$, it suffices to show that
	\begin{align}\label{linear equation p not divides d}
		\begin{split}
			&\sum_{l=0}^{i-2} A_l \sum_{\substack{i \leq j < r-i+1 \\  j \equiv i ~\mathrm{mod}~ (p-1)}} \binom{r-l}{j}  \binom{j}{m} - \frac{d}{p} \sum_{\substack{i <    j < r-i+1 \\  j \equiv i ~\mathrm{mod}~ (p-1)}} \sum_{l=0}^{i-2}  \gamma_l \binom{r-l}{j}  \binom{j}{m}  \equiv 0 \mod p ~\text{ for } 0 \leq m \leq i-2 \\
			&\sum_{l=0}^{i-2} A_l \sum_{\substack{i \leq j < r-i+1 \\  j \equiv i ~\mathrm{mod}~ (p-1)}} \binom{r-l}{j}  \binom{j}{i-1} - \frac{d}{p} \sum_{\substack{i <    j < r-i+1 \\  j \equiv i ~\mathrm{mod}~ (p-1)}}  \sum_{l=0}^{i-2}  \gamma_l \binom{r-l}{j}  \binom{j}{i-1} = (-1)^i  i d^2 \in \mathbb{Z}_{p}^\times.
		\end{split}	
	\end{align}
	This will be done in the next two lemmas.  Thus from the above claim and \cite[Lemma 2.12]{GR19} we get the image of $(T-a_p) (f'_3 + f'_{21} + f'_{22} +f'_1)$ under $\ind_{KZ}^{G} (V_r^{(i-1)}/V_r^{(i)}) \twoheadrightarrow \ind_{KZ}^{G}(V_{p-2} \otimes D^i) $  equals $-(-1)^i  i d^2 \sum_{\lambda \in \mathbb{F}_p}\left[ g_{2,p[\lambda]}^{0},  X^{p-2}\right]$. We  now conclude as in the proof of Theorem~\ref{a=2i-1, left and right half, irred}.
\end{proof}   

The following lemma generalises  Lemma~\ref{gamma l for a = 2i-1, n=1}.
\begin{lemma}\label{gamma l mod p^2 for a = 2i-1, n=1}
	Let  $r \equiv 2i-1 \mod{(p-1)}$ and $r \equiv i \mod p$ with $ 2 \leq i \leq (p-1)/2 $ and $r \geq p $. Let $\gamma_l =(-1)^{l} \binom{r}{l}$ for $l=0, \ldots, i-2$. Then we have the following
	\begin{align*}
		\sum\limits_{l=0}^{i-2} \gamma_l  \sum\limits_{\substack{i < j < r-i+1 \\ j \equiv i \mod (p-1)  }} \binom{r-l}{j} \binom{j}{m}  & \equiv  -(r-i) \binom{i}{m} \left(H_{i-m} + \frac{1}{i-1}\right)+ p \binom{i}{m} \\
		&\qquad \qquad + 
		\begin{cases}
			p \left( (-1)^i (i-1) + \frac{(-1)^i}{i} \right) & \text{ if } m=0, \\
			ip (-1)^i   & \text{ if } m=1, \\
			0  & \text{ if } 2 \leq m \leq i-1
		\end{cases}
		\mod p^2.
	\end{align*} 
	As a consequence, we obtain $\sum\limits_{l=0}^{i-2} \gamma_l  \sum\limits_{\substack{i < j < r-i+1 \\ j \equiv i \mod (p-1)  }} \binom{r-l}{j} \binom{j}{m}$ vanishes modulo $p$ for $m=0,\ldots, i-1$.
\end{lemma}
\begin{proof}
	We first claim that for $0\leq l \leq i-2$ and $0 \leq m \leq i-1$, we have
	\begin{align}
		\begin{split}	
			\sum\limits_{\substack{i \leq j \leq  r-i+1 \\ j \equiv i \mod (p-1)  }} \binom{r-l}{j} \binom{j}{m}  \equiv  & \binom{r-l}{m}  \binom{2i-1-l-m}{i-m} \\
			& \quad + \frac{r-2i+1}{p-1} \binom{r-l}{m} \Bigg( \binom{2i-1-l-m+p-1}{i-m}  \\ 
			& \quad \quad + \binom{2i-1-l-m+p-1}{i-m+p-1} - \binom{2i-1-l-m}{i-m} \Bigg) \mod p^2.
		\end{split}	
	\end{align}
	Indeed, observe that
	\begin{align}\label{binomial sum and roots of unity}
		\begin{split}
			\sum\limits_{\substack{i \leq j \leq  r-i+1 \\ j \equiv i \mod (p-1)  }} \binom{r-l}{j} \binom{j}{m}  &=   \binom{r-l}{m}  \sum\limits_{\substack{i-m \leq j \leq  r-i+1-m \\ j \equiv i-m \mod (p-1)  }} \binom{r-l-m}{j} \\ &=  \binom{r-l}{m}  \frac{1}{p-1}  \sum_{\xi \in \mu_{p-1}} \xi^{-(i-m)} (1+\xi)^{r-l-m}. 
		\end{split}
	\end{align}      	 
	Let $(1+\xi)^{p-1} = 1+p z_{\xi} $, for some $z_{\xi} \in \mathbb{Z}_p$, for every $\xi \in \mu_{p-1} \setminus\{-1\}$. Thus for $n \geq 0$, $0 \leq l \leq i-2$ and $0 \leq m \leq i-1$ we have
	\begin{align}\label{roots of unity a=2i-1}
		\begin{split}	
			\frac{1}{p-1} & \sum_{\xi \in \mu_{p-1}} \xi^{-(i-m)} (1+\xi)^{2i-1-l-m+n(p-1)} \\ 
			& \qquad \qquad \equiv \frac{1}{p-1} \sum_{\xi \in \mu_{p-1}} \xi^{-(i-m)} (1+\xi)^{2i-1-l-m} +  \frac{np}{p-1} \sum_{\xi \in \mu_{p-1}} \xi^{-(i-m)} (1+\xi)^{2i-1-l-m} z_{\xi}  \\
			&\qquad \qquad \equiv  \binom{2i-1-l-m}{i-m} + \frac{np}{p-1} \sum_{\xi \in \mu_{p-1}} \xi^{-(i-m)} (1+\xi)^{2i-1-l-m} z_{\xi} \mod p^2,
		\end{split}	
	\end{align}
	where we have used \eqref{sums of roots of unity} and $2 \leq 2i-1-l-m \leq  2i-1 \leq p-2$ in the last step.
	Taking $n=1$ in \eqref{roots of unity a=2i-1} and using we obtain
	\begin{align*}
		\frac{p}{p-1} \sum_{\xi \in \mu_{p-1}} \xi^{-(i-m)} (1+\xi)^{2i-1-l-m} z_{\xi} & \equiv \binom{2i-1-l-m+p-1}{i-m}+ \binom{2i-1-l-m+p-1}{i-m+p-1} \\
		& \qquad \qquad -\binom{2i-1-l-m}{i-m} \mod p^2. 
	\end{align*}
	Taking $n=(r-2i+1)/(p-1)$ in \eqref{roots of unity a=2i-1} we obtain
	\begin{align*}
		\frac{1}{p-1}  \sum_{\xi \in \mu_{p-1}} \xi^{-(i-m)} (1+\xi)^{r-l-m} & \equiv  \binom{2i-1-l-m}{i-m} + \frac{r-2i+1}{p-1} \Bigg( \binom{2i-1-l-m+p-1}{i-m}  \\ 
		& \quad \quad + \binom{2i-1-l-m+p-1}{i-m+p-1} - \binom{2i-1-l-m}{i-m} \Bigg) \mod p^2.
	\end{align*}
	Now the claim follows from the above congruence and \eqref{binomial sum and roots of unity}. 
	
	From the claim it follows that 
	\begin{align}\label{gamma sum mod p^2 with end points}
		\begin{split}
			\sum\limits_{l=0}^{i-2} \gamma_l  & \sum\limits_{\substack{i \leq j \leq  r-i+1 \\ j \equiv i \mod (p-1)  }}  \binom{r-l}{j} \binom{j}{m} \\
			&\equiv \sum\limits_{l=0}^{i-2} (-1)^l \binom{r}{l} \binom{r-l}{m}  \binom{2i-1-l-m}{i-m} + \frac{r-2i+1}{p-1} \sum\limits_{l=0}^{i-2} (-1)^l \binom{r}{l} \binom{r-l}{m} \Bigg( \binom{2i-1-l-m+p-1}{i-m}  \\ 
			& \qquad \qquad + \binom{2i-1-l-m+p-1}{i-m+p-1} - \binom{2i-1-l-m}{i-m} \Bigg) \mod p^2.
		\end{split}      
	\end{align}
	We now compute the first sum in the above congruence. Note that
	\begin{align}\label{first summand gamma a=2i-1, n=1}
		\begin{split}
			\sum\limits_{l=0}^{i-2} (-1)^l \binom{r}{l} \binom{r-l}{m}  \binom{2i-1-l-m}{i-m} &= \binom{r}{m} \sum\limits_{l=0}^{i-2} (-1)^l  \binom{r-m}{l}  \binom{2i-1-l-m}{i-m} \\
			&=\binom{r}{m} \sum\limits_{l=0}^{i-2} (-1)^l  \binom{r-m}{l}  \binom{2i-1-l-m}{i-1-l} \\
			&=\binom{r}{m} \sum\limits_{l=0}^{i-1} (-1)^l  \binom{r-m}{l}  \binom{2i-1-l-m}{i-1-l}  - (-1)^{i-1}  \binom{r}{m} \binom{r-m}{i-1} . 
		\end{split} 
	\end{align}
	Note that 
	\begin{align*}
		\sum\limits_{l=0}^{i-1} (-1)^l  \binom{r-m}{l}  \binom{2i-1-l-m}{i-1-l}  &= \text{coeff. of } x^{i-1} \text{ in } \sum\limits_{l=0}^{r-m} (-1)^l  \binom{r-m}{l} x^l (1+x)^{2i-1-l-m} \\
		& = \text{coeff. of } x^{i-1} \text{ in } (1+x)^{2i-r-1} \sum\limits_{l=0}^{r-m} (-1)^l  \binom{r-m}{l} x^l (1+x)^{r-m-l} \\
		&=   \text{coeff. of } x^{i-1} \text{ in } ((1+x)^{2i-r-1} \cdot 1) \\
		&= \frac{(2i-r-1)\cdots(i-r+1)}{(i-1)!}\\
		&\equiv 1 - (r-i) H_{i-1} \mod p^2,
	\end{align*}
	since $r \equiv i \mod p$. Substituting this in \eqref{first summand gamma a=2i-1, n=1} we obtain 
	\begin{align}\label{first summand gamma a=2i-1, n=1 final}
		\sum\limits_{l=0}^{i-2} (-1)^l \binom{r}{l} \binom{r-l}{m}  \binom{2i-1-l-m}{i-m} = \binom{r}{m} \left( 1 - (r-i) H_{i-1} + (-1)^{i}   \binom{r-m}{i-1} \right) \mod p^2.
	\end{align}
	We now compute the second sum in \eqref{gamma sum mod p^2 with end points}. Note that by Lucas' theorem, we have 
	\begin{align*}
		\binom{2i-1-l-m+p-1}{i-m} + \binom{2i-1-l-m+p-1}{i-m+p-1} &\equiv \binom{2i-2-l-m}{i-m} + \binom{2i-2-l-m}{i-m-1} \\ &\equiv  \binom{2i-1-l-m}{i-m} \mod  p .
	\end{align*}
	Thus, again by Lucas' theorem we have
	\begin{align*}
		&\sum\limits_{l=0}^{i-2} (-1)^l \binom{r}{l} \binom{r-l}{m} \Bigg( \binom{2i-1-l-m+p-1}{i-m} + \binom{2i-1-l-m+p-1}{i-m+p-1} - \binom{2i-1-l-m}{i-m} \Bigg) \\
		&\equiv \sum\limits_{l=0}^{i-2} (-1)^l \binom{i}{l} \binom{i-l}{m} \Bigg( \binom{2i-1-l-m+p-1}{i-m} + \binom{2i-1-l-m+p-1}{i-m+p-1} - \binom{2i-1-l-m}{i-m} \Bigg)  \\
		&\equiv \sum\limits_{l=0}^{i-2} (-1)^l \binom{i}{m} \binom{i-m}{l} \Bigg( \binom{2i-1-l-m+p-1}{i-m} + \binom{2i-1-l-m+p-1}{i-m+p-1} - \binom{2i-1-l-m}{i-m} \Bigg)  \\
		&\equiv \sum\limits_{l=0}^{i-2} (-1)^l \binom{i}{m} \binom{i-m}{l} \Bigg( \binom{2i-1-l-m+p-1}{i-1-l+p-1} + \binom{2i-1-l-m+p-1}{i-1-l} - \binom{2i-1-l-m}{i-1-l} \Bigg) \\
		&\equiv \sum\limits_{l=0}^{i} (-1)^l \binom{i}{m} \binom{i-m}{l} \Bigg( \binom{2i-1-l-m+p-1}{i-1-l+p-1} + \binom{2i-1-l-m+p-1}{i-1-l} - \binom{2i-1-l-m}{i-1-l} \Bigg)  \\
		&\qquad\qquad -(-1)^{i-1}\binom{i}{m} \binom{i-m}{i-1}\binom{i-m+p-1}{p-1} -(-1)^i \binom{i}{m}\binom{i-m}{i} \binom{i-m+p-2}{p-2}  \\
		&\equiv \sum\limits_{l=0}^{i}  \binom{i}{m} \text{ coefficient of } x^{i-1+p-1} \text{ in } \binom{i-m}{l} (-x)^l (1+x)^{2i-1-l-m+p-1} \\ &\qquad\qquad\qquad\qquad+ \sum\limits_{l=0}^{i}  \binom{i}{m} \text{ coefficient of } x^{i-1} \text{ in } \binom{i-m}{l} (-x)^l (1+x)^{2i-1-l-m+p-1} \\
		&\qquad\qquad\qquad\qquad - \sum\limits_{l=0}^{i}  \binom{i}{m} \text{ coefficient of } x^{i-1} \text{ in } \binom{i-m}{l} (-x)^l (1+x)^{2i-1-l-m} \\
		&\qquad\qquad\qquad\qquad  -(-1)^{i-1}\binom{i}{m} \binom{i-m}{i-1}\binom{i-m+p-1}{p-1} -(-1)^i \binom{i}{m}\binom{i-m}{i} \binom{i-m+p-2}{p-2} \\
		&\equiv \binom{i}{m} \text{ coefficient of } x^{i-1+p-1} \text{ in } (1+x)^{i-1+p-1} + \binom{i}{m} \text{ coefficient of } x^{i-1} \text{ in } (1+x)^{i-1+p-1} \\
		& \qquad\qquad\qquad\qquad  -   \binom{i}{m} \text{ coefficient of } x^{i-1} \text{ in } (1+x)^{i-1}  -(-1)^{i-1}\binom{i}{m} \binom{i-m}{i-1}\binom{i-m+p-1}{p-1} \\
		& \qquad\qquad\qquad \qquad -(-1)^i \binom{i}{m}\binom{i-m}{i} \binom{i-m+p-2}{p-2} \\
		& \equiv \binom{i}{m} \left(  \binom{i+p-2}{i-1}-(-1)^{i-1} \binom{i-m}{i-1}\binom{i-m+p-1}{p-1} -(-1)^i \binom{i-m}{i} \binom{i-m+p-2}{p-2} \right) \mod p^2.
	\end{align*}
	Note that  $\binom{i-m+p-1}{p-1}  \equiv p/(i-m) \mod p^2$ and hence $\binom{i-m+p-2}{p-2} \equiv -p/(i-m)(i-m-1) \mod p^2$. Thus 
	\begin{align*}
		&\sum\limits_{l=0}^{i-2} (-1)^l \binom{r}{l} \binom{r-l}{m} \Bigg( \binom{2i-1-l-m+p-1}{i-m} + \binom{2i-1-l-m+p-1}{i-m+p-1} - \binom{2i-1-l-m}{i-m} \Bigg) \\
		& \qquad \qquad \equiv \binom{i}{m} \left( \frac{p}{i-1} -(-1)^{i-1}  \binom{i-m}{i-1} \frac{p}{i-m} + (-1)^{i}  \binom{i-m}{i} \frac{p}{(i-m)(i-m-1)} \right) \\
		& \qquad \qquad \equiv \frac{p}{i-1}\begin{cases}  
			\left( 1 +(-1)^i (i-1) + (-1)^{i} \frac{1}{i} \right)&\text{ if } m =0,\\
			i (1+(-1)^i) &\text{ if } m =1,   \\ 		
			\binom{i}{m} &\text{ if } m \geq 2,
		\end{cases} \mod p^2.
	\end{align*}
	Note that $(r-2i+1)/(p-1) \equiv i-1 \mod p$. Using the above congruence and \eqref{first summand gamma a=2i-1, n=1 final} in \eqref{gamma sum mod p^2 with end points} we obtain
	\begin{align}\label{gamma sum mod p^2 with end points final}
		\sum\limits_{l=0}^{i-2} &\gamma_l  \sum\limits_{\substack{i \leq j \leq  r-i+1 \\ j \equiv i \mod (p-1)  }}  \binom{r-l}{j} \binom{j}{m}  \\ \nonumber & \equiv 
		\begin{cases}  
			\left( 1 - (r-i) H_{i-1} + (-1)^{i}   \binom{r}{i-1} \right)+ p\left( 1 +(-1)^i (i-1) + (-1)^{i} \frac{1}{i} \right) &\text{ if }  m=0,\\
			\binom{r}{1} \left( 1 - (r-i) H_{i-1} + (-1)^{i}   \binom{r-1}{i-1} \right) + ip  (1+(-1)^i) &\text{ if }  m=1,\\
			\binom{r}{m} \left( 1 - (r-i) H_{i-1} + (-1)^{i}   \binom{r-m}{i-1} \right) + p \binom{i}{m} &\text{ if } 2 \leq m \leq i-1
		\end{cases} \mod p^2.
	\end{align}
	
	We claim that 
	\begin{align}
		\sum\limits_{l=0}^{i-2} \gamma_l \binom{r-l}{i} \binom{i}{m} &\equiv \binom{i}{m} \left( 1+ \frac{r-i}{i-1} + \frac{r-i}{i} \right) \mod p^2 \label{left end point sum gamma a=2i-1, n=1}\\
		\sum\limits_{l=0}^{i-2} \gamma_l \binom{r-l}{r-i+1} \binom{r-i+1}{m} & \equiv (-1)^i \binom{r}{m}\binom{r-m}{i-1} \mod p^2 \label{right end point sum gamma a=2i-1, n=1}.
	\end{align}
	Subtracting \eqref{left end point sum gamma a=2i-1, n=1} and \eqref{right end point sum gamma a=2i-1, n=1} from \eqref{gamma sum mod p^2 with end points final},  we obtain 
	\begin{align*}
		\sum\limits_{l=0}^{i-2} &\gamma_l \sum\limits_{\substack{i < j <  r-i+1 \\ j \equiv i \mod (p-1)  }}  \binom{r-l}{j} \binom{j}{m}  \\ \nonumber & \equiv 
		\begin{cases}  
			- (r-i) H_{i-1}  -\frac{r-i}{i-1}-\frac{r-i}{i}  + p\left( 1 +(-1)^i (i-1) + (-1)^{i} \frac{1}{i} \right) &\text{ if }  m=0,\\
			- r(r-i) H_{i-1}  - i\frac{(r-i)}{i-1} + ip  (1+(-1)^i) &\text{ if }  m=1,\\
			\binom{r}{m} - \binom{i}{m} - \binom{r}{m}(r-i) H_{i-1}   -\frac{r-i}{i-1} \binom{i}{m} -\frac{r-i}{i} \binom{i}{m}  + p \binom{i}{m}  &\text{ if } 2 \leq m \leq i-1
		\end{cases} \mod p^2.
	\end{align*}
	This proves the lemma for $m=0$. Noting that $r(r-i) \equiv i(r-i) \mod p^2$ we obtain lemma for $m=1$. Applying Lemma~\ref{binomial coefficient under congruences} $(ii)$ (with $s$ and $t$ there equal to $i$ and $1$ respectively), we get $\binom{r}{m} - \binom{i}{m} \equiv (r-i) \binom{i}{m} (H_{i}-H_{i-m})$ mod $p^2$. Substituting this  above  and using $(r-i)\binom{r}{m} H_{i-1} \equiv (r-i)\binom{i}{m} H_{i-1}\mod p^2$ we obtain the lemma for $2 \leq m \leq i-1$. 
	
	Thus it remains to prove \eqref{left end point sum gamma a=2i-1, n=1} and \eqref{right end point sum gamma a=2i-1, n=1}. We first show \eqref{left end point sum gamma a=2i-1, n=1} holds. Observe that
	\begin{align*}
		\sum\limits_{l=0}^{i-2} \gamma_l \binom{r-l}{i}  &= \sum\limits_{l=0}^{i-2} (-1)^l \binom{r}{l} \binom{r-l}{i}  \\
		&= \sum\limits_{l=0}^{i} (-1)^l \binom{r}{l} \binom{r-l}{i}  -(-1)^{i-1} \binom{r}{i-1} \binom{r-i+1}{i}   -(-1)^{i} \binom{r}{i} \binom{r-i}{i} . 
	\end{align*}
	Note that by Lucas' theorem we have $\binom{r}{l} \binom{r-l}{i} \equiv \binom{i}{l} \binom{r-l}{i} \mod p^2$ for $l \geq 1$. For $l=0$ we have $\binom{r}{l} \binom{r-l}{i} = \binom{r}{i} = \binom{i}{l} \binom{r-l}{i} $. Thus
	\begin{align*}
		\sum\limits_{l=0}^{i-2} \gamma_l \binom{r-l}{i}  & \equiv \sum\limits_{l=0}^{i} (-1)^l \binom{i}{l} \binom{r-l}{i}  -(-1)^{i-1} \binom{r}{i-1} \binom{r-i+1}{i} -(-1)^{i} \binom{r}{i} \binom{r-i}{i}  \mod p^2. 
	\end{align*}	
	Note that 
	\begin{align*}
		\sum\limits_{l=0}^{i} (-1)^l \binom{i}{l} \binom{r-l}{i} 
		&= \sum\limits_{l=0}^{i} (-1)^l \binom{i}{l} \binom{r-l}{r-l-i} \\
		& = \text{ coefficient of } x^{r-i} \text{ in }  \sum\limits_{l=0}^{i}  \binom{i}{l} (-x)^l (1+x)^{r-l} \\
		& =\text{ coefficient of } x^{r-i} \text{ in } (1+x)^{r-i} \cdot 1 =1.
	\end{align*}
	Also by Lucas' theorem, we have 
	$$
	(-1)^i \binom{r}{i-1} \binom{r-i+1}{i} = (-1)^i \frac{(r-i+1)(r-i)}{i(i-1)} \binom{r}{i-1} \binom{r-i-1}{i-2} \equiv \frac{r-i}{i-1} \mod p^2.
	$$
	A similar computation shows that 
	$
	(-1)^{i+1} \binom{r}{i} \binom{r-i}{i} \equiv \frac{r-i}{i} \mod p^2
	$.
	Putting all these together we obtain  \eqref{left end point sum gamma a=2i-1, n=1}.
	
	We next show \eqref{right end point sum gamma a=2i-1, n=1} holds. Indeed, we have  
	\begin{align*}
		\sum\limits_{l=0}^{i-2} \gamma_l \binom{r-l}{r-i+1}   &= \sum\limits_{l=0}^{i-2} (-1)^l \binom{r}{l} \binom{r-l}{r-i+1}   \\
		& = \sum\limits_{l=0}^{i-2} (-1)^l \binom{r}{i-1} \binom{i-1}{l}   \\
		&=\binom{r}{i-1}   \sum\limits_{l=0}^{i-1} (-1)^l \binom{i-1}{l} -(-1)^{i-1} \binom{r}{i-1} = (-1)^{i} \binom{r}{i-1}.
	\end{align*}
	Multiplying both sides by $\binom{r-i+1}{m}$ and  noting that $ \binom{r-i+1}{m} \binom{r}{i-1} = \binom{r}{m}\binom{r-m}{i-1}$, we obtain \eqref{right end point sum gamma a=2i-1, n=1}. This proves  \eqref{left end point sum gamma a=2i-1, n=1}, \eqref{right end point sum gamma a=2i-1, n=1} and completes the proof of the lemma.
\end{proof}
We are now in a position to solve the congruences \eqref{linear equation p not divides d} appearing in Theorem~\ref{a=2i-1, right half, p not divides d}.
\begin{lemma}
	Let  $r \equiv 2i-1 \mod{(p-1)}$ and $r \equiv i \mod p$ with $ 2 \leq i \leq (p-1)/2 $ and $r \geq p $. Let $\gamma_l =(-1)^{l} \binom{r}{l}$ for $l=0, \ldots, i-2$. Then there exists $A_0, A_1, \ldots, A_{i-2} \in \mathbb{Z}_p$ such that
	\begin{align}
		&\sum_{l=0}^{i-2} A_l \sum_{\substack{i \leq j < r-i+1 \\  j \equiv i ~\mathrm{mod}~ (p-1)}} \binom{r-l}{j}  \binom{j}{m} - \frac{d}{p}  \sum_{l=0}^{i-2}  \gamma_l \sum_{\substack{i <    j < r-i+1 \\  j \equiv i ~\mathrm{mod}~ (p-1)}} \binom{r-l}{j}  \binom{j}{m}  \equiv 0 \mod p ~\text{ for } 0 \leq m \leq i-2 \label{linear eqns p not divides d small m}\\
		&\sum_{l=0}^{i-2} A_l \sum_{\substack{i \leq j < r-i+1 \\  j \equiv i ~\mathrm{mod}~ (p-1)}} \binom{r-l}{j}  \binom{j}{i-1} - \frac{d}{p}  \sum_{l=0}^{i-2} \gamma_l \sum_{\substack{i <    j < r-i+1 \\  j \equiv i ~\mathrm{mod}~ (p-1)}}  \binom{r-l}{j}  \binom{j}{i-1} \equiv (-1)^i i d^2 \mod p.\label{linear eqns p not divides d m=i-1}
	\end{align}            
\end{lemma}
\begin{proof}
	Let $s=i+(i-1)p$. Note that for $0\leq l \leq i-2$ and $0 \leq m \leq i-1$ we have   $\binom{r-l}{r-i+1} = \binom{r-l}{i-1-l} \equiv \binom{s-l}{i-1-l} = \binom{s-l}{s-i+1} \mod p$  and $\binom{r-i+1}{m} \equiv \binom{s-i+1}{m} \mod p$. Thus by  Corollary~\ref{cor: binomial sums under congruences 2}, we get 
	\begin{align}\label{binomial sum for a=2i-1, n=1}
		\sum_{\substack{i \leq j < r-i+1 \\  j \equiv i ~\mathrm{mod}~ (p-1)}} \binom{r-l}{j}  \binom{j}{m} \equiv \sum_{\substack{i \leq j < s-i+1 \\  j \equiv i ~\mathrm{mod}~ (p-1)}} \binom{s-l}{j}  \binom{j}{m}  \mod p  
	\end{align}
	$\text{ for } 0 \leq l \leq i-2 \text{ and } 0\leq m \leq i-1$. We first solve  \eqref{linear eqns p not divides d small m} and later show that \eqref{linear eqns p not divides d m=i-1} holds for the same values of $A_0, \ldots, A_{i-2}$.  By \eqref{binomial sum for a=2i-1, n=1} it is enough to solve 
	\begin{align}\label{reduction A_l a=2i-1, n=1}
		\sum_{l=0}^{i-2} A_l \sum_{\substack{i \leq j < s-i+1 \\  j \equiv i ~\mathrm{mod}~ (p-1)}} \binom{s-l}{j}  \binom{j}{m} \equiv  \frac{d}{p}  \sum_{l=0}^{i-2}  \gamma_l \sum_{\substack{i <    j < r-i+1 \\  j \equiv i ~\mathrm{mod}~ (p-1)}} \binom{r-l}{j}  \binom{j}{m}   \mod p ~\text{ for } 0 \leq m \leq i-2. 
	\end{align}
	Writing these congruence in matrix form we get 
	\begin{align*}
		A (A_0, \ldots, A_{i-2})^t \equiv \left( \frac{d}{p} \sum_{l=0}^{i-2} \sum_{\substack{i <    j < r-i+1 \\  j \equiv i ~\mathrm{mod}~ (p-1)}}  \gamma_l \binom{r-l}{j}  \binom{j}{0}, \ldots, \frac{d}{p}  \sum_{l=0}^{i-2} \sum_{\substack{i <    j < r-i+1 \\  j \equiv i ~\mathrm{mod}~ (p-1)}} \gamma_l \binom{r-l}{j}  \binom{j}{i-2} \right)^t \mod p, \\
	\end{align*}
	where $$ A = \left( \sum\limits_{\substack{i \leq    j < s-i+1 \\  j \equiv i ~\mathrm{mod}~ (p-1)}}\binom{s-l}{j}  \binom{j}{m} \right)_{0 \leq m,l \leq i-2}.$$
	Note that  every $i \leq  j < s-i+1 $ with $ j \equiv i ~\mathrm{mod}~ (p-1)$ can be expressed as $j=i+k(p-1)$ for some $0 \leq k \leq i-2$. Thus $A =BC$ with 
	\begin{align*}
		B= \left( \binom{i+k(p-1)}{m} \right)_{0 \leq m,k \leq i-2} \text{ and } C = \left( \binom{s-l}{i+k(p-1)} \right)_{0 \leq k,l \leq i-2}.
	\end{align*}
	By Corollary~\ref{cor: GV det} $(i)$, we get $B$ is invertible. By Lucas' theorem, we have $\binom{s-l}{i+k(p-1)} \equiv \binom{i-1}{k} \binom{i-l}{i-k}  \mod p$. Thus $\det(C) \equiv \prod_{k=0}^{i-2} \binom{i-1}{k} \det\left( \binom{i-l}{i-k}\right)_{0 \leq k,l \leq i-2} = \prod_{k=0}^{i-2} \binom{i-1}{k} \det\left( \binom{l}{k}\right)_{2 \leq k,l \leq i} $. As $\left( \binom{l}{k}\right)_{2 \leq k,l \leq i}$ is upper triangular with all the diagonal entries equal to $1$, we get $C$ is invertible. This shows that $A$ is invertible modulo $p$ and hence the congruences \eqref{reduction A_l a=2i-1, n=1} are solvable with $A_l$ in $\Z_p$. 
	
	We now determine the value $A_0$ modulo $p$ which we will need later. Let $A^{-1} = (a'_{lm})_{0\leq l,m \leq i-2} $. Then we have $$A_0 = \sum_{m=0}^{i-2} a'_{0m} \frac{d}{p} \sum_{l=0}^{i-2}  \gamma_l \sum_{\substack{i <    j < r-i+1 \\  j \equiv i ~\mathrm{mod}~ (p-1)}} \binom{r-l}{j}  \binom{j}{m}.$$
	Thus it is enough to determine the first row of $A^{-1}$, i.e., $a'_{0m}$. We claim that
	\begin{align}
		a'_{0m} \equiv (-1)^{i-m} (i-1-m) \mod p.
	\end{align}
	To prove this it is enough to show that 
	\begin{align}\label{A0 reduction a=2i-1, n=1}
		\sum_{m=0}^{i-2} (-1)^{i-m} (i-1-m) \sum_{\substack{i \leq j < s-i+1 \\  j \equiv i ~\mathrm{mod}~ (p-1)}} \binom{s-l}{j}  \binom{j}{m} \equiv \begin{cases}
			1 & \text{ if } l = 0, \\
			0 & \text{ if } 1 \leq l \leq i-2,
		\end{cases}\mod p.
	\end{align}
	Note that for $0\leq m \leq i-2$, we have
	\begin{align*}
		\sum_{m=0}^{i-2} (-1)^{i-m} (i-1-m) & \sum_{\substack{i \leq j < s-i+1 \\  j \equiv i ~\mathrm{mod}~ (p-1)}} \binom{s-l}{j}  \binom{j}{m} \\
		& =  \sum_{\substack{i \leq j < s-i+1 \\  j \equiv i ~\mathrm{mod}~ (p-1)}} \binom{s-l}{j}  \sum_{m=0}^{i-2} (-1)^{i-m} (i-1-m) \binom{j}{m} \\
		&=\sum_{k=0}^{i-2} \binom{s-l}{i+k(p-1)} \sum_{m=0}^{i-2} (-1)^{i-m} (i-1-m) \binom{i+k(p-1)}{m} \\
		& \equiv \sum_{k=0}^{i-2} \binom{s-l}{i+k(p-1)} \sum_{m=0}^{i-2} (-1)^{i-m} (i-1-m) \binom{i-k}{m} \\
		& \equiv \sum_{k=0}^{i-2} \binom{s-l}{i+k(p-1)} \sum_{m=0}^{i-1} (-1)^{i-m} (i-1-m) \binom{i-k}{m} \mod p,
	\end{align*}
	where in the penultimate step we used Lucas' theorem. Further note that for $0 \leq k \leq i-2$, we have 
	\begin{align*}
		\sum_{m=0}^{i-1} (-1)^{i-m} &(i-1-m) \binom{i-k}{m} \\
		&= \binom{i-k}{i} + \sum_{m=0}^{i} (-1)^{i-m} (i-1-m) \binom{i-k}{m} \\
		&=\binom{i-k}{i}+(-1)^{i} (i-1)  \sum_{m=0}^{i} (-1)^{m} \binom{i-k}{m} - (-1)^{i} \sum_{m=0}^{i} (-1)^{m} m \binom{i-k}{m} \\
		&= \binom{i-k}{i} + (-1)^{i} (i-1)  \sum_{m=0}^{i} (-1)^{m} \binom{i-k}{m} - (-1)^{i-1} (i-k) \sum_{m=0}^{i} (-1)^{m-1}  \binom{i-k-1}{m-1} \\
		& = \binom{i-k}{i} + (-1)^{i} (i-1)  \sum_{m=0}^{i-k} (-1)^{m} \binom{i-k}{m} - (-1)^{i-1} (i-k) \sum_{m=1}^{i-k} (-1)^{m-1}  \binom{i-k-1}{m-1} \\
		&= \binom{i-k}{i}.
	\end{align*}
	Substituting this above we get
	\begin{align*}
		\sum_{m=0}^{i-2} (-1)^{i-m} (i-1-m)  \sum_{\substack{i \leq j < s-i+1 \\  j \equiv i ~\mathrm{mod}~ (p-1)}} \binom{s-l}{j}  \binom{j}{m} & \equiv  \sum_{k=0}^{i-2} \binom{s-l}{i+k(p-1)} \binom{i-k}{i} \\
		& \equiv \binom{s-l}{i} \equiv \binom{i-l}{i} \mod p
	\end{align*}
	from which  \eqref{A0 reduction a=2i-1, n=1} follows. We obtain
	\begin{align}\label{A0 expression final a=2i-1}
		A_0 \equiv \sum_{m=0}^{i-2} (-1)^{i-m} (i-1-m)\cdot \frac{d}{p} \sum_{l=0}^{i-2}  \gamma_l \sum_{\substack{i <    j < r-i+1 \\  j \equiv i ~\mathrm{mod}~ (p-1)}} \binom{r-l}{j}  \binom{j}{m} \mod p.
	\end{align}
	
	We now prove \eqref{linear eqns p not divides d m=i-1}. To do this we express the LHS of \eqref{linear eqns p not divides d m=i-1} in terms of $A_0$ and the LHSs of the congruences \eqref{linear eqns p not divides d small m}. Note that for $0 \leq k \leq i-1$, we have
	\begin{align*}
		\sum_{m=0}^{i-1} (-1)^{m} \binom{i-k}{m} & = -(-1)^i \binom{i-k}{i} + \sum_{m=0}^{i} (-1)^{m} \binom{i-k}{m} \\
		& = -(-1)^i \binom{i-k}{i} + \sum_{m=0}^{i-k} (-1)^{m} \binom{i-k}{m} = -(-1)^i \binom{i-k}{i} .
	\end{align*}
	Hence
	\[
	(-1)^{i-1}\binom{i-k}{i-1} = (-1)^{i-1} \binom{i-k}{i} -\sum_{m=0}^{i-2} (-1)^{m} \binom{i-k}{m}.
	\]
	Thus for $0 \leq k \leq i-1$, by Lucas' theorem  and the above we have 
	\begin{align*}
		\binom{i+k(p-1)}{i-1} \equiv \binom{i-k}{i-1} & \equiv   \binom{i-k}{i} - \sum_{m=0}^{i-2} (-1)^{i-1-m} \binom{i-k}{m} \\ & \equiv \binom{i-k}{i} - \sum_{m=0}^{i-2} (-1)^{i-1-m} \binom{i+k(p-1)}{m} \mod p.
	\end{align*}
	Note that $\binom{i-k}{i} = \delta_{k0}$. Thus multiplying the above congruence with $\binom{s-l}{i+k(p-1)}$ and then taking the sum over $k=0,\ldots,i-2$, we have
	$$\sum_{k=0}^{i-2} \binom{s-l}{i+k(p-1)}\binom{i+k(p-1)}{i-1} \equiv \binom{s-l}{i} - \sum_{m=0}^{i-2} (-1)^{i-1-m} \sum_{k=0}^{i-2} \binom{s-l}{i+k(p-1)} \binom{i+k(p-1)}{m} \mod p.$$
	Multiplying the above congruence with $A_l$ and then taking the sum over $l=0,\ldots,i-2$, we get 
	\begin{align*}
		\sum_{l=0}^{i-2} A_{l} \sum\limits_{\substack{i\leq j < s-i+1 \\  j \equiv i ~\mathrm{mod}~ (p-1)}} \binom{s-l}{j}\binom{j}{i-1} & \equiv  \sum_{l=0}^{i-2} A_{l} \binom{s-l}{i} - \sum_{m=0}^{i-2} (-1)^{i-1-m} \sum_{l=0}^{i-2} A_{l} \sum\limits_{\substack{i\leq j < s-i+1 \\  j \equiv i ~\mathrm{mod}~ (p-1)}} \binom{s-l}{j} \binom{j}{m} \\
		& \equiv   A_0 - \sum_{m=0}^{i-2} (-1)^{i-1-m} \sum_{l=0}^{i-2} A_{l} \sum\limits_{\substack{i\leq j < s-i+1 \\  j \equiv i ~\mathrm{mod}~ (p-1)}} \binom{s-l}{j} \binom{j}{m} \mod p,
	\end{align*}
	where in the last step we used $\binom{s-l}{i} \equiv \binom{i-l}{i} \equiv \delta_{l0} \mod p$.
	Thus from \eqref{A0 expression final a=2i-1} and \eqref{reduction A_l a=2i-1, n=1} we get
	\begin{align*}
		\sum_{l=0}^{i-2} A_{l} \sum\limits_{\substack{i\leq j < s-i+1 \\  j \equiv i ~\mathrm{mod}~ (p-1)}} \binom{s-l}{j}\binom{j}{i-1} \equiv \sum_{m=0}^{i-2} (-1)^{i-m}  (i-m) \frac{d}{p} \sum_{l=0}^{i-2} \gamma_l \sum\limits_{\substack{i< j < r-i+1 \\  j \equiv i ~\mathrm{mod}~ (p-1)}} \binom{r-l}{j} \binom{j}{m}  \mod p.      
	\end{align*}
	By Lemma~\ref{gamma l mod p^2 for a = 2i-1, n=1} and \eqref{d mod p}, we have 
	\begin{align}\label{conclusion of d mod p and gamma sum mod p^2}
		\begin{split}
			\sum\limits_{l=0}^{i-2} \gamma_l  \sum\limits_{\substack{i < j < r-i+1 \\ j \equiv i \mod (p-1)  }} \binom{r-l}{j} \binom{j}{m}  & \equiv  -(r-i) \binom{i}{m} H_{i-m} + (-1)^{i} ipd \binom{i}{m} \\
			&\qquad \qquad + 
			\begin{cases}
				p \left( (-1)^i (i-1) + \frac{(-1)^i}{i} \right) & \text{ if } m=0, \\
				ip (-1)^i   & \text{ if } m=1, \\
				0  & \text{ if } 2 \leq m \leq i-1
			\end{cases}
			\mod p^2.
		\end{split}   
	\end{align} 
	Thus from above
	\begin{align*}
		\sum_{l=0}^{i-2} A_{l} & \sum\limits_{\substack{i\leq j < s-i+1 \\  j \equiv i ~\mathrm{mod}~ (p-1)}} \binom{s-l}{j}\binom{j}{i-1}\\
		& \equiv \sum_{m=0}^{i-2} (-1)^{i-m} (i-m)\binom{i}{m} \left\{ -\frac{d(r-i)}{p}  H_{i-m} +  (-1)^i id^2 \right\} + d \\
		& \equiv  -\sum_{m=0}^{i-2} (-1)^{i-m} (i-m)\binom{i}{m} \frac{d(r-i)}{p}  H_{i-m} + id^2 \sum_{m=0}^{i-2} (-1)^{m} (i-m)\binom{i}{m} +d \\
		&\equiv  -\frac{d(r-i)}{p} \sum_{m=0}^{i-2} (-1)^{i-m} (i-m)\binom{i}{m}   H_{i-m} + i^2d^2 \sum_{m=0}^{i-2} (-1)^{m} \binom{i-1}{m}+d \\
		&\equiv  -\frac{d(r-i)}{p}\sum_{m=0}^{i-2} (-1)^{i-m} (i-m)\binom{i}{m}   H_{i-m} +(-1)^{i} i^2d^2  +d   \mod p.
	\end{align*}
	By \eqref{conclusion of d mod p and gamma sum mod p^2} with $ 2 \leq m =i-1 \leq i-1$ we have 
	\begin{align*}
		\frac{d}{p} \sum_{l=0}^{i-2} \gamma_l \sum\limits_{\substack{i< j < r-i+1 \\  j \equiv i ~\mathrm{mod}~ (p-1)}} \binom{r-l}{j} \binom{j}{i-1}  \equiv -i\frac{d(r-i)}{p}  H_{1} + (-1)^i i^2d^2  \mod p.
	\end{align*}
	Subtracting the above two congruences we have
	\begin{align*}
		\sum_{l=0}^{i-2} A_{l} & \sum\limits_{\substack{i\leq j < r-i+1 \\  j \equiv i ~\mathrm{mod}~ (p-1)}} \binom{r-l}{j}\binom{j}{i-1} - \frac{d}{p} \sum_{l=0}^{i-2} \gamma_l \sum\limits_{\substack{i< j < r-i+1 \\  j \equiv i ~\mathrm{mod}~ (p-1)}} \binom{r-l}{j} \binom{j}{i-1}  \\
		& \equiv d- \frac{d(r-i)}{p} \sum_{m=0}^{i-1} (-1)^{i-m} (i-m)\binom{i}{m}   H_{i-m} \\
		& \equiv d- \frac{d(r-i)}{p} \sum_{m=1}^{i} (-1)^{m} m \binom{i}{m}   H_{m}\\
		& \equiv d- \frac{id(r-i)}{p} \sum_{m=1}^{i} (-1)^{m} \binom{i-1}{m-1}   H_{m} \\
		& \equiv d- \frac{d(r-i)}{(i-1)p} \equiv (-1)^i id^2\mod p.
	\end{align*}
	where in the penultimate step we used the mathematica program fastZeil to show that $\sum\limits_{m=1}^{i} (-1)^{m} \binom{i-1}{m-1}   H_{m} = \frac{1}{i(i-1)}$ and in last step we used \eqref{d mod p}. Now \eqref{linear eqns p not divides d m=i-1} follows from above and \eqref{binomial sum for a=2i-1, n=1} when $i \geq 3$.
	
	We now show \eqref{linear eqns p not divides d m=i-1} also holds for the case $i=2$. By \eqref{A0 expression final a=2i-1}
	\begin{align*}
		A_{0}   \equiv \frac{d}{p} \sum_{l=0}^{i-2} \gamma_l \sum\limits_{\substack{i< j < r-i+1 \\  j \equiv i ~\mathrm{mod}~ (p-1)}} \binom{r-l}{j} \mod p.
	\end{align*}
	Noting that $i=2$ it follows from Lemma~\ref{gamma l mod p^2 for a = 2i-1, n=1} that
	\[
	A_0 = \frac{d}{p} \left( -\frac{5}{2}(r-i) +  \frac{5}{2}p\right).
	\]
	By \cite[Lemma 2.14]{GR19} we have 
	\[
	\sum\limits_{\substack{i\leq j < r-i+1 \\  j \equiv i ~\mathrm{mod}~ (p-1)}} \binom{r}{j} \binom{j}{1} \equiv r \binom{2}{1} - \binom{r}{r-i+1} \binom{r-i+1}{1} = 2r-r(r-1) \equiv 4 - 2 \equiv 2 \mod p.
	\]
	Thus
	\[
	A_{0}  \sum\limits_{\substack{i\leq j < r-i+1 \\  j \equiv i ~\mathrm{mod}~ (p-1)}} \binom{r}{j} \binom{j}{1} \equiv
	\frac{d}{p} \left( -5(r-i) +  5p \right) \mod p.
	\]
	By Lemma~\ref{gamma l mod p^2 for a = 2i-1, n=1}, we have
	\[
	\frac{d}{p} \sum_{l=0}^{i-2} \gamma_l \sum\limits_{\substack{i< j < r-i+1 \\  j \equiv i ~\mathrm{mod}~ (p-1)}} \binom{r-l}{j} \binom{j}{1} \equiv \frac{d}{p} (-4(r-i)+4p)\mod p
	\]
	Hence 
	\[
	\sum_{l=0}^{i-2} A_{l}  \sum\limits_{\substack{i\leq j < r-i+1 \\  j \equiv i ~\mathrm{mod}~ (p-1)}} \binom{r-l}{j}\binom{j}{1} - \frac{d}{p} \sum_{l=0}^{i-2} \gamma_l \sum\limits_{\substack{i< j < r-i+1 \\  j \equiv i ~\mathrm{mod}~ (p-1)}} \binom{r-l}{j} \binom{j}{1} \equiv -\frac{d}{p} (r-i - p) \equiv (-1)^i id^2 \mod p,
	\]
	where we have used \eqref{d mod p} in the last step. This proves \eqref{linear eqns p not divides d m=i-1} for $i=2$ also.
\end{proof}
We now show that $\ind_{KZ}^G(V_{p-2} \otimes D^i) \twoheadrightarrow \bar{\Theta}_{k,a_p}$ factors through $T$ if $v(a_p^2 -  i d^2 p^{2i+1}) = 2i+1$ and $v(a_p^2) =2i+1$.
\begin{theorem}\label{a=2i-1, mid-point, irred}
	Let $r \geq i(p+1)+p $,  $r \equiv 2i-1 \mod{(p-1)}$ and $r\equiv i \mod p$ with $ 2 \leq i \leq \frac{p-1}{2} $ with $ v(a_{p})  = i + \frac{1}{2} $. Let $-d = \frac{1}{p} \binom{r-i+1}{i} + \frac{(-1)^{i+1}}{i}$.  If $v(a_p^2 - i d^2 p^{2i+1}) = 2i+1$, then
	\[ \frac{\ind_{KZ}^{G}(V_{p-2} \otimes D^{i})}{T} \twoheadrightarrow \bar{\Theta}_{k,a_p}.\]
\end{theorem}
\begin{proof} 
	We consider the functions as in Theorem~\ref{a=2i-1, right half, p not divides d}. It follows that $T^+f'_3, -a_p f'_3,T^+f'_{21}$ and $T^+f'_{22}$ all vanish modulo $p$ as observed in Theorem~\ref{a=2i-1, right half, p not divides d}. Again as in the proof of Theorem~\ref{a=2i-1, right half, p not divides d}, it follows $T^-f'_{21}+T^{-}f'_{22} -a_p f'_{1}  $ and $T^{-}f'_1$ both vanish modulo $p$.  A check similar to that in the proof of Theorem~\ref{a=2i-1, right half, p not divides d} shows that
	\begin{align*}
		T^-f'_3 -a_p f'_{22} +T^+f'_1 \equiv \sum_{  \lambda \in \mathbb{F}_{p}} \left[ g_{2,p[\lambda]}^{0}, \frac{d}{p} \sum_{\substack{i <    j < r-i+1 \\  j \equiv i ~\mathrm{mod}~ (p-1)}} \left( \sum_{l=0}^{i-2}  \gamma_l \binom{r-l}{j} - \alpha'_j \right) X^{r-j}Y^{j} \right]  \mod p.
	\end{align*}
	Thus we get
	\begin{align*}
		(T-a_p) (f'_3 + f'_{21} + f'_{22} +f'_1) & \equiv  T^-f'_3 -a_p f'_{22} +T^+f'_1  -a_p f'_{21}\\
		&  \equiv \sum_{  \lambda \in \mathbb{F}_{p}} \left[ g_{2,p[\lambda]}^{0}, F(X,Y) \right]  
		+ \sum_{  \lambda \in \mathbb{F}_{p}} \left[ g_{2,p[\lambda]}^{0},  \frac{a_p^2}{p^{2i+1}}    (-\theta)^{i+1} X^{-2}Y^{r-(i+1)(p+1)+2} \right] \\
		& \qquad \qquad \qquad \qquad\qquad \qquad \mod (p\Z_p+\langle X^{r-i}Y^i\rangle),              
	\end{align*}
	where $F(X,Y)$ is as in \eqref{polynomial modification a=2i-1, n=1}.
	Arguing as in Theorem~\ref{a=2i-1, left and right half, irred} and Theorem~\ref{a=2i-1, right half, p not divides d} we get $(T-a_p) (f'_3 + f'_{21} + f'_{22} +f'_1)$ maps to 
	\[
	\left( (-1)^i \frac{a_p^2}{p^{2i+1}} - (-1)^i i d^2\right)\sum_{  \lambda \in \mathbb{F}_{p}} \left[ g_{2,p[\lambda]}^{0}, X^{p-2} \right] 
	\]
	under the composition of the following maps $ \ind_{KZ}^{G}(V_{r}^{(i-1)}/V_{r}^{(i)}) \twoheadrightarrow \ind_{KZ}^{G}(V_{p-2} \otimes D^i)$. Under the hypothesis of the theorem we see that $(-1)^i \left( \frac{a_p^2}{p^{2i+1}} - i d^2\right)$ is a $p$-adic unit. This shows that $\ind_{KZ}^G(V_{p-2} \otimes D^i) \twoheadrightarrow \bar{\Theta}_{k,a_p}$ factors through $T$.
\end{proof}
\subsubsection{Above the diagonal}
We next turn our attention to eliminating JH factors  above the diagonal (i.e on and above the superdiagonal).                  
        \begin{lemma}\label{choice beta hybrid ad1}
            Let  $r \equiv a \mod{(p-1)}$ with $ 1 \leq a \leq p-1 $ and $r \geq i(p+1)+p $ with $ v(a_{p})  \in (i,i+1) $. Let $ s =  a-i+n+(i-n)p$ and $v(r-s) =t $  and $ 1 \leq t \leq  T < n < i < a$. If  $2i-a \leq 2T$ and $a <2i$, then there exist $\beta_0, \ldots, \beta_{i-T-1} \in \Z_p$ and $\gamma \in \Z_p$ such that 
        	    \begin{enumerate}
        	    	\item[$(i)$] $\sum\limits_{l=0}^{i-T-1} \beta_l  \sum\limits_{ \substack{i-T <  j < r-(a-i+T) \\  j \equiv i-T ~\mathrm{mod}~ (p-1)} }   \binom{r-l}{j} \binom{j}{m}  + \gamma p^t  \binom {i-T}{m}\equiv 
        	    	0  \mod p^{t+1} ~ \mathrm{for}~ m=0,\ldots, i-T-1$ 
        	    	\item[$(ii)$] $\sum\limits_{l=0}^{i-T-1} \beta_l  \sum\limits_{ \substack{i-T <  j < r-(a-i+T) \\  j \equiv i-T ~\mathrm{mod}~ (p-1)} }   \binom{r-l}{j} \binom{j}{i-T}   \equiv 
        	    	 p^t  \mod p^{t+1}$ 
        	    	\item[$(iii)$] $\sum\limits_{l=0}^{i-T-1} \beta_l \sum\limits_{ \substack{i-T <  j < r-(a-i+T) \\  j \equiv i-T ~\mathrm{mod}~ (p-1)} }   \binom{r-l}{j} \binom{j}{m} \equiv  0 \mod p^{t-v(m!)} ~  \mathrm{for}~ m=i-T+1,\ldots, i+t$.
        	    \end{enumerate}
        \end{lemma}
        \begin{proof}
        	     First we  prove  $(i)$ and $(ii)$. We now compute some binomial sums. By Corollary~\ref{cor: binomial sums under congruences 1}, we have
        	    \begin{align*}
        	    	\sum\limits_{ \substack{i-T <  j < r-(a-i+T) \\  j \equiv i-T ~\mathrm{mod}~ (p-1)} }   \binom{r-l}{j} \binom{j}{m} &\equiv \left( \binom{r-l}{m} - \binom{s-l}{m}\right)\left(\binom{[a-l-m]}{[i-T-m]}+\delta_{[i-T-m],p-1} \right) \\ &
        	    	\qquad \qquad +\sum_{\substack{ i-T < j < s-(a-i+T) \\ j \equiv i-T ~\mathrm{mod}~(p-1)}} \binom{s-l}{j} \binom{j}{m}  
        	    	+ \left(\binom{s-l}{i-T} - \binom{r-l}{i-T}\right) \binom{i-T}{m}  \\
        	    	&\qquad \qquad + \Bigg(\binom{s-l}{s-(a-i+T)} \binom{s-(a-i+T)}{m}  \\ & \qquad \qquad \qquad \qquad - \binom{r-l}{r-(a-i-T)} \binom{r-(a-i+T)}{m}\Bigg)   \mod p^{t+1}. 
        	    \end{align*}
        	      For $ 0 \leq  m \leq i-T$ and $0\leq l \leq i-T-1$, we have $l+m \leq 2i-2T-1 < a$, since $2i-a \leq 2T$. Thus  $[a-l-m] = a-l-m $.
                 If $m=i-T$, then we have  $\binom{[a-l-m]}{[i-T-m]} = \binom{a-i+T-l}{p-1} = 0$. Also, for $m=i-T$, then $\delta_{[i-T-m],p-1} = 1= \binom{a-l-m}{0} =\binom{a-l-m}{i-T-m}$.  Thus for $0 \leq m \leq i-T$, we have
                \begin{align}\label{binomial sum  hybrid ad}
                \begin{split}
        	    	\sum\limits_{ \substack{i-T <  j < r-(a-i+T) \\  j \equiv i-T ~\mathrm{mod}~ (p-1)} }   \binom{r-l}{j} \binom{j}{m} &\equiv \left( \binom{r-l}{m} - \binom{s-l}{m}\right)\binom{a-l-m}{i-T-m}  \\ &
        	    	\qquad \qquad +\sum_{\substack{ i-T < j < s-(a-i+T) \\ j \equiv i-T ~\mathrm{mod}~(p-1)}} \binom{s-l}{j} \binom{j}{m}  
        	    	+ \left(\binom{s-l}{i-T} - \binom{r-l}{i-T}\right) \binom{i-T}{m}  \\
        	    	&\qquad \qquad + \Bigg(\binom{s-l}{s-(a-i+T)} \binom{s-(a-i+T)}{m}  \\ & \qquad \qquad \qquad \qquad - \binom{r-l}{r-(a-i-T)} \binom{r-(a-i+T)}{m}\Bigg)   \mod p^{t+1}. 
              \end{split}      
        	\end{align}
               This finishes the computation of the binomial sums that we need. 
               
               We now derive two identities involving binomial coefficients which will be needed later while performing some row operations. 
               Note that 
               \begin{align*}
                   \binom{r-l}{m-k} &\binom{a-l-m}{i-T-m} - \frac{i-T-m+1+kp}{m-k} \binom{r-l}{(m-1)-k} \binom{a-l-(m-1)}{i-T-(m-1)} \\
                   &= \binom{r-l}{m-k} \binom{a-l-m}{i-T-m}- \frac{i-T-m+1}{m-k} \binom{r-l}{m-1-k}\binom{a-l-(m-1)}{i-T-(m-1)} \\ 
                   & \qquad \qquad \qquad \qquad - \frac{kp}{m-k} \binom{r-l}{m-1-k} \binom{a-l-(m-1)}{i-T-(m-1)} \\
                   &= \frac{r-a+k}{m-k} \binom{r-l}{m-1-k}\binom{a-l-m}{i-T-m}- \frac{kp}{m-k} \binom{r-l}{m-1-k} \binom{a-l-(m-1)}{i-T-(m-1)}.
               \end{align*}
               Multiplying $\prod_{k'=0}^{k-1} \frac{(r-a+k')}{m-k'}$ on both sides we get
               \begin{align*}
                   \prod_{k'=0}^{k-1} \frac{(r-a+k')}{m-k'} &\left(\binom{r-l}{m-k} \binom{a-l-m}{i-T-m} - \frac{i-T-m+1+kp}{m-k} \binom{r-l}{(m-1)-k} \binom{a-l-(m-1)}{i-T-(m-1)}\right) \\
                   &= \prod_{k'=0}^{k} \frac{(r-a+k')}{m-k'} \binom{r-l}{m-1-k}\binom{a-l-m}{i-T-m} \\
                   &\qquad \qquad \qquad - \frac{kp}{m-k} \prod_{k'=0}^{k-1} \frac{(r-a+k')}{m-k'}  \binom{r-l}{m-1-k} \binom{a-l-(m-1)}{i-T-(m-1)}.
               \end{align*}
               By the same computation as above, we get
               \begin{align*}
                   \prod_{k'=0}^{k-1} \frac{(s-a+k')}{m-k'} &\left(\binom{s-l}{m-k} \binom{a-l-m}{i-T-m} - \frac{i-T-(m-1)+kp}{m-k} \binom{s-l}{(m-1)-k} \binom{a-l-(m-1)}{i-T-(m-1)}\right) \\
                   &= \prod_{k'=0}^{k} \frac{(s-a+k')}{m-k'} \binom{s-l}{m-1-k}\binom{a-l-m}{i-T-m} \\
                   &\qquad \qquad \qquad - \frac{kp}{m-k} \prod_{k'=0}^{k-1} \frac{(s-a+k')}{m-k'}  \binom{s-l}{m-1-k} \binom{a-l-(m-1)}{i-T-(m-1)}.
               \end{align*}
               Subtracting, we see that
               \begin{align}\label{identity 2 adhc}
               \begin{split}
                   &\frac{1}{m \cdots (m-k+1)} \Bigg( \left(\prod_{k'=0}^{k-1} (r-a+k') \right) \binom{r-l}{m-k} -\left(\prod_{k'=0}^{k-1} (s-a+k')\right)\binom{s-l}{m-k}  \Bigg) \binom{a-l-m}{i-T-m} \\ 
                  & \qquad - \frac{i-T-(m-1)+kp}{m} \frac{1}{(m-1) \cdots (m-k)}\Bigg( \left(\prod_{k'=0}^{k-1} (r-a+k')\right) \binom{r-l}{m-1-k} \\ & \qquad\qquad \qquad \qquad \qquad \qquad \qquad \qquad \qquad \qquad \qquad \qquad \qquad-\left(\prod_{k'=0}^{k-1} (s-a+k')\right) \binom{s-l}{m-1-k}  \Bigg) \binom{a-l-(m-1)}{i-T-(m-1)} \\
                   &\qquad = \frac{1}{m \cdots (m-k)} \Bigg( \left(\prod_{k'=0}^{k} (r-a+k')\right) \binom{r-l}{m-1-k} -\left(\prod_{k'=0}^{k} (s-a+k') \right)\binom{s-l}{m-1-k}  \Bigg) \binom{a-l-m}{i-T-m} 
                   \\ &\qquad\qquad\qquad\qquad\qquad\qquad\qquad\qquad \qquad\qquad\qquad\qquad\qquad\qquad\qquad\qquad\qquad\qquad\qquad\qquad+O(p(r-s)).
               \end{split}    
               \end{align}
               Also, note that for $j' \geq m \geq 0$ and $j \geq 0$ we have
                \begin{align}\label{identity 1 adhc}
                    \binom{j'}{m} - \frac{j-(m-1)}{m} \binom{j'}{m-1}  = \frac{j'-j}{m} \binom{j'}{m-1}.
                \end{align}
                
               We now begin solving for $\beta_l$ and $\gamma$ satisfying the congruences $(i)$ and $(ii)$.  Let
               \begin{align}\label{A matrix green right}
               A= \left(\begin{array}{@{}c|c@{}}
               \begin{matrix}
                \left( \binom{r-l}{m} - \binom{s-l}{m}\right)\binom{a-l-m}{i-T-m} +   \sum\limits_{k=1}^{i-n-1} \binom{s-l}{i-T+k(p-1)} \binom{i-T+k(p-1)}{m}   \\ 
        	    	+ \left(\binom{s-l}{i-T} - \binom{r-l}{i-T}\right) \binom{i-T}{m}   \\
        	    	+ \binom{s-l}{s-(a-i+T)} \binom{s-(a-i+T)}{m}   - \binom{r-l}{r-(a-i+T)} \binom{r-(a-i+T)}{m}  
                \end{matrix}    & p^t\binom{i-T}{m} -p^t \delta_{i-T,m} 
               \end{array}
               \right)_{m=0,\ldots, i-T},
               \end{align}
               where $l$ ranges from $0$ to $i-T-1$. By \eqref{binomial sum  hybrid ad}, it is enough to solve for $\beta_0, \ldots, \beta_{i-T-1}$ and $\gamma$ in $\Z_p$ satisfying
               \begin{align}\label{eq: linear eqns adhc green}
                   A (\beta_0, \ldots, \beta_{i-T-1}, \gamma)^{\text{tr}} = (0,\ldots, 0,p^t)^{\text{tr}}.
                \end{align}
                We will use Cramer's rule. We first show that the determinant of $A$ is non-zero. To show that $\beta_l \in \Z_p$, we also need to determine the power of $p$ dividing $\det(A)$. To this end we show $p^{t(n-T+2)} \parallel \det(A)$. To achieve this we will perform several row operations so that the first row and the last $(n-T+1)$ rows are multiples of $p^t$. Apply the following row operations to $A$:
            \begin{alignat*}{4}
                R_{i-T-1} &\rightarrow R_{i-T-1}~- ~\frac{i-T-(i-T-2)}{i-T-1} R_{i-T-2} \\
                  & ~~\vdots  \\
                R_{m} &\rightarrow ~ ~~R_{m}~ - ~\frac{i-T-(m-1)}{m} R_{m-1} \\
                 & ~~\vdots  \\
                R_{1} &\rightarrow ~~ R_{1}~- ~\frac{i-T}{1} R_{0}.
            \end{alignat*}   
            Using \eqref{identity 2 adhc} with $k$ there equal to $0$ for the first term,  \eqref{identity 1 adhc} with $j$  equal to $i-T$ and $j'= i-T+k(p-1)$ for $k=1,\ldots,i-n-1$ for the second term, $j'=i-T$ for the third term, $j'=s-(a-i+T)$ for the fourth term, $j'=r-(a-i+T)$ for the fifth term, for all columns but the last, and $j'=i-T$ for the last column,  we get 
            \begin{align*}
                \det(A)= \left\lvert\begin{array}{c|c}
                    \sum\limits_{k=1}^{i-n-1} \binom{s-l}{i-T+k(p-1)}  + \binom{s-l}{i-T} - \binom{r-l}{i-T} + \binom{s-l}{s-(a-i+T)} - \binom{r-l}{r-(a-i+T)}  & p^t  \\ 
                    \hline \\
                     \left( \binom{r-l}{m-1} \frac{r-a}{m}- \binom{s-l}{m-1} \frac{s-a}{m} \right) \binom{a-l-m}{i-T-m}\\
                    + \sum\limits_{k=1}^{i-n-1} \binom{s-l}{i-T+k(p-1)} \binom{i-T+k(p-1)}{m-1} \frac{k(p-1)}{m}  & 0 \\
                    + \binom{s-l}{s-(a-i+T)} \binom{s-(a-i+T)}{m-1} \frac{s-a}{m} - \binom{r-l}{r-(a-i+T)}\binom{r-(a-i+T)}{m-1} \frac{r-a}{m} +O(p(r-s)) \\ 
                    \hline \\
                   \sum\limits_{k=1}^{i-n-1} \binom{s-l}{i-T+k(p-1)}\binom{i-T+k(p-1)}{i-T}+ \binom{s-l}{s-(a-i+T)} \binom{s-(a-i+T)}{i-T}- \binom{r-l}{r-(a-i+T)}\binom{r-(a-i+T)}{i-T} & 0
                \end{array}
                \right\rvert
            \end{align*}
            where the range of $m$ in the upper blocks is $0$, the range of $m$ in the middle blocks is $1,\ldots,i-T-1$ and the range of $m$ in the lower blocks is $i-T$. We now focus on the middle blocks.
            Note that  $\binom{s-l}{i-T}$ and $\binom{r-l}{i-T}$ don't appear from the second row onwards. We now remove the term $\binom{s-l}{i-T+(p-1)}$ and $\binom{r-l}{i-T+(p-1)}$ from the third row onwards. To do this, we apply the following set of row operations: 
            \begin{alignat*}{3}
                R_{i-T-1} &\rightarrow  ~R_{i-T-1} ~-~ \frac{i-T-(i-T-2)+p}{i-T-1} R_{i-T-2} \\
                  & ~~\vdots \\
                R_{m} &\rightarrow ~ R_{m}~-~ \frac{i-T-(m-1)+p}{m} R_{m-1} \\
                & ~~\vdots \\
                R_{2} &\rightarrow~ R_{2} ~-~ \frac{i-T-1+p}{2} R_{1}.
            \end{alignat*}
            Using  \eqref{identity 2 adhc} with $k$ there equal to $1$ for the first term and \eqref{identity 1 adhc}  with $j$ there equal to $i-T+(p-1)$ and $j'=i-T+(p-1), \ldots,  i-T+(i-n)(p-1)=s-(a-i+T)$  for the remaining terms, we get 
            \begin{align*}
                \det(A)= \left\lvert\begin{array}{c|c}                    
                    \left( \prod\limits_{k'=0}^{m-1} (r-a+k') - \prod\limits_{k'=0}^{m-1} (s-a+k') \right) \frac{1}{m!} \binom{a-l-m}{i-T-m} \\
                    +\sum\limits_{k=1}^{i-n-1} \binom{s-l}{i-T+k(p-1)}  \frac{1}{m!} \prod\limits_{k'=0}^{m-1} ((k-k')(p-1)) \\
                    + \left(\binom{s-l}{i-T} - \binom{r-l}{i-T} \right) \frac{1}{m!} \prod\limits_{k'=0}^{m-1} (k'(1-p)) & p^t \delta_{m,0} \\+ \binom{s-l}{s-(a-i+T)} \frac{1}{m!} \prod\limits_{k'=0}^{m-1} (s-a-k'(p-1)) - \binom{r-l}{r-(a-i+T)}\frac{1}{m!} \prod\limits_{k'=0}^{m-1} (r-a-k'(p-1)) +O(p(r-s)) \\ 
                    \hline \\
                    \left( \binom{r-l}{m-2} \prod\limits_{k'=0}^{1} (r-a+k') - \binom{s-l}{m-2} \prod\limits_{k'=0}^{1} (s-a+k') \right) \frac{1}{m(m-1)}\binom{a-l-m}{i-T-m} \\
                    +\sum\limits_{k=1}^{i-n-1} \binom{s-l}{i-T+k(p-1)} \binom{i-T+k(p-1)}{m-2}  \frac{(p-1)^2k(k-1)}{m(m-1)}   & 0 \\
                    + \binom{s-l}{s-(a-i+T)} \binom{s-(a-i+T)}{m-2} \prod\limits_{k'=0}^{1} (s-a-k'(p-1))  \frac{1}{m(m-1)}\\ - \binom{r-l}{r-(a-i+T)}\binom{r-(a-i+T)}{m-2} \prod\limits_{k'=0}^{1} (r-a-k'(p-1))  \frac{1}{m(m-1)} + O(p(r-s))\\ 
                    \hline \\
                   \sum\limits_{k=1}^{i-n-1} \binom{s-l}{i-T+k(p-1)}\binom{i-T+k(p-1)}{i-T}+ \binom{s-l}{s-(a-i+T)} \binom{s-(a-i+T)}{i-T}- \binom{r-l}{r-(a-i+T)}\binom{r-(a-i+T)}{i-T}& 0
                \end{array}
                \right\rvert,
            \end{align*}
            where the range of $m$ in the upper blocks is $0,1$, the range of $m$ in the middle blocks is $2,\ldots,i-T-1$ and the range of $m$ in the lower blocks is $i-T$.
            More generally, we may remove the terms $\binom{s-l}{i-T+k(p-1)}$, $\binom{r-l}{i-T+k(p-1)}$ from the $(k+2)^{\mathrm{th}}$-row onwards for every $k=0,1,\ldots,i-n-1$. To do this, we apply the following set of row operations: 
            \begin{align*}
                R_{i-T-1} &\rightarrow ~R_{i-T-1}- \frac{i-T-(i-T-2)+kp}{i-T} R_{i-T-2} \\
                  &~~\vdots  \\
                R_{m} &\rightarrow ~ R_{m}- \frac{i-T-(m-1)+kp}{m} R_{m-1} \\\
                 &~~\vdots  \\
                R_{k+1} &\rightarrow ~ R_{k+1}- \frac{i-T-k+kp}{k+1} R_{k}
            \end{align*}
            for $k=0,\ldots, i-n-1$. Using \eqref{identity 2 adhc} for $k$ and \eqref{identity 1 adhc} with $j$ there equal to $i-T+k(p-1)$ and $j'= i-T+k(p-1), \ldots, i-T+(i-n)(p-1) = s-(a-i+T)$ we get
            \begin{align*}
            \det(A)= \left\lvert \begin{array}{c|c}
                    \left( \prod\limits_{k'=0}^{m-1} (r-a+k') - \prod\limits_{k'=0}^{m-1} (s-a+k') \right) \frac{1}{m!} \binom{a-l-m}{i-T-m} \\
                   + \sum\limits_{k=1}^{i-n-1} \binom{s-l}{i-T+k(p-1)}  \frac{1}{m!} \prod\limits_{k'=0}^{m-1} ((k-k')(p-1)) \\
                   + \left(\binom{s-l}{i-T} - \binom{r-l}{i-T} \right) \frac{1}{m!} \prod\limits_{k'=0}^{m-1} (k'(1-p))) & p^t \delta_{m,0} \\+ \binom{s-l}{s-(a-i+T)} \frac{1}{m!} \prod\limits_{k'=0}^{m-1} (s-a-k'(p-1)) \\ - \binom{r-l}{r-(a-i+T)} \frac{1}{m!} \prod\limits_{k'=0}^{m-1} (r-a-k'(p-1)) + O(p(r-s)) \\ 
                    \hline \\                   
                    \left( \binom{r-l}{m-(i-n)} \prod\limits_{k'=0}^{i-n-1} (r-a+k') - \binom{s-l}{m-(i-n)} \prod\limits_{k'=0}^{i-n-1} (s-a+k') \right) \frac{(m-(i-n))!}{m!}\binom{a-l-m}{i-T-m} \\
                    +\sum\limits_{k=1}^{i-n-1} \binom{s-l}{i-T+k(p-1)} \binom{i-T+k(p-1)}{m-(i-n)} \frac{(m-(i-n))!}{m!} \prod\limits_{k'=0}^{i-n-1} ((k-k')(p-1)) & 0 \\
                    + \binom{s-l}{s-(a-i+T)} \binom{s-(a-i+T)}{m-(i-n)} \frac{(m-(i-n))!}{m!} \prod\limits_{k'=0}^{i-n-1} (s-a-k'(p-1))  \\ - \binom{r-l}{r-(a-i+T)}\binom{r-(a-i+T)}{m-(i-n)} \frac{(m-(i-n))!}{m!} \prod\limits_{k'=0}^{i-n-1} (r-a-k'(p-1))   + O(p(r-s)) \\ 
                    \hline \\                   
                    \sum\limits_{k=1}^{i-n-1} \binom{s-l}{i-T+k(p-1)}\binom{i-T+k(p-1)}{i-T}+ \binom{s-l}{s-(a-i+T)} \binom{s-(a-i+T)}{i-T}- \binom{r-l}{r-(a-i+T)}\binom{r-(a-i+T)}{i-T}& 0
                \end{array}
                \right\rvert,
            \end{align*}
            where the range of $m$ in the  upper blocks is $0,\ldots,i-n-1$, the range of $m$ in the middle blocks is $i-n,\ldots, i-T-1$ and the range of $m$ in the lower blocks is $i-T$. We now simplify the entries in the middle left block and show that they are multiples of $p^t$. For the first term, we have
            \begin{align}\label{est. error term 1 adhc}
            \begin{split}
                \Bigg( &\binom{r-l}{m-(i-n)}  \prod\limits_{k'=0}^{i-n-1} (r-a+k') - \binom{s-l}{m-(i-n)} \prod\limits_{k'=0}^{i-n-1} (s-a+k') \Bigg) \frac{(m-(i-n))!}{m!}\binom{a-l-m}{i-T-m} \\
                &= \frac{1}{m!} \binom{a-l-m}{i-T-m} \left( \prod_{k'=0}^{i-n-1} (r-a+k') \times \prod_{k'=0}^{m-(i-n)-1} (r-l-k') -\prod_{k'=0}^{i-n-1} (s-a+k') \times \prod_{k'=0}^{m-(i-n)-1} (s-l-k') \right) \\
                &\equiv \frac{r-s}{m!} \binom{a-l-m}{i-T-m} \prod_{k'=0}^{i-n-1} (s-a+k') \times \prod_{k'=0}^{m-(i-n)-1} (s-l-k') 
                \left( \sum_{k'=0}^{i-n-1} \frac{1}{s-a+k'} + \sum_{k'=0}^{m-(i-n)-1} \frac{1}{s-l-k'} \right)  \\
                &\equiv \frac{r-s}{m!} \binom{a-l-m}{i-T-m} (-1)^{i-n} (i-n)! \frac{(a-i+n-l)!}{(a-l-m)!} \left( -H_{i-n} + H_{a-i+n-l} -H_{a-l-m} \right) \\
                & \equiv \frac{r-s}{m!} \cdot \frac{(-1)^{i-n} (i-n)!}{(i-T-m)!} \cdot \frac{(a-i+n-l)!}{(a-i+T-l)!} \left( -H_{i-n} + H_{a-i+n-l} -H_{a-l-m} \right) \mod p(r-s), 
                \end{split}
                \end{align}
                where in the penultimate step we used that $s \equiv a-i+n$ mod $p$.  Note that for $k=1,\ldots, i-n-1$ we have $\prod\limits_{k'=0}^{i-n-1} (k-k') =0$. Hence the second term in the middle left block vanishes.
                For the third and fourth terms in the middle left block, note that
                \begin{align}\label{est. error term 2 adhc}
                \begin{split}
                    &\binom{s-l}{s-(a-i+T)} \binom{s-(a-i+T)}{m-(i-n)} \frac{(m-(i-n))!}{m!} \prod\limits_{k'=0}^{i-n-1} (s-a-k'(p-1)) \\
                    & \qquad \qquad \qquad \qquad \qquad \qquad    - \binom{r-l}{r-(a-i+T)}\binom{r-(a-i+T)}{m-(i-n)} \frac{(m-(i-n))!}{m!} \prod\limits_{k'=0}^{i-n-1} (r-a-k'(p-1)) \\
                    &\qquad \qquad  = \frac{1}{m!(a-i+T-l)!} \Bigg( \prod_{k'=0}^{i-n-1} (s-a-k'(p-1)) \times \prod_{k'=l}^{a-i+T+m-(i-n)-1} (s-k')\\
                    & \qquad \qquad \qquad \qquad  \qquad \qquad \qquad \qquad - \prod_{k'=0}^{i-n-1} (r-a-k'(p-1)) \times \prod_{k'=l}^{a-i+T+m-(i-n)-1} (r-k') \Bigg)  \\
                    & \qquad \qquad  \equiv \frac{s-r}{m!(a-i+T-l)!}  \prod_{k'=0}^{i-n-1} (s-a-k'(p-1)) \times \prod_{k'=l}^{a-i+T+m-(i-n)-1} (s-k') \\
                    & \qquad \qquad \qquad \qquad  \qquad \qquad \qquad \qquad  \times \left( \sum_{k'=0}^{i-n-1} \frac{1}{s-a-k'(p-1)} + \sum_{k'=l}^{a-i+T+m-(i-n)-1} \frac{1}{s-k'}\right) \\
                   & \qquad \qquad  \equiv \frac{s-r}{m!(a-i+T-l)!} \cdot (-1)^{i-n}(i-n)! \cdot \frac{(a-i+n-l)!}{(i-T-m)!} \cdot \left( - H_{i-n} + H_{a-i+n-l} -H_{i-T-m}\right) \\
                   & \qquad \qquad  \equiv \frac{s-r}{m!} \cdot \frac{(-1)^{i-n}(i-n)!}{(i-T-m)!} \cdot \frac{(a-i+n-l)!}{(a-i+T-l)!} \cdot \left(  - H_{i-n} + H_{a-i+n-l}-H_{i-T-m} \right) \mod p(r-s).
                \end{split}  
                \end{align}
                Combining all the above, we see that each entry in the middle left block modulo $p(r-s)$ equals
                \begin{align}\label{est. error term  adhc1}
                \frac{r-s}{m!} \cdot \frac{(-1)^{i-n}(i-n)!}{(i+T-m)!} \cdot \frac{(a-i+n-l)!}{(a-i+T-l)!} \cdot \left( H_{i-T-m} - H_{a-l-m} \right).
                \end{align}
                Thus 
                \begin{align*}
            \det(A)= \left\lvert \begin{array}{c|c}
                    \left( \prod\limits_{k'=0}^{m-1} (r-a+k') - \prod\limits_{k'=0}^{m-1} (s-a+k') \right) \frac{1}{m!} \binom{a-l-m}{i-T-m}+\sum\limits_{k=1}^{i-n-1} (p-1)^m \binom{s-l}{i-T+k(p-1)}   \binom{k}{m}  \\
                   + \left(\binom{s-l}{i-T} - \binom{r-l}{i-T} \right) \binom{0}{m} +(p-1)^m \left( \binom{s-l}{s-(a-i+T)} \binom{\frac{s-a}{p-1}}{m} -  \binom{r-l}{r-(a-i+T)} \binom{\frac{r-a}{p-1}}{m} \right) + O(p(r-s)) & p^t \delta_{m,0}  \\ 
                    \hline \\
                     \frac{r-s}{m!} \cdot \frac{(-1)^{i-n}(i-n)!}{(i-T-m)!} \cdot \frac{(a-i+n-l)!}{(a-i+T-l)!} \cdot \left( H_{i-T-m} - H_{a-l-m} \right) + O(p(r-s)) & 0 \\ 
                    \hline \\                   
                    \sum\limits_{k=1}^{i-n-1} \binom{s-l}{i-T+k(p-1)}\binom{i-T+k(p-1)}{i-T}+ \binom{s-l}{s-(a-i+T)} \binom{s-(a-i+T)}{i-T}- \binom{r-l}{r-(a-i+T)}\binom{r-(a-i+T)}{i-T}& 0
                \end{array}
                \right\rvert,
            \end{align*}
            where $m$ ranges through $0, \ldots, i-n-1$ in the upper blocks, $m$ ranges through $i-n, \ldots, i-T-1$ in the middle blocks, $m$ ranges through $i-T$ in the lower blocks. As $p^t\mid (r-s)$, we see that all the entries in the middle blocks are multiples of $p^t$. Expanding the determinant using the last column we obtain
            \begin{align*}
             \det(A)= (-1)^{i-T} p^t \left\lvert \begin{array}{c}  
              \left( \prod\limits_{k'=0}^{m-1} (r-a+k') - \prod\limits_{k=0}^{m-1} (s-a+k) \right) \frac{1}{m!} \binom{a-l-m}{i-T-m} +\sum\limits_{k=1}^{i-n-1}  (p-1)^m  \binom{s-l}{i-T+k(p-1)} \binom{k}{m} \\
                   + (p-1)^m \Big( \binom{s-l}{s-(a-i+T)} \binom{i-n}{m} - \binom{r-l}{r-(a-i+T)} \binom{\frac{r-a}{p-1}}{m} \Big)  + O(p(r-s))   \\ 
                    \hline \\
                     \frac{r-s}{m!} \cdot \frac{(-1)^{i-n}(i-n)!}{(i-T-m)!} \cdot \frac{(a-i+n-l)!}{(a-i+T-l)!} \cdot \left( H_{i-T-m} - H_{a-l-m} \right) + O(p(r-s))   \\ 
                    \hline \\                   
                    \sum\limits_{k=1}^{i-n-1} \binom{s-l}{i-T+k(p-1)}\binom{i-T+k(p-1)}{i-T}+ \binom{s-l}{s-(a-i+T)} \binom{s-(a-i+T)}{i-T}- \binom{r-l}{r-(a-i+T)}\binom{r-(a-i+T)}{i-T}                   
                \end{array}
                \right\rvert ,
            \end{align*}
            where the range of $m$ in the upper block is $1, \ldots, i-n-1$, the range of $m$ in the middle block is $i-n, \ldots, i-T-1$ and the range of $m$ in the lower block $i-T$. Going modulo $O(r-s)$ instead of $O(p(r-s))$ in the upper block we obtain
            \begin{align*}
             \det(A) &= (-1)^{i-T} p^t \left\lvert \begin{array}{c}  
             \sum\limits_{k=1}^{i-n-1}(p-1)^m \binom{s-l}{i-T+k(p-1)}   \binom{k}{m} +O(r-s)  \\ 
                    \hline \\
                     \frac{r-s}{m!} \cdot \frac{(-1)^{i-n}(i-n)!}{(i-T-m)!} \cdot \frac{(a-i+n-l)!}{(a-i+T-l)!} \cdot \left( H_{i-T-m} - H_{a-l-m} \right) + O(p(r-s))  \\ 
                    \hline \\                   
                    \sum\limits_{k=1}^{i-n-1} \binom{s-l}{i-T+k(p-1)}\binom{i-T+k(p-1)}{i-T}+ \binom{s-l}{s-(a-i+T)} \binom{s-(a-i+T)}{i-T}- \binom{r-l}{r-(a-i+T)}\binom{r-(a-i+T)}{i-T}
                \end{array}
                \right\rvert .
            \end{align*}
             Applying the following sequence of row operations:
            \begin{alignat*}{3}
                R_{i-n-2} &\rightarrow  ~R_{i-n-2}-~  \frac{\binom{i-n-1}{i-n-2}}{(p-1)}R_{i-n-1} \\
                 & ~~\vdots  \\
                R_{m} &\rightarrow ~ R_{m}-~ \sum_{k=m+1}^{i-n-1}\frac{\binom{k}{m}}{(p-1)^{k-m}} R_{k} \\
                 & ~~\vdots  \\
                R_{1} &\rightarrow ~ R_{1}- \sum_{k=2}^{i-n-1}\frac{\binom{k}{1}}{(p-1)^{k-1}} R_{k},  
            \end{alignat*}
            we see that in the summation over $k$ only the $k=m$ term survives in the $m^{\mathrm{th}}$-row of the upper block and we get 
            \begin{equation}\label{det A in terms A1,A2}
                \det(A) = (-1)^{i-T}p^t  \left\lvert \begin{array}{c}  
             (p-1)^m \binom{s-l}{i-T+m(p-1)}    + O(r-s) \\ 
                    \hline \\
                     \frac{r-s}{m!} \cdot \frac{(-1)^{i-n}(i-n)!}{(i-T-m)!} \cdot \frac{(a-i+n-l)!}{(a-i+T-l)!} \cdot \left( H_{i-T-m} - H_{a-l-m} \right) + O(p(r-s))  \\ \hline \\
                     \sum\limits_{k=1}^{i-n-1} \binom{s-l}{i-T+k(p-1)}\binom{i-T+k(p-1)}{i-T}+ \binom{s-l}{s-(a-i+T)} \binom{s-(a-i+T)}{i-T}- \binom{r-l}{r-(a-i+T)}\binom{r-(a-i+T)}{i-T}
                \end{array}
                \right\rvert.
            \end{equation}
            
                To remove the sum appearing in the lower block, we apply the following row operation: 
            \begin{align*}
            	R_{i-T} \rightarrow R_{i-T} - \sum_{m=1}^{i-n-1} (p-1)^{-m} \binom{i-T+m(p-1)}{i-T} R_m.
            \end{align*}
            By Lucas' theorem, $\binom{i-T+m(p-1)}{i-T} \equiv 0 \mod p$ for $m=1,\ldots,i-n-1$, so $\binom{i-T+m(p-1)}{i-T} O(r-s) = O(p(r-s))$. Thus we obtain
             \begin{equation*}
            		\det(A) = p^t \left\lvert\begin{array}{c}
            			(p-1)^m \binom{s-l}{i-T+m(p-1)} +O(r-s) \\
            			\hline \\
            			\frac{r-s}{m!} \times \frac{(-1)^{i-n} (i-n)!}{(i-T-m)!} \times  \frac{(a-i+n-l)!}{(a-i+T-l)!} \times (H_{i-T-m} - H_{a-l-m}) + + O(p(r-s))
            			\\ \hline \\
            			 \binom{s-l}{s-(a-i+T)}\binom{s-(a-i+T)}{i-T}-\binom{r-l}{r-(a-i+T)}\binom{r-(a-i+T)}{i-T} + O(p(r-s)) 
            		\end{array} \right\rvert.
            \end{equation*}
           We now simplify the entries in the lower block. As $T < n<i$ and $s \equiv a-i+n$ mod $p$, we have $0 < n-T < i-T$ and $p\mid (s-(a-i+T)-(n-T))$, respectively. Thus  $  \binom{s-(a-i+T)}{i-T} = (s-(a-i+n))\Z_p^\times$. Hence for $l=0,\ldots,a-i+T-1$, we have
            \begin{align*}
            	\binom{s-l}{s-(a-i+T)}&\binom{s-(a-i+T)}{i-T} - \binom{r-l}{r-(a-i+T)}\binom{r-(a-i+T)}{i-T} \\ 
            	& =   \left( \binom{s-l}{s-(a-i+T)}-\binom{r-l}{r-(a-i+T)} \right) \binom{s-(a-i+T)}{i-T} \\
            	& \qquad +  \left( \binom{s-(a-i+T)}{i-T}-\binom{r-(a-i+T)}{i-T} \right) \binom{r-l}{r-(a-i+T)}  \\
            	& =  O(p(r-s)) + \frac{s-r}{s-(a-i+n)}\binom{s-l}{s-(a-i+T)} \binom{s-(a-i+T)}{i-T} ,
            \end{align*}
            where we have used  Lemma~\ref{binomial coefficient under congruences} $(i)$ and $(iii)$ in the last step.
            Substituting this above, we get
            \begin{align*}
                \det(A) = p^t \times \left\lvert\begin{array}{c}
            			(p-1)^m \binom{s-l}{i-T+m(p-1)}  +  O(r-s) \\
            			\hline \\
            			\frac{r-s}{m!} \times \frac{(-1)^{i-n} (i-n)!}{(i-T-m)!} \times  \frac{(a-i+n-l)!}{(a-i+T-l)!} \times (H_{i-T-m} - H_{a-l-m}) +O(p(r-s))
            			\\ \hline \\
            			\frac{s-r}{s-(a-i+n)}\binom{s-l}{s-(a-i+T)} \binom{s-(a-i+T)}{i-T} +O(p(r-s))
            		\end{array}\right\rvert,
            \end{align*}
            where the range of $m$ in the upper block, resp. middle block, resp. lower block is $1,\ldots,i-n-1$, resp. $i-n,\ldots,i-T-1$, resp. $i-T$  and the range of $l$  is $0,\ldots,i-T-1$. We now use Lucas' theorem to simplify the upper block and lower block entries. Note that for $0 \leq l \leq i-T-1$ and $1 \leq m \leq i-n $, by Lucas' theorem we have
            \begin{align*}
                \binom{s-l}{i-T+m(p-1)} &\equiv \binom{i-n}{m} \binom{a-i+n-l}{i-T-m} \\
                &\equiv  \frac{(a-i+n-l)!}{(a-i+T-l)!}\frac{(i-n-m)!}{(i-T-m)!}\binom{i-n}{m}\binom{a-i+T-l}{i-n-m}\mod p.
            \end{align*}
            Also note $s-(a-i+T) = i-T+(i-n)(p-1)$. Thus $\binom{s-l}{i-T+m(p-1)}$ also appears in the the lower block with $m=i-n$. Hence we obtain 
            \begin{align*}
                \det(A) =p^t \times \left\lvert\begin{array}{c}
            			(p-1)^m \frac{(a-i+n-l)!}{(a-i+T-l)!}\frac{(i-n-m)!}{(i-T-m)!}\binom{i-n}{m}\binom{a-i+T-l}{i-n-m}  +  O(p) \\
            			\hline \\
            			\frac{r-s}{m!} \times \frac{(-1)^{i-n} (i-n)!}{(i-T-m)!} \times  \frac{(a-i+n-l)!}{(a-i+T-l)!} \times (H_{i-T-m} - H_{a-l-m}) +O(p(r-s))
            			\\ \hline \\
            			\frac{s-r}{s-(a-i+n)}\frac{(a-i+n-l)!}{(a-i+T-l)!}\frac{1}{(n-T)!}  \binom{s-(a-i+T)}{i-T} + O(p(r-s))
            		\end{array}\right\rvert.
            \end{align*}
            Pulling out $(p-1)^m \frac{(i-n-m)!}{(i-T-m)!}\binom{i-n}{m}$ from every row in the upper block, 
            $\frac{r-s}{m!} \times \frac{(-1)^{i-n} (i-n)!}{(i-T-m)!}$ from every row in the middle block 
            and $\frac{s-r}{s-(a-i+n)}\frac{1}{(n-T)!} \binom{s-(a-i+T)}{i-T}$ from the last row, and $\frac{(a-i+n-l)!}{(a-i+T-l)!}$ from the $l^{\mathrm{th}}$-column we get
            \begin{align}\label{detA interms of detB right bold green}
            \begin{split}
                \det(A) &= p^t \times \prod_{m=1}^{i-n-1} (p-1)^m \frac{(i-n-m)!}{(i-T-m)!}\binom{i-n}{m} \times  \prod_{m=i-n}^{i-T-1} \frac{r-s}{m!} \times \frac{(-1)^{i-n} (i-n)!}{(i-T-m)!} \\
                & \qquad \qquad \times \frac{s-r}{s-(a-i+n)} \frac{1}{(n-T)!}\binom{s-(a-i+T)}{i-T}  \times \prod_{l=0}^{i-T-1} \frac{(a-i+n-l)!}{(a-i+T-l)!} \times \det(B),
            \end{split}    
            \end{align}
            where
            \begin{align*}
                B =  \left[\begin{array}{c}
            			\binom{a-i+T-l}{i-n-m}  +  O(p) \\
            			\hline \\
            			(H_{i-T-m} - H_{a-l-m}) +O(p)
            			\\ \hline \\
            			1  + O(p)
            		\end{array}\right],
            \end{align*}
            with the range of $m$ in the upper block, resp. middle block, resp. lower block is $1,\ldots,i-n-1$, resp. $i-n,\ldots,i-T-1$, resp. $i-T$,  and the range of $l$  is $0,\ldots,i-T-1$.

            We now compute $\det(B)$ modulo $p$. To kill all but the last entry  in the last row, we apply the following column operations 
            \begin{align*}
                    C_0 &\rightarrow\; C_0 - C_{1} \\
                    &~~\vdots \\
                    C_l &\rightarrow\; C_l - C_{l+1}, \\
                    &~~\vdots \\
                    C_{i-T-2} &\rightarrow\; C_{i-T-2} - C_{i-T-1}.
            \end{align*}
            By Pascal's identity, for $0\le m\le i-n-1$ and $0\le l\le i-T-1$ we have
            \[
               \binom{a-i+T-l}{i-n-m}-\binom{a-i+T-l-1}{i-n-m} = \binom{a-i+T-l-1}{i-n-m-1},
            \]
            so we obtain 
            \begin{align*}
	                \det(B) \equiv  \det \left[
	                \begin{array}{c}
		                    \left( \binom{a-i+T-1-l}{i-n-1-m} \right)_{\substack{m=1,\ldots,i-n-1 \\ l=0,\dots,i-T-2}}
                             \\[10pt]\hline\\[-5pt]
		                    \left(-\frac{1}{a-l-m}\right)_{\substack{m=i-n,\ldots,i-T-1 \\ l=0,\dots,i-T-2}}
	                  \end{array}
	                \right] \mod p.
            \end{align*}
            We now apply the column operations  
            \begin{align*}
	                C_{0} & \rightarrow C_{0} + \sum_{l'=1}^{i-T-2} (-1)^{l'} \binom{i-T-2}{l'}C_{l'}, \\
	                  & ~~ \vdots \\
	                C_{l} & \rightarrow C_{l} + \sum_{l'=l+1}^{i-T-2} (-1)^{l'-l} \binom{i-T-2-l}{l'-l}C_{l'} \\
	                  & ~~\vdots \\
	                C_{i-T-3} &\rightarrow C_{i-T-3} - C_{i-T-2}
            \end{align*}
            to simplify the upper block entries. More precisely, we would like to kill all entries except the last in the last row of the upper block, kill all entries except the last two in the second last row of the upper block and so on.
            Observe that
            \begin{align}\label{eq: comb. id 1 green right}
            \begin{split}	
	                \binom{a-i+T-1-l}{i-n-1-m} &+ \sum_{l'=l+1}^{i-T-2} (-1)^{l'-l}\binom{a-i+T-1-l'}{i-n-1-m}\binom{i-T-2-l}{l'-l} \\
	                &= \sum_{l'=0}^{i-T-2-l} (-1)^{l'}\binom{a-i+T-1-l-l'}{i-n-1-m}\binom{i-T-2-l}{l'} \\
	                &= \binom{a-2i+2T+1}{a-2i+m+n+T-l} \\
                    &=\binom{a-2i+2T+1}{l-m-n+T+1},
            \end{split}
            \end{align}	
            where the penultimate step follows from Lemma~\ref{lem: combinatorial id for det(B)} $(v)$ applied with $M=a-i+T-1-l$, $N=i-T-2-l$ and $k=a-2i+m+n+T-l$.  
            Also, taking $\alpha = a-l-m$ and $N=i-T-2-l$ in Lemma~\ref{lem: combinatorial id for det(B)} $(ii)$ we get
            \begin{align}\label{eq: comb. id 2 green right}
            \begin{split}
	             \frac{1}{a-l-m}+ \sum_{l'=l+1}^{i-T-2} (-1)^{l'-l} \binom{i-T-2-l}{l'-l} \frac{1}{a-l'-m} &= \sum_{l'=0}^{i-T-2-l} (-1)^{l'} \binom{i-T-2-l}{l'} \frac{1}{a-l-m-l'}   \\
                 &=(-1)^{i-T-2-l} \frac{(i-T-2-l)!}{(a-i+T+2-m)\cdots(a-l-m)}.
            \end{split}     
            \end{align}
             Using \eqref{eq: comb. id 1 green right} for the upper  block and \eqref{eq: comb. id 2 green right} for the lower  block  we get 
            \begin{align*}
	                  \det(B) 
                    \equiv  \det \left[
	                \begin{array}{c}
		                \left( \binom{a-2i+2T+1}{l-m-n+T+1}\right)_{\substack{m=1,\ldots,i-n-1 \\ l=0,\dots,i-T-2}} \\[10pt]\hline\\[-5pt]
		                \left(\frac{(-1)^{i-T-1-l}(i-T-2-l)!}{(a-i+T+2-m)\cdots(a-l-m)}\right)_{\substack{m=i-n,\ldots,i-T-1 \\ l=0,\dots,i-T-2}} 
	                \end{array} \right] \mod p.
            \end{align*}
            Note that $\binom{a-2i+2T+1}{l-m-n+T+1} = 0$ if $l<m+n-T-1$. Expanding the determinant first using row $m=i-n-1$, second using row $m=i-n-2$ and so on till the first row, we get
            \begin{align*}
	                  \det(B) 
                    &\equiv   (-1)^{(n-T)(i-n-1)}\det 
	                \begin{bmatrix}
		                \frac{(-1)^{i-T-1-l}(i-T-2-l)!}{(a-i+T+2-m)\cdots(a-l-m)}
	                \end{bmatrix}_{\substack{m=i-n,\ldots,i-T-1 \\ l=0,\dots,n-T-1}} \\
                    &\equiv   \prod_{l=0}^{n-T-1} (-1)^{i-T-1-l}\,l!\,(i-T-2-l)!\times (-1)^{(n-T)(i-n-1)}\det 
	                \begin{bmatrix}
		                \frac{(a-i+T+1-m)!}{(a-l-m)!l!}
	                \end{bmatrix}_{\substack{m=i-n,\ldots,i-T-1 \\ l=0,\dots,n-T-1}} \\
                    & \equiv   \prod_{l=0}^{n-T-1} (-1)^{n-T-l}\,l!\,(i-T-2-l)!\times \det 
	                \begin{bmatrix}
		                \frac{(a-i+T+1-m)!}{(a-l-m)!\,l!}
	                \end{bmatrix}_{\substack{m=i-n,\ldots,i-T-1 \\ l=0,\dots,n-T-1}} \mod p
            \end{align*}
            where in the second last step we pulled out $(-1)^{i-T-1-l}\,l!\,(i-T-2-l)!$ from the $l^{\mathrm{th}}$-column. Pulling out $\frac{(a-i+T+1-m)!}{(a-m)!}$ from the $m^{\mathrm{th}}$-row by  for $m=i-n,\ldots,i-T-1$,  we  obtain
            \begin{align*}
	                  \det(B)            
                    &\equiv   \prod_{l=0}^{n-T-1} (-1)^{n-T-l}\,l!\,(i-T-2-l)! \times \prod_{m=i-n}^{i-T-1}\frac{(a-i+T+1-m)!}{(a-m)!} \\
                    & \qquad \qquad \times \det 
	                \begin{bmatrix}
		                \binom{a-m}{l}
	                \end{bmatrix}_{\substack{m=i-n,\ldots,i-T-1 \\ l=0,\dots,n-T-1}} \mod p.
            \end{align*}
            Observe that
            \[
             \det \begin{bmatrix}\binom{a-m}{l}\end{bmatrix}_{\substack{m=i-n,\ldots,i-T-1 \\ l=0,\dots,n-T-1}} = \det\begin{bmatrix}\binom{a-i+n-m}{l}\end{bmatrix}_{\substack{m=0,\ldots,n-T-1 \\ l=0,\dots,n-T-1}} =  \pm 1,
            \]
            where we have used Lemma~\ref{cor: GV det} in the last step. Substituting this in the above, we get  $\det(B) \in \Z_p^\times$. Hence it follows from \eqref{detA interms of detB right bold green} that $\det(A) \in p^t(r-s)^{n-T+1} \Z_p^\times$. Thus by Cramer's rule, the linear equations \eqref{eq: linear eqns adhc green} has a solution in $\Q_p$.
        
            By Cramer's rule, we have $\beta_l = \pm p^t \det(A_{i-T,l})/\det(A)$ and $\gamma =  p^t \det(A_{i-T,i-T})/\det(A)$. To show $\beta_0,\ldots,\beta_{i-T-1}$ and $\gamma$ belong to $\Z_p$ it is enough to show $\det(A_{i-T,l}) \in p^{t(n-T+1)}\Z_p$. This follows from Lemma~\ref{minor trick 2}. This proves $(i)$ and $(ii)$.

             We now prove $(iii)$. This is similar to the argument given at the end of Lemma~\ref{choice beta hybrid green}. As $\beta_l \in \Z_p$, by Corollary~\ref{cor: binomial sums under congruences 2}, for all $m\geq 0$ we have
            \begin{align}\label{eq: binomial sum without end points green}
                \sum_{l=0}^{i-T-1} \beta_l \sum\limits_{ \substack{i-T <  j < r-(a-i+T) \\  j \equiv i-T ~\mathrm{mod}~(p-1)} }   \binom{r-l}{j} \binom{j}{m} \equiv \sum_{l=0}^{i-T-1} \beta_l \sum\limits_{ \substack{i-T <  j < s-(a-i+T) \\  j \equiv i-T ~\mathrm{mod}~(p-1)} }   \binom{s-l}{j} \binom{j}{m} \mod p^{t-v(m!)} .
            \end{align}
            Here we have omitted the end points as they are congruent modulo $p^t$. Since many terms in the matrix $A$ vanish modulo $p^t$ in \eqref{A matrix green right}, it follows from \eqref{eq: linear eqns adhc green} that
            \begin{align}\label{eq: (i) and (ii) mod p green}
            \begin{split}
               \sum_{k=1}^{i-n-1} \left( \sum_{l=0}^{i-T-1} \beta_l \binom{s-l}{i-T+k(p-1)}\right) \binom{i-T+k(p-1)}{m} 
               \equiv 0 \mod p^t
            \end{split}   
            \end{align}
            for $m=0,\ldots, i-T$. By Lemma~\ref{dmm' trick} (applied with $N=i-n-2$ $(\leq i-T)$ and $c=i-T+(p-1)$), the congruence \eqref{eq: (i) and (ii) mod p green} holds for all $m \geq 0$. Substituting this in \eqref{eq: binomial sum without end points green}, we obtain $(iii)$.
        \end{proof}  

        We now show that the conclusion of the previous lemma holds even if $T=n$. 
        \begin{lemma}\label{choice beta hybrid ad last column} 
            Let  $r \equiv a \mod{(p-1)}$ with $ 1 \leq a \leq p-1 $ and $r \geq i(p+1)+p $ with $ v(a_{p})  \in (i,i+1) $. Let $ s =  a-i+n+(i-n)p$ and $v(r-s) =t $  and $ 1 \leq t < T =n< i < a$. If  $2i-a \leq 2T$, then there exist $\beta_0, \ldots, \beta_{i-n-1} \in \Z_p$ and $\gamma \in \Z_p$ such that 
        	    \begin{enumerate}
        	    	\item[$(i)$] $\sum\limits_{l=0}^{i-T-1} \beta_l  \sum\limits_{ \substack{i-T <  j < r-(a-i+T) \\  j \equiv i-T ~\mathrm{mod}~ (p-1)} }   \binom{r-l}{j} \binom{j}{m}  + \gamma p^t  \binom {i-T}{m}\equiv 
        	    	0  \mod p^{t+1} ~ \mathrm{for}~ m=0,\ldots, i-T-1$. 
        	    	\item[$(ii)$] $\sum\limits_{l=0}^{i-T-1} \beta_l  \sum\limits_{ \substack{i-T <  j < r-(a-i+T) \\  j \equiv i-T ~\mathrm{mod}~ (p-1)} }   \binom{r-l}{j} \binom{j}{i-T}   \equiv 
        	    	 p^t  \mod p^{t+1}$. 
        	    	\item[$(iii)$] $\sum\limits_{l=0}^{i-T-1} \beta_l \sum\limits_{ \substack{i-T <  j < r-(a-i+T) \\  j \equiv i-T ~\mathrm{mod}~ (p-1)} }   \binom{r-l}{j} \binom{j}{m} \equiv  0 \mod p^{t-v(m!)} ~  \mathrm{for}~ m=i-T+1,\ldots, i+t$.
        	    \end{enumerate}
        \end{lemma}
        \begin{proof}
        	     First we  prove  $(i)$ and $(ii)$. As in Lemma~\ref{choice beta hybrid ad1}, using Corollary~\ref{cor: binomial sums under congruences 1}, for $0 \leq m \leq i-T$ and $0\leq l \leq i-T-1$,  we have
        	    \begin{align*}
        	    	\sum\limits_{\substack{i-T <  j < r-(a-i+T) \\  j \equiv i-T ~\mathrm{mod}~ (p-1)} }   \binom{r-l}{j} \binom{j}{m} &\equiv \left( \binom{r-l}{m} - \binom{s-l}{m} \right) \binom{a-l-m}{i-T-m} \\ &
        	    	\qquad \qquad +\sum_{\substack{ i-T < j < s-(a-i+T) \\ j \equiv i-T ~\mathrm{mod}~(p-1)}} \binom{s-l}{j} \binom{j}{m}  
        	    	+ \left(\binom{s-l}{i-T} - \binom{r-l}{i-T}\right) \binom{i-T}{m}  \\
        	    	&\qquad \qquad + \Bigg(\binom{s-l}{s-(a-i+T)} \binom{s-(a-i+T)}{m}  \\ & \qquad \qquad \qquad \qquad - \binom{r-l}{r-(a-i+T)} \binom{r-(a-i+T)}{m}\Bigg)   \mod p^{t+1}. 
        	    \end{align*}
        	     We now begin solving for $\beta_l$ and $\gamma$ satisfying the congruences $(i)$ and $(ii)$.  Let
               \begin{align*}
               A= \left(\begin{array}{@{}c|c@{}}
               \begin{matrix}
                \left( \binom{r-l}{m} - \binom{s-l}{m}\right)\binom{a-l-m}{i-T-m} +   \sum\limits_{k=1}^{i-n-1} \binom{s-l}{i-T+k(p-1)} \binom{i-T+k(p-1)}{m}   \\ 
        	    	+ \left(\binom{s-l}{i-T} - \binom{r-l}{i-T}\right) \binom{i-T}{m}   \\
        	    	+ \binom{s-l}{s-(a-i+T)} \binom{s-(a-i+T)}{m}   - \binom{r-l}{r-(a-i+T)} \binom{r-(a-i+T)}{m}  
                \end{matrix}    & p^t\binom{i-T}{m} -p^t \delta_{i-T,m} 
               \end{array}
               \right)_{m=0,\ldots, i-T},
               \end{align*}
               where $l$ ranges from $0$ to $i-T-1$. By \eqref{binomial sum  hybrid ad}, it is enough to solve for $\beta_0, \ldots, \beta_{i-T-1}$ and $\gamma$ in $\Z_p$ satisfying
               \begin{align}\label{eq: linear eqns adhc green last column}
                   A (\beta_0, \ldots, \beta_{i-T-1}, \gamma)^{\text{tr}} = (0,\ldots, 0,p^t)^{\text{tr}}.
                \end{align}

                If $T=n=i-1$, then using the above expression, we have
                \begin{align*}
                    A & = \begin{pmatrix}
                        s-r+ \binom{s}{a-i+n} -\binom{r}{a-i+n} & p^t \\ \binom{s}{a-i+n}(s-(a-i+n))   - \binom{r}{a-i+n} (r-(a-i+n)) & 0
                    \end{pmatrix} \\
                    & \equiv \begin{pmatrix}
                        s-r+ \binom{s}{a-i+n} -\binom{r}{a-i+n} & p^t \\ (s-r)\binom{s}{a-i+n}   & 0
                    \end{pmatrix} \mod p^{t+1}.
                \end{align*}
                By Lucas' theorem $\binom{s}{a-i+n} \equiv 1 \mod p$. Take $\beta_0 \in \Z_p^\times$ such that $\beta_0 (s-r)\binom{s}{a-i+n} = p^t $. Take $\gamma = -\left((s-r)+ \binom{s}{a-i+n} -\binom{r}{a-i+n}\right) \beta_0/p^t \in \Z_p$.  This shows that $(i)$ and $(ii)$ are solvable if $n=i-1$.
                
                Assume $1 \leq n \leq i-2$. We will use Cramer's rule. We first show that the determinant of $A$ is non-zero. To show that $\beta_l, \gamma \in \Z_p$, we also need to determine the power of $p$ dividing $\det(A)$. To this end we show $p^{2t} \parallel \det(A)$. To achieve this we will perform  row operations. Apply the following row operations to $A$:
            \begin{alignat*}{4}
                  R_{i-T-1} &\rightarrow R_{i-T-1}~- ~\frac{i-T-(i-T-2)}{i-T-1} R_{i-T-2} \\
                  & ~~\vdots  \\
                R_{m} &\rightarrow ~ ~~R_{m}~ - ~\frac{i-T-(m-1)}{m} R_{m-1} \\
                 & ~~\vdots  \\
                R_{1} &\rightarrow ~~ R_{1}~- ~\frac{i-T}{1} R_{0}.
            \end{alignat*}   
            Using \eqref{identity 1 adhc} with $j$  equal to $i-T$ and $j'= i-T+k(p-1)$ for $k=1,\ldots,i-n-1$  we get 
            \begin{align*}
                \det(A)= \left\lvert\begin{array}{c|c}
                    \sum\limits_{k=1}^{i-n-1} \binom{s-l}{i-T+k(p-1)}  + O(r-s)  & p^t  \\ 
                    \hline \\                     
                     \sum\limits_{k=1}^{i-n-1} \binom{s-l}{i-T+k(p-1)} \binom{i-T+k(p-1)}{m-1} \frac{k(p-1)}{m}  +O(r-s) & 0 \\
                    \hline \\
                   \sum\limits_{k=1}^{i-n-1} \binom{s-l}{i-T+k(p-1)}\binom{i-T+k(p-1)}{i-T}+ \binom{s-l}{s-(a-i+T)} \binom{s-(a-i+T)}{i-T}- \binom{r-l}{r-(a-i+T)}\binom{r-(a-i+T)}{i-T} & 0
                \end{array}
                \right\rvert
            \end{align*}
            where the range of $m$ in the upper blocks is $0$, the range of $m$ in the middle blocks is $1,\ldots,i-T-1$ and the range of $m$ in the lower blocks is $i-T$. Expanding the determinant using the last column, we obtain
            \begin{align*}
                \det(A)= (-1)^{i-T} p^t \left\lvert\begin{array}{c}                    
                     \sum\limits_{k=1}^{i-n-1} \binom{s-l}{i-T+k(p-1)} \binom{i-T+k(p-1)}{m-1} \frac{k(p-1)}{m}  +O(r-s) \\
                    \hline \\
                   \sum\limits_{k=1}^{i-n-1} \binom{s-l}{i-T+k(p-1)}\binom{i-T+k(p-1)}{i-T}+ \binom{s-l}{s-(a-i+T)} \binom{s-(a-i+T)}{i-T}- \binom{r-l}{r-(a-i+T)}\binom{r-(a-i+T)}{i-T}
                \end{array}
                \right\rvert
            \end{align*}
            where the range of $m$ in the upper blocks is $1,\ldots,i-T-1$ and the range of $m$ in the lower blocks is $i-T$.  As $p \mid (s-(a-i+T)), (r-(a-i+T))$ and $1 \leq i-T < p-1$, by Lemma~\ref{binomial coefficient under congruences} $(i)$ and $(iii)$, we have 
            \begin{align*}
                \binom{s-l}{s-(a-i+T)} \binom{s-(a-i+T)}{i-T}&- \binom{r-l}{r-(a-i+T)}\binom{r-(a-i+T)}{i-T} \\
                &\equiv \binom{s-l}{s-(a-i+T)} \left(\binom{s-(a-i+T)}{i-T}- \binom{r-(a-i+T)}{i-T}\right) \\
                &\equiv \frac{s-r}{s-(a-i+n)}\binom{s-l}{s-(a-i+T)}\binom{s-(a-i+T)}{i-T}  \mod p^{t+1}.
            \end{align*}
            Hence  
            \begin{align*}
                \det(A)= (-1)^{i-T} p^t \left\lvert\begin{array}{c}                    
                     \sum\limits_{k=1}^{i-n-1} \binom{s-l}{i-T+k(p-1)} \binom{i-T+k(p-1)}{m-1} \frac{k(p-1)}{m}  +O(r-s) \\
                    \hline \\
                   \sum\limits_{k=1}^{i-n-1} \binom{s-l}{i-T+k(p-1)}\binom{i-T+k(p-1)}{i-T}+ \frac{s-r}{s-(a-i+n)} \binom{s-l}{s-(a-i+T)}\binom{s-(a-i+T)}{i-T} + O(p(r-s))
                \end{array}
                \right\rvert
            \end{align*}
            By Lucas theorem, we have $\binom{i-T+k(p-1)}{i-T} \equiv 0 \mod p$ for $k=1,\ldots,i-n-1$. Thus every entry in the last row is a multiple of $p$. Hence 
            \begin{align}\label{A=BC last column}
            \begin{split}
                \det(A)& \equiv (-1)^{i-T} p^t \left\lvert\begin{array}{c}                    
                     \sum\limits_{k=1}^{i-n-1} \binom{s-l}{i-T+k(p-1)} \binom{i-T+k(p-1)}{m-1} \frac{k(p-1)}{m}  \\
                    \hline \\
                   \sum\limits_{k=1}^{i-n-1} \binom{s-l}{i-T+k(p-1)}\binom{i-T+k(p-1)}{i-T}+ \frac{s-r}{s-(a-i+n)} \binom{s-l}{s-(a-i+T)}\binom{s-(a-i+T)}{i-T}
                \end{array}
                \right\rvert \mod p^{2t+1} \Z_p \\
                & = (-1)^{i-T} p^t \det(BC) \mod p^{2t+1} \Z_p,
            \end{split}    
            \end{align}
            where
            \begin{align*}
                B &= \left[\begin{array}{c|c} \left(\binom{i-T+k(p-1)}{m-1} \frac{k(p-1)}{m} \right)_{\substack{m =1,\ldots,i-T-1 \\ k=1,\ldots,i-n-1}} &  0 \\ \hline \\ \left(\binom{i-T+k(p-1)}{i-T}\right)_{k=1,\ldots,i-n-1}  & \frac{s-r}{s-(a-i+n)} \binom{s-(a-i+T)}{i-T} \end{array}\right] \\ C &= \left[ \binom{s-l}{i-T+k(p-1)} \right]_{\substack{k =1,\ldots,i-n \\ l=0,\ldots,i-T-1}} .
            \end{align*}
            
            We now compute $\det(B)$. Note that 
            \begin{align*}
                \det(B) &= \frac{s-r}{s-(a-i+n)} \binom{s-(a-i+T)}{i-T} \det\left[\binom{i-T+k(p-1)}{m-1} \frac{k(p-1)}{m}\right]_{\substack{m=1,\ldots,i-T-1 \\ k=1,\ldots,i-n-1}} \\
                &= \frac{s-r}{s-(a-i+n)} \binom{s-(a-i+T)}{i-T}\times (p-1)^{i-n-1} \times\det\left[ \binom{i-T+k(p-1)}{m-1}\right]_{\substack{m =1,\ldots,i-T-1\\ k=1,\ldots,i-n-1}} ,
            \end{align*}
            where in the last step we pulled out $(p-1)/m$ from every row and $k$ from every column. By Corollary~\ref{cor: GV det} $(i)$, we have $\det(B) \in \frac{s-r}{s-(a-i+n)} \binom{s-(a-i+T)}{i-T}\Z_p^\times$. Since $T=n$ and $\binom{s-(a-i+n)}{i-n} = \frac{s-(a-i+n)}{i-n} \binom{s-(a-i+n)-1}{i-n-1}$, we obtain $\det(B) \in p^t \Z_p^\times$ by Lucas' theorem. 
            
            We now compute $\det(C)$. By Lucas' theorem, we have $\binom{s-l}{i-T+k(p-1)} \equiv \binom{i-n}{k} \binom{a-i+n-l}{i-T-k} \mod p$ for $k=1,\ldots,i-n$ and $l=0,\ldots,i-T-1$. Thus 
            \begin{align*}
                \det(C) &\equiv \det\left[\binom{i-n}{k} \binom{a-i+n-l}{i-T-k}\right]_{\substack{k =1,\ldots,i-n \\ l=0,\ldots,i-T-1}} \\
                &= \prod_{k=1}^{i-n} \binom{i-n}{k} \times \det\left[ \binom{a-i+n-l}{i-T-k}\right]_{\substack{k =1,\ldots,i-n \\ l=0,\ldots,i-T-1}} \\
                & = \prod_{k=1}^{i-n} \binom{i-n}{k} \times \det\left[ \binom{a-2i+n+T+1+l}{k}\right]_{\substack{k =0,\ldots,i-n-1 \\ l=0,\ldots,i-T-1}} \mod p,
            \end{align*}
            where in the penultimate step we pulled out $\binom{i-n}{k}$ from every row, and in the last step we reversed the order of rows and columns. Now by Corollary~\ref{cor: GV det} $(i)$, we have $\det(C) \in \Z_p^\times$.  Thus it follows from $\det(B) \in p^t\Z_p^\times$ and \eqref{A=BC last column} that $p^{2t} \parallel \det(A)$. Hence \eqref{eq: linear eqns adhc green last column} has a solution in $\Q_p$. As explained in the proof of Lemma~\ref{choice beta hybrid ad1},  $\beta_0,\ldots,\beta_{i-n-1}$ and $\gamma$ belong to $\Z_p$. Finally, $(iii)$ follows from $(i)$ and $(ii)$ by an argument similar to the one used at the end of the proof of Lemma~\ref{choice beta hybrid ad1}.
        \end{proof}          
        \begin{theorem}\label{thm: hybrid ad}
        	Let  $r \equiv a \mod{(p-1)}$ with $ 1 \leq a \leq p-1 $ and $r \geq i(p+1)+p $ with $ v(a_{p})  \in (i,i+1) $. Let $ s =  a-i+n+(i-n)p$ and $v(r-s) =t $  and $ 1 \leq t \leq T$  $\leq  n < i < a$. If  $T+t \leq 2i-a \leq 2T$, then the image of $ \mathrm{ind}_{KZ}^{G}(V_{r}^{(i-T)}) $ is the same as the image of $ \mathrm{ind}_{KZ}^{G}(V_{r}^{(i-T+1)}) $ in $ \bar{\Theta}_{k,a_{p}} $.
        \end{theorem}
        \begin{proof}
        	Let $\beta_l$ and $\gamma$ be the $p$-adic integers  chosen in Lemma~\ref{choice beta hybrid ad1} (if $T<n$) and Lemma~\ref{choice beta hybrid ad last column} (if $T=n$). Then by Lemma~\ref{lem:choice of beta}, there exist $\alpha_{j} \in \mathbb{Z}_{p}$ satisfying
        	\begin{enumerate}
        		\item[$(1)$]  $\alpha_j \equiv \sum\limits_{l=0}^{i-T-1} \beta_{l}  \binom{r-l}{j} ~\mathrm{mod}~p^{t} $, for all $ i-T< j < r-(a-i+T)$  with $j \equiv i-T $ mod $(p-1)$
        		\item[$(2)$]   $\sum\limits_{\substack{ i-T < j < r-(a-i+T) \\ j \equiv i-T ~\mathrm{mod}~(p-1)}}^{} \alpha_j \binom{j}{m} \equiv 0 $ mod $p^{i+t+1-m}$ for $m=0,\ldots, \min\{i+t, p-1\}$.
        	\end{enumerate}
        	From the congruence condition $(1)$ and Lemma~\ref{choice beta hybrid ad1} $(iii)$ (if $T<n$) and Lemma~\ref{choice beta hybrid ad last column} $(iii)$ (if $T=n$), we also have 
        	\begin{enumerate}
        		\item[$(2')$] $\sum\limits_{\substack{ i-T < j < r-(a-i+T) \\ j \equiv i-T ~\mathrm{mod}~(p-1)}}^{} \alpha_j \binom{j}{m} \equiv 0 $ mod $p^{t-v(m!)}$ for $m=p,\ldots, i+t$. 
        	\end{enumerate}
                 Note that for $m=p,\ldots, i+t$ we have $t-v(m!) = t-1 \geq i+t+1-p \geq i+t+1-m$  as $t < i \leq p-2$.
        	Let
        	\begin{align*}
        		f_2 &=  
        		\sum_{  \lambda \in \mathbb{F}_{p}^{\times}} 
        		\Bigg[ g_{2,p[\lambda]}^{0}, \sum_{l=0}^{i-T-1} \frac{[\lambda]^{l-(a-i+T)}}{p^{l+t}} 
        		\beta_{l} (-\theta)^{l+t+1} X^{-t-1}Y^{r-(l+t+1)(p+1)+t+1} \Bigg]  \\
        		& \quad + \left[ g_{2,0}^{0}, \frac{1-p}{p^{a-i+T+t}}  \sum_{l=0}^{i-T-1} \beta_{l}  \binom{r-l}{r-(a-i+T)}  (-\theta)^{a-i+T+t+1} X^{-t-1}Y^{r-(a-i+T+t+1)(p+1)+t+1} \right]  \\
        		f_1 &=  \left[g_{1,0}^0, \frac{p-1}{p^t a_p} \sum_{ \substack{i-T <  j < r-(a-i+T) \\  j \equiv i-T ~\mathrm{mod}~ (p-1)} } \alpha_{j} X^{r-j} Y^{j} \right] \\
        		f_{0} &= \left[ \mathrm{id}, \frac{1-p}{p^{i-T+t}}  \left(\sum_{l=0}^{i-T-1} \beta_{l} \binom{r-l}{i-T} \right)  \theta^{i-T+t+1} X^{r-(i-T+t+1)(p+1)+t+1} Y^{-t-1}\right].
        	\end{align*}
        	From Lemma~\ref{theta and T plus} and the condition $a-i+T+t+1 \leq i+1$, it follows that $T^+ f_2$ vanishes modulo $p$. It easy to see that $-a_p f_2$,  $-a_p f_0 $ and $T^{-} f_0$ all vanish modulo $p$ using $t\leq T$ and $T+t \leq 2i-a$. From $(2)$ and (the discussion below) $(2')$ and the condition $T+t \leq 2i-a \leq a-2 \leq p-3$, it follows that $T^{+}f_1$ vanishes modulo $p$. Using $t \leq T $ and $2i-a \leq p-3$ one checks that $T^{-}f_1$ also vanishes modulo $p$. It can be checked that 
        	\begin{align*}
        		T^{-} f_2 -a_p f_1 + T^+ f_0 \equiv   \left[g_{1,0}^0, \frac{p-1}{p^t} \sum_{ \substack{i-T <  j < r-(a-i+T) \\  j \equiv i-T ~\mathrm{mod}~ (p-1)} } \left( \sum_{l=0}^{i-T-1}  \beta_{l} \binom{r-l}{j}- \alpha_{j} \right) X^{r-j} Y^{j} \right] ~\mathrm{mod}~ p.
        	\end{align*}
        	Let
        	\begin{align*}
        		F(X,Y) = \frac{p-1}{p^t} \sum_{ \substack{i-T <  j < r-(a-i+T) \\  j \equiv i-T ~\mathrm{mod}~ (p-1)} } \left( \sum_{l=0}^{i-T-1}  \beta_{l} \binom{r-l}{j}- \alpha_{j} \right) X^{r-j} Y^{j} + (p-1)\gamma X^{r-i+T}Y^{i-T}.
        	\end{align*}
        	By $(1)$ above, we have $F(X,Y) \in \mathbb{Z}_p[X,Y]$.  To prove the theorem, it is enough to show that $\overline{F(X,Y)}$ generates $V_{r}^{(i-T)}/V_r^{(i-T+1)}$. Using the properties of $\beta_l$ and $\gamma$ from Lemma~\ref{choice beta hybrid ad1} (if $T<n$) and Lemma~\ref{choice beta hybrid ad last column}  (if $T=n$), and the choice of $\alpha_{j}$, for $m=0,1,\ldots, i-T$, we have 		     
        	\begin{align}\label{derivative condition left green}
        		\sum_{ \substack{i-T <  j < r-(a-i+T) \\  j \equiv i-T ~\mathrm{mod}~ (p-1)} } &\binom{j}{m} \left( \sum_{l=0}^{i-T-1}  \beta_{l}\binom{r-l}{j}- \alpha_{j} \right) + p^t \gamma \left( \binom{i-T}{m} - \delta_{i-T,m}\right)\equiv \delta_{i-T,m}p^t ~\mathrm{mod}~p^{t+1}.
        	\end{align}		     
        	Note that  the coefficients of $X^{r}, \ldots, X^{r-(i-T-1)}Y^{i-T-1}$ in $F(X,Y)$ are zero. Since $i-T< a-i+T+ p-1 $, it follows that the coefficients of $X^{i-T}Y^{r-i+T}, \ldots, Y^r$ in $F(X,Y)$ are zero. By \eqref{derivative condition left green} and \cite[Lemma 2.8]{GR19}, we have $\theta^{i-T} \mid \overline{F(X,Y)}$. Applying \eqref{derivative condition left green} and \cite[Lemma 2.12]{GR19} with $m,l$ there equal $i-T$,  we obtain 
        	\begin{align*}
        		\overline{F(X,Y)} & \equiv (p-1) \theta^{i-T} X^{r-(i-T)(p+1)-(p-1)}Y^{p-1} + (p-1) \gamma \theta^{i-T} X^{r-(i-T)(p+1)} ~\mathrm{mod}~V_r^{(i-T+1)}. 
        	\end{align*}  
        	Applying Lemma~\ref{generating polynomial quotient} with $m$ there equal to $i-T$, it follows that $\overline{F(X,Y)} $ generates $V_r^{(i-T)}/V_r^{(i-T+1)}$. This finishes the proof of the theorem.  
        \end{proof}

        Let the notation be as in the previous theorem. We now consider the case $T+t \geq 2i-a+2$.

       \begin{lemma}\label{lem: choice beta asd a>2i}
             Let $1 \leq i \leq p-2$. Let $r \geq i(p+1)+p $, $r \equiv a \mod{(p-1)}$ with $ 1 \leq a \leq p-1 $ and $r \equiv a-i+n \mod p$ with $1 \leq n < i$. Let $ s =  a-i+n+(i-n)p$ and $v(r-s) =t $  with $t\geq 1$. Fix an integer $1  \leq T \leq n$. Assume $a>2i-2T$.  There exists $p$-adic integers $\beta_0, \ldots, \beta_{i-T} \in \mathbb{Z}_p$ satisfying 	
	        \begin{enumerate}
		        \item[$(i)$] For $m=0,\ldots,i-T$, we have
		        $$
		           \sum\limits_{\substack{a - i + T \leq j < r-i+T \\ j \equiv a-i+T \mod(p-1)}} \left( p \sum\limits_{l = 0}^{i-T-1} \beta_l \binom{r-l}{j} + \beta_{i-T} \binom{r-i+T}{j} \right) \binom{j}{m}  \equiv p^{t+1} \delta_{i-T,m} \mod p^{t+2}.
		        $$
		        \item[$(ii)$] For $m=i-T+1,\ldots,i+t+2$ and $m \not\equiv a-i+T \mod (p-1)$, we have
		          $$
		               \sum\limits_{\substack{a - i + T \leq j < r-i + T \\ j \equiv a - i + T \mod(p-1)}} \left( p \sum\limits_{l = 0}^{i-T-1} \beta_l \binom{r - l}{j} + \beta_{i-T} \binom{r - i + T}{j} \right) \binom{j}{m}  \equiv 0 \mod p^{t+1 - v(m!)}.
                $$ 
                \item[$(iii)$] For $m \equiv a-i+T \mod (p-1)$ , we have 
                $$
		               \sum\limits_{\substack{m < j < r-i + T \\ j \equiv a - i + T \mod(p-1)}} \left( p \sum\limits_{l = 0}^{i-T-1} \beta_l \binom{r - l}{j} + \beta_{i-T} \binom{r - i + T}{j} \right) \binom{j}{m}  \equiv 0 \mod p^{t+1-v(m!)}.
                $$ 
		        \item[$(iv)$] $ p \sum\limits_{l = 0}^{i-T-1} \beta_l \binom{r-l}{a-i + T} + \beta_{i-T} \binom{r-i + T}{a-i+T} \equiv  0 \mod p^{t}$.
		
		        \item[$(v)$] $p \sum\limits_{l = 0}^{i-T-1} \beta_l \binom{r - l}{r - i + T - (p-1)} + \beta_{i-T} \binom{r - i + T}{r - i + T - (p-1)} \equiv 0 ~\mathrm{mod}~ p^{t-1}$.
	        \end{enumerate}
        \end{lemma}
        \begin{proof}
            We will first show $(i)$ holds. Later we will show $(i)$ implies $(ii), (iii),(iv)$ and $(v)$. To solve the congruences $(i)$, we now compute the coefficients of $\beta_l$  modulo $p^{t+2}$. By Lemma~\ref{cor: binomial sums under congruences 1}, for $0 \leq l,m \leq i-T$ we have
            \begin{align*}
	            \sum_{\substack{0 \leq j < r - i + T \\ j \equiv a - i + T ~\mathrm{mod}~ (p-1)}}\binom{r- l}{j} \binom{j}{m}\equiv & \left\{\binom{r - l}{m} - \binom{s - l}{m}\right\}\binom{a - l - m}{a - i + T - m} \\
	            & \quad + \sum_{\substack{0 \le j < s - i + T \\ j \equiv a - i + T ~\mathrm{mod}~(p-1)}} \binom{s - l}{j} \binom{j}{m}\\
                & \quad + \binom{s - l}{s - i + T} \binom{s - i + T}{m}- \binom{r - l}{r - i + T} \binom{r - i + T}{m}  \mod p^{t+1},
	        \end{align*}
            where we have used $a-i+T > i-T$.  
            
        If $l=i-T$ and $m=0,\ldots,i-T$ , then by Lemma~\ref{lem: binomial sum Doc. Math general}, we even have 
            \begin{align*}
                \sum_{\substack{0 < j < r - i + T \\ j \equiv a - i + T \mod (p-1)}} \binom{r - i + T}{j} \binom{j}{m} \equiv & p \left\{\binom{r - i+T}{m} - \binom{s - i+T}{m}\right\} \frac{a - s}{a - i + T - m}\\
	            & \quad +\sum_{\substack{0 < j < s - i + T \\ j \equiv a - i + T \mod (p-1)}} \binom{s-i+ T}{j} \binom{j}{m} \\
	            & \quad + p \frac{s - r}{a - i + T - m} \binom{s - i + T}{m} \mod p^{t+2}.
	        \end{align*}
            Let 
            \begin{equation}\label{eq: A sdad}
	            A = \left[ \begin{array}{c|c}
                    p \left\{\binom{r - l}{m} - \binom{s - l}{m}\right\} \binom{a - l - m}{a - i + T - m} &  p \left\{ \binom{r - i+T}{m} - \binom{s - i+T}{m} \right\} \frac{a - s}{a - i + T - m} \\ \\
		            + p \sum\limits_{\substack{0 \leq j< s - i + T \\ j \equiv a - i + T ~\mathrm{mod}~{p-1}}} \binom{s - l}{j} \binom{j}{m} &  +\sum\limits_{\substack{0 < j < s - i + T \\ j \equiv a - i + T ~\mathrm{mod}~ (p-1)}} \binom{s - i + T}{j} \binom{j}{m}\\ \\
		            + p \binom{s - l}{s - i + T} \binom{s - i + T}{m}- p \binom{r - l}{r - i + T} \binom{r - i + T}{m} & + p \frac{s - r}{a - i + T - m} \binom{s - i + T}{m} 
	              \end{array}\right]_{m=0,\ldots,i-T},
            \end{equation}
            where the range of $l$ in the left block is $0,\ldots,i-T-1$ and the range of $l$ in the right most column is $i-T$.

            To solve congruences $(i)$, it is enough to show that the following congruence
            \[
                A \begin{bmatrix} \beta_0 \\ \vdots \\ \beta_{i-T-1} \\ \beta_{i-T} \end{bmatrix} \equiv  \begin{bmatrix} 0\\ \vdots \\ 0 \\ p^{t+1} \end{bmatrix} \mod p^{t+2} \mathbb{Z}_{p}
            \] 
            has a solution in $\mathbb{Z}_{p}$. 

            To solve the above congruence, it is enough to show the following equation
            \begin{equation}\label{eq: linear eq. in matrix form a>2i}
                A \begin{bmatrix} \beta_0 \\ \vdots \\ \beta_{i-T-1} \\ \beta_{i-T} \end{bmatrix} = \begin{bmatrix} 0\\ \vdots \\ 0 \\ p^{t+1} \end{bmatrix}
            \end{equation}
            has a solution over $\mathbb{Q}$ with $\beta_l \in \mathbb{Z}_{p}$. To show this, we use Cramer's rule.

            First, we will  perform a sequence of row operations on $A$ so that the bottom $(n-T+1)$ rows are multiples of $p(r-s)$ and every entry in the top $(i-n-1)$ rows is a binomial coefficient up to an error term.

            Recall that for a non-negative integer $k $, $(n)_k:= n \cdots (n-k+1)$ denotes the falling factorial. Let $A(0)=A$. We define matrices $A(0), \ldots, A(i-n)$ recursively, where  $A(k'+1)$ is obtained from $A(k')$ by performing the following row operations
            \begin{alignat*}{4}
	              R_{i-T} &\rightarrow R_{i-T}~- ~\frac{a-i+T-(i-T-1)+k'p}{i-T} R_{i-T-1} \\
	            & ~~\vdots  \\
	            R_{m} &\rightarrow ~ ~~R_{m}~ - ~\frac{a-i+T-(m-1)+k'p}{m} R_{m-1} \\
	              & ~~\vdots  \\
	            R_{k'+1} &\rightarrow ~~ R_{k'+1}~- ~\frac{a-i+T-k'+k'p}{k'+1} R_{k'}
            \end{alignat*} 
            on $A(k')$. We claim that for $k'=0,\ldots,i-n$, we have
            \begin{equation}\label{eq: matrix formula recursion a>2i}
            \resizebox{0.999\hsize}{!}{%
                $A(k')=\left[\begin{array}{c|c}
		           p \sum\limits_{k=m}^{i-n-1} (p-1)^m \binom{s-l}{a-i+T+k(p-1)} \binom{k}{m} +  O(p(r-s))
		           &\sum\limits_{k=m}^{i-n-1} (p-1)^m\binom{s-i+T}{a-i+T+k(p-1)} \binom{k}{m} +  O(p(r-s))\\
		           \hline \\
		           p \left( \binom{r-l}{m-k'} \frac{(r-a+k'-1)_{k'}}{(m)_{k'}}- \binom{s-l}{m-k'} \frac{(s-a+k'-1)_{k'}}{(m)_{k'}} \right) \binom{a-l-m}{a-i+T-m}
                  & p\left( \binom{r-i+T}{m-k'} \frac{(r-a+k')_{k'}}{(m)_{k'}}- \binom{s-i+T}{m-k'}\frac{(s-a+k')_{k'}}{(m)_{k'}}\right) \frac{a-s}{a-i+T-m} \\  \\
		           + p(p-1)^{k'} \sum\limits_{k=k'}^{i-n-1} \binom{s-l}{a-i+T+k(p-1)} \binom{a-i+T+k(p-1)}{m-k'} \frac{(k)_{k'}}{(m)_{k'}}  
                  &  (p-1)^{k'}\sum\limits_{k=k'}^{i-n-1} \binom{s-i+T}{a-i+T+k(p-1)} \binom{a-i+T+k(p-1)}{m-k'} \frac{(k)_{k'}}{(m)_{k'}}   \\  \\
		           + p \binom{s-l}{s-i+T} \binom{s-i+T}{m-k'} \frac{\prod\limits_{k=0}^{k'-1}(s-a-k(p-1))}{(m)_{k'}} - p \binom{r-l}{r-i+T}\binom{r-i+T}{m-k'} \frac{\prod\limits_{k=0}^{k'-1}(r-a-k(p-1))}{(m)_{k'}} 
                 &+  p \frac{s - r}{a - i + T - m} \binom{s - i + T}{m-k'} \frac{(s-a+k')_{k'}}{(m)_{k'}} \\ \\
                 + O(p^2(r-s)) & + O(p^2(r-s))\\ 
	          \end{array} \right],$
             } %
            \end{equation}
            where the range of $m$ in the upper block (resp. lower block) is $0,\ldots,k'-1$ (resp. $k',\ldots,i-T$) and the range of $l$ in the left block (resp. right block) is $0,\ldots,i-T-1$ (resp. $i-T$). We prove this by recursion on $k'$. Note that every $0\leq j < s-i+T$ with $j\equiv a-i+T$ can be expressed as $a-i+T+k(p-1)$ for some $0\leq k \leq i-n-1$. This implies $A(0)=A$. Assume that the above formula holds for $k'$. We need to prove $A(k'+1)$ is also given by the above expression.

            Note that going modulo $p(r-s)$ in $k'^{\mathrm{th}}$ row of $A(k')$ gives $k'^{\text{th}}$ row of $A(k'+1)$. To obtain the first term in the bottom left block of $A(k'+1)$ we use the following identity
            \begin{equation}\label{identity 2}
	          \begin{split}
		          &\frac{1}{(m)_{k'}} \Bigg( (r-a+k'-1)_{k'}  \binom{r-l}{m-k'} - (s-a+k'-1)_{k'} \binom{s-l}{m-k'}  \Bigg) \binom{a-l-m}{a-i+T-m} \\ 
		          & \qquad\qquad - \frac{a-i+T-(m-1)+k'p}{m} \times   \\ 
		          & \qquad \qquad \qquad  \frac{1}{(m-1)_{k'}} \Bigg( (r-a+k'-1)_{k'} \binom{r-l}{m-1-k'} - (s-a+k'-1)_{k'}  \binom{s-l}{m-1-k'}  \Bigg) \binom{a-l-(m-1)}{a-i+T-(m-1)} \\
		        &\quad = \frac{1}{(m)_{k'+1}} \Bigg( (r-a+k')_{k'+1} \binom{r-l}{m-1-k'} -(s-a+k')_{k'+1} \binom{s-l}{m-1-k'}  \Bigg) \binom{a-l-m}{a-i+T-m} +O(p(r-s)).
	        \end{split}    
            \end{equation}           
           Note that for $j' > m \geq k' \geq 0$ and $j \geq0$ we have
           \begin{align}\label{identity 1}
	               \frac{1}{(m)_{k'}}\binom{j'}{m-k'} - \frac{j-(m-1-k')}{m} \frac{1}{(m-1)_{k'}}\binom{j'}{m-1-k'}  = \frac{j'-j}{(m)_{k'+1}} \binom{j'}{m-k'-1}. 
            \end{align}
             Taking $j=a-i+T+k'(p-1)$ and $j'=a-i+T+k'(p-1), \ldots, a-i+T+(i-n-1)(p-1)$ in the above equation we obtain the second term in the bottom left block of $A(k'+1)$. Taking $j=a-i+T+k'(p-1)$ and $j'=s-i+T$ (resp. $j'=r-i+T$) in the above equation we obtain the third and fourth term in the bottom left block of $A(k'+1)$. 

             To obtain the third term in the bottom right block of $A(k'+1)$ note that \begin{equation}\label{identity 3}
	        \begin{split}
                 &\frac{1}{a-i+T-m} \frac{(s-a+k')_{k'}}{(m)_{k'}}  \binom{s-i+T}{m-k'} \\
                 & \qquad \qquad \qquad  -  \frac{a-i+T-(m-1)+k'p}{m} \frac{1}{a-i+T-(m-1)} \frac{(s-a+k')_{k'}}{(m-1)_{k'}}  \binom{s-i+T}{m-1-k'}\\
                 & \qquad \qquad = \frac{1}{a-i+T-m} \frac{(s-a+k'+1)_{k'+1}}{(m)_{k'+1}}  \binom{s-i+T}{m-(k'+1)} \\
                 & \qquad \qquad \qquad \qquad \qquad - \frac{k'p}{m} \frac{1}{a-i+T-(m-1)} \frac{(s-a+k')_{k'}}{(m-1)_{k'}}   \binom{s-i+T}{m-1-k'}.
            \end{split}
            \end{equation}
            Note that the above identity holds when $s$ is replaced by $r$. To obtain the second term in the bottom right block of $A(k'+1)$, we use \eqref{identity 1} with $l=i-T$, $j=a-i+T+k'(p-1)$ and $j'=a-i+T+k'(p-1), \ldots, a-i+T+(i-n-1)(p-1)$. Using \eqref{identity 3} for $r$ and $s$, we see that
            \begin{equation}\label{identity 4}
	        \begin{split}
		        &\frac{1}{a-i+T-m} \Bigg( \frac{(r-a+k')_{k'}}{(m)_{k'}}  \binom{r-i+T}{m-k'} - \frac{(s-a+k')_{k'}}{(m)_{k'}} \binom{s-i+T}{m-k'}  \Bigg) \\
		        & \qquad \qquad \qquad  -  \frac{a-i+T-(m-1)+k'p}{m} \times \\
		        & \qquad \qquad \qquad \qquad  \frac{1}{a-i+T-(m-1)} \Bigg( \frac{(r-a+k')_{k'}}{(m-1)_{k'}} \binom{r-i+T}{m-1-k'} - \frac{(s-a+k')_{k'}}{(m-1)_{k'}} \binom{s-i+T}{m-1-k'}  \Bigg)\\
		        & \qquad \qquad = \frac{1}{a-i+T-m} \Bigg( \frac{(r-a+k'+1)_{k'+1}}{(m)_{k'+1}}  \binom{r-i+T}{m-(k'+1)} - \frac{(s-a+k'+1)_{k'+1}}{(m)_{k'+1}} \binom{s-i+T}{m-(k'+1)}  \Bigg) \\
		        & \qquad \qquad \qquad \qquad \qquad +O(p(r-s)).
	        \end{split}
            \end{equation}
            This gives the first term in the bottom right block of $A(k'+1)$. Hence we obtain the claim.

           From the above claim, it follows that $\det(A) = \det(A(i-n))$ and $\det(A(i-n))$ equals
           \begin{equation*}
            \resizebox{0.999\hsize}{!}{%
                $\left\lvert\begin{array}{c|c}
		           p \sum\limits_{k=m}^{i-n-1} (p-1)^m \binom{s-l}{a-i+T+k(p-1)} \binom{k}{m} +  O(p(r-s))
		           &\sum\limits_{k=m}^{i-n-1} (p-1)^m\binom{s-i+T}{a-i+T+k(p-1)} \binom{k}{m} +  O(p(r-s))\\
		           \hline \\
		           p \left( \binom{r-l}{m-(i-n)} \frac{(r-a+i-n-1)_{i-n}}{(m)_{i-n}}- \binom{s-l}{m-(i-n)} \frac{(s-a+i-n-1)_{i-n}}{(m)_{i-n}} \right) \binom{a-l-m}{a-i+T-m}
                 & p\left( \binom{r-i+T}{m-(i-n)} \frac{(r-a+i-n)_{i-n}}{(m)_{i-n}}- \binom{s-i+T}{m-(i-n)}\frac{(s-a+i-n)_{i-n}}{(m)_{i-n}}\right) \frac{a-s}{a-i+T-m} \\  \\
		           + p \binom{s-l}{s-i+T} \binom{s-i+T}{m-(i-n)} \frac{\prod\limits_{k=0}^{i-n-1}(s-a-k(p-1))}{(m)_{i-n}} - p \binom{r-l}{r-i+T}\binom{r-i+T}{m-(i-n)} \frac{\prod\limits_{k=0}^{i-n-1}(r-a-k(p-1))}{(m)_{i-n}} 
                 &+  p \frac{s - r}{a - i + T - m} \binom{s - i + T}{m-(i-n)} \frac{(s-a+i-n)_{i-n}}{(m)_{i-n}} \\ \\
                 + O(p^2(r-s)) & + O(p^2(r-s))\\ 
	          \end{array} \right\rvert$,
             } %
            \end{equation*}
            where the range of $m$ in the upper block (resp. lower block) is $0,\ldots,i-n-1$ (resp. $i-n,\ldots,i-T$) and the range of $l$ in the left block (resp. right block) is $0,\ldots,i-T-1$ (resp. $i-T$). We now simplify every term in the lower blocks up to $O(p^2(r-s))$.  Using \eqref{est. error term 1 adhc} and \eqref{est. error term 2 adhc} (with $i-T$ there replaced by $a-i+T$) we get the lower left block entry equals
            \begin{equation}\label{est. error term 1 a>2i}
                p\frac{(r-s)}{m!}  \times  \frac{(-1)^{i-n} (i-n)!}{(a-i+T-m)!} \times \frac{(a-i+n-l)!}{(i-T-l)!}\times \bigl(H_{a-i+T-m}-H_{a-l-m}\bigr)+O(p^2(r-s)).
            \end{equation}
            Observe that the second term in the lower right block is a multiple of $p^2(r-s)$ as $(s-a+i-n)_{(i-n)}$ is a multiple of $p$. For the first term in the lower right block, note that
            \begin{align*}
                &\binom{r-i+T}{m-(i-n)} \frac{(r-a+i-n)_{(i-n)}}{(m)_{(i-n)}} - \binom{s-a+i+T}{m-(i-n)}\frac{(s-a+i+n)_{(i-n)}}{(m)_{(i-n)}} \\
                & \qquad  = \frac{1}{m!} \Big((r-i+T)\cdots(r-n-m+T+1)\cdot(r-a+i-n)\cdots(r-a+1)\\
                & \qquad \qquad \qquad -(s-i+T)\cdots(s-n-m+T+1) \cdot (s-a+i-n)\cdots(s-a+1) \Big) \\
                & \qquad \equiv \frac{r-s}{m!} \Big((a-2i+n+T)\cdots(a-i-m+T+1) \Big) \cdot \Big((-1)\cdots(-(i-n-1)) \Big)  \\
                & \qquad \equiv (r-s) \cdot \frac{(-1)^{i-n-1}}{m!} \cdot  \frac{(a-2i+n+T)!}{(a-i+T-m)!} \cdot  \frac{(i-n)!}{(i-n)} \mod p(r-s).
            \end{align*}
            Since $a-s \equiv (i-n) \mod p$, we get
            \begin{align}\label{est. error term 2 a>2i}
            \begin{split}
	            &\left\{ \binom{r-i+T}{m-(i-n)} \frac{(r-a+i-n)_{(i-n)}}{(m)_{(i-n)}} -\binom{s-a+i+T}{m-(i-n)}\frac{(s-a+i+n)_{(i-n)}}{(m)_{(i-n)}} \right\} \frac{a-s}{a-i+T-m} \\
	            & \qquad \equiv -\frac{r-s}{a-i+T-m} \cdot  \frac{(-1)^{i-n}}{m!} \cdot  \frac{(a-2i+n+T)!}{(a-i+T-m)!}  \cdot (i-n)! \mod p(r-s).
            \end{split}    
            \end{align}
            Substituting \eqref{est. error term 1 a>2i} and \eqref{est. error term 2 a>2i} in the above expression of $\det(A(i-n))$ we get $\det(A)$ equals
            \begin{equation*}
            \resizebox{0.999\hsize}{!}{%
                $\left\lvert\begin{array}{c|c}
		           p \sum\limits_{k=m}^{i-n-1} (p-1)^m \binom{s-l}{a-i+T+k(p-1)} \binom{k}{m} +  O(p(r-s))
		           &\sum\limits_{k=m}^{i-n-1} (p-1)^m\binom{s-i+T}{a-i+T+k(p-1)} \binom{k}{m} +  O(p(r-s))\\
		           \hline \\
		            p(r-s) \times \frac{(-1)^{i-n} (i-n)!}{m!(a-i+T-m)!} \times \frac{(a-i+n-l)!}{(i-T-l)!} \times   \bigl(H_{a-i+T-m}-H_{a-l-m}\bigr) 
                 & -p\frac{r-s}{a-i+T-m} \cdot  \frac{(-1)^{i-n}}{m!} \cdot  \frac{(a-2i+n+T)!}{(a-i+T-m)!}  \cdot (i-n)! \\ \\
                 + O(p^2(r-s)) & + O(p^2(r-s))\\ 
	          \end{array} \right\rvert$,
             } %
            \end{equation*}
            where the range of $m$ in the upper block (resp. lower block) is $0,\ldots,i-n-1$ (resp. $i-n,\ldots,i-T$) and the range of $l$ in left block (resp. right block) is $0,\ldots,i-T-1$ (resp. $i-T$). Pulling out $(p-1)^m$ from every row in the upper block and $p(r-s) \times \frac{(-1)^{i-n} (i-n)!}{(m!(i-T-m)!)}$ for every row in the lower block we get
           \begin{align*}
	            \det(A) &= p^{n-T+1} (r-s)^{n-T+1} \prod_{m=0}^{i-n-1} (p-1)^m  \times \prod_{m=i-n}^{i-T} \frac{(-1)^{i-n} (i-n)!}{m!(a-i+T-m)!} \\
	              &\qquad \times \det 
                \left[\begin{array}{c|c}
		        p\sum\limits_{k=0}^{i-n-1} \binom{s-l}{a-i+T+k(p-1)} \binom{k}{m} &
		        \sum\limits_{k=0}^{i-n-1} \binom{s-i+T}{a-i+T+k(p-1)} \binom{k}{m} 
		        \\ +O(p(r-s)) & +O(p(r-s)) \\ \hline \\
		          \frac{(a-i+n-l)!}{(i-T-l)!} \times  (H_{a-i+T-m}-H_{a-l-m}) +O(p) &  - \frac{(a-2i+n+T)!}{a-i+T-m} +O(p)
	            \end{array}\right].
            \end{align*}
            We now perform the following row operations: 
            \begin{align*}
                R_{i-n-2} &\rightarrow R_{i-n-2}- (i-n-1)R_{i-n-1} \\ 
                &~~\vdots \\
                R_m &\rightarrow R_m - \sum_{k=m+1}^{i-n-1}\binom{k}{m} R_k \\
                &~~\vdots \\
                R_0 &\rightarrow R_0 - \sum_{k=1}^{i-n-1}\binom{k}{0}R_k
            \end{align*}
            so that every entry in the upper block is given by a binomial coefficient up to an error term. Thus 
            \begin{align*}
            \begin{split}
	                \det(A) &= p^{n-T+1}(r-s)^{n-T+1} \prod_{m=0}^{i-n-1} (p-1)^m  \times \prod_{m=i-n}^{i-T} \frac{(-1)^{i-n} (i-n)!}{m!(a-i+T-m)!} \\
	                &\qquad\times  \det \left(
	                \begin{array}{c|c}
		            p \binom{s-l}{a-i+T+m(p-1)} +O(p(r-s)) 
		            & \binom{s-i+T}{a-i+T+m(p-1)} +O(p(r-s)) \\ \hline \\
		            \frac{(a-i+n-l)!}{(i-T-l)!} \times  (H_{a-i+T-m}-H_{a-l-m}) +O(p) &  - \frac{(a-2i+n+T)!}{a-i+T-m} +O(p)
	                \end{array}\right).
            \end{split}        
            \end{align*}
            Note that by Lucas theorem
           \begin{align*}
	                \binom{s-l}{a-i+T+m(p-1)} \equiv \binom{i-n}{m}\binom{a-i+n-l}{a-i+T-m} \mod p
            \end{align*}
            for $0 \leq l \leq i-T-1$ and $0\leq m \leq i-n-1$. Also, for $0\leq m \leq i-n-1$ we have
            \begin{align*}
                    \binom{s-i+T}{a-i+T+m(p-1)} & = \frac{(s-i+T)\cdots (s-a+i-n)}{(a-i+T+m(p-1))\cdots (i-n+m(p-1))} \cdot
                    \binom{s-a+i-n-1}{i-n-1+m(p-1)} \\
                    & \equiv (i-n)p \times  \frac{(a-2i+n+T)!(i-n-1-m)!}{(a-i+T-m)!} \\
                    &\qquad \qquad \qquad \times \binom{i-n-1}{m} \binom{p-1}{i-n-1-m} \\
                    &\equiv p \times  \frac{(a-2i+n+T)!(i-n-m)!}{(a-i+T-m)!} (-1)^{i-n-1-m} \binom{i-n}{m} \mod {p^{2}},
            \end{align*}
            where the penultimate step follows from Lucas' theorem and the last step uses $\binom{p-1}{i-n-1-m} \equiv (-1)^{i-n-1-m} \mod p$. Substituting this in the above expression for $\det(A)$ and pulling out  $(a-2i+n+T)!$ from  the last column and $p\binom{i-n}{m}$ from $m^{\text{th}}$ row for $m=0,\ldots,i-n-1$ we get 
            \begin{align}\label{eq: det(A) interms of det(B) a>2i}
            \begin{split}
	                \det(A) &= p^{i-T+1} \cdot (r-s)^{n-T+1} \prod_{m=0}^{i-n-1} (p-1)^m  \binom{i-n}{m} \times \prod_{m=i-n}^{i-T} \frac{(-1)^{i-n} (i-n)!}{m!(a-i+T-m)!} \\
	                  &\qquad \times \det \left(
	                \begin{array}{c|c}
		            \binom{a-i+n-l}{a-i+T-m} +O(r-s) & \frac{(-1)^{i-n-1-m}(i-n-m)!}{(a-i+T-m)!} +O(r-s) \\ \hline \\
		            \frac{(a-i+n-l)!}{(i-T-l)!} \times  (H_{a-i+T-m}-H_{a-l-m}) +O(p) &  -\frac{1}{a-i+T-m} +O(p)
	                \end{array}\right). \\
	                & \equiv p^{i-T+1} \cdot (r-s)^{n-T+1} \prod_{m=0}^{i-n-1} (p-1)^m  \binom{i-n}{m} \times \prod_{m=i-n}^{i-T} \frac{(-1)^{i-n} (i-n)!}{m!(a-i+T-m)!} \\
	                &\qquad \qquad\qquad \times \det(B) \qquad \mod(p^{i-T+2}(r-s)^{(n-T+1)}),
            \end{split}        
            \end{align}
            where
           \begin{align*}
	                B = \left(\begin{array}{c|c} 
	   	            \left(\binom{a-i+n-l}{a-i+T-m}\right)_{\substack{0 \leq m \leq i-n-1 \\[1pt] 0 \le l \le i-T-1}}
	   	            & \left(\frac{(-1)^{i-n-1-m}(i-n-m)!}{(a-i+T-m)!}\right)_{\substack{0 \leq m \leq i-n-1 \\[1pt] l = i-T}} \\  \\ \hline \\
	   	              \Big(\frac{(a-i+n-l)!(H_{a-i+T-m}-H_{a-l-m})}{(i-T-l)!} \Big)_{\substack{i-n \leq m \le i-T \\[1pt] 0 \leq l \leq i-T-1}}
	                &
	   	            \left( -\frac{1}{a-i+T-m}\right)_{\substack{i-n \leq m \le i-T \\[1pt] l = i-T}}
	        \end{array}\right).
            \end{align*}
            Below we show that $\det(B) \in \mathbb{Z}_p^\times$. To do this we perform column operations on $B$. We first show that after performing column operations on $B$ we can reduce to a matrix with every row in the upper block containing exactly one non-zero entry. 

            Note that
            \[ \binom{a-i+n-l}{a-i+T-m} = \frac{(a-i+n-l)!}{(i-T-l)!} \cdot \frac{(i-n-m)!}{(a-i+T-m)!} \binom{i-T-l}{i-n-m}.
            \]
            Thus pulling out $\frac{(a-i+n-l)!}{(i-T-l)!}$ from $l^{\mathrm{th}}$ column for $l=0,\ldots, i-T-1$ and $\frac{(i-n-m)!}{(a-i+T-m)!}$ from $m^{\mathrm{th}}$ row for $m=0,\ldots, i-n-1$, we see that
            \begin{align*}
	                  \det(B) 
                    &= \prod_{l=0}^{i-T-1} \frac{(a-i+n-l)!}{(i-T-l)!} \times \prod_{m=0}^{i-n-1} \frac{(i-n-m)!}{(a-i+T-m)!} \\
	                & \qquad \qquad \times 
                    \det \left[
	                \begin{array}{c|c}
		            \left( \binom{i-T-l}{i-n-m} \right)_{\substack{m=0,\dots,i-n-1 \\ l=0,\dots,i-T-1}}
		            & \left((-1)^{i-n-1-m}\right)_{m=0,\ldots,i-n-1 } \\[10pt]\hline\\[-5pt]
		            \left(H_{a-i+T-m}-H_{a-l-m}\right)_{\substack{m=i-n,\dots,i-T \\ l=0,\dots,i-T-1}}
		            &\left(-\frac{1}{a-i+T-m}\right)_{m=i-n,\dots,i-T }
	                \end{array}\right].
            \end{align*}
            Using Pascal's identity, for $0\le m\le i-n-1$ and $0\le l\le i-T-1$ we have
            \[
               \binom{i-T-l}{i-n-m}-\binom{i-T-l-1}{i-n-m} = \binom{i-T-l-1}{i-n-m-1}.
            \]
            Applying the column operations 
            \begin{align*}
                    C_0 &\rightarrow\; C_0 - C_{1} \\
                    &~~\vdots \\
                    C_l &\rightarrow\; C_l - C_{l+1}, \\
                    &~~\vdots \\
                    C_{i-T-2} &\rightarrow\; C_{i-T-2} - C_{i-T-1}.
            \end{align*}
            and using the above identity we obtain 
            \begin{align*}
	                \det(B) &=\prod_{l=0}^{i-T-1} \frac{(a-i+n-l)!}{(i-T-l)!} \times \prod_{m=0}^{i-n-1} \frac{(i-n-m)!}{(a-i+T-m)!}\\
	                &\qquad\qquad \times \det \left[
	                \begin{array}{c|c}
		                    \left( \binom{i-T-1-l}{i-n-1-m} \right)_{\substack{m=0,\ldots,i-n-1 \\ l=0,\dots,i-T-1}}
                            & \left((-1)^{i-n-1-m}\right)_{m=0,\ldots,i-n-1 } \\[10pt]\hline\\[-5pt]
		                    \left(-\frac{1}{a-l-m}\right)_{\substack{m=i-n,\ldots,i-T \\ l=0,\dots,i-T-1}}
		                      &\left(-\frac{1}{a-i+T-m}\right)_{m=i-n,\ldots,i-T}
	                  \end{array}
	                \right].
            \end{align*}
            Taking $N=i-T-1-l$ and $k'=i-n-1-m$ in Lemma~\ref{lem: combinatorial id for det(B)} $(i)$ we get
            \begin{align}\label{eq: comb. id 1}
            \begin{split}	
	                \binom{i-T-1-l}{i-n-1-m} &+ \sum_{l'=l+1}^{i-T-1} (-1)^{l'-l}\binom{i-T-1-l}{l'-l}\binom{i-T-1-l'}{i-n-1-m} \\
	                &= \sum_{l'=0}^{i-T-1-l} (-1)^{l'}\binom{i-T-1-l}{l'}\binom{i-T-1-l-l'}{i-n-1-m} \\
	                &= \begin{cases}
		               1, & \text{if } l=n-T+m,\\
		               0, & \text{otherwise}.
	                  \end{cases}
            \end{split}
            \end{align}	
            Also, taking $\alpha = a-l-m$ and $N=i-T-1-l$ in Lemma~\ref{lem: combinatorial id for det(B)} $(ii)$ we get
            \begin{align}\label{eq: comb. id 2}
            \begin{split}
	             \frac{1}{a-l-m}+ \sum_{l'=l+1}^{i-T-1} (-1)^{l'-l} \binom{i-T-1-l}{l'-l} \frac{1}{a-l'-m} &= \sum_{l'=0}^{i-T-1-l} (-1)^{l'} \binom{i-T-1-l}{l'} \frac{1}{a-l-m-l'}   \\
                 &=(-1)^{i-T-1-l} \frac{(i-T-1-l)!}{(a-i+T+1-m)\cdots(a-l-m)}.
            \end{split}     
            \end{align}
            Applying the column operations 
            \begin{align*}
	                C_{0} & \rightarrow C_{0} + \sum_{l'=1}^{i-T-1} (-1)^{l'} \binom{i-T-1}{l'}C_{l'}, \\
	                  & ~~ \vdots \\
	                C_{l} & \rightarrow C_{l} + \sum_{l'=l+1}^{i-T-1} (-1)^{l'-l} \binom{i-T-1-l}{l'-l}C_{l'} \\
	                  & ~~\vdots \\
	                C_{i-T-2} &\rightarrow C_{i-T-2} - C_{i-T-1}
            \end{align*}
            and using \eqref{eq: comb. id 1} for the upper left block and \eqref{eq: comb. id 2} for the lower left block  we get 
            \begin{align*}
	                  \det(B) 
                    &=\prod_{l=0}^{i-T-1} \frac{(a-i+n-l)!}{(i-T-l)!} \times \prod_{m=0}^{i-n-1} \frac{(i-n-m)!}{(a-i+T-m)!}\\
	                &\qquad\qquad \times \det \left[
	                \begin{array}{c|c}
		                \left( \delta_{l,n-T+m} \right)_{\substack{m=0,\ldots,i-n-1 \\ l=0,\dots,i-T-1}} 
                         & \left((-1)^{i-n-1-m}\right)_{m=0,\ldots,i-n-1 } \\[10pt]\hline\\[-5pt]
		                \left(\frac{(-1)^{i-T-l}(i-T-1-l)!}{(a-i+T+1-m)\cdots(a-l-m)}\right)_{\substack{m=i-n,\ldots,i-T \\ l=0,\dots,i-T-1}}
		                &\left(-\frac{1}{a-i+T-m}\right)_{m=i-n,\ldots,i-T}
	                \end{array} \right].
            \end{align*}
            This matrix is close to what we wanted to achieve except that the upper right block should be zero.  Applying the following column operation
            \[ C_{i-T} \rightarrow C_{i-T} + \sum_{l'=n-T}^{i-T-1} (-1)^{i-T-l'} C_{l'} \]
            we get 
            \begin{align*}
	            \det(B) 
                &=\prod_{l=0}^{i-T-1} \frac{(a-i+n-l)!}{(i-T-l)!} \times \prod_{m=0}^{i-n-1} \frac{(i-n-m)!}{(a-i+T-m)!}\\
	            &\qquad\qquad \times \det \left[
	              \begin{array}{c|c}
		        \left( \delta_{l,n-T+m} \right)_{\substack{m=0,\ldots,i-n-1 \\ l=0,\dots,i-T-1}} 
                &\left(0\right)_{m=0,\ldots,i-n-1 }\\[10pt]\hline\\[-5pt]
		        \left(\frac{(-1)^{i-T-l}(i-T-1-l)!}{(a-i+T+1-m)\cdots (a-l-m)}\right)_{\substack{m=i-n,\ldots,i-T \\ l=0,\dots,i-T-1}}
		        & \left(\frac{-(i-n)!}{(a-i+T-m)\cdots(a-n+T-m)}\right)_{m=i-n,\ldots,i-T}
	            \end{array}\right],
            \end{align*}
            where we have used Lemma~\ref{lem: combinatorial id for det(B)} $(iii)$ with $\alpha=a-i+T-m$ and $N=i-n$ to obtain the expression for the lower right block.  Thus we obtain the form of the matrix we claimed earlier. 

        Expanding the determinant using the first row, then the second row, and so on, we get
        \begin{align*}
	        \det(B) & = \prod_{l=0}^{i-T-1} \frac{(a-i+n-l)!}{(i-T-l)!} \prod_{m=0}^{i-n-1} \frac{(i-n-m)!}{(a-i+T-m)!}\times (-1)^{(n-T)(i-n)} \times \\
	          & \qquad \det
	        \left[\begin{array}{c|c}
		      \left(\frac{(-1)^{i-T-l}(i-T-1-l)!}{(a-i+T+1-m)\cdots(a-l-m)}\right)_{\substack{m=i-n,\ldots,i-T \\ l=0,\dots,n-T-1}}
		    &\left(\frac{-(i-n)!}{(a-i+T-m)\cdots(a-n+T-m)}\right)_{m=i-n,\ldots,i-T}
	        \end{array} \right] \\
	        & = \prod_{l=0}^{i-T-1} \frac{(-1)^{(i-T-l)}(a-i+n-l)!}{(i-T-l)} \prod_{m=0}^{i-n-1} \frac{(i-n-m)!}{(a-i+T-m)!} \times (-1)^{(n-T)(i-n)} \times (i-n-1)! \\
	          & \qquad \times  \det \left[
	        \begin{array}{c|c}
		        \left(\frac{(a-i+T-m)!}{(a-l-m)!}\right)_{\substack{m=i-n,\ldots,i-T \\ l=0,\dots,n-T-1}}
		          & \left(-(i-n)\frac{(a-i+T-m-1)!}{(a-n+T-m)!}\right)_{\substack{m=i-n,\ldots,i-T \\ l=n-T }}
	        \end{array}
	        \right],
        \end{align*}
        where in the last step we have pulled out $(-1)^{i-T-l}(i-T-1-l)!$ from the $l^{\text{th}}$ column for $l=0,\ldots, n-T-1$ and $(i-n-1)!$ from the $(n-T)^\mathrm{th}$ column. We now perform one final set of column operations so that every entry of the last matrix can be realised as a binomial coefficient. Applying the column operations
        \begin{align*}
	        C_{n-T-1} &\rightarrow C_{n-T-1} + \frac{1}{i-n} C_{n-T} \\
	        &~~ \vdots \\
	        C_{l} &\rightarrow C_{l} + \frac{1}{i-T-1-l} C_{l+1} \\
	        &~~ \vdots \\
	        C_{0} &\rightarrow C_{0} + \frac{1}{i-T-1} C_{1},
        \end{align*}
         we get
        \begin{align*}
	            \det(B) & =  \prod_{l=0}^{i-T-1} \frac{(-1)^{(i-T-l)}(a-i+n-l)!}{(i-T-l)} \prod_{m=0}^{i-n-1} \frac{(i-n-m)!}{(a-i+T-m)!}\times (-1)^{(n-T)(i-n)}\\ 
                & \qquad  \times (i-n-1)! \times \det \left[
		          -(i-T-l)\frac{(a-i+T-m-1)!}{(a-l-m)!} \right]_{\substack{m=i-n,\ldots,i-T \\ l=0,\ldots,n-T}},
        \end{align*}
        where we have used the identity 
        \begin{align*}
	           \frac{(a-i+T-m)!}{(a-l-m)!} -\frac{(a-i+T-m-1)!}{(a-l-m-1)!} = -(i-T-l) \frac{(a-i+T-m-1)!}{(a-l-m)!}.
        \end{align*}
        Pulling out $(a-i+T-1-m)!/(a-m)!$ from the $m^{\text{th}}$ row and $l!(l-i+T)$ from the $l^{\text{th}}$ column we get
        \begin{align*}
	            \det(B) & = \prod_{l=0}^{i-T-1} \frac{(-1)^{(i-T-l)}(a-i+n-l)!}{(i-T-l)} \times \prod_{m=0}^{i-n-1} \frac{(i-n-m)!}{(a-i+T-m)!}\times (-1)^{(n-T)(i-n)}\times (i-n-1)! \\ 
                & \qquad  \times \prod_{m=i-n}^{i-T}\frac{(a-i+T-1-m)!}{(a-m)!}   \times \prod_{l=0}^{n-T} (l-i+T)  
                \times \det_{\substack{m=i-n,\ldots,i-T \\ l=0,\ldots,n-T}}
	            \left( \frac{(a-m)!}{(a-l-m)!l!} \right).
        \end{align*}
        Note that
        \begin{align*}
	            \det_{\substack{m=i-n,\ldots,i-T \\ l=0,\ldots,n-T}} \left(\frac{(a-m)!}{(a-l-m)!l!}\right) 
                &= \det_{\substack{m=i-n,\ldots,i-T \\ l=0,\ldots,n-T}}\left( \binom{a-m}{l} \right) \\ 
                &=\det_{\substack{m=0,\ldots,n-T \\ l=0,\ldots,n-T}}\left( \binom{a-i+n-m}{l} \right). \\
	             &= (-1)^{\binom{n-T+1}{2}} \det_{\substack{ m=0,\ldots,n-T \\ l=0,\ldots,n-T}}\left( \binom{a-i+T+m}{l} \right).  
        \end{align*}
        Applying Lemma~\ref{cor: GV det} $(i)$ with $d=1$, we see that the last determinant is a $p$-adic unit. Hence  $\det(B) \in \mathbb{Z}_{p}^\times$. Substituting this in \eqref{eq: det(A) interms of det(B) a>2i}, we get $\det(A) \in p^{i-T+1} (r-s)^{n-T+1} \mathbb{Z}_p^\times$. In particular, the system of linear equations \eqref{eq: linear eq. in matrix form a>2i} has a solution in $\Q$ by Cramer's rule.

        We now show that \eqref{eq: linear eq. in matrix form a>2i} has a solution in $\mathbb{Z}_p$. Again by Cramer's rule, we must show $p^{t+1}\det(A_{i-T,l})/\det(A) \in \Z_p$ for all $0 \leq l \leq i-T$, where $A_{i-T,l}$ is the minor of $(i-T,l)$ entry of $A$. As $\det(A) \in p^{i-T+1} (r-s)^{n-T+1} \mathbb{Z}_p^\times$ and $v(r-s)=t$, it is enough to show $\det(A_{i-T,l}) \in p^{i-T} (r-s)^{n-T} \mathbb{Z}_p$. 


      Note that by \cite[Lemma 2.5]{BG15} and the condition $a-i+T>i-T$, we have 
      \begin{align*}
          \sum_{\substack{0<j<s-i+T \\ j \equiv a-i+T \mod (p-1)}}^{} \binom{s-i+T}{j} \binom{j}{m} = \binom{s-i+T}{m} \sum_{\substack{0<j<s-i+T \\ j \equiv a-i+T \mod (p-1)}}^{} \binom{s-i+T-m}{j-m} \equiv 0 \mod p
      \end{align*}
      for $m=0,\ldots,i-T$.
      Thus every entry of $A$ in \eqref{eq: A sdad} is a multiple of $p$.  We define matrices $A_{i-T,l}(0), \ldots,$ $A_{i-T,l}(i-n)$ recursively. Let $A_{i-T,l}(0):=A_{i-T,l}$. Having defined $A_{i-T,l}(k')$, we define $A_{i-T,l}(k'+1)$  as the matrix obtained by applying the following row operations 
    \begin{align*}
	      R_{i-T-1} &\rightarrow ~R_{i-T-1}- \frac{a-i+T-(i-T-2)+k'p}{i-T-1} R_{i-T-2} \\
	      &~~\vdots  \\
	    R_{m} &\rightarrow ~ R_{m}- \frac{a-i+T-(m-1)+k'p}{m} R_{m-1} \\\
	    &~~\vdots  \\
	    R_{k'+1} &\rightarrow ~ R_{k'+1}- \frac{a-i+T-k'+k'p}{k'+1} R_{k'}
    \end{align*}
    on $A_{i-T,l}(k')$ for $k'=0,\ldots, i-n-1$. By a similar argument used in the earlier computation at the beginning of the proof we can show that $A_{i-T,l}(k')$ is the same as the $(i-T,l)$-th minor of $A(k')$. This gives that the last $(n-T)$ rows of $A_{i-T,l}(i-n)$ are multiples of $p(r-s)$. Pulling out $p$ from the first $i-n$ rows and $p(r-s)$ from the last $(n-T)$ rows we see that $\det(A_{i-T,l}) \in p^{i-T} (r-s)^{n-T} \mathbb{Z}_p$, as desired.  Hence the system of linear equations given by \eqref{eq: linear eq. in matrix form a>2i} has a solution in $\Z_p$. This proves $(i)$.
    
    Reducing the equations \eqref{eq: linear eq. in matrix form a>2i} modulo $p^{t+1}$ we get
     \begin{equation}\label{eq: congruences for binomial sums for s a>2i}
         \sum\limits_{\substack{a-i+T \leq j < s-i+T \\ j \equiv a-i+T ~\mathrm{mod}~(p-1)}} \left(\sum_{l=0}^{i-T-1}p \beta_{l} \binom{s-l}{j} + \beta_{i-T}\binom{s-i+T}{j} \right) \binom{j}{m} \equiv 0 \mod p^{t+1}.
     \end{equation}
    for $m=0,\ldots,i-T$.
    
    Using Corollary~\ref{cor: GV det} $(i)$, we see that 
     \begin{align}\label{eq:matrix invertible GV det ad a>2i}
         \left( \binom{a-i+T+k(p-1)}{m}\right)_{\substack{0\leq m \leq i-n-1 \\ 0\leq k \leq i-n-1}} \in \gl_{i-n}(\Z_p)
     \end{align}
     is invertible. We now show that $(ii)$ holds. Observe that by Corollary~\ref{cor: binomial sums under congruences 2}, for  $l=0,\ldots, i-T-1$ we have 
     \begin{align*}
         p\sum\limits_{\substack{0 \leq j < r-i+T \\ j \equiv a-i+T ~\mathrm{mod}~(p-1)}} \binom{r-l}{j}  \binom{j}{m} &\equiv p\sum\limits_{\substack{0 \leq j < s-i+T \\ j \equiv a-i+T ~\mathrm{mod}~(p-1)}} \binom{s-l}{j}  \binom{j}{m} \mod p^{t+1-v(m!)} 
     \end{align*}
     where we have used $\binom{r-l}{r-i+T}\binom{r-i+T}{m} \equiv \binom{s-l}{s-i+T}\binom{s-i+T}{m} \mod p^{t-v(m!)}$ (cf. Lemma~\ref{binomial coefficient under congruences} $(i)$). Again  by Corollary~\ref{cor: binomial sums under congruences 2}, for $m \not \equiv a-i+T \mod (p-1)$ we have
     \begin{align*}
         \sum\limits_{\substack{0 \leq j < r-i+T \\ j \equiv a-i+T ~\mathrm{mod}~(p-1)}} \binom{r-i+T}{j}  \binom{j}{m} \equiv \sum\limits_{\substack{0 \leq j < s-i+T \\ j \equiv a-i+T ~\mathrm{mod}~(p-1)}} \binom{s-i+T}{j} \binom{j}{m} \mod p^{t+1-v(m!)}
     \end{align*}
     where we have used $\delta_{p-1,[a-i+T-m]}=0$. From the above congruences it is enough to show $(ii)$ holds with $r$  replaced by $s$. Now part $(ii)$ follows  by a similar argument as in  Lemma~\ref{choice beta dcbd} using \eqref{eq: congruences for binomial sums for s a>2i} and \eqref{eq:matrix invertible GV det ad a>2i}.

     We now prove $(iii)$. Again by Corollary~\ref{cor: binomial sums under congruences 2}, for $l=0,\ldots, i-T-1$  we have 
     \begin{align*}
         p\sum\limits_{\substack{a-i+T < j < r-i+T \\ j \equiv a-i+T ~\mathrm{mod}~(p-1)}} \binom{r-l}{j}  \binom{j}{m} &\equiv p\sum\limits_{\substack{a-i+T < j < s-i+T \\ j \equiv a-i+T ~\mathrm{mod}~(p-1)}} \binom{s-l}{j}  \binom{j}{m} \mod p^{t+1-v(m!)} .
     \end{align*}
     Further,  by Corollary~\ref{cor: binomial sums under congruences 2}, for  $m  \equiv a-i+T \mod (p-1)$ we have
     \begin{align*}
         \sum\limits_{\substack{m < j < r-i+T \\ j \equiv a-i+T ~\mathrm{mod}~(p-1)}} \binom{r-i+T}{j}  \binom{j}{m} \equiv \sum\limits_{\substack{m < j < s-i+T \\ j \equiv a-i+T ~\mathrm{mod}~(p-1)}} \binom{s-i+T}{j} \binom{j}{m} \mod p^{t+1-v(m!)},  
     \end{align*}
     where we have used $\binom{j}{m} = 0$ for $j<m$ and $\delta_{p-1,[a-i+T-m]}=1$. Now $(iii)$ follows by a similar argument as used in the proof of part $(ii)$.

     We now prove $(iv)$. By Lemma~\ref{binomial coefficient under congruences} $(i)$ we have $\binom{r-l}{a-i+T} \equiv \binom{s-l}{a-i+T} \mod p^{t}$. Thus it is enough to show $(iv)$  holds when $r$ is replaced by $s$. It follows from \eqref{eq:matrix invertible GV det ad a>2i} that there exist $d_0,\ldots,d_{i-n-1} \in \Z_p$ such that 
     \begin{align*}
     \sum_{m=0}^{i-n-1} d_{m} \binom{a-i+T+k(p-1)}{m} =  
     \begin{cases}
         1 & \text{ if } k=0,\\
         0 & \text{ if } 1 \leq k \leq i-n-1.
     \end{cases}
     \end{align*}
     Thus we have
     \begin{align*}
         p & \sum_{l=0}^{i-T-1} \beta_{l} \binom{s-l}{a-i+T} + \beta_{i-T} \binom{s-i+T}{a-i+T} \\
         &=   \sum\limits_{\substack{0 \leq j < s-i+T \\ j \equiv a-i+T ~\mathrm{mod}~(p-1)}} \sum_{m=0}^{i-n-1} d_{m} \left( \sum_{l=0}^{i-T-1} p\beta_{l}  \binom{s-l}{j}   + \beta_{i-T}  \binom{s-i+T}{j}  \right)    \binom{j}{m} \\
         &\equiv 0 \mod p^{t},
     \end{align*}
     where the last step follows from \eqref{eq: congruences for binomial sums for s a>2i}. 

     We now prove $(v)$. Note that by Lemma~\ref{binomial coefficient under congruences} $(i)$ we have
     \[
     \binom{r-l}{r-i+T-(p-1)} = \binom{r-l}{i-T-l+p-1} \equiv \binom{s-l}{i-T-l+p-1}  \equiv \binom{s-l}{s-i-T-(p-1)}\mod p^{t-1}.
     \]
     Now using a similar argument as in the proof of $(iv)$ we obtain $(v)$.
    \end{proof}
    \begin{lemma}\label{lem: choice gamma asd a>2i}
             Let $1 \leq i \leq p-2$. Let $r \geq i(p+1)+p $, $r \equiv a \mod{(p-1)}$ with $ 1 \leq a \leq p-1 $ and $r \equiv a-i+n \mod p$ with $1 \leq n < i$. Let $ s =  a-i+n+(i-n)p$ and $v(r-s) =t $. Fix an integer $1 \leq T \leq n$. Assume $a>2i-2T$ (this can be relaxed to $a>2i-n-T$).  There exists $p$-adic integers $\gamma_0, \ldots, \gamma_{i-n} \in \mathbb{Z}_p$ satisfying 	
	        \begin{enumerate}
		        \item[$(i)$] For $m=0,\ldots,2p-1$, we have
		        $$
		           \sum\limits_{\substack{a - i + T < j < r-i+T \\ j \equiv a-i+T \mod(p-1)}}  \sum\limits_{l = 0}^{i-n-1} \gamma_l \binom{r-l}{j}  \binom{j}{m}  \equiv 0 \mod p^{t-v(m!)}.
		        $$
		        \item[$(ii)$] $ \sum\limits_{l = 0}^{i-n-1} \gamma_l  \binom{r-l}{a-i+T} \equiv 1  \mod p^{t}$.		
	        \end{enumerate}
        \end{lemma}
        \begin{proof}
            The proof is similar to Lemma~\ref{choice beta sdcbd}. Using $r\equiv s \mod p^{t}$, Lemma~\ref{binomial coefficient under congruences} $(i)$ and Corollary~\ref{cor: binomial sums under congruences 1}, it is enough to show that the lemma holds when $r=s$. We first show the existence of $\gamma_{0}, \ldots, \gamma_{i-n} \in \mathbb{Z}_{p}$  satisfying
		\begin{align}
			\sum\limits_{l=0}^{i-n-1} \gamma_{l} \sum\limits_{ \substack{a-i+T <   j < s-i+T \\  j \equiv a-i+T ~\mathrm{mod}~ (p-1)} }\binom{s-l}{j} \binom{j}{m} &= 0 ~\text{ for } m=0,\ldots, i-n-2 \label{red eq 1 sdcad a>2i}\\
                \sum\limits_{l=0}^{i-n-1} \gamma_{l} \binom{s-l}{a-i+T} &= 1. \label{red eq 3 sdcad a>2i}
		\end{align}
        Note that if $n=i-1$, then \eqref{red eq 1 sdcad a>2i} is vacuously true. Then the above system of equations  can be written as 
		\begin{align*}
			A \begin{bmatrix}
				\gamma_{0} \\ \vdots \\ \gamma_{i-n-1}
			\end{bmatrix} = \begin{bmatrix}
				1 \\ 0 \\ \vdots \\ 0 
			\end{bmatrix}.
		\end{align*}
		where 
		\begin{align*}
			A = \left[ \begin{array}{c}
                    \left( \binom{s-l}{a-i+T} \right)_{l=0, \ldots, i-n} \\ \hline
				\left(\sum\limits_{ \substack{a-i+T <   j < s-i+T \\  j \equiv a-i+T ~\mathrm{mod}~ (p-1)} }\binom{s-l}{j} \binom{j}{m}\right)_{\substack{m=0,\ldots, i-n-2\\ l=0, \ldots, i-n}} 
			\end{array}\right]. 
		\end{align*}
		By a similar argument as in  Lemma~\ref{choice beta sdcbd}, we have $\det(A) \in \mathbb{Z}_p^\times$. This shows that \eqref{red eq 1 sdcad a>2i} and \eqref{red eq 3 sdcad a>2i} have a solution in $\mathbb{Z}_p$.  As  observed towards the end of the proof of Lemma~\ref{choice beta dcbd}, it can be shown that \eqref{red eq 1 sdcad a>2i} implies $(i)$ and $(ii)$.
        \end{proof}
        Under the assumptions of the previous lemma, we show that certain functions taking non-integral values vanish in $\bar{\Theta}_{k,a_p}$.
        \begin{lemma}\label{lem: killing non-integral polynomial a>2i-2T}
            Let $1 \leq i \leq p-2$. Let $r \geq i(p+1)+p $, $r \equiv a \mod{(p-1)}$ with $ 1 \leq a \leq p-1 $ and $r \equiv a-i+n \mod p$ with $1 \leq n < i$. Let $ s =  a-i+n+(i-n)p$ and $v(r-s) =t $. Fix an integer $t \leq T \leq n$. If $i-T<a-i+T \leq i+1$, $2i-a+2 \leq t+T$ and $2i+3 \leq a+p-1$. Then for every $m \geq 1$ and $\lambda \in I_{m-1}$,  there exists $f \in \mathrm{ind}_{KZ}^{G}(\mathrm{Sym}^r\bar{\Q}^2_p)$ such that
            \begin{align*}
                (T-a_p)f \equiv \left[ g^{0}_{m,\lambda}, \frac{p^{a-i+T}}{p a_{p} } X^{r-(a-i+T)}Y^{a-i+T}\right] \mod p.
            \end{align*}
        \end{lemma}
        \begin{proof}
            Since $\gl_2(\Q_p)$ acts transitively on the vertices of the Bruhat-Tits tree, it is enough to show that the above conclusion holds when $m=0$.
            Let $\gamma_0, \ldots, \gamma_{i-n-1}$ be $p$-adic integers as in Lemma~\ref{lem: choice gamma asd a>2i}. Let 
	        \begin{align*}		
		          f_1 &= \sum_{\lambda \in \mathbb{F}_p^\times} \left[ g^{0}_{1,[\lambda]}, \sum_{l = 0}^{i-n-1} \frac{p^{a-i+T}}{p^{l+1} a_{p}} \gamma_l [\lambda]^{l-i+T} (-\theta)^{i+l+2}X^{-i-2} Y^{r-(i+l+2)(p+1) +i+2} \right] \\
		        & \qquad \qquad - \left[ g^{0}_{1,0}, \left( \sum_{l = 0}^{i-n-1} (p-1) \gamma_l \binom{r-l}{r-i+T}  \right) \frac{p^{a-i+T}}{p^{i-T+1} a_{p} } (-\theta)^{ i+1} X^{-T-1} Y^{r - (i+1) (p+1)+T+1} \right]
	        \end{align*}
	        Clearly $T^{+}f_1$ vanishes modulo $p$ by Lemma~\ref{theta and T plus}. Note that $-a_p f_1$ vanishes if $a-i+T$ is strictly greater than $i-n$ and $i-T+1$. As $t< T \leq  n$ and $a-i+T\geq i-t+2$ it can be checked that these inequalities hold. Applying Lemma~\ref{lem:choice of beta}, there exist $\alpha_j'$ such that
	        \begin{enumerate}
		        \item[$(1)$] $ \alpha_j' \equiv \sum\limits_{l = 0}^{i-n-1} \gamma_l \binom{r - l}{j} \mod{p^{t}}$ for $a - i + T < j < r - i + T$ with $j \equiv a - i + T \mod (p-1)$,
		        \item[$(2)$] For  $m = 0, \ldots, \min\{ p-1,i+t+1\}$, we have
		          \[
		              \sum_{\substack{a - i + T < j < r-i + T \\ j \equiv a - i + T \pmod{p-1}}} \alpha_j' \binom{j}{m} \equiv 0 \mod p^{i+t+1-m}.
		        \]
	        \end{enumerate}
            Note that for $m=p,\ldots, i+t+1$ we have $t-v(m!) = t-1 \geq i+t+1-p \geq i+t+1-m$  as  $i \leq p-2$. Thus it follows from  Lemma~\ref{lem: choice gamma asd a>2i}~$(i)$ and the congruence condition $(1)$ that
            \begin{enumerate}
                \item [(2$'$)] For  $m = p, \ldots,i+t+1$, we have
		          \[
		              \sum_{\substack{a - i + T < j < r-i + T \\ j \equiv a - i + T \mod(p-1)}} \alpha_j' \binom{j}{m} \equiv 0 \mod p^{i+t+1-m}.
		        \]
	        \end{enumerate}
            Let
            \begin{align*}
                f_0 = \left[\mathrm{Id}, \frac{(p-1)p^{a-i+T}}{p a_{p}^2 } \sum_{\substack{a - i + T < j < r - i + T \\ j \equiv a - i + T \mod{(p-1)}}} \alpha'_j X^{r-j} Y^j \right]
            \end{align*}
             From  $(2)$, $(2')$ and $a+T+t \leq 2i+2$, we see that $T^{+}f_0$ vanishes modulo $p$ at the vertices $g^{0}_{1,[\mu]}$ for $\mu \neq 0$. Since $i<a$, $i< p-1$ and $a-2i+2T \geq  1$ we get $2(a-i+T) + p-1 \geq 2i+3$. Thus $T^{+}f_0$ also vanishes modulo $p$ at the vertex $g^{0}_{1,0}$. So $T^{+}f_0$ vanishes modulo $p$. As $2i+3 \leq a+p-1$, we obtain $T^{-}f_0$ vanishes modulo $p$.  Observe that
	        \begin{align*}
		        T^{-}f_1 - a_pf_0  &\equiv \left[ \mathrm{Id}, \frac{(p-1)p^{a-i+T}}{p a_{p} } \sum_{l = 0}^{i-n-1}  \gamma_l \binom{r-l}{a-i+T} X^{r-(a-i+T)}Y^{a-i+T}\right] \\
                & \qquad +\left[ \mathrm{Id}, \frac{(p-1)p^{a-i+T}}{p a_{p} } \sum_{\substack{a - i + T < j < r-i + T \\ j \equiv a - i + T \mod(p-1)}} \left(\sum_{l = 0}^{i-n-1}  \gamma_l \binom{r-l}{j} - \alpha_j'\right) X^{r-j}Y^{j} \right] \mod p. 
	        \end{align*}
	        From $(1)$ and $2i-a+2 \leq t+T$ it follows that the second term vanishes modulo $p$. Using this inequality and Lemma~\ref{lem: choice gamma asd a>2i} ($ii$) we get 
            \[
                  (T - a_p)(f_1 + f_0 )  \equiv T^{-} f_1 - a_p f_0 \equiv  \left[ \mathrm{Id}, \frac{(p-1)p^{a-i+T}}{p a_{p} } X^{r-(a-i+T)}Y^{a-i+T}\right] \mod p.
            \]
            This finishes the proof.
        \end{proof}
        
        The following theorem allows us to treat the triangular region in the hybrid conjecture (see Figure~\ref{fig: hybrid a odd}, \ref{fig: hybrid a even}.
        
        \begin{theorem}\label{thm: above super diagonal a large}
        	Let  $a_p \in \bar{\Q}_p$ with $ v(a_{p})  \in (i,i+1) $ and $1 \leq i \leq p-2$. Let $r \geq i(p+1)+p $, $r \equiv a \mod{(p-1)}$ with $ 1 \leq a \leq p-1 $ and $r \equiv a-i+n \mod p$ with $1 \leq n < i$. Let $ s =  a-i+n+(i-n)p$ and $v(r-s) =t $. Fix an integer $1 \leq T \leq n$. Assume that the following holds:
            \begin{enumerate}
                \item[$(i)$] $ 1 \leq t < T \leq n < i < a$
                \item[$(ii)$] $a > 2i-2T$
                \item[$(iii)$] if $a-i+T \leq i+1$, then we further assume that
                \begin{enumerate}
                    \item[$(a)$] $2i-a+2 \leq T+t$
                    \item[$(b)$] $i = p-2 \Longrightarrow a \neq p-1$.
                \end{enumerate}
            \end{enumerate}
            Then the image of $ \mathrm{ind}_{KZ}^{G}(V_{r}^{(i-T)}) $ is the same as the image of $ \mathrm{ind}_{KZ}^{G}(V_{r}^{(i-T+1)}) $ in $ \bar{\Theta}_{k,a_{p}} $.
        \end{theorem}
        \begin{proof}
            Let $\beta_l$ be integers as chosen in Lemma~\ref{lem: choice beta asd a>2i}. Consider 
            \begin{align*}
	            f_2 &= \sum_{\lambda \in \mathbb{F}_p^\times} \left[ g^{0}_{2,p[\lambda]}, \sum_{l = 0}^{i-T-1} \frac{\beta_l}{p^{l+t}} [\lambda]^{l-i+T} (-\theta)^{l+t+1} X^{-t-1} Y^{r - (l+t+1)(p+1)+t+1} \right] \\
                &\qquad \qquad + \sum_{\lambda \in \mathbb{F}_p^\times} \left[ g^{0}_{2,p[\lambda]},\frac{\beta_{i-T}}{p^{i-T+t+1}} (-\theta)^{i -T+t+ 2} X^{-t-2} Y^{r - (i-T+t+2)(p+1)+t+2}\right] \\
	            & \qquad \qquad  + \left[ g^{0}_{2,0}, -\left( \sum_{l = 0}^{i-T-1} p \beta_l \binom{r-l}{r-i+T} + \beta_{i-T} \right) \frac{p-1}{p^{i-T+t+1}}(-\theta)^{i+1} X^{-T-1} Y^{r-(i+1)(p+1)+T+1} \right].
	        \end{align*}
            Note that $T^{+} f_2$  vanishes since $ t < T$. Also	$-a_p f_2$ vanishes since  $i - T + t + 1 \leq i <  v(a_p)$. Observe that
            \begin{align*}
	            T^{-} f_2 \equiv \left[ g^{0}_{1,0}, \frac{p-1}{p^{t+1}} \sum_{\substack{a - i + T \leq j < r- i + T \\ j \equiv a-i+T \mod(p-1)}} \left( p \sum_{l=0}^{i-T-1} \beta_l \binom{r-l}{j} +\beta_{i-T} \binom{r-i+T}{j}\right) X^{r - j} Y^j \right].
	        \end{align*}
            By  Lemma~\ref{lem: choice beta asd a>2i} $(i)-(ii)$ and Lemma~\ref{lem:choice of beta}, there exist $\alpha_j$ such that
	          \begin{enumerate}
		        \item[(1)] \[ \alpha_j \equiv p \sum_{l = 0}^{i-T-1} \beta_l \binom{r-l}{j} + \beta_{i-T} \binom{r - i + T}{j} \mod{p^{t+1}}\]
		        for $a - i + T \le j < r - i + T$ with $j \equiv a - i + T \mod (p-1)$, and
                
		          \item[$(2a)$]  \[ \sum_{\substack{a - i + T \leq j < r-i + T \\ j \equiv a - i + T \mod(p-1)}} \alpha_j \binom{j}{m} \equiv 0 \mod{p^{i+t+2-m}} \]
                  for  $m = 0, \ldots, \min\{i +t + 2,p-1\}$ and $m \neq a-i+T$.
                  \item[$(2b)$] 
                  Suppose $m=a-i+T$ and $a-i+T \leq i+t+2$. Adding $j=a-i+T$ term on both sides of Lemma~\ref{lem: choice beta asd a>2i} $(iii)$ and Lemma~\ref{lem:choice of beta},  we obtain
                  \[
                \sum_{\substack{a - i + T \leq j < r-i + T \\ j \equiv a - i + T \mod(p-1)}} \alpha_j \binom{j}{a-i+T} \equiv  p \sum_{l = 0}^{i-T-1} \beta_l \binom{r-l}{a-i+T} + \beta_{i-T} \binom{r - i + T}{a-i+T} \mod{p^{i+t+2-(a-i+T)}}. 
                \]
		      \item[($2'$)]
              From the hypothesis, we have $a-i+T < a \leq p-1$ and $i+t+2 < p-1+T+1 \leq p-1+a-i+T$. Hence $p,\ldots, i+t+2 \not \equiv a-i+T \mod (p-1)$. Thus, from the congruence condition $(1)$ and Lemma~\ref{lem: choice beta asd a>2i} $(ii)$, we also have 
        	 \[\sum_{\substack{a - i + T \leq j < r-i + T \\ j \equiv a - i + T \mod(p-1)}} \alpha_j \binom{j}{m} \equiv 0  \mod p^{t+1-v(m!)} \text{ for } m=p,\ldots, i+t+2.\] 
             Note that for $m=p,\ldots, i+t+2$, we have $t-v(m!) = t-1 \geq i+t+1-p \geq i+t+1-m$  as $i \leq p-2$.
        	\end{enumerate}
                  Further, from Lemma~\ref{lem: choice beta asd a>2i} $(iv)$, $(v)$ and the congruence condition $(1)$ we have
            \begin{enumerate}
		          \item[(3)] $ \alpha_{a-i+T} \equiv 0 \mod{p^{t}}$.
	              \item[(4)] $ \alpha_{r-i+T-(p-1)} \equiv 0 \mod{p^{t-1}}$.
	        \end{enumerate}
            Let 
            \begin{align*}
	            f_1 = \left[g^{0}_{1,0}, \frac{p-1}{ p^{t+1} a_p} \sum_{\substack{a - i + T \leq j < r - i + T \\ j \equiv a - i + T \mod{(p-1)}}} \alpha_j X^{r-j} Y^j \right].
	        \end{align*}
	        Since $\alpha_{r-i+T-(p-1)} \equiv 0 \mod{p^{t}}$ and $v(pa_p) < i+2 \leq p \leq i-T+p-1$, it follows that $T^{-} f_1$ vanishes. Thus
	          \begin{align*}
	            (T - a_p)(f_2 + f_1) \equiv T^{+} f_1 + T^{-} f_2 - a_p f_1 \equiv T^{+}f_1 +  \left[g^{0}_{1,0}, F(X,Y) \right] \mod p,
	        \end{align*}
	        where
	          \begin{align*}
		          F(X,Y) = \frac{p-1}{p^{t+1}}\sum_{\substack{a - i + T \leq j < r- i + T \\ j \equiv a - i + T \mod(p-1)}} \left( p \sum_{l=0}^{i-T-1} \beta_l \binom{r-l}{j}+ \beta_{i-T} \binom{r-i+T}{j}- \alpha_j \right) X^{r-j} Y^j.
	          \end{align*}
            From $(2a)$, Lemma~\ref{lem: choice beta asd a>2i} $(i)$ and the assumption $a-i+T > i-T$, it follows that $F(X,Y) \in V_{r}^{(i-T)}$. Again, from the assumption $a-i+T > i-T$ and  \cite[Lemma 2.13]{GR19}, we have 
	        \[
                F(X,Y) \equiv (p-1) \theta^{i-T} X^{r-(i-T)(p+1)-(a-2i+2T)}Y^{a-2i+2T} \mod V_{r}^{(i-T+1)}.
	          \]
            By Lemma~\ref{Glover-Brueil map image} $(ii)$, it follows that the image of $F(X,Y)$ under the projection $V_{r}^{(i-T)}/V_r^{(i-T+1)} \twoheadrightarrow V_{p-1-(a-2i+2T)} \otimes D^{(a-i+T)}$ is non-zero. From the assumption $a-i+T > i-T$ it follows that \eqref{exact seq. singular quotients} is non-split for $m=i-T$. Hence, $F(X,Y)$ generates the quotient $V_{r}^{(i-T)}/V_r^{(i-T+1)}$. So we are done if $T^{+}f_1$ can be killed. 

            We now show $T^{+}f_1$ can be eliminated. It follows from $(2a)$, $(2b)$ and $(2')$ that  
	          \begin{align*}
		          T^{+}f_1 \equiv & \sum_{\lambda \in \mathbb{F}_p^\times} \left[ g^{0}_{2,p[\lambda]}, \frac{(p-1)p^{a-i+T}}{p a_p } \frac{p\sum\limits_{l=0}^{i-T-1} \beta_l \binom{r-l}{a-i+T}+ \beta_{i-T}\binom{r-i+T}{a-i+T}}{p^t} X^{r-(a-i+T)}Y^{a-i+T} \right]\\
                  & \qquad \quad + \left[ g^{0}_{2,0}, \frac{(p-1)p^{a-i+T}}{p a_p } \frac{\alpha_{a - i + T}}{p^t} X^{r-(a-i+T)}Y^{a-i+T} \right]  \mod p.
	        \end{align*}
            From $(3)$, we have $\alpha_{a-i+T}  \in p^t\mathbb{Z}_{p}$. By Lemma~\ref{lem: choice beta asd a>2i} $(iv)$, we have $p\sum\limits_{l=0}^{i-T-1} \beta_l \binom{r-l}{a-i+T}+ \beta_{i-T}\binom{r-i+T}{a-i+T} \in p^t \Z_p$.
            If $v(a_p)+1 < a-i+T$, then we are done. Assume $v(a_p)+1 > a-i+T$.  Since $i<a$ and $i< p-1$ we have $2i+2 \leq a+p-1$ and equality holds only if $a=i+1=p-1$. Thus, $2i+2 < a+p-1$ as $a\neq p-1$ when $i=p-2$. Hence, $T^{+}f_1$ can be eliminated using Lemma~\ref{lem: killing non-integral polynomial a>2i-2T}. Thus, in both cases $T^{+}f_1$ can be eliminated. Hence, we are done. 
        \end{proof}
        
        We now treat the columns
        $2i-n-a \leq T \leq i-a/2$ in the hybrid conjecture.
        The argument is similar to the one used to prove  Theorem~\ref{thm: hybrid ad}. 
        \begin{lemma}\label{choice beta hybrid green}
            Let  $r \equiv a \mod{(p-1)}$ with $ 1 \leq a \leq p-1 $ and $r \geq i(p+1)+p $ with $ v(a_{p})  \in (i,i+1) $. Let $ s =  a-i+n+(i-n)p$ and $v(r-s) =t $  and $ 1 \leq t \leq  T < n < i < a$. If  $2i-n-T \leq a \leq  2i-2T $, then there exist $\beta_0, \ldots,  \beta_{i-T} \in \Z_p$ and $\gamma \in \Z_p$ with $\beta_{a-i+T}=0$ such that 
        	    \begin{enumerate}
        	    	\item[$(i)$] $\sum\limits_{l=0}^{i-T} \beta_l  \sum\limits_{ \substack{i-T <  j < r-(a-i+T) \\  j \equiv i-T ~\mathrm{mod}~ (p-1)} }   \binom{r-l}{j} \binom{j}{m}  + \gamma p^t  \binom {i-T}{m}\equiv 
        	    	0  \mod p^{t+1} ~ \mathrm{for}~ m=0,\ldots, i-T-1$ 
        	    	\item[$(ii)$] $\sum\limits_{l=0}^{i-T} \beta_l  \sum\limits_{ \substack{i-T <  j < r-(a-i+T) \\  j \equiv i-T ~\mathrm{mod}~ (p-1)} }   \binom{r-l}{j} \binom{j}{i-T}   \equiv 
        	    	 p^t  \mod p^{t+1}$ 
        	    	\item[$(iii)$] $\sum\limits_{l=0}^{i-T} \beta_l \sum\limits_{ \substack{i-T <  j < r-(a-i+T) \\  j \equiv i-T ~\mathrm{mod}~ (p-1)} }   \binom{r-l}{j} \binom{j}{m} \equiv  0 \mod p^{t-v(m!)} ~  \mathrm{for}~ m=i-T+1,\ldots, i+t$.
        	    \end{enumerate}
        \end{lemma}
        \begin{proof}
            First we  prove  $(i)$ and $(ii)$. We now compute some binomial sums. By Corollary~\ref{cor: binomial sums under congruences 1}, for $m=0,\ldots, i-T$ we have
        	    \begin{align*}
        	    	\sum\limits_{ \substack{i-T \leq j \leq r-(a-i+T) \\  j \equiv i-T ~\mathrm{mod}~ (p-1)} }   \binom{r-l}{j} \binom{j}{m} &\equiv \left( \binom{r-l}{m} - \binom{s-l}{m}\right)\left(\binom{[a-l-m]}{[i-T-m]}+\delta_{[i-T-m],p-1} \right) \\ &
        	    	\qquad \qquad +\sum_{\substack{ i-T \leq j \leq s-(a-i+T) \\ j \equiv i-T ~\mathrm{mod}~(p-1)}} \binom{s-l}{j} \binom{j}{m}   \mod p^{t+1}. 
        	    \end{align*}
                We now claim that for $m=0,\ldots,i-T$, we have
                \begin{align*}
                    \binom{[a-l-m]}{[i-T-m]}+\delta_{[i-T-m],p-1} \equiv 
                    \begin{cases}
                        \binom{a-l-m}{i-T-m} ~&\mathrm{if}~l=0,\ldots,a-i+T-1, \\
                        \binom{p-1+a-l-m}{i-T-m} ~&\mathrm{if}~l=a-i+T+1,\ldots,i-T,
                    \end{cases}
                    \mod p.
                \end{align*}
        	      We prove the claim by treating various cases. First, suppose that $m=0,\ldots,i-T-1$. Then $[i-T-m]=i-T-m<p-1$ and $\delta_{[i-T-m],p-1}=0$. If $0\leq l \leq a-i+T-1$, then  $l+m \leq a-2 < a$ and $[a-l-m]=a-l-m$, so the above claim follows for $l=0,\ldots,a-i+T-1$. If $l=a-i+T+1,\ldots,i-T$ and $l+m \geq a$, then $[a-l-m]=p-1+a-l-m$ and the claim is true. If $l=a-i+T+1,\ldots,i-T$ and $l+m < a$, then $ \binom{a-l-m}{i-T-m} =0$. Furthermore, by Lucas' theorem $\binom{p-1+a-l-m}{i-T-m} \equiv \binom{a-l-m-1}{i-T-m}  \equiv 0 \mod p$. This proves the claim for $m=0,\ldots,i-T-1$. Now assume $m=i-T$. Then $[i-T-m]=p-1$ and $\delta_{[i-T-m],p-1}=1$. As $l \neq a-i+T$, we have $[a-l-m] \neq p-1$ and $\binom{[a-l-m]}{[i-T-m]} = \binom{[a-l-m]}{p-1}=0$. Now the case $m=i-T$ follows. This proves the claim in all possible cases. 
                  Thus for $m=0,\ldots, i-T$ and $l=0,\ldots,a-i+T-1$ we have
        	    \begin{align*}
        	    	\sum\limits_{ \substack{i-T < j < r-(a-i+T) \\  j \equiv i-T ~\mathrm{mod}~ (p-1)} }   \binom{r-l}{j} \binom{j}{m} &\equiv \left( \binom{r-l}{m} - \binom{s-l}{m}\right) \binom{a-l-m}{i-T-m} \\ &
        	    	\qquad  +\sum_{\substack{ i-T < j < s-(a-i+T) \\ j \equiv i-T ~\mathrm{mod}~(p-1)}} \binom{s-l}{j} \binom{j}{m}   \\
                    & \qquad + \left( \binom{s-l}{i-T} - \binom{r-l}{i-T}\right) \binom{i-T}{m} \\
                    & \qquad + \binom{s-l}{s-(a-i+T)} \binom{s-(a-i+T)}{m} \\
                    &\qquad - \binom{r-l}{r-(a-i+T)} \binom{r-(a-i+T)}{m} \mod p^{t+1}. 
        	    \end{align*}
                Also,  for $m=0,\ldots, i-T$ and $l=a-i+T+1,\ldots,i-T$ we have
                \begin{align*}
        	    	\sum\limits_{ \substack{i-T < j < r-(a-i+T) \\  j \equiv i-T ~\mathrm{mod}~ (p-1)} }   \binom{r-l}{j} \binom{j}{m} &\equiv \left( \binom{r-l}{m} - \binom{s-l}{m}\right) \binom{p-1+a-l-m}{i-T-m} \\ &
        	    	\qquad  +\sum_{\substack{ i-T < j < s-(a-i+T) \\ j \equiv i-T ~\mathrm{mod}~(p-1)}} \binom{s-l}{j} \binom{j}{m}   \\
                    & \qquad + \left( \binom{s-l}{i-T} - \binom{r-l}{i-T}\right) \binom{i-T}{m}  \mod p^{t+1}. 
        	    \end{align*}
                Note that every $0< j < s-(a-i+T)$ with $j \equiv i-T$ mod $(p-1)$ can be expressed as $i-T+k(p-1)$ for some $1\leq k \leq i-n-1$. Consider the following matrix
                \begin{equation}\label{eq: A green}
	                   A = \left[ \begin{array}{c|c|c}
                     \left\{\binom{r - l}{m} - \binom{s - l}{m}\right\} \binom{a - l - m}{i - T - m}+ &   \left\{\binom{r - l}{m} - \binom{s - l}{m}\right\} \binom{p-1+a - l - m}{i - T - m} + & {} \\[7pt]
		             \sum\limits_{k=1}^{i-n-1} \binom{s - l}{i-T+k(p-1)} \binom{i-T+k(p-1)}{m}+ & \sum\limits_{k=1}^{i-n-1} \binom{s - l}{i-T+k(p-1)} \binom{i-T+k(p-1)}{m}+  & p^t \binom{i-T}{m}-  \\[7pt]
                     \left( \binom{s-l}{i-T} - \binom{r-l}{i-T}\right) \binom{i-T}{m}+ & \left( \binom{s-l}{i-T} - \binom{r-l}{i-T}\right) \binom{i-T}{m}  & p^t \delta_{i-T,m} \\[7pt] 
		             \binom{s - l}{s -(a-i + T)} \binom{s -(a-i + T)}{m}- \binom{r-l}{r-(a-i+T)} \binom{r-(a-i + T)}{m} &{} & {}
	              \end{array}\right]_{m=0,\ldots,i-T},
            \end{equation}
            where the range of $l$ in the left block is $0,\ldots,a-i+T-1$, range of $l$ in the middle block is $a-i+T+1,\ldots,i-T$  and the range of $l$ in the rightmost column is $i-T+1$. Note that the above matrix is similar to the one appearing in the proof of Lemma~\ref{choice beta hybrid ad1}.

             To solve the congruences $(i)$ and $(ii)$, it is enough to show that the following congruence
            \[
                A \begin{bmatrix} \beta_0 \\ \vdots \\ \beta_{a-i+T-1} \\ \beta_{a-i+T+1} \\ \vdots \\ \beta_{i-T} \\ \gamma \end{bmatrix} \equiv  \begin{bmatrix} 0\\ \vdots \\ 0 \\ p^{t} \end{bmatrix} \mod p^{t+1} \mathbb{Z}_{p}
            \] 
            has a solution in $\mathbb{Z}_{p}$. To solve the above congruence, it is enough to show the following equation
            \begin{equation}\label{eq: linear eq. in matrix form green}
                A \begin{bmatrix}\beta_0 \\ \vdots \\ \beta_{a-i+T-1} \\ \beta_{a-i+T+1} \\ \vdots \\ \beta_{i-T} \\ \gamma \end{bmatrix} = \begin{bmatrix} 0\\ \vdots \\ \vdots \\ 0 \\ p^{t} \end{bmatrix}
            \end{equation}
            has a solution over $\mathbb{Q}_p$ with $\beta_l, \gamma \in \mathbb{Z}_{p}$. To show this, we use Cramer's rule.

            We show $p^{(n-T+2)t} \parallel \det(A)$. To achieve this, we will  perform a sequence of row operations on $A$ so that the $m=i-n,\ldots,i-T-1, i-T$ rows are multiples of $(r-s)$ and the last column  becomes $(p^t,0,\ldots,0 )$. 
            Apply the following row operations (these are the same as the ones used  in Lemma~\ref{choice beta hybrid ad1}): 
            \begin{alignat*}{4}
                R_{i-T-1} &\rightarrow R_{i-T-1}~- ~\frac{i-T-(i-T-2)}{i-T-1} R_{i-T-2} \\
                  & ~~\vdots  \\
                R_{m} &\rightarrow ~ ~~R_{m}~ - ~\frac{i-T-(m-1)}{m} R_{m-1} \\
                 & ~~\vdots  \\
                R_{1} &\rightarrow ~~ R_{1}~- ~\frac{i-T}{1} R_{0}.
            \end{alignat*}   
            Then expand using the last column to obtain $\det(A) = (-1)^{i-T}p^t\times \det(A_1)$, where
            \begin{align*}
                A_1 = \left[ \begin{array}{c|c}
                     \left\{\binom{r - l}{m-1} \frac{r-a}{m} - \binom{s - l}{m-1} \frac{s-a}{m} \right\} \binom{a - l - m}{i - T - m} &   \left\{\binom{r - l}{m-1} \frac{r-a+1}{m}- \binom{s - l}{m-1} \frac{s-a+1}{m}\right\} \binom{p-1+a - l - m}{i - T - m}  \\ \\
		            + \sum\limits_{k=1}^{i-n-1} \binom{s - l}{i-T+k(p-1)} \binom{i-T+k(p-1)}{m-1}\frac{k(p-1)}{m} &  +  \sum\limits_{k=1}^{i-n-1} \binom{s - l}{i-T+k(p-1)} \binom{i-T+k(p-1)}{m-1}\frac{k(p-1)}{m}  \\ \\
		            + \binom{s - l}{s - (a-i + T)} \binom{s - (a-i + T)}{m-1} \frac{s-a}{m}- \binom{r - l}{r - (a-i + T)} \binom{r - (a-i + T)}{m-1} \frac{r-a}{m} & +O(p(r-s))\\ \\+O(p(r-s)) & \\ \\ \hline \\    \sum\limits_{k=1}^{i-n-1} \binom{s - l}{i-T+k(p-1)} \binom{i-T+k(p-1)}{i-T}  & \sum\limits_{k=1}^{i-n-1} \binom{s - l}{i-T+k(p-1)} \binom{i-T+k(p-1)}{i-T}   \\
                    + \binom{s - l}{s - (a-i + T)} \binom{s - (a-i + T)}{i-T}- \binom{r - l}{r - (a-i + T)} \binom{r - (a-i + T)}{i-T} & 
	              \end{array}\right],
            \end{align*}
            where the range of $l$ in the left blocks is $0,\ldots,a-i+T-1$ and the range of $l$ in the right blocks is $a-i+T+1,\ldots,i-T$ , and where the range of $m$ in the upper blocks is $1,\ldots,i-T-1$ and the range of $m$ in the lower blocks is $i-T$. The proof of the above statement is similar to the one given in Lemma~\ref{choice beta hybrid ad1} using  \eqref{identity 2 adhc} and \eqref{identity 1 adhc} for $k'$ there equal to zero. 
            
            We now  compute $\det(A_1)$.  We define matrices $A(1), \ldots, A(i-n)$ recursively, by setting  $A(1)=A_1$ and by letting $A(k'+1)$ to be the matrix obtained from $A(k')$ by performing the following row operations:
            \begin{alignat*}{4}
            	R_{i-T-1} &\rightarrow R_{i-T-1}~- ~\frac{i-T-(i-T-2)+k'p}{i-T-1} R_{i-T-2} \\
            	& ~~\vdots  \\
            	R_{m} &\rightarrow ~ ~~R_{m}~ - ~\frac{i-T-(m-1)+k'p}{m} R_{m-1} \\
            	& ~~\vdots  \\
            	R_{k'+1} &\rightarrow ~~ R_{k'+1}~ - ~\frac{i-T-k'+k'p}{k'+1} R_{k'}
            \end{alignat*} 
             for $k'=1,\ldots,i-n-1$. We claim that for $k'=1,\ldots,i-n$, we have
            \begin{equation}\label{eq: matrix formula recursion green}
            	\resizebox{0.999\hsize}{!}{%
            		$A(k')=\left[\begin{array}{c|c}
            			 \sum\limits_{k=m}^{i-n-1} (p-1)^m \binom{s-l}{i-T+k(p-1)} \binom{k}{m} +  O(r-s)
            			&\sum\limits_{k=m}^{i-n-1} (p-1)^m\binom{s-l}{i-T+k(p-1)} \binom{k}{m} +  O(r-s)\\
            			\hline \\
            			 \left( \binom{r-l}{m-k'} \frac{(r-a+k'-1)_{k'}}{(m)_{k'}}- \binom{s-l}{m-k'} \frac{(s-a+k'-1)_{k'}}{(m)_{k'}} \right) \binom{a-l-m}{i-T-m}
            			& \left( \binom{r-l}{m-k'} \frac{(r-a+k')_{k'}}{(m)_{k'}}- \binom{s-l}{m-k'}\frac{(s-a+k')_{k'}}{(m)_{k'}}\right) \binom{p-1+a-l-m}{i-T-m} \\  \\
            			+ (p-1)^{k'} \sum\limits_{k=k'}^{i-n-1} \binom{s-l}{i-T+k(p-1)} \binom{i-T+k(p-1)}{m-k'} \frac{(k)_{k'}}{(m)_{k'}}  
            			&  + (p-1)^{k'}\sum\limits_{k=k'}^{i-n-1} \binom{s-l}{i-T+k(p-1)} \binom{i-T+k(p-1)}{m-k'} \frac{(k)_{k'}}{(m)_{k'}}   \\  \\
            			+  \binom{s-l}{s-(a-i+T)} \binom{s-(a-i+T)}{m-k'} \frac{\prod\limits_{k=0}^{k'-1}(s-a-k(p-1))}{(m)_{k'}} -  \binom{r-l}{r-(a-i+T)}\binom{r-(a-i+T)}{m-k'} \frac{\prod\limits_{k=0}^{k'-1}(r-a-k(p-1))}{(m)_{k'}} 
            			&+O(p(r-s))\\ \\
            			+ O(p(r-s)) & \\ \hline \\
            			 \sum\limits_{k=1}^{i-n-1} \binom{s - l}{i-T+k(p-1)} \binom{i-T+k(p-1)}{i-T}  & \sum\limits_{k=1}^{i-n-1} \binom{s - l}{i-T+k(p-1)} \binom{i-T+k(p-1)}{i-T} \\
            			 + \binom{s-l}{s-(a-i+T)}\binom{s-(a-i+T)}{i-T}-\binom{r-l}{r-(a-i+T)}\binom{r-(a-i+T)}{i-T}  &
            		\end{array} \right],$
            	} %
            \end{equation}
            where the range of $m$ in the upper blocks  is $1,\ldots,k'-1$, the range of $m$ in the middle blocks is $k',\ldots,i-T-1$ and the range of $m$ in the lower blocks is $ i-T $, and the range of $l$ in the left blocks (resp. right blocks) is $0,\ldots,a-i+T-1$ (resp. $a-i+T+1, \ldots, i-T$). This can be proved by recursion as in Lemma~\ref{lem: choice beta asd a>2i} using the identities \eqref{identity 2 adhc} and \eqref{identity 1 adhc}.
            
            From the above claim, it follows that $\det(A_1) = \det(A(i-n))$ and $\det(A(i-n))$ equals
            \begin{equation*}
            	\resizebox{0.999\hsize}{!}{%
            		$\left\lvert\begin{array}{c|c}
            			 \sum\limits_{k=m}^{i-n-1} (p-1)^m \binom{s-l}{i-T+k(p-1)} \binom{k}{m} +  O(r-s)
            			&\sum\limits_{k=m}^{i-n-1} (p-1)^m\binom{s-l}{i-T+k(p-1)} \binom{k}{m} +  O(r-s) \\
            			\hline \\
            			 \left( \binom{r-l}{m-(i-n)} \frac{(r-a+i-n-1)_{i-n}}{(m)_{i-n}}- \binom{s-l}{m-(i-n)} \frac{(s-a+i-n-1)_{i-n}}{(m)_{i-n}} \right) \binom{a-l-m}{i-T-m}
            			& \left( \binom{r-l}{m-(i-n)} \frac{(r-a+i-n)_{i-n}}{(m)_{i-n}}- \binom{s-l}{m-(i-n)} \frac{(s-a+i-n)_{i-n}}{(m)_{i-n}} \right) \binom{p-1+a-l-m}{i-T-m} \\  \\
            			+  \binom{s-l}{s-(a-i+T)} \binom{s-(a-i+T)}{m-(i-n)} \frac{\prod\limits_{k=0}^{i-n-1}(s-a-k(p-1))}{(m)_{i-n}} -  \binom{r-l}{r-(a-i+T)}\binom{r-(a-i+T)}{m-(i-n)} \frac{\prod\limits_{k=0}^{i-n-1}(r-a-k(p-1))}{(m)_{i-n}} 
            			&+ O(p(r-s)) \\ \\
            			+ O(p(r-s)) & \\ \hline \\
            			\sum\limits_{k=1}^{i-n-1} \binom{s - l}{i-T+k(p-1)} \binom{i-T+k(p-1)}{i-T}  & \sum\limits_{k=1}^{i-n-1} \binom{s - l}{i-T+k(p-1)} \binom{i-T+k(p-1)}{i-T} \\
            			+ \binom{s-l}{s-(a-i+T)}\binom{s-(a-i+T)}{i-T}-\binom{r-l}{r-(a-i+T)}\binom{r-(a-i+T)}{i-T}  & 
            		\end{array} \right\rvert$,
            	} %
            \end{equation*}
            where the range of $m$ in the upper blocks  is $1,\ldots,i-n-1$, the range of $m$ in the middle blocks is $i-n,\ldots,i-T-1$ and the range of $m$ in the lower blocks is $ i-T $, and the range of $l$ in the left blocks (resp. right blocks) is $0,\ldots,a-i+T-1$ (resp. $a-i+T+1, \ldots, i-T$). 
            
            We now simplify every term in the middle blocks up to $O(p(r-s))$. We first simplify the terms appearing in the middle left block. From \eqref{est. error term adhc1} we have 
            \begin{align}\label{eq: error term left green}
            	\begin{split}
            	 &\left( \binom{r-l}{m-(i-n)} \frac{(r-a+i-n-1)_{i-n}}{(m)_{i-n}}- \binom{s-l}{m-(i-n)} \frac{(s-a+i-n-1)_{i-n}}{(m)_{i-n}} \right) \binom{a-l-m}{i-T-m} \\
            	 & \qquad \qquad \qquad + \binom{s-l}{s-(a-i+T)} \binom{s-(a-i+T)}{m-(i-n)} \frac{\prod\limits_{k=0}^{i-n-1}(s-a-k(p-1))}{(m)_{i-n}}  \\ 
            	 & \qquad \qquad \qquad -  \binom{r-l}{r-(a-i+T)}\binom{r-(a-i+T)}{m-(i-n)} \frac{\prod\limits_{k=0}^{i-n-1}(r-a-k(p-1))}{(m)_{i-n}}  \\
            	 & = \frac{r-s}{m!} \times \frac{(-1)^{i-n} (i-n)!}{(i-T-m)!} \times  \frac{(a-i+n-l)!}{(a-i+T-l)!} \times (H_{i-T-m} - H_{a-l-m}) + O(p(r-s)),
            	\end{split}
            \end{align}
            for  $m=i-n,\ldots,i-T-1$ and $l=0,\ldots,a-i+T-1$. We now estimate the terms in the middle right block. For $m=i-n,\ldots,i-T-1$ and $l=a-i+T+1,\ldots,i-T$ we have
            \begin{align*}
            	\binom{r-l}{m-(i-n)} & \frac{(r-a+i-n)_{i-n}}{(m)_{i-n}} - \binom{s-l}{m-(i-n)} \frac{(s-a+i-n)_{i-n}}{(m)_{i-n}} \\
            	&=  \frac{(r-a+i-n)_{i-n}}{(m)_{i-n}} \left( \binom{r-l}{m-(i-n)} - \binom{s-l}{m-(i-n)} \right)\\
            	&\qquad + \binom{s-l}{m-(i-n)} \left( \frac{(r-a+i-n)_{i-n}}{(m)_{i-n}} - \frac{(s-a+i-n)_{i-n}}{(m)_{i-n}} \right) \\
            	& = \binom{s-l}{m-(i-n)} \binom{m}{i-n}^{-1} \left( \binom{r-a+i-n}{i-n} - \binom{s-a+i-n}{i-n} \right) + O(p(r-s)) \\
            	& =  \frac{r-s}{s-(a-i+n)} \binom{m}{i-n}^{-1}  \binom{s-l}{m-(i-n)}   \binom{s-a+i-n}{i-n} + O(p(r-s)) \\
            	& = (r-s) \frac{(-1)^{i-n-1}}{i-n} \binom{m}{i-n}^{-1} \binom{s-l}{m-(i-n)} + O(p(r-s)) ,
            \end{align*}
            where we used $p\mid(r-(a-i+n))$ and Lemma~\ref{binomial coefficient under congruences} $(i)$ in the second equality, $p\mid(s-(a-i+n))$ and Lemma~\ref{binomial coefficient under congruences} $(iii)$ in the penultimate step and $s=a-i+n+(i-n)p$ and Lucas' theorem in the last step. Thus for $m=i-n,\ldots,i-T-1$ and $l=a-i+T+1,\ldots,i-T$, we have
            \begin{align*}
            	   \Bigg(\binom{r-l}{m-(i-n)} & \frac{(r-a+i-n)_{i-n}}{(m)_{i-n}} - \binom{s-l}{m-(i-n)} \frac{(s-a+i-n)_{i-n}}{(m)_{i-n}}\Bigg) \binom{p-1+a-l-m}{i-T-m} \\
            	   & = (r-s)  \frac{(-1)^{i-n-1}}{i-n} \binom{m}{i-n}^{-1} \binom{s-l}{m-(i-n)} \binom{p-1+a-l-m}{i-T-m} + O(p(r-s)). 
            \end{align*}
            We now further simplify the above expression. By Lucas' theorem, we have $\binom{p-1+a-l-m}{i-T-m} \equiv 0$ mod $p$ if $l+m<a$ and $\binom{s-l}{m-(i-n)} \equiv \binom{a-i+n-l}{m-(i-n)} \equiv 0$ mod $p$ if $l+m>a$. Thus for $m=i-n,\ldots,i-T-1$ and $l=a-i+T+1,\ldots,i-T$, we have
               \begin{align*}
               	\binom{s-l}{m-(i-n)} \binom{p-1+a-l-m}{i-T-m} \equiv \begin{cases}
               		\binom{p-1}{i-T-m} & \text{ if } l+m =a,\\
               		0 & \text{ if } l+m \neq a
               	\end{cases} \mod p
               \end{align*}         
               Substituting this in the previous congruence, for $m=i-n,\ldots,i-T-1$ and $l=a-i+T+1,\ldots,i-T$, we have
               \begin{align}\label{eq: error term right green}
               	\begin{split}
               	\Bigg(\binom{r-l}{m-(i-n)} & \frac{(r-a+i-n)_{i-n}}{(m)_{i-n}} - \binom{s-l}{m-(i-n)} \frac{(s-a+i-n)_{i-n}}{(m)_{i-n}}\Bigg) \binom{p-1+a-l-m}{i-T-m} \\
               	& = (r-s)  \frac{(-1)^{i-n-1}}{i-n} \binom{m}{i-n}^{-1}  \binom{p-1}{i-T-m} \delta_{l+m,a}+ O(p(r-s)). 
               	\end{split}
               \end{align}
             Thus by \eqref{eq: error term left green} and \eqref{eq: error term right green}, we have $\det(A_1)$ equals
            \begin{equation*}
            	\resizebox{0.999\hsize}{!}{%
            		$\left\lvert\begin{array}{c|c}
            			\sum\limits_{k=m}^{i-n-1} (p-1)^m \binom{s-l}{i-T+k(p-1)} \binom{k}{m} +  O(r-s)
            			&\sum\limits_{k=m}^{i-n-1} (p-1)^m\binom{s-l}{i-T+k(p-1)} \binom{k}{m} +  O(r-s) \\
            			\hline \\
            			\frac{r-s}{m!} \times \frac{(-1)^{i-n} (i-n)!}{(i-T-m)!} \times  \frac{(a-i+n-l)!}{(a-i+T-l)!} \times (H_{i-T-m} - H_{a-l-m}) 
            			& (r-s)  \frac{(-1)^{i-n-1}}{i-n} \binom{m}{i-n}^{-1}  \binom{p-1}{i-T-m} \delta_{l+m,a}\\  \\
            			+ O(p(r-s)) & + O(p(r-s))\\ \hline \\
            			\sum\limits_{k=1}^{i-n-1} \binom{s - l}{i-T+k(p-1)} \binom{i-T+k(p-1)}{i-T}  & \sum\limits_{k=1}^{i-n-1} \binom{s - l}{i-T+k(p-1)} \binom{i-T+k(p-1)}{i-T} \\
            			+ \binom{s-l}{s-(a-i+T)}\binom{s-(a-i+T)}{i-T}-\binom{r-l}{r-(a-i+T)}\binom{r-(a-i+T)}{i-T}  & 
            		\end{array} \right\rvert$
            	} %
            \end{equation*}
            where the range of $m$ and $l$ in every block is as before. We now want to perform row operations so that every entry in the upper blocks is given by a single binomial coefficient. Apply 
            \begin{align*}
            	R_{i-n-2} &\rightarrow R_{i-n-2}- \frac{(i-n-1)}{p-1} R_{i-n-1} \\ 
            	&~~\vdots \\
            	R_m &\rightarrow R_m - \sum_{k=m+1}^{i-n-1}\frac{\binom{k}{m}}{(p-1)^{k-m}} R_k \\
            	&~~\vdots \\
            	R_1 &\rightarrow R_1 - \sum_{k=2}^{i-n-1}\frac{\binom{k}{1}}{(p-1)^{k-1}}R_k.
            \end{align*}
             Thus $\det(A_1)$ equals
             \begin{equation*}
            	\resizebox{0.999\hsize}{!}{%
            		$\left\lvert\begin{array}{c|c}
            			(p-1)^m \binom{s-l}{i-T+m(p-1)}  +  O(r-s)
            			& (p-1)^m\binom{s-l}{i-T+m(p-1)} +  O(r-s) \\
            			\hline \\
            			\frac{r-s}{m!} \times \frac{(-1)^{i-n} (i-n)!}{(i-T-m)!} \times  \frac{(a-i+n-l)!}{(a-i+T-l)!} \times (H_{i-T-m} - H_{a-l-m}) 
            			& (r-s) \frac{(-1)^{i-n}(i-n)!}{(m)_{i-n}}  \binom{p-1}{i-T-m} \delta_{l+m,a}\\  \\
            			+ O(p(r-s)) & + O(p(r-s))\\ \hline \\
            			\sum\limits_{k=1}^{i-n-1} \binom{s - l}{i-T+k(p-1)} \binom{i-T+k(p-1)}{i-T}  & \sum\limits_{k=1}^{i-n-1} \binom{s - l}{i-T+k(p-1)} \binom{i-T+k(p-1)}{i-T} \\
            			+ \binom{s-l}{s-(a-i+T)}\binom{s-(a-i+T)}{i-T}-\binom{r-l}{r-(a-i+T)}\binom{r-(a-i+T)}{i-T}  & 
            		\end{array} \right\rvert$.
            	} %
            \end{equation*}
            To remove the sum appearing in the lower blocks, we apply the following row operation: 
            \begin{align*}
            	R_{i-T} \rightarrow R_{i-T} - \sum_{m=1}^{i-n-1} (p-1)^{-m} \binom{i-T+m(p-1)}{i-T} R_m.
            \end{align*}
            Noting that $\binom{i-T+m(p-1)}{i-T} \equiv 0 \mod p$ for $m=1,\ldots,i-n-1$, we see that $\det(A_1)$ equals
             \begin{equation*}
            	\resizebox{0.999\hsize}{!}{%
            		$\left\lvert\begin{array}{c|c}
            			(p-1)^m \binom{s-l}{i-T+m(p-1)}  +  O(r-s)
            			& (p-1)^m\binom{s-l}{i-T+m(p-1)} +  O(r-s) \\
            			\hline \\
            			\frac{r-s}{m!} \times \frac{(-1)^{i-n} (i-n)!}{(i-T-m)!} \times  \frac{(a-i+n-l)!}{(a-i+T-l)!} \times (H_{i-T-m} - H_{a-l-m}) 
            			& (r-s) \frac{(-1)^{i-n}(i-n)!}{(m)_{i-n}}  \binom{p-1}{i-T-m} \delta_{l+m,a} \\  \\
            			+ O(p(r-s)) & + O(p(r-s))\\ \hline \\
            			 \binom{s-l}{s-(a-i+T)}\binom{s-(a-i+T)}{i-T}-\binom{r-l}{r-(a-i+T)}\binom{r-(a-i+T)}{i-T} + O(p(r-s)) & O(p(r-s))
            		\end{array} \right\rvert$.
            	} %
            \end{equation*}
           We now simplify the entries in the lower left block. As $T \leq n<i$ and $s \equiv a-i+n$ mod $p$, we have $0 \leq n-T < i-T$ and $p\mid (s-(a-i+T)-(n-T))$. Thus for $l=0,\ldots,a-i+T-1$, we have
            \begin{align*}
            	\binom{s-l}{s-(a-i+T)}&\binom{s-(a-i+T)}{i-T} - \binom{r-l}{r-(a-i+T)}\binom{r-(a-i+T)}{i-T} \\ 
            	& =   \left( \binom{s-l}{s-(a-i+T)}-\binom{r-l}{r-(a-i+T)} \right) \binom{s-(a-i+T)}{i-T} \\
            	& \qquad +  \left( \binom{s-(a-i+T)}{i-T}-\binom{r-(a-i+T)}{i-T} \right) \binom{r-l}{r-(a-i+T)}  \\
            	& = O(p(r-s)) + \frac{r-s}{s-(a-i+n)}\binom{s-l}{s-(a-i+T)} \binom{s-(a-i+T)}{i-T},
            \end{align*}
            where we have used  Lemma~\ref{binomial coefficient under congruences} $(i)$ and $(iii)$ in the last step.
            Substituting this above, we get $\det(A_1)$ equals
            \begin{equation*}
            	\resizebox{0.999\hsize}{!}{%
            		$\left\lvert\begin{array}{c|c}
            			(p-1)^m \binom{s-l}{i-T+m(p-1)}  +  O(r-s)
            			& (p-1)^m\binom{s-l}{i-T+m(p-1)} +  O(r-s) \\
            			\hline \\
            			\frac{r-s}{m!} \times \frac{(-1)^{i-n} (i-n)!}{(i-T-m)!} \times  \frac{(a-i+n-l)!}{(a-i+T-l)!} \times (H_{i-T-m} - H_{a-l-m}) 
            			& (r-s) \frac{(-1)^{i-n}(i-n)!}{(m)_{i-n}}  \binom{p-1}{i-T-m} \delta_{l+m,a} \\  \\
            			+ O(p(r-s)) & + O(p(r-s))\\ \hline \\
            			\frac{r-s}{s-(a-i+n)}\binom{s-l}{s-(a-i+T)} \binom{s-(a-i+T)}{i-T} + O(p(r-s)) & O(p(r-s))
            		\end{array} \right\rvert$.
            	} %
            \end{equation*}
            Pulling out $(p-1)^m$ from every row in the upper blocks, $(-1)^{i-n}(i-n)!(r-s)$ from every row in the middle blocks and $(r-s)$ from the last row, we see that 
            \begin{align*}
            	\det(A_1) & = (p-1)^{(i-n)(i-n-1)/2} \times (-1)^{(i-n)(n-T)}  ((i-n)!)^{n-T} (r-s)^{n-T+1}   \times \det(B_1),
           \end{align*}
           where
           \begin{align*}	 
            	B_1 = \left[\begin{array}{c|c}
            		 \binom{s-l}{i-T+m(p-1)}  +  O(r-s)
            		& \binom{s-l}{i-T+m(p-1)} +  O(r-s) \\
            		\hline \\
            	 \frac{1}{(i-T-m)!m!} \times  \frac{(a-i+n-l)!}{(a-i+T-l)!} \times (H_{i-T-m} - H_{a-l-m}) + O(p)
            		&  \frac{1}{(m)_{i-n}}  \binom{p-1}{i-T-m} \delta_{l+m,a} + O(p) \\  \hline \\
            		\frac{1}{s-(a-i+n)}\binom{s-l}{s-(a-i+T)} \binom{s-(a-i+T)}{i-T} + O(p) & O(p)
            	\end{array} \right],
            \end{align*}
            and the range of $m$ in the upper blocks (resp. middle blocks and lower blocks) is $1,\ldots,i-n-1$ (resp. $i-n,\ldots,i-T-1$ and $i-T$) and the range of $l$ in the left blocks (resp. right blocks) is $0,\ldots,a-i+T-1$ (resp. $a-i+T+1,\ldots,i-T$). As $\det(A)=(-1)^{i-T}p^t \det(A_1)$, this implies  $p^{t(n-T+2)} \parallel \det(A)$, if $\det(B_1)$ is a $p$-adic unit.
             
            We claim that $\det(B_1)$ is indeed a $p$-adic unit. We determine each entry of $B_1$ modulo $p$. For $l=0,\ldots,i-T$ and $m=1,\ldots,i-n-1$, by Lucas' theorem, we have 
            \begin{align*}
            	\binom{s-l}{i-T+m(p-1)} \equiv \binom{i-n}{m} \binom{a-i+n-l}{i-T-m} \equiv 
            	\begin{cases}
            		   \binom{i-n}{m} \binom{a-i+n-l}{i-T-m} &\text{ if } l < a-i+T, \\
            		   0 &\text{ if } l > a-i+T,
            	\end{cases} \mod p,
            \end{align*}
            where we used $a\geq 2i-n-T$.
            This determines entries in the upper blocks of $B_1$ modulo $p$. For the entries in the last row, again by Lucas' theorem,  for $l=0,\ldots,a-i+T-1$ we have
            \begin{align*}
            	  \binom{s-l}{s-(a-i+T)} = \binom{s-l}{(a-i+T)-l} \equiv \binom{a-i+n-l}{a-i+T-l} \mod p.
            \end{align*}
            Furthermore, we have
            \begin{align}
               \begin{split}
                \frac{1}{s-(a-i+n)}\binom{s-(a-i+T)}{i-T}  &= \frac{(s-(a-i+T)) \cdots (s-(a-i+n)+1)}{(i-T)\cdots(i-n)} \binom{s-(a-i+n)-1}{i-n-1} \\
                &\equiv \frac{(n-T) \cdots 1}{(i-T)\cdots(i-n)} \binom{p-1}{i-n-1} \\
                & \equiv  (-1)^{i-n-1} \frac{(n-T)!(i-n-1)!}{(i-T)!} \mod p,
               \end{split} 
            \end{align}
            where in the penultimate step we used $s\equiv a-i+n \mod p$ and Lucas' theorem. Thus
            \begin{align*}	 
            	B_1 \equiv \left[\begin{array}{c|c}
            		 \binom{i-n}{m}\binom{a-i+n-l}{i-T-m}  
            		& 0 \\
            		\hline \\
            	 \frac{1}{(i-T-m)!m!} \times  \frac{(a-i+n-l)!}{(a-i+T-l)!} \times (H_{i-T-m} - H_{a-l-m}) 
            		&  \frac{1}{(m)_{i-n}}  \binom{p-1}{i-T-m} \delta_{l+m,a} \\  \hline \\
            		(-1)^{i-n-1} \binom{a-i+n-l}{a-i+T-l}  \frac{(n-T)! (i-n-1)!}{(i-T)!} & 0
            	\end{array} \right] \mod p.
            \end{align*}
            We observe that every column in the right half has exactly one non-zero entry. Expanding the determinant along the column $l=i-T$, then $l=i-T-1$ and so on till $l=a-i+T+1$, we get
            \begin{align*}
                \det(B_1) &=   \prod_{m=a-i+T}^{i-T-1} \frac{(-1)^{a}}{(m)_{i-n}}  \binom{p-1}{i-T-m} \\
                &\qquad \qquad \qquad \times \left\lvert \begin{array}{c}
            		 \binom{i-n}{m}\binom{a-i+n-l}{i-T-m} \\
            		\hline \\
            	 \frac{1}{(i-T-m)!m!} \times  \frac{(a-i+n-l)!}{(a-i+T-l)!} \times (H_{i-T-m} - H_{a-l-m}) 
            		\\ \hline \\
            		(-1)^{i-n-1} \binom{a-i+n-l}{a-i+T-l}  \frac{(n-T)! (i-n-1)!}{(i-T)!}    
            	\end{array} \right\rvert
            \end{align*}
             where the range of $m$ in the upper block (resp. middle block and lower block) is $1,\ldots,i-n-1$ (resp. $i-n,\ldots,a-i+T-1$ and $i-T$) and the range of $l$  is $0,\ldots,a-i+T-1$. Pulling out $\binom{i-n}{m} \frac{(i-n-m)!}{(i-T-m)!}$ from every row in the upper block,  $\frac{1}{m!(i-T-m)!}$ from every row in the middle block and $(-1)^{i-n-1}  \frac{ (i-n-1)!}{(i-T)!}$ from the last row, and $\frac{(a-i+n-l)!}{(a-i+T-l)!}$ from the $l^{\mathrm{th}}$ column, we see that
             \begin{align*}
                \det(B_1) &=   \prod_{m=a-i+T}^{i-T-1} \frac{(-1)^{a}}{(m)_{i-n}}  \binom{p-1}{i-T-m} \times  \prod_{m=1}^{i-n-1}\binom{i-n}{m} \frac{(i-n-m)!}{(i-T-m)!}  \times \prod_{m=i-n}^{a-i+T-1} \frac{1}{m!(i-T-m)!}  \\
                &\qquad \qquad \qquad \times (-1)^{i-n-1}   \frac{ (i-n-1)!}{(i-T)!} \times \prod_{l=0}^{a-i+T-1} \frac{(a-i+n-l)!}{(a-i+T-l)!} \times \det(B_2),
            \end{align*}
            where
            \begin{align*}
                B_2 = \left[ \begin{array}{c}
            		 \binom{a-i+T-l}{i-n-m} \\ \\
            		\hline \\
            	       H_{i-T-m} - H_{a-l-m} \\
            		\\ \hline \\
            		1   
            	\end{array} \right].
            \end{align*}
            It suffices to show $\det(B_2)$ is a $p$-adic unit. The computation below is similar to the one given in Lemma~\ref{choice beta hybrid ad1}. To kill every entry in the last row except the last, we  apply the column operations: 
            \begin{align*}
                    C_0 &\rightarrow\; C_0 - C_{1} \\
                    &~~\vdots \\
                    C_l &\rightarrow\; C_l - C_{l+1}, \\
                    &~~\vdots \\
                    C_{a-i+T-2} &\rightarrow\; C_{a-i+T-2} - C_{a-i+T-1}.
            \end{align*}
             Applying Pascal's identity for the binomial coefficients and expanding the determinant using the last row, we obtain 
            \begin{align*}
                \det(B_2) =  \det \left[ \begin{array}{c} \binom{a-i+T-l-1}{i-n-1-m}  \\ \\
            		\hline \\
            	      -\frac{1}{a-l-m} \\ \end{array} \right]
            \end{align*}
             where the range of $m$ in the upper block (resp. bottom block) is $1,\ldots,i-n-1$ (resp. $i-n,\ldots,a-i+T-1$) and the range of $l$  is $0,\ldots,a-i+T-2$.  
             
             We now derive some identities that we need for the next set of column operations. Note that 
            \begin{align}\label{eq: col comb. id 1}
            \begin{split}	
	                \binom{a-i+T-1-l}{i-n-1-m} &+ \sum_{l'=l+1}^{a-i+T-2} (-1)^{l'-l}\binom{a-i+T-2-l}{l'-l}\binom{a-i+T-1-l'}{i-n-1-m} \\
	                &= \sum_{l'=0}^{a-i+T-2-l} (-1)^{l'}\binom{a-i+T-2-l}{l'}\binom{a-i+T-1-l-l'}{i-n-1-m} \\
	                &= \begin{cases}
		               1, & \text{if } l=a-2i+n+T+m-1,~a-2i+n+T+m,\\
		               0, & \text{otherwise},
	                  \end{cases}
            \end{split}
            \end{align}	
            by Lemma~\ref{lem: combinatorial id for det(B)} $(iv)$ with $N=a-i+T-2-l$ and $k'=i-n-1-m$.
            Also, observe that
            \begin{align}\label{eq: col comb. id 2}
            \begin{split}
	             \frac{1}{a-l-m}+ \sum_{l'=l+1}^{a-i+T-2} (-1)^{l'-l} & \binom{a-i+T-2-l}{l'-l} \frac{1}{a-l'-m} \\
                 &= \sum_{l'=0}^{a-i+T-2-l} (-1)^{l'} \binom{a-i+T-2-l}{l'} \frac{1}{a-l-m-l'}   \\
                 &=(-1)^{a-i+T-l} \frac{(a-i+T-2-l)!}{(a-l-m) \cdots (i-T+2-m)}
            \end{split}     
            \end{align}             
            by Lemma~\ref{lem: combinatorial id for det(B)} $(ii)$ with $\alpha = a-l-m$ and $N=a-i+T-2-l$. 
            
             Applying the column operations 
            \begin{align*}
	                C_{0} & \rightarrow C_{0} + \sum_{l'=1}^{a-i+T-2} (-1)^{l'} \binom{a-i+T-2}{l'}C_{l'}, \\
	                  & ~~ \vdots \\
	                C_{l} & \rightarrow C_{l} + \sum_{l'=l+1}^{a-i+T-2} (-1)^{l'-l} \binom{a-i+T-2-l}{l'-l}C_{l'} \\
	                  & ~~\vdots \\
	                C_{a-i+T-3} &\rightarrow C_{a-i+T-3} - C_{a-i+T-2}
            \end{align*}
            and using \eqref{eq: col comb. id 1} for the upper left block and \eqref{eq: col comb. id 2} for the lower left block,  we get 
            \begin{align*}
	                  \det(B_2) 
                    &=   \det \left[
	                \begin{array}{c}
		                  \delta_{l,a-2i+n+T+m-1}+\delta_{l,a-2i+n+T+m}
                         \\[10pt]\hline\\[-5pt]
		                \frac{(-1)^{a-i+T-l+1}(a-i+T-2-l)!}{(a-l-m)\cdots(i-T+2-m)}
	                \end{array} \right] 
            \end{align*}
            Expanding the determinant first using the row $m=i-n-1$, second using  the row $m=i-n-2$ and so on we see that 
            \begin{align*}
            	\det(B_2) 
            	&=   (-1)^{(a+n+T)(i-n-1)} \times \det \left[
            	\begin{array}{c}
            		\frac{(-1)^{a-i+T-l+1}(a-i+T-2-l)!}{(a-l-m)\cdots (i-T+2-m)}
            	\end{array} \right],
            \end{align*}
            where the range of $m$ is $i-n, \ldots, a-i+T-1$ and the range of $l$ is $0,\ldots,a-2i+n+T-1$. 
            Re-indexing the rows we get
            \begin{align*}
            	\det(B_2) 
            	&=  (-1)^{(a+n+T)(i-n-1)} \times \det \left[
            	\begin{array}{c}
            		\frac{(-1)^{a-i+T-l+1}(a-i+T-2-l)!}{(a-i+n-l-m)\cdots(n-T+2-m)}
            	\end{array} \right],
            \end{align*}
            where the  range of $m,l$ is $0,\ldots,a-2i+n+T-1$. Pulling out $\frac{(n-T+1-m)!}{(a-i+n-m)!}$ from the $m^{\mathrm{th}}$-row and $(-1)^{a-i+T-l+1}(a-i+T-2-l)! l!$ from the $l^{\mathrm{th}}$-column, we get 
              \begin{align*}
              	\det(B_2) 
              	&=   (-1)^{(a+n+T)(i-n-1)} \times \prod_{l=0}^{a-2i+n+T-1} (-1)^{a-i+T-l+1}(a-i+T-2-l)!l! \\
              	&\qquad \times \prod_{m=0}^{a-2i+n+T-1}\frac{(n-T+1-m)!}{(a-i+n-m)!} \times 
              	\det \left[
              	\begin{array}{c}
              		\binom{a-i+n-m}{l}
              	\end{array} \right].
            \end{align*}
            From Lemma~\ref{cor: GV det}, we have the determinant on the second line is $\pm 1$. Thus $\det(B_2)$ is a $p$-adic unit, and so is $\det(B_1)$. Hence $p^{t(n-T+2)} \parallel \det(A)$. 
            Thus, there exist $\beta_1,\ldots,\beta_{a-i+T-1}$, $\beta_{a-i+T+1},\ldots,\beta_{i-T} \in \Q_p$ and $\gamma \in \Q_p$ satisfying \eqref{eq: linear eq. in matrix form green}. 
            
            We now show that they are $p$-adic integers. By Cramer's rule, we have
            \begin{align*}
                 \beta_l = \pm p^t \frac{\det(A_{i-T,l})}{\det(A)} \quad \forall ~l \neq {a-i+T}.
            \end{align*}
            By Lemma~\ref{minor trick 2}, we have $\det(A_{i-T,l}) \in p^{t(n-T+1)}\Z_p$ thus $\beta_l \in \Z_p$. A similar argument shows $\gamma \in \Z_p$. This proves $(i)$ and  $(ii)$.  

            We now prove $(iii)$. As $\beta_l \in \Z_p$, by Corollary~\ref{cor: binomial sums under congruences 2}, we have
            \begin{align*}
                \sum_{l=0}^{i-T} \beta_l \sum\limits_{ \substack{i-T <  j < r-(a-i+T) \\  j \equiv i-T ~\mathrm{mod}~(p-1)} }   \binom{r-l}{j} \binom{j}{m} \equiv \sum_{l=0}^{i-T} \beta_l \sum\limits_{ \substack{i-T <  j < s-(a-i+T) \\  j \equiv i-T ~\mathrm{mod}~(p-1)} }   \binom{s-l}{j} \binom{j}{m} \mod p^{t-v(m!)} \quad \text{for all } m.
            \end{align*}
            Here we have omitted the end points as they are congruent modulo $p^t$. Since many terms in the matrix $A$ vanish modulo $p^t$, it follows from \eqref{eq: linear eq. in matrix form green} then
            \begin{align}\label{eq: (i) mod p green}
            \begin{split}
               \sum_{k=1}^{i-n-1} \left( \sum_{l=0}^{i-T} \beta_l \binom{s-l}{i-T+k(p-1)}\right) \binom{i-T+k(p-1)}{m} 
               \equiv 0 \mod p^t
            \end{split}   
            \end{align}
            for $m=0,\ldots, i-T$. By Lemma~\ref{dmm' trick} (applied with $N=i-n-2$ $(\leq i-T)$ and $c=i-T+(p-1)$), the congruence \eqref{eq: (i) mod p green} holds for all $m \geq 0$. Part $(iii)$ follows.
            %
        \end{proof}
        \begin{theorem}\label{thm: hybrid green}
        	Let  $r \equiv a \mod{(p-1)}$ with $ 1 \leq a \leq p-1 $ and $r \geq i(p+1)+p $ with $ v(a_{p})  \in (i,i+1) $. Let $ s =  a-i+n+(i-n)p$ and $v(r-s) =t $  and $ 1 \leq t \leq T < n < i < a$. If  $2T \leq 2i-a $ and $2i-a-n \leq T$, then the image of $ \mathrm{ind}_{KZ}^{G}(V_{r}^{(i-T)}) $ is the same as the image of $ \mathrm{ind}_{KZ}^{G}(V_{r}^{(i-T+1)}) $ in $ \bar{\Theta}_{k,a_{p}} $.
        \end{theorem}
        \begin{proof}
        	Let $\beta_l$ and $\gamma$ be the $p$-adic integers  chosen in Lemma~\ref{choice beta hybrid green}. Then we have
            \begin{align*}
         \sum\limits_{l=0}^{i-T} \beta_l  \sum\limits_{ \substack{i-T <  j < r-(a-i+T) \\  j \equiv i-T ~\mathrm{mod}~ (p-1)} }   \binom{r-l}{j} \binom{j}{m}  + \gamma p^t  \binom {i-T}{m}\equiv 
        	    	0  \mod p^{t+1} ~ \mathrm{for}~ m=0,\ldots, i-T-1
            \end{align*}
            Then by Lemma~\ref{lem:choice of beta} there exist $\alpha_{j} \in \mathbb{Z}_{p}$ satisfying
        	\begin{enumerate}
        		\item[$(1)$]  $\alpha_j \equiv \sum\limits_{l=0}^{i-T} \beta_{l}  \binom{r-l}{j} ~\mathrm{mod}~p^{t} $, for all $ i-T< j < r-(a-i+T)$  with $j \equiv i-T $ mod $(p-1)$
        		\item[$(2)$]   $\sum\limits_{\substack{ i-T < j < r-(a-i+T) \\ j \equiv i-T ~\mathrm{mod}~(p-1)}}^{} \alpha_j \binom{j}{m} \equiv 0 $ mod $p^{i+t+1-m}$ for $m=0,\ldots, \min\{i+t, p-1\}$.
        	\end{enumerate}
        	From the congruence condition $(1)$ and Lemma~\ref{choice beta hybrid green} $(iii)$, we also have 
        	\begin{enumerate}
        		\item[$(2')$] $\sum\limits_{\substack{ i-T < j < r-(a-i+T) \\ j \equiv i-T ~\mathrm{mod}~(p-1)}}^{} \alpha_j \binom{j}{m} \equiv 0 $ mod $p^{t-v(m!)}$ for $m=p,\ldots, i+t$. 
        	\end{enumerate}
                 Note that for $m=p,\ldots, i+t$ we have $t-v(m!) = t-1 \geq i+t+1-p \geq i+t+1-m$  as $t < i \leq p-2$.
        	Let
        	\begin{align*}
        		f_2 &=  
        		\sum_{  \lambda \in \mathbb{F}_{p}^{\times}} 
        		\Bigg[ g_{2,p[\lambda]}^{0}, \sum_{\substack{l=0}}^{i-T} \frac{[\lambda]^{l-(a-i+T)}}{p^{l+t}} 
        		\beta_{l} (-\theta)^{l+t+1} X^{-t-1}Y^{r-(l+t+1)(p+1)+t+1} \Bigg]  \\
        		& \quad + \left[ g_{2,0}^{0}, \frac{1-p}{p^{a-i+T+t}}  \sum_{\substack{l=0}}^{i-T} \beta_{l}  \binom{r-l}{r-(a-i+T)}  (-\theta)^{a-i+T+t+1} X^{-t-1}Y^{r-(a-i+T+t+1)(p+1)+t+1} \right]  \\
        		f_1 &=  \left[g_{1,0}^0, \frac{p-1}{p^t a_p} \sum_{ \substack{i-T <  j < r-(a-i+T) \\  j \equiv i-T ~\mathrm{mod}~ (p-1)} } \alpha_{j} X^{r-j} Y^{j} \right] \\
        		f_{0} &= \left[ \mathrm{id}, \frac{1-p}{p^{i-T+t}}  \left(\sum_{\substack{l=0}}^{i-T} \beta_{l} \binom{r-l}{i-T} \right)  \theta^{i-T+t+1} X^{r-(i-T+t+1)(p+1)+t+1} Y^{-t-1}\right].
        	\end{align*}
        	By Lemma~\ref{theta and T plus} it follows that $T^+ f_2$ vanishes modulo $p$. It is easy to see that $-a_p f_2$,  $-a_p f_0 $ and $T^{-} f_0$ all vanish modulo $p$ using $t\leq T$ and $a \leq 2i-2T \leq 2i-T-t$. From $(2)$ and (the discussion below) $(2')$  it follows $T^{+}f_1$ vanishes modulo $p$. Using $t \leq T $ and $2i-a \leq i-1 \leq p-3$ one checks that $T^{-}f_1$ also vanishes modulo $p$. It can be checked that 
        	\begin{align*}
        		T^{-} f_2 -a_p f_1 + T^+ f_0 \equiv   \left[g_{1,0}^0, \frac{p-1}{p^t} \sum_{ \substack{i-T <  j < r-(a-i+T) \\  j \equiv i-T ~\mathrm{mod}~ (p-1)} } \left( \sum_{\substack{l=0}}^{i-T}  \beta_{l} \binom{r-l}{j}- \alpha_{j} \right) X^{r-j} Y^{j} \right] ~\mathrm{mod}~ p.
        	\end{align*}
        	From here on the argument is similar to the one given in the proof of Theorem~\ref{thm: hybrid ad}. Let
        	\begin{align*}
        		F(X,Y) = \frac{p-1}{p^t} \sum_{ \substack{i-T <  j < r-(a-i+T) \\  j \equiv i-T ~\mathrm{mod}~ (p-1)} } \left( \sum_{l=0 }^{i-T}  \beta_{l} \binom{r-l}{j}- \alpha_{j} \right) X^{r-j} Y^{j} + (p-1)\gamma X^{r-i+T}Y^{i-T}.
        	\end{align*}
        	By $(1)$ above, we have $F(X,Y) \in \mathbb{Z}_p[X,Y]$.  To prove the theorem, it is enough to show that $\overline{F(X,Y)}$ generates $V_{r}^{(i-T)}/V_r^{(i-T+1)}$. 
        	Note that  the coefficients of $X^{r}, \ldots, X^{r-(i-T-1)}Y^{i-T-1}$ in $F(X,Y)$ are zero. Since $i-T< a-i+T+ p-1 $, it follows that the coefficients of $X^{i-T}Y^{r-i+T}, \ldots, Y^r$ in $F(X,Y)$ are zero. By Lemma~\ref{choice beta hybrid green} $(i)$, $(2)$ and \cite[Lemma 2.8]{GR19}, we have $\theta^{i-T} \mid \overline{F(X,Y)}$. Applying Lemma~\ref{choice beta hybrid green} $(ii)$, $(2)$ and \cite[Lemma 2.12]{GR19}, with $m,l$ there equal $i-T$,  we obtain 
        	\begin{align*}
        		\overline{F(X,Y)} & \equiv (p-1) \theta^{i-T} X^{r-(i-T)(p+1)-(p-1)}Y^{p-1} + (p-1) \bar{\gamma} \theta^{i-T} X^{r-(i-T)(p+1)} ~\mathrm{mod}~V_r^{(i-T+1)}. 
        	\end{align*}  
        	Applying Lemma~\ref{generating polynomial quotient} with $m$ there equal to $i-T$, it follows that the Weil involution of $\overline{F(X,Y)} $ generates $V_r^{(i-T)}/V_r^{(i-T+1)}$. This finishes the proof of the theorem.
        \end{proof}
    We use the results obtained so far to describe $\bar{\Theta}_{k,a_{p}}$ when $a < 2i$ and $r \equiv a-i+n+(i-n)p \mod p(p-1)$ for some $1 \leq n \leq i-1$. We first consider the case $2i-2n-1 < a < 2i-n$.
    
    \begin{theorem}[Hybrid Conjecture]\label{thm: Shape theta a<2i-n}
		Let  $r \equiv a \mod{(p-1)}$ with $ 1 \leq a \leq p-1 $ and $r \geq i(p+1)+p $ with $ v(a_{p})  \in (i,i+1) $. Let $ s =  a-i+n+(i-n)p$ and $v(r-s) =t $  with $ 1 \leq  n \leq i-1 $. Assume that $i<a$  and  $2i-2n-1 < a < 2i-n$. Then
		\begin{enumerate}
			 \item[$(i)$]  $\mathrm{ind}_{KZ}^{G}(V_r^{(i-t+1)}/V_r^{(i-t+2)}) \twoheadrightarrow \bar{\Theta}_{k,a_p}$ if $  t \leq 2i-n-a$
          \item[$(ii)$] $\mathrm{ind}_{KZ}^{G}(V_r^{(a+t-i-1)}/V_r^{(a+t-i)}) \twoheadrightarrow \bar{\Theta}_{k,a_p}$ if $ 2i-n-a < t \leq i- \lfloor \frac{a}{2} \rfloor$  
          \item[$(iii)$] $\mathrm{ind}_{KZ}^{G}(V_r^{(i-t)}/V_r^{(i-t+1)}) \twoheadrightarrow \bar{\Theta}_{k,a_p}$ if $  i- \lfloor \frac{a}{2} \rfloor < t \leq n $ 
			 \item[$(iv)$]  $\mathrm{ind}_{KZ}^{G}(V_r^{(i-n)}/V_r^{(i-n+1)}) \twoheadrightarrow \bar{\Theta}_{k,a_p}$ if $ t \geq n+1$, in fact, we have  $$\mathrm{ind}_{KZ}^{G}( V_{p-1-(a-2i+2n)} \otimes D^{a-i+n}) \twoheadrightarrow \bar{\Theta}_{k,a_p}$$ if $ t \geq n+1$,
		\end{enumerate}
        where we assume $i=p-2 \Longrightarrow a \neq p-1$ if $2i-n-a < t \leq n$.
	\end{theorem}
    \begin{proof} From the hypotheses, we have $a-i+1 \leq a-i+n < i$. Also note that $a < 2(a-i+n)+1$. Thus by Lemma~\ref{JH factor Q 2i > a} $(ii)$ part $(a)$ applied with $r_0$ there equal to $a-i+n$, we have JH factors of $Q(i)$ are 
    \begin{align*}
    \{ V_{p-1-(a-2l)} &\otimes D^{a-l}: 0 \leq l \leq a-i-1\} \cup \{ V_{p-1-(a-2i+2n)} \otimes D^{a-i+n}\} \\ &\cup \text{ JH factors of } \{ V_{r}^{(l)}/V_{r}^{(l+1)} : i-n <l \leq i \}.
	\end{align*}
    Here the middle JH factor occurs in the cosocle of $V_{r}^{(i-n)}/V_{r}^{(i-n+1)}$. By Theorem~\ref{Elimination i < a and not in interval}, we have that the JH factors of the quotients $V_{r}^{(l)}/V_{r}^{(l+1)}$ for $0 \leq l \leq a-i-1$ vanish in $\bar{\Theta}_{k,a_{p}}$. Thus it remains to determine which of the following quotients $\{ V_{r}^{(i-T)}/V_{r}^{(i-T+1)}: 0   \leq T  \leq n \}$ survive in $\bar{\Theta}_{k,a_{p}}$. Since $a<2i-n$ and $i<a$, we get $n<i-1$.
        \begin{enumerate}
            \item[$(i)$] We will show the images of all JH factors except for those in $V_{r}^{(i-t+1)}/V_{r}^{(i-t+2)}$  vanish in  $\bar{\Theta}_{k,a_{p}}$.
            \begin{itemize}
                \item If $T<t-1$, then $n+T < 2i-a-1$. Thus by Theorem~\ref{eliminating JH dcbd}, we get  the  $\mathrm{ind}_{KZ}^{G}(V_{r}^{(i-T)}/V_{r}^{(i-T+1)})$  vanishes in  $\bar{\Theta}_{k,a_{p}}$ if $T < t-1$.
                \item If  $t\leq T \leq 2i-n-a-1$, then by the hypothesis $a \geq 2i-2n$ we have $T\leq n-1$. Hence by Theorem~\ref{eliminating JH dcad}, we obtain $\mathrm{ind}_{KZ}^{G}(V_{r}^{(i-T)}/V_{r}^{(i-T+1)})$  vanishes in  $\bar{\Theta}_{k,a_{p}}$ if $t\leq T \leq 2i-n-a-1$.
                \item If $2i-n-a \leq T < i - \frac{a}{2}$, then by Theorem~\ref{thm: hybrid green}, we have  $\mathrm{ind}_{KZ}^{G}(V_{r}^{(i-T)}/V_{r}^{(i-T+1)})$  vanishes in  $\bar{\Theta}_{k,a_{p}}$. 
                
                \item If $i - \frac{a}{2} \leq T \leq n$, then $T+t \leq 2i-a$. Hence by Theorem~\ref{thm: hybrid ad}, we have $\mathrm{ind}_{KZ}^{G}(V_{r}^{(i-T)}/V_{r}^{(i-T+1)})$  vanishes in  $\bar{\Theta}_{k,a_{p}}$.
            \end{itemize}
             This proves $(i)$.
            
            \item[$(ii)$]  We will  show the images of all JH factors coming from $\mathrm{ind}_{KZ}^{G}(V_{r}^{(i-T)}/V_{r}^{(i-T+1)})$  vanish in  $\bar{\Theta}_{k,a_{p}}$ except for $T = 2i-a-t+1$. Since $2i-2n-1 <a$, it follows that $\lfloor \frac{a}{2} \rfloor \geq i-n$. As $t\leq i - \lfloor \frac{a}{2}\rfloor$, we get $t\leq n$.
            \begin{itemize}
                \item If $T\leq 2i-a-n-1$, then $n+T \leq 2i-a-1$, $T<n$ and $T+2 \leq 2i-a-n+1 \leq t$. Thus by Theorem~\ref{eliminating JH hcbd}, we get   $\mathrm{ind}_{KZ}^{G}(V_{r}^{(i-T)}/V_{r}^{(i-T+1)})$  vanishes in  $\bar{\Theta}_{k,a_{p}}$ if $T \leq 2i-a-n-1$.
                \item  If $2i-a-n-1<T<t$, then $T<t\leq n$ and it follows from Theorem~\ref{eliminating JH sdcbd} that  $\mathrm{ind}_{KZ}^{G}(V_{r}^{(i-T)}/V_{r}^{(i-T+1)})$  vanishes in  $\bar{\Theta}_{k,a_{p}}$.
                \item If $t \leq T < i - \frac{a}{2}$, then  we have $2i-a-n <  t\leq T<i - \frac{a}{2} \leq n$. Thus by Theorem~\ref{thm: hybrid green}, we have  $\mathrm{ind}_{KZ}^{G}(V_{r}^{(i-T)}/V_{r}^{(i-T+1)})$  vanishes in  $\bar{\Theta}_{k,a_{p}}$ for $t \leq T \leq i - \frac{a}{2}$.
                \item If $i - \frac{a}{2} \leq T \leq 2i-a-t$, then  by Theorem~\ref{thm: hybrid ad}  we have  $\mathrm{ind}_{KZ}^{G}(V_{r}^{(i-T)}/V_{r}^{(i-T+1)})$  vanishes in  $\bar{\Theta}_{k,a_{p}}$ as $2i-a-t<n$.
                \item If $ 2i-a-t+2 \leq T \leq n$, then we check that the hypothesis of Theorem~\ref{thm: above super diagonal a large} hold. It follows from the condition $t \leq i- \lfloor \frac{a}{2} \rfloor$ that $2t \leq 2i-a+1 $. Thus $T \geq 2i-a-t+2 = t+1 +(2i-a+1-2t) \geq t+1$ and $a \geq 2i-T-t+2 > 2i-2T$. This checks the hypotheses $(i)$ and $(ii)$ of Theorem~\ref{thm: above super diagonal a large}. Since $T \leq n \leq 2i-a-1$, we need to verify the hypothesis $(iii)$ of Theorem~\ref{thm: above super diagonal a large}. This clearly follows from the assumption $ 2i-a-t+2 \leq T$.
                Hence by Theorem~\ref{thm: above super diagonal a large}  we have  $\mathrm{ind}_{KZ}^{G}(V_{r}^{(i-T)}/V_{r}^{(i-T+1)})$  vanishes in  $\bar{\Theta}_{k,a_{p}}$.
            \end{itemize}
             This proves $(ii)$.
            
            \item[$(iii)$] We show  $\mathrm{ind}_{KZ}^{G}(V_{r}^{(i-T)}/V_{r}^{(i-T+1)})$  vanishes in  $\bar{\Theta}_{k,a_{p}}$ except for $T = t$.
            \begin{itemize}
                \item If $T \leq 2i-a-n-1$, then $n+T \leq 2i-a-1$, $T\leq 2i-a-n-1<n$ and $T+2 \leq 2i-a-n+1 \leq i -  \lfloor \frac{a}{2} \rfloor \leq t$. Thus by Theorem~\ref{eliminating JH hcbd}, we get   $\mathrm{ind}_{KZ}^{G}(V_{r}^{(i-T)}/V_{r}^{(i-T+1)})$  vanishes in  $\bar{\Theta}_{k,a_{p}}$ if $T < 2i-a-n-1$.
                \item If $2i-a-n \leq T \leq t-1 $, then by Theorem~\ref{eliminating JH sdcbd}, we see that  $\mathrm{ind}_{KZ}^{G}(V_{r}^{(i-T)}/V_{r}^{(i-T+1)})$  vanishes in  $\bar{\Theta}_{k,a_{p}}$ as $T<t\leq n$.
                \item If $ t+1 \leq T \leq n$, then we check the hypothesis of Theorem~\ref{thm: above super diagonal a large}. Clearly, the hypothesis $(i)$ holds. It can be checked that the hypothesis $(ii)$ and $(iii)$ of Theorem~\ref{thm: above super diagonal a large} hold using the string of inequalities $i - \frac{a}{2} \leq i - \lfloor \frac{a}{2} \rfloor < t < T $. Hence it follows from  Theorem~\ref{thm: above super diagonal a large} that  the quotient $\mathrm{ind}_{KZ}^{G}(V_{r}^{(i-T)}/V_{r}^{(i-T+1)})$  vanishes in  $\bar{\Theta}_{k,a_{p}}$ if $t+1 \leq T \leq n$.
            \end{itemize}
            This proves $(iii)$.  
            
            \item[$(iv)$] We show  $\mathrm{ind}_{KZ}^{G}(V_{r}^{(i-T)}/V_{r}^{(i-T+1)})$  vanishes in  $\bar{\Theta}_{k,a_{p}}$ except for $T = n$.
            \begin{itemize}
                \item If $T \leq 2i-a-n-1$, then $n+T \leq 2i-a-1$ and $T+2 \leq 2i-a-n+1 \leq n+1 \leq t$. Thus by Theorem~\ref{eliminating JH hcbd}, we get   $\mathrm{ind}_{KZ}^{G}(V_{r}^{(i-T)}/V_{r}^{(i-T+1)})$  vanishes in  $\bar{\Theta}_{k,a_{p}}$ if $T \leq 2i-a-n-1$.
                \item If $2i-a-n \leq T \leq n-1$, then by Theorem~\ref{eliminating JH sdcbd}, we see that  $\mathrm{ind}_{KZ}^{G}(V_{r}^{(i-T)}/V_{r}^{(i-T+1)})$  vanishes in  $\bar{\Theta}_{k,a_{p}}$ as $t \geq n+1$.
            \end{itemize}
             This proves $(iv)$. \qedhere
        \end{enumerate}
    \end{proof}
    The conclusion of the  theorem can be summarized by the following diagrams:  
    	\begin{figure}[hbt!]
             \centering
			\resizebox{.8\textwidth}{!}{
			\begin{tikzpicture}[font=\small]
				
				\def\Nr{13}        
				\def\Nc{13}        
				\def\cwA{2}
				\def\cw{1}
				\def\rhA{1.2}
				\def\rh{1}
				
				\draw (0,0) -- (0,-\rhA-\Nr*\rh);
				\draw (\cwA,0) -- (\cwA,-\rhA-\Nr*\rh);
				\foreach \i in {1,...,\Nc}{
					\draw (\cwA+\i*\cw,0) -- (\cwA+\i*\cw,-\rhA-\Nr*\rh);
				}
				
				\draw (0,0) -- (\cwA+\Nc*\cw,0);
				\draw (0,-\rhA) -- (\cwA+\Nc*\cw,-\rhA);
				\foreach \j in {1,...,\Nr}{
					\draw (0,-\rhA-\j*\rh) -- (\cwA+\Nc*\cw,-\rhA-\j*\rh);
				}
				
				\draw (0,0) rectangle (\cwA,-\rhA);
				\draw (0,0) -- (\cwA,-\rhA);
				
				\node at (\cwA/3,-2*\rhA/3) {$t$};
				\node at (2*\cwA/3,-\rhA/3) {$T$};
				
				\node at (\cwA+0.5,-0.5*\rhA) {$0$};
				\node at (\cwA+1.5,-0.5*\rhA) {$1$};
				\node at (\cwA+2.5,-0.5*\rhA) {$\cdot$};
				\node at (\cwA+3.5,-0.5*\rhA) {$\cdot$};
				\node at (\cwA+4.5,-0.5*\rhA) {$k-1$};
				\node at (\cwA+5.5,-0.5*\rhA) {$k$};
				\node at (\cwA+6.5,-0.5*\rhA) {$\cdot$};
				\node at (\cwA+7.5,-0.5*\rhA) {$\cdot$};
				\node at (\cwA+8.5,-0.5*\rhA) {$\cdot$};
				\node at (\cwA+9.5,-0.5*\rhA) {{\scriptsize $i-\lfloor\frac{a}{2}\rfloor$}};
				\node at (\cwA+10.5,-0.5*\rhA) {$\cdot$};
				\node at (\cwA+11.5,-0.5*\rhA) {$\cdot$};
				\node at (\cwA+12.5,-0.5*\rhA) {$n$};
				
				\node at (0.5*\cwA,-\rhA-0.5) {$t=1$};
				\node at (0.5*\cwA,-\rhA-1.5) {$t=2$};
				\node at (0.5*\cwA,-\rhA-2.5) {$\cdot$};
				\node at (0.5*\cwA,-\rhA-3.5) {$\cdot$};
				\node at (0.5*\cwA,-\rhA-4.5) {$t=k$};
				\node at (0.5*\cwA,-\rhA-5.5) {$t=k+1$};
				\node at (0.5*\cwA,-\rhA-6.5) {$\cdot$};
				\node at (0.5*\cwA,-\rhA-7.5) {$\cdot$};
				\node at (0.5*\cwA,-\rhA-8.5) {$t=i-\lfloor \frac{a}{2}\rfloor$};
				\node at (0.5*\cwA,-\rhA-9.5) {$\cdot$};
				\node at (0.5*\cwA,-\rhA-10.5) {$\cdot$};
				\node at (0.5*\cwA,-\rhA-11.5) {$t=n$};
				\node at (0.5*\cwA,-\rhA-12.5) {$t \geq n+1$};
				
				\foreach \i in {1,...,\Nc}{
					\foreach \j in {1,...,\Nr}{
						\pgfmathsetmacro{\x}{\cwA+(\i-0.5)}
						\pgfmathsetmacro{\y}{-\rhA-(\j-0.5)}
						
						\def\putmark{\xmark}
						
						\ifnum\j<6
						\ifnum\i=\j
						\def\putmark{\cmark}
						\fi
						\fi
						
						\ifnum\j=13 \ifnum\i=13 \def\putmark{\cmark}\fi\fi
						\ifnum\j=12 \ifnum\i=13 \def\putmark{\cmark}\fi\fi
						\ifnum\j=11 \ifnum\i=12 \def\putmark{\cmark}\fi\fi
						\ifnum\j=10 \ifnum\i=11 \def\putmark{\cmark}\fi\fi
						\ifnum\j=9  \ifnum\i=10 \def\putmark{\cmark}\fi\fi
						\ifnum\j=8  \ifnum\i=11 \def\putmark{\cmark}\fi\fi
						\ifnum\j=7  \ifnum\i=12 \def\putmark{\cmark}\fi\fi
						\ifnum\j=6  \ifnum\i=13 \def\putmark{\cmark}\fi\fi
						
						\node at (\x,\y) {\putmark};
					}
				}
				
			\end{tikzpicture}
			}
            \caption{\label{fig: hybrid a odd} Here $2i-2n-1<a<2i$ and $a$ is odd, and $k=2i-a-n$.}
		\end{figure}

        	\begin{figure}[hbt!]
             \centering
			\resizebox{.8\textwidth}{!}{
			\begin{tikzpicture}[font=\small]
				
				\def\Nr{14}        
				\def\Nc{14}        
				\def\cwA{2}
				\def\cw{1}
				\def\rhA{1.2}
				\def\rh{1}
				
				\draw (0,0) -- (0,-\rhA-\Nr*\rh);
				\draw (\cwA,0) -- (\cwA,-\rhA-\Nr*\rh);
				\foreach \i in {1,...,\Nc}{
					\draw (\cwA+\i*\cw,0) -- (\cwA+\i*\cw,-\rhA-\Nr*\rh);
				}
				
				\draw (0,0) -- (\cwA+\Nc*\cw,0);
				\draw (0,-\rhA) -- (\cwA+\Nc*\cw,-\rhA);
				\foreach \j in {1,...,\Nr}{
					\draw (0,-\rhA-\j*\rh) -- (\cwA+\Nc*\cw,-\rhA-\j*\rh);
				}
				
				\draw (0,0) rectangle (\cwA,-\rhA);
				\draw (0,0) -- (\cwA,-\rhA);
				
				\node at (\cwA/3,-2*\rhA/3) {$t$};
				\node at (2*\cwA/3,-\rhA/3) {$T$};
				
				\node at (\cwA+0.5,-0.5*\rhA) {$0$};
				\node at (\cwA+1.5,-0.5*\rhA) {$1$};
				\node at (\cwA+2.5,-0.5*\rhA) {$\cdot$};
				\node at (\cwA+3.5,-0.5*\rhA) {$\cdot$};
				\node at (\cwA+4.5,-0.5*\rhA) {$k-1$};
				\node at (\cwA+5.5,-0.5*\rhA) {$k$};
				\node at (\cwA+6.5,-0.5*\rhA) {$\cdot$};
				\node at (\cwA+7.5,-0.5*\rhA) {$\cdot$};
                \node at (\cwA+8.5,-0.5*\rhA) {$\cdot$};
				\node at (\cwA+9.5,-0.5*\rhA) {{\scriptsize $i-\lfloor\frac{a}{2}\rfloor$}};
				\node at (\cwA+10.5,-0.5*\rhA) {{\small {$i-\lfloor \frac{a}{2}\rfloor \atop +1$}}};
				\node at (\cwA+11.5,-0.5*\rhA) {$\cdot$};
				\node at (\cwA+12.5,-0.5*\rhA) {$\cdot$};
				\node at (\cwA+13.5,-0.5*\rhA) {$n$};
				
				\node at (0.5*\cwA,-\rhA-0.5) {$t=1$};
				\node at (0.5*\cwA,-\rhA-1.5) {$t=2$};
				\node at (0.5*\cwA,-\rhA-2.5) {$\cdot$};
				\node at (0.5*\cwA,-\rhA-3.5) {$\cdot$};
				\node at (0.5*\cwA,-\rhA-4.5) {$t=k$};
				\node at (0.5*\cwA,-\rhA-5.5) {$t=k+1$};
				\node at (0.5*\cwA,-\rhA-6.5) {$\cdot$};
				\node at (0.5*\cwA,-\rhA-7.5) {$\cdot$};
				\node at (0.5*\cwA,-\rhA-8.5) {$t=i-\lfloor \frac{a}{2}\rfloor$};
				\node at (0.5*\cwA,-\rhA-9.5) {{\scriptsize $t=i-\lfloor \frac{a}{2}\rfloor+1$}};
				\node at (0.5*\cwA,-\rhA-10.5) {$\cdot$};
				\node at (0.5*\cwA,-\rhA-11.5) {$\cdot$};
				\node at (0.5*\cwA,-\rhA-12.5) {$t= n$};
				\node at (0.5*\cwA,-\rhA-13.5) {$t\geq n+1$};
				\foreach \i in {1,...,\Nc}{
					\foreach \j in {1,...,\Nr}{
						\pgfmathsetmacro{\x}{\cwA+(\i-0.5)}
						\pgfmathsetmacro{\y}{-\rhA-(\j-0.5)}
						
						\def\putmark{\xmark}
						
						\ifnum\j<6
						\ifnum\i=\j
						\def\putmark{\cmark}
						\fi
						\fi
						
                       \ifnum\j=14 \ifnum\i=14 \def\putmark{\cmark}\fi\fi
						\ifnum\j=13 \ifnum\i=14 \def\putmark{\cmark}\fi\fi
						\ifnum\j=12 \ifnum\i=13 \def\putmark{\cmark}\fi\fi
						\ifnum\j=11 \ifnum\i=12 \def\putmark{\cmark}\fi\fi
						\ifnum\j=10 \ifnum\i=11 \def\putmark{\cmark}\fi\fi
						\ifnum\j=9  \ifnum\i=11 \def\putmark{\cmark}\fi\fi
						\ifnum\j=8  \ifnum\i=12 \def\putmark{\cmark}\fi\fi
						\ifnum\j=7  \ifnum\i=13 \def\putmark{\cmark}\fi\fi
						\ifnum\j=6  \ifnum\i=14 \def\putmark{\cmark}\fi\fi
						
						\node at (\x,\y) {\putmark};
					}
				}				
			\end{tikzpicture}
			}
            \caption{\label{fig: hybrid a even} Here $2i-2n-1<a<2i$ and $a$ is even, and $k=2i-a-n$.}
		\end{figure}
\newpage
    We next consider the case $2i-n \leq a < 2i$.
    \begin{theorem}[Hybrid Conjecture]\label{Shape theta 2i-n<a}
    Let  $r \equiv a \mod{(p-1)}$ with $ 1 \leq a \leq p-1 $ and $r \geq i(p+1)+p $ with $ v(a_{p})  \in (i,i+1) $. Let $ s =  a-i+n+(i-n)p$ and $v(r-s) =t $  with $ 1 \leq  n \leq i-1 $ and $t \geq 1$. Assume that $i<a$  and $2i-n \leq a < 2i$. Then
		\begin{enumerate}
			 \item[$(i)$] $\mathrm{ind}_{KZ}^{G}(V_r^{(a+t-i-1)}/V_r^{(a+t-i)}) \twoheadrightarrow \bar{\Theta}_{k,a_p}$ if $ 0 < t \leq i- \lfloor \frac{a}{2} \rfloor$,  in fact, if $a=2i-n$ and $t=1$ we have $$\mathrm{ind}_{KZ}^{G}( V_{p-1-(a-2i+2n)} \otimes D^{a-i+n}) \twoheadrightarrow \bar{\Theta}_{k,a_p}$$
          \item[$(ii)$] $\mathrm{ind}_{KZ}^{G}(V_r^{(i-t)}/V_r^{(i-t+1)}) \twoheadrightarrow \bar{\Theta}_{k,a_p}$ if $  i- \lfloor \frac{a}{2} \rfloor < t \leq n $   
			 \item[$(iii)$]  $\mathrm{ind}_{KZ}^{G}(V_r^{(i-n)}/V_r^{(i-n+1)}) \twoheadrightarrow \bar{\Theta}_{k,a_p}$ if $ t \geq n+1$, in fact,  we have  $$\mathrm{ind}_{KZ}^{G}( V_{p-1-(a-2i+2n)} \otimes D^{a-i+n}) \twoheadrightarrow \bar{\Theta}_{k,a_p}$$ if $ t \geq n+1$,
		\end{enumerate}
        where we assume $i=p-2 \Longrightarrow a \neq p-1$ if $0 < t \leq n$.
	\end{theorem}
    \begin{proof}
        We first determine the JH factors of $Q(i)$ using Lemma~\ref{JH factor Q 2i > a}. From the hypotheses, we have $r \equiv a-i+n \mod p$ and $i \leq a-i+n \leq a-1$. Thus  we have $r \equiv i,\ldots, a-1 \mod p$. If $r \equiv i \mod p$, then by Lemma~\ref{JH factor Q 2i > a} $(ii)$ part $(a)$ applied with $r_0=i$ we get the
        JH factors of $Q(i)$ are 
       \begin{align}\label{JH factors a<2i, r_0=i} 
          \{ V_{p-1-(a-2l)} \otimes D^{a-l}: 0 \leq l \leq a-i-1\} \cup\{V_{p-1+a-2i} \otimes D^{i}\}  \cup \text{ JH factors of } \{ V_{r}^{(l)}/V_{r}^{(l+1)} : a-i <l \leq i \}.
	\end{align}
     If $r \equiv i+1,\ldots,a-1 \mod p$, then  by Lemma~\ref{JH factor Q 2i > a} $(iii)$ with $j$ there equal to $i-n+1$, we have the JH factors of $Q(i)$ are 
    \begin{align}\label{JH factors a<2i, r_0 > i}
    \{ V_{p-1-(a-2l)} \otimes D^{a-l}: 0 \leq l \leq i-n\}  \cup \text{ JH factors of } \{ V_{r}^{(l)}/V_{r}^{(l+1)} : i-n <l \leq i \}.
	\end{align}
    Combining \eqref{JH factors a<2i, r_0=i} and \eqref{JH factors a<2i, r_0 > i}, we get JH factors of $Q(i)$ are
    \begin{align}\label{JH factors a<2i, r_0 ge i}
    \{ V_{p-1-(a-2l)} \otimes D^{a-l}: 0 \leq l \leq i-n\}  \cup \text{ JH factors of } \{ V_{r}^{(l)}/V_{r}^{(l+1)} : i-n <l \leq i \}
	\end{align}
    whenever $r \equiv a-i+n \mod p$ with $2i-a \leq n \leq i-1$.
    
     We now eliminate the first few JH factors appearing in \eqref{JH factors a<2i, r_0 ge i} using  Theorem~\ref{Elimination i < a and not in interval}. For $0\leq l \leq i-n-1$, we have $l \leq a-i-1$ since $2i-a \leq n$. Also, $l=a-i-1$ if only if $l=i-n-1$ and $n=2i-a$. Thus by Theorem~\ref{Elimination i < a and not in interval}, we have that the JH factors of the quotients $V_{r}^{(l)}/V_{r}^{(l+1)}$ for $0 \leq l \leq i-n-1$ vanish in $\bar{\Theta}_{k,a_{p}}$. Thus it remains to determine which of the following quotients $\{ V_{r}^{(i-T)}/V_{r}^{(i-T+1)}: 0   \leq T  \leq n \}$ survive in $\bar{\Theta}_{k,a_{p}}$. Now the proof is similar to Theorem~\ref{thm: Shape theta a<2i-n}. 
     \begin{enumerate}
		    \item[$(i)$]  We will  show the images of all the JH factors coming from $\mathrm{ind}_{KZ}^{G}(V_{r}^{(i-T)}/V_{r}^{(i-T+1)})$  vanish in  $\bar{\Theta}_{k,a_{p}}$ except for $T = 2i-a-t+1$.  Since $2i-2n \leq a$, it follows that $\lfloor \frac{a}{2} \rfloor \geq i-n$. As $t\leq i - \lfloor \frac{a}{2}\rfloor$, we get $t\leq n$.
            \begin{itemize}
                \item  If $0 \leq T <t$, then $T<t\leq n$ and it follows from Theorem~\ref{eliminating JH sdcbd} that  $\mathrm{ind}_{KZ}^{G}(V_{r}^{(i-T)}/V_{r}^{(i-T+1)})$  vanishes in  $\bar{\Theta}_{k,a_{p}}$.
                \item If $t \leq T \leq i - \frac{a}{2}$, then  we have $2i-a-n \leq 0 <  t\leq T<i - \frac{a}{2} \leq n$. Thus by Theorem~\ref{thm: hybrid green}, we have  $\mathrm{ind}_{KZ}^{G}(V_{r}^{(i-T)}/V_{r}^{(i-T+1)})$  vanishes in  $\bar{\Theta}_{k,a_{p}}$ for $t \leq T \leq i - \frac{a}{2}$.
                \item If $i - \frac{a}{2} \leq T \leq 2i-a-t$, then  by Theorem~\ref{thm: hybrid ad}  we have  $\mathrm{ind}_{KZ}^{G}(V_{r}^{(i-T)}/V_{r}^{(i-T+1)})$  vanishes in  $\bar{\Theta}_{k,a_{p}}$.
                \item If $ 2i-a-t+2 \leq T \leq n$, then we check that the hypothesis of Theorem~\ref{thm: above super diagonal a large} hold. It follows from the condition $t \leq i- \lfloor \frac{a}{2} \rfloor$ that $2t \leq 2i-a+1 $. Thus $T \geq 2i-a-t+2 = t+1 +(2i-a+1-2t) \geq t+1$ and $a \geq 2i-T-t+2 > 2i-2T$. This checks the hypotheses $(i)$ and $(ii)$ of Theorem~\ref{thm: above super diagonal a large}. If $a-i+T \leq i+1$, we need to verify the hypothesis $(iii)$ of Theorem~\ref{thm: above super diagonal a large}. This clearly follows from the assumption $ 2i-a-t+2 \leq T$. Hence by Theorem~\ref{thm: above super diagonal a large}  we have  $\mathrm{ind}_{KZ}^{G}(V_{r}^{(i-T)}/V_{r}^{(i-T+1)})$  vanishes in  $\bar{\Theta}_{k,a_{p}}$.
            \end{itemize}
            Using $a\geq 2i-n$, it can be checked that $l=i-(2i-a-t+1) \leq i-n$ if and only if $t=1$, $a=2i-n$ and $l = i-n$. Now the second statement follows from \eqref{JH factors a<2i, r_0 ge i}.
            \item[$(ii)$]   We show  $\mathrm{ind}_{KZ}^{G}(V_{r}^{(i-T)}/V_{r}^{(i-T+1)})$  vanishes in  $\bar{\Theta}_{k,a_{p}}$ except for $T = t$.
            \begin{itemize}
                \item {If $0 \leq T \leq t-1 $, then $2i-a-n \leq 0 \leq T$. Thus by Theorem~\ref{eliminating JH sdcbd}, we see that  $\mathrm{ind}_{KZ}^{G}(V_{r}^{(i-T)}/V_{r}^{(i-T+1)})$  vanishes in  $\bar{\Theta}_{k,a_{p}}$ as $T<t\leq n$.}
                \item If $ t+1 \leq T \leq n$, then we check the hypothesis of Theorem~\ref{thm: above super diagonal a large}. Clearly, the hypothesis $(i)$ holds. It can be checked that the hypothesis $(ii)$ and $(iii)$ of Theorem~\ref{thm: above super diagonal a large} hold using the string of inequalities $i - \frac{a}{2} \leq i - \lfloor \frac{a}{2} \rfloor < t < T $. Hence it follows from  Theorem~\ref{thm: above super diagonal a large} that  the quotient $\mathrm{ind}_{KZ}^{G}(V_{r}^{(i-T)}/V_{r}^{(i-T+1)})$  vanishes in  $\bar{\Theta}_{k,a_{p}}$ if $t+1 \leq T \leq n$.             
            \end{itemize}
            \item[$(iii)$] Clearly $2i-a-n -T\leq 0 $ for any $T\geq 0$. If $0 \leq T \leq n-1$, then by Theorem~\ref{eliminating JH sdcbd}, we see that  $\mathrm{ind}_{KZ}^{G}(V_{r}^{(i-T)}/V_{r}^{(i-T+1)})$  vanishes in  $\bar{\Theta}_{k,a_{p}}$ as $t \geq n+1$.
            \qedhere
		\end{enumerate}
    \end{proof}

 \begin{remark}\label{r = 27}
        In 
        part $(i)$
        of Theorem~\ref{Shape theta 2i-n<a}, if $n = 1$,  then $$\mathrm{ind}_{KZ}^{G}( V_{p-2} \otimes D^{i}) \twoheadrightarrow \bar{\Theta}_{k,a_p}$$
        and the question arises whether
        this map factors through $T$ or a
        quadratic polynomial in $T$.
        When $p = 5$, $r = 27$ and $i = 2$,
        so that $a = 3$, $n = 1$, $s = 7$ and $t = 1$, we can show that the above map 
        factors through $T$ under the additional assumption
        $$v(a_p^2-3p^5) = 2i+1 \quad \text{if} \quad v(a_p) = i + \frac{1}{2}.$$
        On the other hand, \cite{Roz} shows that the map
        factors through $T^2+1$ if, for instance, $a_p^2 = 3 p^5$. 
        Thus both the reducible and irreducible possibilities for $\bar{V}_{k,a_p}$ arise.
        \end{remark}

       \subsection{The case \texorpdfstring{$a\geq 2i$}{}}\label{sec: bad large}
       In this section we prove a \textquote{superdiagonal conjecture}. 
       
       We first use the results of the previous section to describe $\bar{\Theta}_{k,a_{p}}$ when $a \geq 2i$ and $r \equiv a-i+n+(i-n)p \mod p(p-1)$ for some $1 \leq n \leq i-1$.
	
	\begin{theorem}[Superdiagonal Conjecture] \label{Shape theta a>2i bad congruence class}
		Let  $r \equiv a \mod{(p-1)}$ with $ 1 \leq a \leq p-1 $ and $r \geq i(p+1)+p $ with $ v(a_{p})  \in (i,i+1) $. Let $ s =  a-i+n+(i-n)p$ and $v(r-s) =t $  with $ 1 \leq  n \leq i-1 $. Assume that $a\geq 2i$. Then
		\begin{enumerate}
			 \item[$(i)$]  $\mathrm{ind}_{KZ}^{G}(V_r^{(i-t)}/V_r^{(i-t+1)}) \twoheadrightarrow \bar{\Theta}_{k,a_p}$ if $  t < n$. 
			 \item[$(ii)$]  $\mathrm{ind}_{KZ}^{G}(V_r^{(i-n)}/V_r^{(i-n+1)}) \twoheadrightarrow \bar{\Theta}_{k,a_p}$ if $ t \geq n$.
		\end{enumerate}
	\end{theorem}
	\begin{proof}
		By Lemma~\ref{JH factor Q 2i < a} $(ii)$ and Lemma~\ref{JH factor Q 2i = a} $(iv)$, we have		 
		\[ \text{JH factors of } Q(i)  = \{ V_{p-1-a+2l} \otimes D^{a-l}: 0 \leq l < i-n+1 \} \cup \text{JH factors of } \{ V_{r}^{(l)}/V_{r}^{(l+1)}: i-n+1     \leq l  \leq i \}.\]
		By Theorem~\ref{Elimination i < a and not in interval}, we have that the JH factors of the quotients $V_{r}^{(l)}/V_{r}^{(l+1)}$ for $0 \leq l < i-n$ vanish in $\bar{\Theta}_{k,a_{p}}$. Thus it remains to determine which of the following quotients $\{ V_{r}^{(i-T)}/V_{r}^{(i-T+1)}: 0     \leq T  \leq n \}$ survive in $\bar{\Theta}_{k,a_{p}}$. 
		
		\begin{enumerate}
		    \item[$(i)$] We first eliminate the quotients $V_{r}^{(i-T)}/V_{r}^{(i-T+1)}$ for $t<T$. By Theorem~\ref{thm: above super diagonal a large} we see that the quotient $\mathrm{ind}_{KZ}^{G}(V_{r}^{(i-T)}/V_{r}^{(i-T+1)})$  vanishes in $\bar{\Theta}_{k,a_{p}}$ if $t < T \leq n$. Indeed, the hypotheses $(i)$ and $(ii)$ of Theorem~\ref{thm: above super diagonal a large} hold, and the hypothesis $(iii)$ is vacuously true. Since $t< n$, by Theorem~\ref{eliminating JH sdcbd}, we have $\mathrm{ind}_{KZ}^{G}(V_{r}^{(i-T)}/V_{r}^{(i-T+1)})$ dies in $\bar{\Theta}_{k,a_{p}}$ for $t>T$ also vanish in $\bar{\Theta}_{k,a_{p}}$. 
            \item[$(ii)$]  Note that again by Theorem~\ref{eliminating JH sdcbd}, we have $\mathrm{ind}_{KZ}^{G}(V_{r}^{(i-T)}/V_{r}^{(i-T+1)})$ dies in $\bar{\Theta}_{k,a_{p}}$ for $0 \leq T<n$. Indeed $T < n \leq t$. \qedhere
		\end{enumerate}
	\end{proof}
    The above theorem can be illustrated pictorially as follows:
    		\begin{figure}[hbt!]
             \centering
			\resizebox{.6\textwidth}{!}{
				\begin{tikzpicture}[font=\small]
					\def\rows{10}
					\def\cols{11}
					\def\cwA{2}   
					\def\cw{1}    
					\def\rhA{1.2} 
					\def\rh{1}    
					
					\draw (0,0) -- (0,-\rhA-\rows*\rh);
					\draw (\cwA,0) -- (\cwA,-\rhA-\rows*\rh);
					\foreach \j in {1,...,\cols}{ \draw (\cwA+\j*\cw,0) -- (\cwA+\j*\cw,-\rhA-\rows*\rh); }
					\draw (0,0) -- (\cwA+\cols*\cw,0);
					\draw (0,-\rhA) -- (\cwA+\cols*\cw,-\rhA);
					\foreach \i in {1,...,\rows}{ \draw (0,-\rhA-\i*\rh) -- (\cwA+\cols*\cw,-\rhA-\i*\rh); }
					
					\draw (0,0) rectangle (\cwA,-\rhA);
					\draw (0,0) -- (\cwA,-\rhA); 
					
					\pgfmathsetmacro{\tx}{\cwA/3}
					\pgfmathsetmacro{\ty}{-2*\rhA/3}
					\pgfmathsetmacro{\Tx}{2*\cwA/3}
					\pgfmathsetmacro{\Ty}{-1*\rhA/3}
					
					\node at (\tx,\ty) {\(t\)};  
					\node at (\Tx,\Ty) {\(T\)};  
					
					\node at (\cwA+0.5*\cw,-0.5*\rhA) {$0$};
					\node at (\cwA+1.5*\cw,-0.5*\rhA) {$1$};
					\node at (\cwA+2.5*\cw,-0.5*\rhA) {$2$};
					\node at (\cwA+3.5*\cw,-0.5*\rhA) {$3$};
					\node at (\cwA+4.5*\cw,-0.5*\rhA) {$\cdot$};
					\node at (\cwA+5.5*\cw,-0.5*\rhA) {$\cdot$};
					\node at (\cwA+6.5*\cw,-0.5*\rhA) {$\cdot$};
					\node at (\cwA+7.5*\cw,-0.5*\rhA) {$\cdot$};
					\node at (\cwA+8.5*\cw,-0.5*\rhA) {$\cdot$};
					\node at (\cwA+9.5*\cw,-0.5*\rhA) {$n-1$};
					\node at (\cwA+10.5*\cw,-0.5*\rhA) {$n$};
					
					\node at (0.5*\cwA,-\rhA-0.5*\rh) {$t=1$};
					\node at (0.5*\cwA,-\rhA-1.5*\rh) {$t=2$};
					\node at (0.5*\cwA,-\rhA-2.5*\rh) {$t=3$};
					\node at (0.5*\cwA,-\rhA-3.5*\rh) {$t=4$};
					\node at (0.5*\cwA,-\rhA-4.5*\rh) {$\cdot$};
					\node at (0.5*\cwA,-\rhA-5.5*\rh) {$\cdot$};
					\node at (0.5*\cwA,-\rhA-6.5*\rh) {$\cdot$};
					\node at (0.5*\cwA,-\rhA-7.5*\rh) {$\cdot$};
					\node at (0.5*\cwA,-\rhA-8.5*\rh) {$t=n-1$};
					\node at (0.5*\cwA,-\rhA-9.5*\rh) {$t\geq n$};
					
					\foreach \i in {1,...,\cols}{
						\foreach \j in {1,...,\rows}{
							\pgfmathsetmacro{\cx}{\cwA+(\i-0.5)*\cw}
							\pgfmathsetmacro{\cy}{-\rhA-(\j-0.5)*\rh}
							
							\ifnum\i>4
							\ifnum\i<9
							\ifnum\j>4
							\ifnum\j<9
							\else
							\ifnum\i=\numexpr\j+1\relax
							\node at (\cx,\cy){\cmark};
							\else
							\node at (\cx,\cy){\xmark};
							\fi
							\fi
							\else
							\ifnum\i=\numexpr\j+1\relax
							\node at (\cx,\cy){\cmark};
							\else
							\node at (\cx,\cy){\xmark};
							\fi
							\fi
							\else
							\ifnum\i=\numexpr\j+1\relax
							\node at (\cx,\cy){\cmark};
							\else
							\node at (\cx,\cy){\xmark};
							\fi
							\fi
							\else
							\ifnum\i=\numexpr\j+1\relax
							\node at (\cx,\cy){\cmark};
							\else
							\node at (\cx,\cy){\xmark};
							\fi
							\fi
						}
					}
					
					\foreach \d in {-2,...,4}{
						\pgfmathtruncatemacro{\imin}{max(5,5+\d)}
						\pgfmathtruncatemacro{\imax}{min(8,8+\d)}
						\ifnum\imin<\imax
						\pgfmathsetmacro{\xstart}{\cwA+(\imin-0.5)*\cw}
						\pgfmathsetmacro{\ystart}{-\rhA-(\imin-\d-0.5)*\rh}
						\pgfmathsetmacro{\xend}{\cwA+(\imax-0.5)*\cw}
						\pgfmathsetmacro{\yend}{-\rhA-(\imax-\d-0.5)*\rh}
						\ifnum\d=1
						\draw[red, dash pattern=on 1.5pt off 3pt]
						(\xstart,\ystart) -- (\xend,\yend);
						\else
						\draw[dash pattern=on 1.2pt off 3pt]
						(\xstart,\ystart) -- (\xend,\yend);
						\fi
						\fi
					}
				\end{tikzpicture}
			}
           \caption{\label{fig:super diagonal conj. a>2i} Contribution of  JH factors when $a\geq 2i$ and $r\equiv a-i+n \mod p$.} 
		\end{figure}
        
        In the above figure, $\times$  inside a grid at position $(t,T)$ indicates that, for a given value of $t$ the image of $\mathrm{ind}_{KZ}^{G}(V_r^{(i-T)}/V_r^{(i-T+1)})$ vanishes  in $\bar{\Theta}_{k,a_p}$. Thus for a given value of $t$, the sub-quotient that possibly survives in $\bar{\Theta}_{k,a_p}$ is marked by \textcolor{red}{$\checkmark$}. In the above picture, we notice that \textcolor{red}{$\checkmark$} always appears along the diagonal ($t=T+1$). This observation motivates the name \textquote{super diagonal conjecture}.

        \begin{remark}
        When $p = 7$, $r = 54$ and $i = 3$,
        so that $a = 6$, $n = 2$, $s = 12$ and $t = 1$, then 
        \cite{Roz} shows that the cosocle of $V_r^{(1)}/V_r^{(2)}$ contributes,
        corroborating part $(i)$ of  Theorem~\ref{Shape theta a>2i bad congruence class}. 
        \end{remark}

       We now treat the cases that are not covered by Corollary~\ref{Shape of Theta i < a not in interval}~$(ii)$ and Theorem~\ref{Shape theta a>2i bad congruence class} in the  case $a=2i$, namely $r\equiv a-i-1,a-i \mod p$.

       \begin{theorem}\label{Shape theta a=2i bad n=-1,0}
		Let  $r \equiv a \mod{(p-1)}$ with $ 1 \leq a \leq p-1 $ and $r \geq i(p+1)+p $ with $ v(a_{p})  \in (i,i+1) $. Assume that $a = 2i$.
        \begin{enumerate}
            \item[$(i)$] If $r\equiv a-i-1~\mathrm{mod}~p$, then $\mathrm{ind}_{KZ}^{G}(V_{0}\otimes D^i) \twoheadrightarrow \bar{\Theta}_{k,a_p}$.
            \item[$(ii)$] If $r\equiv a-i ~\mathrm{mod}~p$, then $\mathrm{ind}_{KZ}^{G}(V_{p-1}\otimes D^i) \twoheadrightarrow \bar{\Theta}_{k,a_p}$. 
        \end{enumerate}
       \end{theorem} 
       \begin{proof}
           By Theorem~\ref{Elimination i < a and not in interval}, we have the image of  $\mathrm{ind}_{KZ}^{G}(V_r^{(n)})$ is same as the image of  $\mathrm{ind}_{KZ}^{G}(V_r^{(n+1)})$  in $\bar{\Theta}_{k,a_p}$ for $n=0,\ldots,i-1$. 
           Thus the quotient $V_r/V_r^{(i)}$ vanishes in $\bar{\Theta}_{k,a_p}$. Now $(i)$ and $(ii)$ follow from Lemma~\ref{JH factor Q 2i = a} $(ii)$ and $(iii)$ respectively.
       \end{proof}
        \subsection{The case  \texorpdfstring{$n=i$}{}}\label{sec: n=i}
         Finally, we treat the case $r\equiv a \mod p(p-1)$ for all $1 \leq a \leq p-1$. Thus $n=i$ and $s=a-i+n+(i-n)p$ equals $a$. We first consider the case $a>2i+1$. 
        \begin{theorem}\label{eliminating JH n=i bsd}
	        Let  $r \equiv a \mod{p(p-1)}$ with $ 1 \leq a \leq p-1 $ and $r \geq i(p+1)+p $ with $ v(a_{p})  \in (i,i+1) $. Let $v(r-a) =t $ for some $t\geq 1$ and $ 0 \leq T < i < a$. If  $t \geq T+1$ and $a >2i+1$, then the image of $ \mathrm{ind}_{KZ}^{G}(V_{r}^{(i-T)}) $ is the same as the image of $\mathrm{ind}_{KZ}^{G} (V_{r}^{(i-T+1)}) $ in $\bar{\Theta}_{k,a_{p}} $.
        \end{theorem}
        \begin{proof}
	        Consider the function
	          \begin{align*}
		        f_{2} = \sum_{\lambda \in \mathbb{F}_{p}^{\times}} \Bigg[ g_{2,p[\lambda]}^{0}, & \frac{1}{a_{p}} \bigg( \frac{p}{[\lambda]}\bigg)^{i-T}   (-\theta)^{T+1} X^{-T-1}Y^{r-(T+1)(p+1)+T+1} \Bigg] \\
		          & + \frac{(1-p)}{a_{p}}  \Bigg[ g_{2,0}^{0}, \binom{r}{r-i+T} (-\theta)^{i+1} X^{-T-1}Y^{r-(i+1)(p+1)+T+1} \Bigg].
	        \end{align*}
	          From Lemma~\ref{theta and T plus} it follows that $T^+ f_2$ vanishes modulo $p$.  Further
	          \begin{align}
		          T^{-}f_2 \equiv \left[ g_{1,0}^0, \frac{p^{i-T}(p-1)}{a_p}  \sum_{ \substack{a-i+T \leq   j < r-i+T \\  j \equiv a-i+T ~\mathrm{mod}~ (p-1)} }  \binom{r}{j} X^{r-j} Y^{j} \right]  ~\mathrm{mod}~p.
	        \end{align}
    
	        By the hypotheses $a>2i+1$ and $T<i$, we see that $i+T+1<p-1$. By Corollary~\ref{cor: binomial sums under congruences 2}, for $m=0,\ldots,i+T+1$  we have 
	        \begin{align*}
		        \sum_{ \substack{a-i+T \leq   j < r-i+T \\  j \equiv a-i+T ~\mathrm{mod}~ (p-1)} }  \binom{r}{j} \binom{j}{m} \equiv    \binom{a}{a-i+T} \binom{a-i+T}{m}-\binom{r}{r-i+T}\binom{r-i+T}{m} \equiv 0 \mod p^{t},
	        \end{align*}
	        where the last congruence follows from Lemma~\ref{binomial coefficient under congruences} $(i)$ and $r\equiv a \mod p^t$. Applying Lemma~\ref{lem:choice of beta} with $m$ there equal to $i+T+1$, we obtain $\alpha_j$ satisfying
	        \begin{enumerate}
		        \item[$(1)$] $\alpha_{j} \equiv  \binom{r}{j}$ mod $p^{t}$, for all     $a-i+T \leq  j < r-i+T$ and $j\equiv a-i+T \mod (p-1)$, and 
		        \item[$(2)$] $\sum\limits_{\substack{ a-i+T \leq j < r-i+T \\ j \equiv a-i+T ~\mathrm{mod}~(p-1)}}^{} \alpha_{j} \binom{j}{m} \equiv 0$ mod $p^{i+t+2-m}$ for $m=0, \ldots, i+T+1$.
	          \end{enumerate} 
	        Let 
	        \begin{align*}
		          f_1 =  \left[g_{1,0}^0, \frac{(p-1)p^{i-T}}{a_p^2} \sum_{ \substack{a-i+T \leq   j < r-i+T \\  j \equiv a-i+T ~\mathrm{mod}~ (p-1)} } \alpha_{j} X^{r-j} Y^{j} \right]
	        \end{align*}
epbg	        From $(2)$ and $a>2i+1$ it follows that $T^+ f_1$ vanishes modulo $p$. Using $ 2T< 2i <a < p-1$ we see that $T^{-}f_1$ vanishes modulo $p$. Using $t\geq T+1$ and $(1)$, we see that
	        \begin{align*}
		          T^{-}f_2 -a_p f_1 \equiv \left[ g_{1,0}^0, \frac{p^{i-T}(p-1)}{a_p}  \sum_{ \substack{a-i+T \leq   j < r-i+T \\  j \equiv a-i+T ~\mathrm{mod}~ (p-1)} } \left( \binom{r}{j} - \alpha_j \right) X^{r-j} Y^{j} \right] \equiv  0 ~\mathrm{mod}~p.
	        \end{align*}
	        Thus 
	        $$
               (T-a_p)(f_2 + f_1) \equiv -a_p f_2 \equiv -  \Bigg[ g_{2,0}^{0},  \binom{r}{r-i+T} (-\theta)^{i+1} X^{-T-1}Y^{r-(i+1)(p+1)+T+1} \Bigg] \mod p,  
            $$
	          since $T<i$. Note 
	        \begin{align*}
		          (-\theta)^{i+1} X^{-T-1}Y^{r-(i+1)(p+1)+T+1} &= (-\theta)^{i-T} \left(\sum\limits_{j=0}^{T+1} (-1)^{j} \binom{T+1}{j} X^{j(p-1)}Y^{r-(i-T)(p+1)-j(p-1)}\right)\\
		           & \equiv (-\theta)^{i-T}  (Y^{r-(i-T)(p+1)} - X^{p-1}Y^{r-(i-T)(p+1)-(p-1)}) \mod{V}_r^{(i-T+1)}.
	          \end{align*}
	        By Lemma~\ref{generating polynomial quotient},  the above polynomial generates ${V}_r^{(i-T)}/{V}_r^{(i-T+1)}$. As $\binom{r}{r-i+T} \not \equiv 0 \mod p$ by Lucas' theorem, this completes the proof of the theorem.
        \end{proof}
        \begin{remark}\label{remark a=2i+1 and n=i below super diagonal}
	        The argument of the previous theorem works in the case $a=2i+1$ and $v(a_p^2)<2i+1$. One can show that the conclusion of the previous theorem is true even in the case $a=2i+1$ and $v(a_p^2)>2i+1$ by considering the functions $f_2' = (a_p^2/p^{2i+1}) f_2$ and $f_1' = (a_p^2/p^{2i+1}) f_1$ cf. \cite[Theorem 8.7]{BG15}.  
        \end{remark}
         We now prove an analogue of the previous theorem if $a\leq 2i$ and $r \equiv a \mod p(p-1)$. By Lemmas~\ref{JH factor Q 2i = a} $(i)$ $(a)$ and \ref{JH factor Q 2i > a} $(i)$ $(a)$ the JH factors of $Q(i)$ are
         \begin{align}\label{JH factors n=i and <2i}
           \{ V_{p-1-a} \otimes D^a\} \cup \text{JH factors of }\{ V_{r}^{(i-T)}/V_{r}^{(i-T+1)}: 2i-a+1 \leq T \leq i-1 \}.
         \end{align}
        \begin{theorem}\label{eliminating JH n=i bsd1}
	        Let  $r \equiv a \mod{p(p-1)}$ with $ 1 \leq a \leq p-1 $ and $r \geq i(p+1)+p $ with $ v(a_{p})  \in (i,i+1) $. Let $v(r-a) =t $ for some $t\geq 1$. Let $ a \leq 2i$ and $ 2i-a+1 \leq T < i < a$. If  $t\geq T+a-2i$, then the image of $ \mathrm{ind}_{KZ}^{G}(V_{r}^{(i-T)}) $ is the same as the image of $\mathrm{ind}_{KZ}^{G} (V_{r}^{(i-T+1)}) $ in $\bar{\Theta}_{k,a_{p}} $.
        \end{theorem}
        \begin{proof}
	        Let $f_2$ and $f_1$ be as defined in Theorem~\ref{eliminating JH n=i bsd}. Consider the functions
	        \begin{align*}
		        f'_{2} &= \frac{a_p^2}{p^a} f_2 = \sum_{\lambda \in \mathbb{F}_{p}^{\times}} \Bigg[ g_{2,p[\lambda]}^{0},  \frac{a_p}{p^a} \bigg( \frac{p}{[\lambda]}\bigg)^{i-T}  (-\theta)^{T+1} X^{-T-1}Y^{r-(T+1)(p+1)+T+1} \Bigg] \\
		        & \qquad \qquad \qquad \qquad \qquad + \frac{(1-p)a_p}{p^a}  \Bigg[ g_{2,0}^{0}, \binom{r}{r-i+T}(\theta)^{i+1} X^{-T-1}Y^{r-(i+1)(p+1)+T+1} \Bigg],\\	
		        f'_{1} &= \frac{a_p^2}{p^a} f_1 = \left[g_{1,0}^0, \frac{(p-1)}{p^{a-i+T}} \sum_{ \substack{a-i+T \leq   j < r-i+T \\  j \equiv a-i+T ~\mathrm{mod}~ (p-1)} } \alpha_{j} X^{r-j} Y^{j} \right].
	        \end{align*}
	        Since $a\leq 2i$, we get $v(a_p^2/p^a ) >0$. It follows that $T^+ f'_2$ vanishes modulo $p$ as before.  Also, $-a_p f_2'$ vanishes. Furthermore, we have
	        \begin{align}
		          T^{-}f'_2 - a_p f'_1 \equiv \left[ g_{1,0}^0, \frac{(p-1)a_p}{p^{a-i+T}}  \sum_{ \substack{a-i+T \leq  j < r-i+T \\  j \equiv a-i+T ~\mathrm{mod}~ (p-1)} }  \left(\binom{r}{j} - \alpha_j \right) X^{r-j} Y^{j} \right] .
	        \end{align}
	          Using $\alpha_{j} \equiv  \binom{r}{j}$ mod $p^{t}$ and $t-T \geq a-2i$  we see that $T^{-}f'_2 - a_p f'_1$ vanishes modulo $p$. Clearly $T^{-}f'_1$ vanishes modulo $p$.  Thus
	        \[
	            (T-a_p)(f'_2+f'_1) \equiv T^{+} f'_1 \mod p.
	        \]
	        Using
	        $$
            \sum\limits_{\substack{ a-i+T \leq j < r-i+T \\ j \equiv a-i+T ~\mathrm{mod}~(p-1)}}^{} \alpha_{j} \binom{j}{m} \equiv 0 \mod p^{i+t+2-m}
            $$
	        for $m=0, \ldots, a-i+T$ we see that
	        \begin{align*}
		        T^{+}f'_1 \equiv \left[g_{2,0}^0, (p-1) \alpha_{a-i+T} X^{r-(a-i+T)} Y^{a-i+T} \right] \mod p.
	          \end{align*}
	        Since $\alpha_{a-i+T} \equiv \binom{r}{a-i+T} \mod p$ and $\binom{r}{a-i+T} \equiv \binom{a}{a-i+T} \not\equiv 0 \mod p$, we get $\alpha_{a - i + T} \not \equiv 0 \mod p$. Note that $X^{r-(a-i+T)} Y^{a-i+T}- X^{i-T} Y^{r-i+T}$ lies in $V_r^{(i-T)} \setminus V_r^{(i-T+1)}$. Further by \cite[Lemma 2.13]{GR19}, it follows that $X^{r-(a-i+T)} Y^{a-i+T}- X^{i-T} Y^{r-i+T}$ generates $V_r^{(i-T)} /V_r^{(i-T+1)}$. This finishes the proof as $X^{i-T} Y^{r-i+T}$ vanishes in $Q(i)$.
        \end{proof}

        Note that 
        \[
            \max\{t,t+2i+1-a\} = \begin{cases}
            t+2i+1-a & \mathrm{if}~a \leq 2i, \\
                t & \mathrm{if}~a\geq 2i+1.               
            \end{cases}
        \]
        The previous two theorems eliminate the JH factors of $V_{r}^{(i-T)}/V_r^{(i-T+1)}$ when $r\equiv a \mod (p-1)$, $p^t\parallel(r-a)$ and $T < \max\{t,t+2i-a+1\}$ (these factors are below the `superdiagonal'). We now eliminate the factors when $T > \max\{t,t+2i-a+1\}$ (above the `superdiagonal').  The result below is a variant of Lemma~\ref{lem: choice beta asd a>2i}. 
        \begin{lemma}\label{lem: choice beta asd n=i}
             Let $1 \leq i \leq p-2$. Let $r \geq i(p+1)+p $, $r \equiv a \mod{p(p-1)}$ with $ 1 \leq a \leq p-1 $. Let $v(r-a) =t $  with $t\geq 1$. Fix an integer $1  \leq T \leq i$. Assume $a>2i-2T$ and $i+t+1<a-i+T$.  There exists $p$-adic integers $\beta_0, \ldots, \beta_{i-T} \in \mathbb{Z}_p$ satisfying 	
	        \begin{enumerate}
		        \item[$(i)$] For $m=0,\ldots,i-T$, we have
		        $$
		           \sum\limits_{\substack{a - i + T \leq j < r-i+T \\ j \equiv a-i+T \mod(p-1)}} \left( p \sum\limits_{l = 0}^{i-T-1} \beta_l \binom{r-l}{j} + \beta_{i-T} \binom{r-i+T}{j} \right) \binom{j}{m}  \equiv p^{t+1} \delta_{i-T,m} \mod p^{t+2}.
		        $$
		        \item[$(ii)$] For $m=i-T+1,\ldots,i+t+1$, we have
		          $$
		               \sum\limits_{\substack{a - i + T \leq j < r-i + T \\ j \equiv a - i + T \mod(p-1)}} \left( p \sum\limits_{l = 0}^{i-T-1} \beta_l \binom{r - l}{j} + \beta_{i-T} \binom{r - i + T}{j} \right) \binom{j}{m}  \equiv 0 \mod p^{t+1}.
                $$ 
	        \end{enumerate}
        \end{lemma}
        \begin{proof}
              We  first show $(i)$ holds. To solve the congruences $(i)$, we now compute the coefficients of $\beta_l$  modulo $p^{t+2}$. By Corollary~\ref{cor: binomial sums under congruences 1}, for $0 \leq m \leq i-T$ and $0 \leq l \leq i-T-1$ we have
            \begin{align*}
	            \sum_{\substack{0 \leq j < r - i + T \\ j \equiv a - i + T ~\mathrm{mod}~ (p-1)}}\binom{r- l}{j} \binom{j}{m}\equiv & \left\{\binom{r - l}{m} - \binom{a - l}{m}\right\}\binom{a - l - m}{a - i + T - m} \\
                & \quad + \binom{a - l}{a - i + T} \binom{a - i + T}{m}- \binom{r - l}{r - i + T} \binom{r - i + T}{m}  \mod p^{t+1}\\
                & \equiv \binom{r - l}{m} \left( \binom{a - l - m}{a - i + T - m} - \binom{r - l - m}{r - i + T - m}\right) \mod p^{t+1},
	        \end{align*}
            where we have used $a-i+T > i-T$ in obtaining the first congruence.   For $0 \leq m \leq i-T$ and $0 \leq l \leq i-T-1$ we have
            \begin{align*}
                \binom{a - l - m}{a - i + T - m} - \binom{r - l - m}{r - i + T - m} &=\binom{a - l - m}{i - T - l} - \binom{r - l - m}{i - T - l} \\
                &\equiv (a-r) \binom{a - l - m}{a - i + T - m} (H_{a-l-m} - H_{a-i+T-m}) \mod p^{t+1},
            \end{align*}
            by Lemma~\ref{binomial coefficient under congruences} $(ii)$.
            Substituting this  above, for $0 \leq m \leq i-T$ and $0 \leq l \leq i-T-1$ and noting $\binom{r-l}{m} \equiv \binom{a-l}{m} \mod p$ by Lucas' theorem, we get
            \begin{align*}
                \sum_{\substack{0 \leq j < r - i + T \\ j \equiv a - i + T ~\mathrm{mod}~ (p-1)}}\binom{r- l}{j} \binom{j}{m} &\equiv (a-r) \binom{r-l}{m}\binom{a - l - m}{a - i + T - m} (H_{a-l-m} - H_{a-i+T-m}) \\
                &\equiv (a-r) \binom{a-l}{m}\binom{a - l - m}{a - i + T - m} (H_{a-l-m} - H_{a-i+T-m}) \\
                &\equiv (a-r)\binom{a-l}{a-i+T} \binom{a-i+T}{m}( H_{a-l-m}-H_{a-i+T-m})  \mod p^{t+1}.
            \end{align*}
            If $l=i-T$ and $m=0,\ldots,i-T$ , then since $i-T< a-i+T$ we see that the numbers $m=0,\ldots,i-T$ are not congruent to $a-i+T$ modulo $(p-1)$. Hence  
            \begin{align*}
                \sum_{\substack{0 \leq j < r - i + T \\ j \equiv a - i + T \mod (p-1)}} \binom{r - i + T}{j} \binom{j}{m} & 
                \equiv \sum_{\substack{m < j < r - i + T \\ j \equiv a - i + T \mod (p-1)}} \binom{r - i + T}{j} \binom{j}{m}\\
                &\equiv  \frac{p(a - r)}{a - i + T - m} \binom{a - i + T}{m} \mod p^{t+2},
	        \end{align*}
            by Lemma~\ref{lem: binomial sum Doc. Math general}. 
            
            Let 
            \begin{equation}\label{eq: A asd n=i}
	            A = \left[ \begin{array}{c|c}
                    p (a-r) \binom{a-l}{a-i+T} \binom{a-i+T}{m}( H_{a-l-m}-H_{a-i+T-m})  
		             &  \frac{p(a- r)}{a - i + T - m} \binom{a - i + T}{m} 
	              \end{array}\right]_{m=0,\ldots,i-T},
            \end{equation}
            where the range of $l$ in the left block is $0,\ldots,i-T-1$ and the range of $l$ in the right most column is $i-T$.

            To solve congruences $(i)$, it is enough to show that the following system
             \begin{equation}\label{eq: matrix n=i}
                A \begin{bmatrix} \beta_0 \\ \vdots \\ \beta_{i-T-1} \\ \beta_{i-T} \end{bmatrix} \equiv \begin{bmatrix} 0\\ \vdots \\ 0 \\ p^{t+1} \end{bmatrix} \mod p^{t+2} \mathbb{Z}_{p}
            \end{equation}
            has a solution in $\mathbb{Z}_{p}$. 

            To solve the above congruence, it is enough to show the following equation
            \begin{equation}\label{eq: matrix linear eq. in matrix form n=i}
                A \begin{bmatrix} \beta_0 \\ \vdots \\ \beta_{i-T-1} \\ \beta_{i-T} \end{bmatrix} = \begin{bmatrix} 0\\ \vdots \\ 0 \\ p^{t+1} \end{bmatrix}
            \end{equation}
            has a solution over $\mathbb{Q}$ with $\beta_l \in \mathbb{Z}_{p}$. To show this we use Cramer's rule.       
            We first compute $\det(A)$. Pulling out $p(a-r)\binom{a-i+T}{m}$ from the $m^{\mathrm{th}}$-row and $\binom{a-l}{a-i+T}$ from the $l^{\mathrm{th}}$-column we get 
            \begin{align*}
                  \det(A) &= \left(\prod_{m=0}^{i-T} p(a-r)\binom{a-i+T}{m} \right) \times \left( \prod_{l=0}^{i-T}\binom{a-l}{a-i+T}\right) \\
                  & \qquad \qquad \times \det\left[ \begin{array}{c|c}
                    ( H_{a-l-m}-H_{a-i+T-m})   & \frac{1}{a - i + T - m} 
	              \end{array}\right]_{m=0,\ldots,i-T}.
            \end{align*}
            Applying the column operations 
            \begin{align*}
                    C_0 &\rightarrow\; C_0 - C_{1} \\
                    &~~\vdots \\
                    C_l &\rightarrow\; C_l - C_{l+1}, \\
                    &~~\vdots \\
                    C_{i-T-2} &\rightarrow\; C_{i-T-2} - C_{i-T-1}
            \end{align*}
            we see that 
            \[
                  \det(A) = \left(\prod_{m=0}^{i-T} p(a-r)\binom{a-i+T}{m} \right) \times \left( \prod_{l=0}^{i-T}\binom{a-l}{a-i+T}\right) \times \det\left[ \begin{array}{c} \frac{1}{a - l - m} 
	              \end{array}\right]_{m,l=0,\ldots,i-T}.
            \]
            From Cauchy’s double alternant (cf. \cite[(2.7)]{Kra99}), we have
            \begin{align*}
                \det_{0\leq m,l \leq i-T}\left( \frac{1}{X_m+Y_l}\right) = \frac{\prod\limits_{0 \leq m < l \leq i-T} (X_{m}-X_{l})(Y_{m}-Y_{l})}{\prod\limits_{0\leq m,l \leq i-T}(X_{m}+Y_{l})}
            \end{align*}
            for all integers $X_{0}, \ldots, X_{i-T}$ and $Y_{0}, \ldots, Y_{i-T}$. Taking $X_{m}=a-m$ and $Y_{l}=-l$ in the above formula we get 
            \begin{align*}
                \det_{m,l=0,\ldots,i-T} \left[ \begin{array}{c} \frac{1}{a - l - m} 
	              \end{array}\right] = \frac{\prod\limits_{0 \leq m < l < i-T} (l-m)^2}{\prod\limits_{0\leq m,l \leq i-T}(a-l-m)} 
                  ,
            \end{align*}
            which belongs to $\Z_p^\times$, since all quantities appearing in the above expression are positive and strictly less than $p$. This shows that $\det(A) \in (p(a-r))^{(i-T+1)} \mathbb{Z}_p^\times$ and \eqref{eq: matrix linear eq. in matrix form n=i} has solutions in $\Q$ by Cramer's rule. Since every entry of $A$ is a multiple of $p(a-r)$, we see that $\det(A_{i-T,l}) \in (p(a-r))^{(i-T)} \mathbb{Z}_p$, where $A_{i-T,l}$ is the minor of the $(i-T,l)$-entry of $A$. Again by Cramer's rule we obtain \eqref{eq: matrix linear eq. in matrix form n=i} has solutions in $\Z_p$.  This proves $(i)$.

            To prove  $(ii)$, we show that the coefficient of each $\beta_l$ vanishes modulo $p^{t+1}$. By  hypothesis, we have $i+t+1< a-i+T \leq a \leq p-1$. By Corollary~\ref{cor: binomial sums under congruences 2}, for $i-T+1 \leq m \leq i+t+1$ and $0 \leq l \leq i-T-1$ we have
            \begin{align*}
                \sum_{\substack{a-i+T \leq j < r - i + T \\ j \equiv a - i + T ~\mathrm{mod}~ (p-1)}}\binom{r- l}{j} \binom{j}{m}\equiv  \binom{a - l}{a- i + T} \binom{a - i + T}{m}- \binom{r - l}{r - i + T} \binom{r - i + T}{m}  \equiv 0 \mod p^{t},
            \end{align*}
            where the last congruence follows from Lemma~\ref{binomial coefficient under congruences} $(i)$. This shows that the coefficient of $\beta_l$ vanishes modulo $p^{t+1}$ for $l=0,\ldots,i-T-1$. The vanishing of the coefficient of $\beta_{i-T}$ modulo $p^{t+1}$ follows from Lemma~\ref{lem: binomial sum Doc. Math general} and the hypotheses $m+1 \leq i+t+2 \leq a-i+T$ and $r\equiv a \mod p^t$. This proves $(ii)$ and the lemma follows.
        \end{proof}
        \begin{theorem}\label{eliminating JH n=i asd}
        	Let  $a_p \in \bar{\Q}_p$ with $ v(a_{p})  \in (i,i+1) $ and $1 \leq i \leq p-2$. Let $r \geq i(p+1)+p $, $r \equiv a \mod{p(p-1)}$ with $ 1 \leq a \leq p-1 $ and $t=v(r-a) $. Fix an integer $1 \leq T \leq i$. Assume that the following hold:
            \begin{enumerate}
                \item[$(i)$] $ i < a$
                \item[$(ii)$] $T > \max\{t,t+2i-a+1\}$.
            \end{enumerate}
            Then the image of $ \mathrm{ind}_{KZ}^{G}(V_{r}^{(i-T)}) $ is the same as the image of $ \mathrm{ind}_{KZ}^{G}(V_{r}^{(i-T+1)}) $ in $ \bar{\Theta}_{k,a_{p}} $.
        \end{theorem}
        \begin{proof}
            Using the above hypotheses, we have $a \geq 2i+t-T+2 >2i-2T$ and $i+t+1<a-i+T$. So the hypotheses of Lemma~\ref{lem: choice beta asd n=i} hold. Let $\beta_l \in \Z_p$ be as in Lemma~\ref{lem: choice beta asd n=i}. Consider the function (as in Theorem~\ref{eliminating JH n=i bsd})
            \begin{align*}
	            f_2 &= \sum_{\lambda \in \mathbb{F}_p^\times} \left[ g^{0}_{2,p[\lambda]}, \sum_{l = 0}^{i-T-1} \frac{\beta_l}{p^{l+t}} [\lambda]^{l-i+T} (-\theta)^{l+t+1} X^{-t-1} Y^{r - (l+t+1)(p+1)+t+1} \right] \\
                &\qquad \qquad + \sum_{\lambda \in \mathbb{F}_p^\times} \left[ g^{0}_{2,p[\lambda]},\frac{\beta_{i-T}}{p^{i-T+t+1}} (-\theta)^{i -T+t+ 2} X^{-t-2} Y^{r - (i-T+t+2)(p+1)+t+2}\right] \\
	            & \qquad \qquad  + \left[ g^{0}_{2,0}, -\left( \sum_{l = 0}^{i-T-1} p \beta_l \binom{r-l}{r-i+T} + \beta_{i-T} \right) \frac{p-1}{p^{i-T+t+1}}(-\theta)^{i+1} X^{-T-1} Y^{r-(i+1)(p+1)+T+1} \right].
	        \end{align*}
            By Lemma~\ref{theta and T plus} and $ t < T$ we have $T^{+} f_2$  vanishes. Also $-a_p f_2$ vanishes since  $i - T + t + 1 \leq i <  v(a_p)$. Observe that
            \begin{align*}
	            T^{-} f_2 \equiv \left[ g^{0}_{1,0}, \frac{p-1}{p^{t+1}} \sum_{\substack{a - i + T \leq j < r- i + T \\ j \equiv a-i+T \mod(p-1)}} \left( p \sum_{l=0}^{i-T-1} \beta_l \binom{r-l}{j} +\beta_{i-T} \binom{r-i+T}{j}\right) X^{r - j} Y^j \right].
	        \end{align*}
             By  Lemma~\ref{lem: choice beta asd n=i} $(i),(ii)$ and Lemma~\ref{lem:choice of beta}, there exist $\alpha_j$ such that
	          \begin{enumerate}
		          \item[(1)] \[ \alpha_j \equiv p \sum_{l = 0}^{i-T-1} \beta_l \binom{r-l}{j} + \beta_{i-T} \binom{r - i + T}{j} \mod{p^{t+1}}\]
		        for $a - i + T \leq j < r - i + T$ with $j \equiv a - i + T \mod (p-1)$, and
                 \item[$(2)$]  \[ \sum_{\substack{a - i + T \leq j < r-i + T \\ j \equiv a - i + T \mod(p-1)}} \alpha_j \binom{j}{m} \equiv 0 \mod{p^{i+t+2-m}} \]
                 for  $m = 0, \ldots, i +t + 1$.
	        \end{enumerate}
            Let 
            \begin{align*}
	            f_1 = \left[g^{0}_{1,0}, \frac{p-1}{ p^{t+1} a_p} \sum_{\substack{a - i + T \leq j < r - i + T \\ j \equiv a - i + T \mod{(p-1)}}} \alpha_j X^{r-j} Y^j \right].
	        \end{align*}
	        It follows from $i+t+1 < a-i+T$ and $(2)$ that $T^{+}f_1$ vanishes. Since $t+T+1 <  i+t+2 \leq  a-i+T \leq p-1$, we get $T^{-} f_1$ vanishes. Thus
	          \begin{align*}
	            (T - a_p)(f_2 + f_1) \equiv  T^{-} f_2 - a_p f_1 \equiv   \left[g^{0}_{1,0}, F(X,Y) \right] \mod p,
	        \end{align*}
	        where
	          \begin{align*}
		          F(X,Y) = \frac{p-1}{p^{t+1}}\sum_{\substack{a - i + T \leq j < r- i + T \\ j \equiv a - i + T \mod(p-1)}} \left( p \sum_{l=0}^{i-T-1} \beta_l \binom{r-l}{j}+ \beta_{i-T} \binom{r-i+T}{j}- \alpha_j \right) X^{r-j} Y^j.
	          \end{align*}
            As observed in the proof of Theorem~\ref{thm: above super diagonal a large} this the above polynomial generates $V_{r}^{(i-T)}/V_r^{(i-T+1)}$.
        \end{proof}

       We now collect the results proved in this section describing $\bar{\Theta}_{k,a_p}$.

       \begin{theorem}\label{Shape theta a>2i bad n=i}
	        Let  $r \equiv a \mod{p(p-1)}$ with $ 1 \leq a \leq p-1 $ and $r \geq i(p+1)+p $ with $ v(a_{p})  \in (i,i+1) $. Assume $a\geq 2i+1$. Let $v(r-a) =t $ for some $t\geq 1$.  If $a=2i+1$, then we assume $v(a_p) \neq i +(1/2)$. Then
            \begin{enumerate}
                \item[$(i)$] If $1 \leq t < i$, then $\mathrm{ind}_{KZ}^{G}(V_{r}^{(i-t)}/V_{r}^{(i-t+1)}) \twoheadrightarrow \bar{\Theta}_{k,a_p}$.
                \item[$(ii)$] If $t \geq i$, then $\mathrm{ind}_{KZ}^{G}(V_{p-a-1} \otimes D^{a}) \twoheadrightarrow \bar{\Theta}_{k,a_p}$.
            \end{enumerate}
        \end{theorem}
        \begin{proof}
            Note that $i< a$. By Lemma~\ref{JH factor Q 2i < a} (ii), we have JH factors of $Q(i)$ are
            \[
               \{V_{p-a-1} \otimes D^{a}\} \cup \text{JH factors of }\{ V_r^{(i-T)}/V_r^{(i-T+1)}: 0 \leq T <i\}.
            \]
            Note that $V_{p-a-1} \otimes D^{a}$ is the cosocle of $V_r/V_r^{(1)}$. We eliminate all but one quotient described above.
            \begin{enumerate}
                \item[$(i)$] If $0\leq T < t$, then $\mathrm{ind}_{KZ}^{G}(V_{r}^{(i-T)})$ and $\mathrm{ind}_{KZ}^{G}(V_{r}^{(i-T+1)})$ have the same image in $\bar{\Theta}_{k,a_p}$ by Theorem~\ref{eliminating JH n=i bsd} and Remark~\ref{remark a=2i+1 and n=i below super diagonal}. If $t < T \leq i$, then $\mathrm{ind}_{KZ}^{G}(V_{r}^{(i-T)})$ and $\mathrm{ind}_{KZ}^{G}(V_{r}^{(i-T+1)})$ have the same image in $\bar{\Theta}_{k,a_p}$ by Theorem~\ref{eliminating JH n=i asd}. Here we used $\max\{t,t+2i-a+1\}=t$ as $a\geq2i+1$. This eliminates all but one JH factor listed above, namely $V_{r}^{(i-t)}/V_{r}^{(i-t+1)}$. 
                \item[$(ii)$] If $0\leq T < i$, then $\mathrm{ind}_{KZ}^{G}(V_{r}^{(i-T)})$ and $\mathrm{ind}_{KZ}^{G}(V_{r}^{(i-T+1)})$ have the same image in $\bar{\Theta}_{k,a_p}$ by Theorem~\ref{eliminating JH n=i bsd} and Remark~\ref{remark a=2i+1 and n=i below super diagonal}. This shows that the only JH factor surviving in $\bar{\Theta}_{k,a_p}$ from the above list is $V_{p-a-1} \otimes D^{a}$. \qedhere
            \end{enumerate}
        \end{proof}

        \begin{theorem}\label{Shape theta a<2i bad n=i}
	        Let  $r \equiv a \mod{p(p-1)}$ with $ 1 \leq a \leq p-1 $ and $r \geq i(p+1)+p $ with $ v(a_{p})  \in (i,i+1) $. Assume  $ i< a \leq 2i$. Let $v(r-a) =t $ for some $t\geq 1$. 
            \begin{enumerate}
                \item[$(i)$] If $1 \leq t < a-i$, then $\mathrm{ind}_{KZ}^{G}(V_{r}^{(a-i-t-1)}/V_{r}^{(a-i-t)}) \twoheadrightarrow \bar{\Theta}_{k,a_p}$.
                \item[$(ii)$] If $t \geq a-i$, then $\mathrm{ind}_{KZ}^{G}(V_{p-a-1} \otimes D^{a}) \twoheadrightarrow \bar{\Theta}_{k,a_p}$.
            \end{enumerate}
        \end{theorem}
        \begin{proof}
            By \eqref{JH factors n=i and <2i} the JH factors of $Q(i)$ are 
            \begin{align*}
                \{ V_{p-1-a} \otimes D^a\} \cup \text{JH factors of }\{ V_{r}^{(i-T)}/V_{r}^{(i-T+1)}: 2i-a+1 \leq T \leq i-1 \}.
            \end{align*}
            Note that $V_{p-a-1} \otimes D^{a}$ is the cosocle of $V_r/V_r^{(1)}$. As in the previous theorem, we eliminate all but one quotient from the above list.
                  \begin{enumerate}
                \item[$(i)$] If $2i-a+1 \leq T \leq t+2i-a$, then $\mathrm{ind}_{KZ}^{G}(V_{r}^{(i-T)})$ and $\mathrm{ind}_{KZ}^{G}(V_{r}^{(i-T+1)})$ have the same image in $\bar{\Theta}_{k,a_p}$ by Theorem~\ref{eliminating JH n=i bsd1}. If $t+2i-a+1 < T \leq i$, then $\mathrm{ind}_{KZ}^{G}(V_{r}^{(i-T)})$ and $\mathrm{ind}_{KZ}^{G}(V_{r}^{(i-T+1)})$ have the same image in $\bar{\Theta}_{k,a_p}$ by Theorem~\ref{eliminating JH n=i asd}. Here we used $\max\{t,t+2i-a+1\}=t+2i-a+1$ as $a\leq 2i$. This eliminates all but one JH factor listed above, namely $V_{r}^{(a-i-t-1)}/V_{r}^{(a-i-t)}$. 
                \item[$(ii)$] If $2i-a+1\leq T < i$, then $\mathrm{ind}_{KZ}^{G}(V_{r}^{(i-T)})$ and $\mathrm{ind}_{KZ}^{G}(V_{r}^{(i-T+1)})$ have the same image in $\bar{\Theta}_{k,a_p}$ by Theorem~\ref{eliminating JH n=i bsd1}. This shows that the only JH factor surviving in $\bar{\Theta}_{k,a_p}$ from the above list is $V_{p-a-1} \otimes D^{a}$. \qedhere
            \end{enumerate}    
        \end{proof}

\vspace{-.3 mm}
\end{document}